\input amstex

\loadeufm
\loadmsbm
\loadeufm

\documentstyle{amsppt}
\input amstex
\catcode `\@=11
\def\logo@{}
\catcode `\@=12
\magnification \magstep1
\NoRunningHeads
\NoBlackBoxes
\TagsOnLeft

\def \={\ = \ }
\def \+{\ +\ }
\def \-{\ - \ }

\def \b|{\big |}

\def \g1{\Gamma_1}

\def \nfp{\demo\nofrills{Proof:\usualspace\usualspace }}

\def\rarr#1#2{\smash{\mathop{\hbox to .5in{\rightarrowfill}}
 	 \limits^{\scriptstyle#1}_{\scriptstyle#2}}}

\def\larr#1#2{\smash{\mathop{\hbox to .5in{\leftarrowfill}}
	  \limits^{\scriptstyle#1}_{\scriptstyle#2}}}

\def\swarr#1#2 {\llap{$\scriptstyle #1$}  \swarrow
  	\vcenter to .5in{}\rlap{$\scriptstyle #2$}}

\topmatter
\title Diophantine Geometry over Groups X:
\centerline{The Elementary Theory of Free Products of Groups}
\endtitle
\author
\centerline{
Z. Sela${}^{1,2}$}
\endauthor
\footnote""{${}^1$Hebrew University, Jerusalem 91904, Israel.}
\footnote""{${}^2$Partially supported by an Israel academy of sciences fellowship.}
\abstract\nofrills{}
This paper is the 10th in a sequence
on the structure of
sets of solutions to systems of equations over groups, projections
of such sets (Diophantine sets), and the structure of definable sets over few classes of groups.
In the 10th paper we study the first order theory of free products of arbitrary groups. With a given (coefficient-free) predicate,
we associate (non-canonically) a finite collection of graded resolutions. These graded resolutions enable us to reduce uniformly
an
arbitrary coefficient-free sentence over a free product to a finite disjunction of conjunctions of sentences over the factors 
of the free product. The graded resolutions enable a uniform quantifier elimination over free products, where a coefficient-free
predicate (over free products) 
is shown to be equivalent to a predicate in an extended language, 
that contains finitely many quantifiers over the factors of the free product, and
only 3 quantifiers over the ambient free product. 

These uniform reductions allow us to answer affirmatively a question of R. Vaught on the elementary equivalence of free products
of pairs of elementarily equivalent groups, and to obtain a generalization of Tarski's problem on the elementary equivalence of 
non-abelian free groups. Finally, we use the resolutions that are associated with a predicate over free products, to
prove that the free product of stable groups is stable, generalizing the main theorem of [Se9] on the stability of free groups.
\endabstract
\endtopmatter

\document

\baselineskip 12pt

In the first 6 papers in the sequence on Diophantine geometry over groups we
studied sets of solutions to systems of equations in a free group, and
developed basic techniques and objects that are required for the analysis of sentences and
elementary sets defined over a free group.  The techniques we developed,
enabled us to present an iterative procedure that analyzes $EAE$ sets defined
over a free group (i.e., sets defined using 3 quantifiers), and shows that
every such set is in the Boolean algebra generated by $AE$ sets
([Se6],41), hence, we obtained a quantifier elimination over a free group.

\noindent
In the 7th paper in the sequence we generalized the techniques and the results from free
groups to torsion-free hyperbolic groups, and in the 8th paper we used the techniques that were
developed for quantifier elimination to prove that the elementary theories of free and torsion-free
hyperbolic groups are stable.
In the 9th paper in the sequence we studied definable equivalence relations over free and hyperbolic groups, proved
that 3 basic families of equivalence relations are imaginaries (non-reals), and finally proved a geometric elimination
of imaginaries when sorts are added for these 3 basic imaginaries.

In this paper we study the first order theory of free products of arbitrary groups. In a joint work with E. Jaligot [Ja-Se] we
started this study, by analyzing the set of solutions to systems of equations over an arbitrary free product. For that purpose,
we used limit groups over free products, and with each system of equations, or alternatively, with each f.p.\ group, we have
associated (non-canonically) a Makanin-Razborov diagram over free products. This Makanin-Razborov diagram encodes the set
of solutions to the finite system of equations over an arbitrary free product, or alternatively describes all the quotients
of a given f.p.\ group that are free products.

We start this paper by studying systems of equations with parameters (section 1). We generalize the notion of a  graded limit group
from free groups to free products, and then define rigid and (weakly) solid limit groups over free products, generalizing the
corresponding notions over free and hyperbolic groups. Unfortunately, the boundedness results that were proved in [Se3] for
the number of rigid and strictly solid families of homomorphisms over free and hyperbolic groups, can no longer be valid over
free products. However, we manage to prove a combinatorial boundedness for rigid and (weakly) strictly solid families (theorems 1.14
and 1.15), that plays an essential role in studying the first order theory of free products, successfully replacing the strong boundedness
results of [Se3]. 

In section 2 we prove a general form of  Merzlyakov theorem (over free groups) on the existence of formal solutions for
sentences and predicates over varieties that are defined over free products. In particular, we show how to associate (non-canonically)
a formal Makanin-Razborov diagram with a given AE sentence or predicate over free products, generalizing the results of [Se2]
over free groups.

In section 3 we start studying sentences and predicates over free products. In section 3 we study AE sentences. The strategy that was
used to study such sentences over free groups in [Se4] can not be applied over free products, hence, we use a modification of
it that uses the tools that were constructed in the first 2 sections and in [Ja-Se]. In section 4 we study AE sets (predicates), and
further apply it to study EAE sets and sentences. In section 5 we study AEAE sets and sentences, and finally, in section 6 we
study general definable sets and sentences over arbitrary free products.

In theorem 6.1, we associate (non-canonically) finitely many graded resolutions with a given coefficient-free 
predicate over free products. This finite collection of resolutions is non-canonical, but it is universal, and it is good for all
non-trivial free products apart from the infinite dihedral group, $D_{\infty}$. In principle, the finite collection of
resolutions enables one to reduce a sentence or a predicate from an ambient free product to its factors. Indeed, in theorem 6.3,
we show that any given coefficient-free sentence over free products is equivalent to a finite disjunction of conjunctions of
(coefficient-free) sentences over the factors of the free product. Furthermore, any given coefficient-free predicate over
free products is equivalent to a coefficient-free predicate in an extended language, that involves finitely many predicates over 
the factors of the free product, and only 3 quantifiers over the ambient free product. Note that since the resolutions that
we associated with a coefficient-free predicate are universal, the reduction of sentences and predicates from the ambient
free product to its factors is uniform, i.e., it is good for all free products, and it does not depend on any particular given one.

The uniform reduction of sentences and predicates, and the resolutions that are associated with a given (coefficient-free) predicate,
enable us to prove some basic results on the first order theory of free products in the next sections. In [Fe-Va] S. Feferman and R. 
Vaught studied the first order properties of certain products of structures. Their methods, that look at the cartesian product of 
given structures, do not cover free products of groups (as they indeed indicated in their paper). 
This and his work with A. Tarski, 
lead R. L.  Vaught to ask the question that we answer in the beginning of section
7:

\vglue 1pc
\proclaim{Theorem 7.1} Let $A_1,B_1,A_2,B_2$ be groups. Suppose that $A_1$ is elementarily equivalent to $A_2$,
and $B_1$ is elementarily equivalent to $B_2$. Then $A_1*B_1$ is elementarily equivalent to $A_2*B_2$.
\endproclaim

The existence of graded resolutions that are associated with a given  sentence over free products enables one to prove the following
theorem, that generalizes Tarski's problem for free groups.
 
\vglue 1pc
\proclaim{Theorem 7.2} Let $A,B$ be non-trivial  groups, and suppose that either $A$ or $B$ is not $Z_2$. Let $F$ be a (possibly cyclic) free group.
Then $A*B$ is elementarily equivalent to $A*B*F$.
\endproclaim

The resolutions that are associated with coefficient-free predicates and sentences over free products, that enable a uniform reduction
from the ambient free product to its factors, allow us to prove other   
uniform properties of sentences over free products.
 
\vglue 1pc
\proclaim{Theorem 7.3} Let $\Phi$ be a coefficient free sentence over groups. There  exists an integer, $k(\Phi)$, so
that for every group, $H$, $\Phi$ is a truth sentence over $H_1*\ldots*H_{k(\Phi)}$, $H_i \simeq H$, if and only if
$\Phi$ is a truth sentence over $H_1*\ldots*H_n$, $H_i \simeq H$, for every $n \geq k(\Phi)$.
\endproclaim

Note that the integer $k(\Phi)$ depends on the coefficient free sentence, $\Phi$, but it does not depend on the group, $H$.
It is easy to see that $k(\Phi)$ can not be chosen to be a universal constant, e.g.,  we can take $\Phi_m$ to be a sentence that 
specifies if the number of conjugacy classes of involutions in the group is at least $m$. For such a sentence, 
$\Phi_m$, $k(\Phi_m)=m$.

Theorem 7.3 can be further strengthened for sequences of groups.
 Let $\Phi$ be a coefficient free sentence over groups. Given any sequence of groups, $G_1,G_2,\ldots$,
we set $M_1=G_1$, $M_2=G_1*G_2$, $M_3=G_1*G_2*G_3$, and so on. The sentence $\Phi$ may be truth or false on any of the groups (free
products) $M_i$, $i=1,\ldots$. Here one can (clearly) not guarantee that the sentence $\Phi$ is constantly truth or 
constantly false 
staring at a bounded index (of the $M_i$'s). However, one can prove the following.

\vglue 1pc
\proclaim{Theorem 7.4} There exists an integer $c(\Phi)$, so that for every sequence of groups,
$G_1,G_2,\ldots$, the sentence $\Phi$ over the sequence of groups, $M_1=G_1,M_2=G_1*G_2,\ldots$ may change signs (from truth to false
or vice versa) at most $c(\Phi)$ times.
\endproclaim

In section 8 we use the resolutions that are associated with a coefficient-free predicate over free products, and combine them with
a modification of the strategy that was applied to prove the stability of free groups in [Se9], to prove that free products of stable
groups is stable.

\vglue 1pc
\proclaim{Theorem 8.1} Let $A$ and $B$ be stable groups. Then $A*B$ is stable.
\endproclaim

This question was brought to our attention by E. Jaligot, and was the main reason for this entire work.
In fact we prove a slightly stronger result, and show that a  free product of a countable collection of groups that are
uniformly stable, is stable.

\vglue 1pc
\proclaim{Theorem 8.7} Let $G_1,G_2,\ldots$ be a sequence of groups. Suppose that every sentence $\Phi$ is uniformly stable 
over the sequence $\{G_i\}$. Then the countable free product, $G_1*G_2*\ldots$, is stable. 
\endproclaim

In section 9, we answer a question that we have learnt from I. Kazachkov [Ca-Ka], and prove that the free product of two 
equationally Noetherian groups is equationally Noetherian (theorem 9.1).

Finally, it is worth noting that our results for free products of groups, or slight strengthenings of them that are still
valid over groups, can be shown to be
false for free products of semigroups, using techniques of Quine [Qu] and Durnev [Du]. e.g., a free product of finite semigroups is
in general unstable (although it is stable if the finite semigroups happen to be groups). Hence, it seems that model theoretic techniques
that handle products of general structures, like the ones that were used by Feferman and Vaught, can not suffice to analyze the
elementary theory 
of free products of groups. 

We would like to thank I. Kazachkov, A. Ould-Houcine, G. Sabbagh, and especially E. Jaligot, 
for presenting to us some of the problems that are 
discussed in this paper. We thank Mladen Bestvina and Mark Feighn for  helpful conversations around Tarski's problem
 that influenced our approach to the problems discussed in this paper, and 
Eliyahu Rips for convincing us to write the paper, and for many discussions on 
further questions it raises.

\vglue 1.5pc
\centerline{\bf{\S1. Graded Limit Groups over Free Products}}
\medskip

In [Ja-Se] we have studied systems of equations over free products. To do that we generalized limit groups over free
groups to limit groups over free products. Limit groups over free products are equipped with
an elliptic structure, i.e., with a finite collection of conjugacy classes of some f.g.\ subgroups, that are
supposed to be mapped into elliptic elements in every homomorphism into a free product under consideration.

With a limit group over free products we have associated a JSJ decomposition over free products (theorem 11 in [Ja-Se]),
that encodes all the abelian decompositions of such limit groups over non-elliptic abelian subgroups. Such a JSJ 
decomposition is non-trivial if the limit group is not (entirely) elliptic, and not abelian nor a 2-orbifold
group. 

We further proved a d.c.c.\ for limit groups over free products (theorem 13  in [Ja-Se]), that  proves that a strictly
decreasing sequence of epimorphisms of limit groups over free products that do not map elliptic elements into
the trivial element, terminate after finitely many steps. This d.c.c.\ which is weaker than the one proved
for limit groups over a free group (theorem 5.1 in [Se1]), is still sufficient for constructing a Makanin-Razborov
diagram for a limit group, or for a system of equations over free products, and such a diagram is the final conclusion
of [Ja-Se].

The Makanin-Razborov diagram over free products encodes all the homomorphism from a given f.p.\ group
into free products, and hence, all the solutions of a given system of equations over arbitrary free products.
However, unlike the construction of Makanin-Razborov diagrams over a free or a hyperbolic group, the construction
of the diagram over free products does not give a canonical diagram, but rather a collection
of (strict) resolutions, that encode the entire set of homomorphisms into free products.

In this section we combine the construction of the Makanin-Razborov diagram over free products, with
the construction of the graded Makanin-Razborov diagram over a free group,  to construct
a graded (relative) diagram over free products.

Let $\Sigma(x,p)=1$ be a system of equations with (a tuple of) variables $x$ and (a tuple of) parameters $p$. 
With $\Sigma(x,p)$ we naturally associate a f.p.\ group, $G(x,p)$, that is generated by copies of the variables $x$
and parameters $p$,
and the relations are words that correspond to the equations of $\Sigma(x,p)$. By theorem 21 in [Ja-Se] with the f.p.\
group $G(x,p)$ we associate (canonically) a finite collection of its maximal limit quotients (over free products),
$(L_1(x,p),E_{L_1}),\ldots,(L_t(x,p),E_{L_t})$. We continue with each of these limit groups in parallel, 
and denote such a limit group
(over free products), $(L(x,p),E_L)$.

With the limit group $(L(x,p),E_L)$ we associate (canonically) a graded (relative) JSJ decomposition over free
products (see theorem 11 in [Ja-Se]), that encodes all the splittings of $(L(x,p),E_L)$ over trivial and
non-elliptic abelian subgroups, in which the 
parameter subgroup $<p>$ is elliptic (i.e., contained in a  vertex group). Like over a free or a hyperbolic group,
and in a similar way to the ungraded case, with the graded JSJ decomposition (over  free products),
we naturally associate a graded modular group (over free products), $GMod(x,p)$, of the limit group 
$(L(x,p),E_L)$.  
Given the graded limit group $(L(x,p),E_L)$ and its graded modular group, $GMod(x,p)$, 
we naturally associate with $(L(x,p),E_L)$ a collection of shortening quotients.

As over free and
hyperbolic groups, and unlike the ungraded
(non-relative) case, it may be that the graded virtually-abelian JSJ decomposition of a graded limit
group (over free products), $(L(x,p),E_L)$, is trivial. By construction, a shortening quotient of a graded
limit group $(L(x,p),E_L)$ is a quotient of $L(x,p)$. In the ungraded case (over free products)
a shortening quotient is always a 
proper quotient, or it is entirely elliptic. As in the free and hyperbolic case a shortening quotient may 
be isomorphic to the 
original limit group $L(x,p)$, even when it is not entirely elliptic.   

\vglue 1pc
\proclaim{Definition 1.1} Let $(L(x,p),E_L)$ be a graded limit group over free products. If the graded (virtually
abelian) JSJ decomposition of $(L(x,p),E_L)$ over free products is trivial, we say that $L(x,p)$ is $rigid$. If
the following 3 conditions hold:
\roster
\item"{(1)}" $(L(x,p),E_L)$ admits no non-trivial 
free decomposition in which the parameter subgroup $<p>$, and the elliptic
subgroups in $(L(x,p),E_L)$, are elliptic. 

\item"{(2)}" the graded JSJ decomposition of $(L(x,p),E_L)$ (over free products)
is non-trivial.  

\item"{(3)}"  
$(L(x,p),E_L)$ has a graded shortening quotient which is isomorphic to $L(x,p)$ (as an abstract group), and this
shortening quotient is not entirely elliptic.
\endroster
We say that $(L(x,p),E_L)$ is $solid$.
\endproclaim

As in the case of free and hyperbolic groups, and unlike the ungraded case, associating shortening quotients
with rigid and solid limit groups (over free products) is not helpful (in order to construct
a diagram that encodes all the homomorphisms of a given graded limit group over  free products). Like in the 
cases of free
and hyperbolic groups, to analyze the collection of homomorphisms that factor through a given rigid or solid
limit group, we need to associate with such a graded limit group a finite collection of $flexible$ $quotients$.

\vglue 1pc
\proclaim{Definition 1.2} Let $(L(x,p),E_L)$ be either a rigid or a solid
graded limit group (over free products). If $L$ is a group of rank at most $d$, then  with $L$ we can associate a sequence of
f.p.\ approximating subgroups $F_d \to G_1 \to G_2 \to \ldots$, so that $G_m$ converges into $L$, and the abelian JSJ decomposition
over free products of $L$ lifts to abelian decompositions of the f.p.\ groups $G_m$. We fix generating sets of the groups $G_m$, and of
the rigid or solid limit group, $(L,E_L)$, that is obtained from a generating set of the free group, $F_d$. 

Given these generating sets, we associate Cayley graphs with each of the f.p.\ groups $G_m$, and with $(L,E_L)$. 
Let $\{h_m: G_m \to A^1_m * \ldots * A^{\ell}_m\}$ be a sequence of homomorphisms 
that converges into a quotient of the limit group over free products,
$(L,E_L)$. With each free product, $A^1_m* \ldots * A^{\ell}_m$, we naturally
associate the pointed Bass-Serre tree, $(T_m,t_m)$,  that is associated with the free product ($(T_m,t_m)$ is dual to a finite tree of groups,
having one vertex with trivial stabilizer, and $\ell$ vertices connected to it with stabilizers, $A^1_m,\ldots,A^{\ell}_m$).
We denote by $d_{T_m}$ the (simplicial)
metric on the tree $T_m$. We say that the sequence of homomorphisms, $\{h_m\}$, is a $flexible$
$sequence$
if one of the following holds:
\roster
\item"{(i)}" each
 homomorphism $h_m$
 can not be shortened (as measured in the trees $T_m$) by an element from the 
 graded modular group (over free products) of the group $G_m$, (that lifts the graded modular group of the 
rigid or solid limit group $L(x,p)$).

\item"{(ii)}"  for each index $m$, let $B_1$ be the ball of radius $1$ in the Cayley graph $X_m$ of the f.p.\ $G_m$. Then:
$$\operatornamewithlimits{\max}_{g \in B_1} \, d_{T_m}(h_m(g)(t_m),t_m) \, > \, 
m \cdot (1+
\operatornamewithlimits{\max}_{1 \leq j \leq u} \, d_{T_m}(h_m(p_j)(t_m),t_m)).$$
where $p_1,\ldots,p_u$ is a fixed generating set of the parameter subgroup $<p>$ in $G_m$, which is the image of this subgroup in $F_d$.
\endroster

A graded limit group (over free products) which is the limit of a flexible sequence is called
a $flexible$ $quotient$ of the rigid or solid graded limit group
$(L(x,p),E_L)$.
\endproclaim

As in the free group case, the following are immediate properties of flexible quotients.  

\vglue 1pc
\proclaim{Lemma 1.3} Let $(Flx(x,p),E_F)$ be a flexible  quotient (over free products)
of the rigid or solid graded limit group (over free products) $(L(x,p),E_L)$. Then:
\roster
\item"{(i)}" $(Flx(x,p),E_F)$ is not a rigid limit group (over free products), and is not entirely elliptic.

\item"{(ii)}" $(Flx(x,p,a),E_F)$ is a proper quotient of $(L(x,p),E_L)$.
\endroster
\endproclaim

\nfp Identical to lemma 10.4 in [Se1].

\line{\hss$\qed$}

\medskip
Let $(L(x,p),E_L)$ be a  rigid or solid limit group over free products. 
As over (ungraded) limit groups over free products, on the set of flexible quotients of $(L(x,p),E_L)$
we can naturally define a partial order and an equivalence relation, similar
to the ones defined on limit groups (over free products) in [Ja-Se]. 
By the same argument that is used to prove proposition 20 in [Ja-Se], the set of flexible 
quotients of $(L(x,p),E_L)$ contains maximal
elements with respect to the partial order. Theorem 1.4 proves that there are at most
finitely many maximal flexible quotients  of a rigid limit group over free products.

\vglue 1pc
\proclaim{Proposition 1.4} Let $(L(x,p),E_L)$ be a  rigid limit group over free products. Then there exist only
finitely many maximal flexible quotients of $(L(x,p),E_L)$ (up to equivalence). 

Furthermore, recall that with the (rigid) limit group over free 
products, $(L,E_L)$, we can naturally associate finitely many elliptic subgroups, $E_1,\ldots,E_r$. Each maximal flexible quotient
of the  rigid limit group, $(L,E_L)$, can be embedded in  an (ungraded)  completion, $Comp$, where the completion, $Comp$, 
is obtained from the elliptic subgroups of the rigid limit group,
$(L,E_L)$, $E_1,\ldots,E_r$, by adding finitely many generators and relations.
\endproclaim

\nfp Identical to the proof of theorem 21 in [Ja-Se].

\line{\hss$\qed$}

Note that in case a rigid limit group, $(L,E_L)$, is finitely presented, theorem 1.4 implies that each maximal flexible
quotient of $(L,E_L)$ can be embedded into a f.p.\ completion (with f.p.\ terminal elliptic subgroups).

Theorem 1.4 proves that rigid limit groups over free products have finitely many maximal flexible quotients. For solid limit groups 
(over free products) we prove a slightly weaker statement. We do not prove the existence of finitely many maximal flexible
quotients of a solid limit group, but instead we prove the existence of finitely many covers of flexible quotients, that cover all
the flexible quotients of the solid limit group (see theorem 24 in [Ja-Se] for a cover of a limit
quotient of a limit group over free products).

\vglue 1pc
\proclaim{Proposition 1.5} Let $(L(x,p),E_L)$ be a  solid limit group over free products. Then there exists a finite
collection of 
covers of maximal flexible quotients of $(L(x,p),E_L)$, so that:
\roster
\item"{(1)}"  each cover is a  proper
quotient of $(L(x,p),E_L)$. 

\item"{(2)}" each cover has a non-trivial graded JSJ decomposition over free products. 

\item"{(3)}" let $E_1,\ldots,E_r$ be the elliptic subgroups of $(L,E_L)$. Then each cover can be embedded in an (ungraded)   completion that
is obtained from the elliptic subgroups of $(L,E_L)$, $E_1,\ldots,E_r$, by adding finitely many generators and relations.

\item"{(4)}" every flexible 
quotient of $(L(x,p),E_L)$ is dominated by at least one of the covers from the finite collection.

\item"{(5)}" let $h_n: L(x,p) \to A^1_n* \ldots * A^{\ell}_n$ be a sequence of homomorphisms 
that converges into a flexible quotient of 
$(L(x,p),E_L)$. Then there exists a subsequence of the homomorphisms $\{h_n\}$ that do all
factor through one of the
covers from the fixed finite collection of covers.
\endroster
\endproclaim

\nfp Identical to the proof of theorem 25 in [Ja-Se].

\line{\hss$\qed$}

As in the case of free groups, having defined flexible quotients over free products, we are
able to define $flexible$ $homomorphisms$, and $rigid$ and $solid$ $homomorphisms$. Since a limit group 
over free products is in general 
f.g.\ and not necessarily finitely presented, for rigid and solid limit groups we define both rigid and solid homomorphisms
that factor through the (rigid or solid) limit group, and $asymptotically$ rigid or solid sequences of homomorphisms that converge to 
limit groups that are dominated by the rigid or solid limit group.

\vglue 1pc
\proclaim{Definition 1.6} Let $(Rgd(x,p),E_R)$ be a  rigid limit group over free products, and let:
$$(Flx_1(x,p),E_{F_1}),\ldots,(Flx_v(x,p),E_{F_v})$$ be the maximal  
flexible quotients of it. A homomorphism
$h: Rgd(x,p) \to A^1* \ldots *A^{\ell}$ that does not factor through any of the
 maximal graded flexible quotients of $(Rgd(x,p),E_R)$ is called
a $rigid$ $homomorphism$  ($specialization$) of the rigid limit group $(Rgd(x,p),E_R)$.
A homomorphism of $(Rgd(x,p),E_R)$ into a free product that does factor through one of the maximal flexible
quotients is called a $flexible$ $homomorphism$ ($specialization$).

Suppose that the rank of the rigid limit group, $(Rgd,E_R)$, is at most $d$. A sequence of homomorphisms, 
$\{h_n: F_d \to A^1_n* \ldots * A^{\ell}_n \}$, 
that converges to a limit quotient $(U,E_U)$ of $(Rgd,E_R)$ is called $asymptotically$ $rigid$ for $(Rgd,E_R)$,
 if the limit group $(U,E_U)$ is not dominated
by any of the flexible quotients of $(Rgd,E_R)$, i.e., if there is no epimorphism $\tau: (Flx,E_F) \to (U,E_U)$, where $(Flx,E_F)$ is one
of the
flexible quotients of $(Rgd,E_R)$. 
A sequence of homomorphisms $\{h_n: F_d \to A^1_n* \ldots * A^{\ell}_n \}$, 
that converges to a limit quotient $(U,E_U)$ of $(Rgd,E_R)$ that is not $asymptotically$ $rigid$ is called
$asymptotically$ $flexible$.
\endproclaim

\vglue 1pc
\proclaim{Definition 1.7} Let $(Sld(x,p),E_S)$ be a  solid  limit group over free products. With 
$(Sld(x,p),E_S)$ we associate a (fixed) finite collection of covers that satisfies the conclusion of
proposition 1.5. 
 A homomorphism
$h: Sld(x,p) \to A^1* \ldots *A^{\ell}$  for which $h=h' \circ \varphi$ where $h'$
factors through one of the covers from the fixed finite collection of covers that is associated with
$(Sld(x,p),E_S)$, and $\varphi$ is a   graded modular automorphism (over free products) of
$(Sld(x,p),E_S)$, is called a $flexible$ $homomorphism$ ($specialization$) of the solid 
limit group $(Sld(x,p),E_S)$ (with respect to the given finite collection of covers). 

\noindent
 A non-flexible homomorphism
$h: Sld(x,p) \to A^1* \ldots * A^{\ell}$ 
is called a $solid$ $homomorphism$ ($specialization$) of the solid limit group (over free products)
$(Sld(x,p),E_S)$ (with respect to the fixed finite collection of covers). 

Suppose that the rank of the solid limit group, $(Sld,E_S)$, is at most $d$. With $(L,E_L)$ we can associate a sequence of f.p.\ groups,
$F_d \to G_1 \to G_2 \ldots$, that converges into $Sld$, so that the graded abelian decomposition of $Sld$ lifts to
abelian decompositions of the f.p.\ groups $\{G_n\}$. A sequence of homomorphisms $\{h_n: G_n \to A^1_n* \ldots * A^{\ell}_n \}$, 
that converges to a limit quotient $(U,E_U)$ of $(Sld,E_S)$ is called $asymptotically$ $solid$ for $(Sld,E_S)$, if there is no
subsequence of homomorphisms (still denoted) $\{h_n\}$, and automorphisms $\{\varphi_n\}$ from the graded modular automorphisms of
the groups $\{g_n\}$, so that the sequence $\{h_n \circ \varphi_n : G_n \to A^1_n*\ldots *A^{\ell}_n\}$ converges into 
a limit quotient of one of the
covers that are associated with $(Sld,E_S)$ according to theorem 1.5.

A sequence of homomorphisms $\{h_n: G_n \to A^1_n* \ldots *A^{\ell}_n \}$, for which there are automorphisms, $\{\varphi_n\}$, so that $h_n \circ \varphi_n$
converges to a limit quotient of one of the covers that are associated with $(Sld,E_S)$ 
is called
$asymptotically$ $flexible$.
\endproclaim

As in the case of free groups, flexible quotients of rigid and solid  limit groups over free products
contain all the "generic infinite families" of
homomorphisms (into free products) of these graded limit groups. Rigid solutions of rigid
 limit groups over free products are the exceptional solutions, and solid solutions of
 solid limit groups
are the exceptional families of solutions.

\vglue 1pc
\proclaim{Proposition 1.8} Let
$(Rgd(x,p),E_R)$ be a  rigid limit group,  and let $(Sld(x,p),E_S)$ be a  solid
  limit group. With $(Sld(x,p),E_S)$ we associate a (fixed) finite collection of covers that satisfies the 
conclusion of proposition 1.5. Let $p_0$ be a specialization of the defining
parameters $p$. Let $x_1,\ldots,x_f$ be a fixed generating set of $Rgd$ or $Sld$,
and let $p_1,\ldots,p_u$ be a fixed generating set for the parameter subgroup $<p>$. 
Then there exists a constant $c(p_0)$ so that:
\roster
\item"{(i)}" let 
$h:Rgd(x,p) \to A^1* \ldots * A^{\ell}$ be a rigid homomorphism for which $h(p)=p_0$. With the free product $A^1* \ldots * A^{\ell}$ we
associate its pointed Bass-Serre tree $(T,t)$. Then:
$$\operatornamewithlimits{\max}_{1 \leq i \leq f} \, d_{T}(h(x_i)(t),t) \, < \, c(p_0) \, \cdot \, 
  \operatornamewithlimits{\max}_{1 \leq j \leq u} \, d_{T}(h(p_j)(t),t)$$ 
Note that the constant $c(p_0)$ depends on the parameter $p_0$, but not on the rigid homomorphism $h$.

\item"{(ii)}" let 
$h:Sld(x,p) \to A^1* \ldots * A^{\ell}$ be a solid homomorphism for which $h(p)=p_0$ with respect to the given finite collection 
of covers that is associated with $(Sld(x,p),E_S)$. 
With the free product $A^1* \ldots * A^{\ell}$ we
associate its pointed Bass-Serre tree $(T,t)$. Then
there exists  a graded modular automorphism (over free products)
of the solid limit group $Sld(x,p)$, $\varphi$, so that $h = h' \circ \varphi$, and: 
$$\operatornamewithlimits{\max}_{1 \leq i \leq f} \, d_{T}(h'(x_i)(t),t) \, < \, c(p_0) \, \cdot \, 
  \operatornamewithlimits{\max}_{1 \leq j \leq u} \, d_{T}(h'(p_j)(t),t)$$ 
\endroster
\endproclaim

\nfp Similar to the proof of proposition 10.7 in [Se1]. 

\line{\hss$\qed$}

\medskip
In studying the first order theory of free and hyperbolic groups, it was essential to strengthen the finiteness
of the number of rigid and solid families of specializations for any given value of the
defining parameters, to a global bound on the number of rigid and strictly solid families (these bounds were
proved in [Se3]).

Over free products the number of rigid or families of solid solutions for a given value of the defining
parameters is not finite in general, hence, we can't expect a strong form of global boundedness. One way
to generalize the boundedness results from free and hyperbolic groups would be to prove that the constant,
$c(p_0)$, that
appears in the formulation of proposition 1.8 can be taken to be uniform, i.e., independent of the specific
value of the parameters $p_0$. This can indeed be done, but it won't suffice to analyze first order predicates
and sentences over free products.

To analyze predicates over free products, we look for a different (stronger) type of generalization of
proposition 1.8. We show that in an appropriate sense, the collection of rigid and solid families of 
homomorphisms, is contained in boundedly many families, where the bound on the number of families does not
depend on the specific value of the defining parameters. Furthermore, these families can be defined using
certain (AE like) predicates, that are crucial in proving a form of quantifier elimination over free products.

Our strategy to define these families and the uniform bound on their number generalizes the argument that
was used in [Se3] to prove a uniform bound on the number of rigid and strictly solid families (over a free group).
We start by proving combinatorial boundedness for rigid and (shortest) solid specializations, which is
similar to the one that was proved in section 1 of [Se3], and then use this combinatorial boundedness
to prove the existence of a uniform bound on the number of rigid and strictly solid families of homomorphisms. 
 
In principle, the combinatorial bound we are looking for is 
showing that there exist finitely many (combinatorial) ways to cut the specializations of the
defining parameters, and finitely many words in these pieces ($fractions$), so that the
specializations of the (given) generating set of a rigid or solid limit group inherited
from any  rigid
or (shortest) solid specialization, can be presented as one of these (finitely many) words in
fractions of the specializations of the defining parameters. 

\noindent
To prove such a combinatorial
bound, we actually prove a stronger result (cf. theorems 1.2 and 1.7 in [Se3]).
To state the stronger statement we need the notion of an
$R$-$P$-covered (graded) homomorphism.
  
Let $X$ be the Cayley graph of a rigid or solid limit group (over free products), $Rgd(x,p)$ or $Sld(x,p)$, with
respect to the generating system $Rgd(x,p)=<x,p>$ (or $Sld(x,p)=<x,p>$). Let $h: Rgd(x,p) \to A^1* \ldots * A^{\ell}$ 
(or $h: Sld(x,p) \to A^1* \ldots * A^{\ell}$)
be a homomorphism,
and and let $(T,t)$ be the  pointed Bass-Serre tree that corresponds to the free product $A^1* \ldots * A^{\ell}$.
Clearly, the  homomorphism $h$ corresponds to a  natural equivariant
map $\tau: X \to T$, where each edge in $X$ is mapped to a (possibly degenerate)
path in $T$.

\vglue 1pc
\proclaim{Definition 1.9} Let $B_R$ be the ball of radius $R$ in the Cayley graph
$X$ of $Rgd(x,p)$ (or $Sld(x,p)$). We say that a homomorphism $h: Rgd(x,p) \to A^1*\ldots *A^{\ell}$ 
(or $h: Sld(x,p) \to A^1* \ldots A^{\ell}$) is
$R$-$P$-$covered$, if the image in $T$ of the union of 
edges labeled by  an
element of $\{p\}$ in $B_R$ covers the entire image in $T$ of the ball $B_1$.
\endproclaim

As over free groups (theorem 1.2 in [Se3]), to control the combinatorial types of rigid solutions over free
products, we need the following
basic
theorem, that bounds combinatorially the structure of rigid solutions (specializations) of
a rigid limit group (over free products).

\vglue 1pc
\proclaim{Theorem 1.10} Let $Rgd(x,p)$ be a  rigid limit group over free products. There exists a
constant $R_0$ so that every rigid homomorphism $h: Rgd(x,p) \to A^1* \ldots *A^{\ell}$ is
$R_0$-$P$-covered.
\endproclaim

\nfp Identical to the proof of theorem 1.2 in [Se3].

\line{\hss$\qed$}

\smallskip
As over free groups, to state a similar theorem for (shortest) solid homomorphisms, we first need to
define $strictly$ $solid$ $homomorphisms$ of a f.p.\ solid limit group (over free products) with respect to a given
finite collection of covers of its flexible quotients (see proposition 1.5). The definition we present is
similar to the definition of strictly solid homomorphisms over free groups that are defined in definition
1.5 in [Se3]. Strictly solid homomorphisms,
and their families, are the ones that needed to be considered in analyzing first order predicates over
free products.

\vglue 1pc
\proclaim{Definition 1.11} Let $(Sld(x,p),E_S)$ be a  solid limit group over free products,
and let: $$(CFlx_1(x,p),E_1),\ldots,(CFlx_s(x,p),E_s)$$ be a finite collection of covers that satisfies the conclusion of 
proposition 1.5.

With each cover, $(CFlx(x,p),E)$, we associate a natural one step resolution given
by the quotient map: $Sld(x,p) \to CFlx(x,p)$, and
the graded abelian JSJ decomposition over free products of the solid limit group $(Sld(x,p),E_S)$, $\Lambda_S$. 
Given the
one step resolution, we associate with it 
the entire collection of
homomorphisms (into free products) that factor through the one step resolution, 
and the limit group corresponding to that
collection that we denote $(Q(x,p),E_Q)$, which is a quotient of the solid limit group $(Sld(x,p),E_S)$.

The limit group $Q(x,p)$ inherits a graded abelian decomposition over free products, $\Delta_Q$, from
the graded abelian decomposition of the solid limit group $(Sld(x,p),E_S)$, 
which is in general either similar to the abelian decomposition 
$\Lambda_S$ (i.e., it has isomorphic $QH$ and abelian vertex and edge groups, and the other (rigid)
vertex groups in $\Delta_Q$ are quotients of the corresponding ones in $\Lambda_S$), or it is a
degeneration of the abelian decomposition $\Lambda_S$, that is associated with $(Sld(x,p),E_S)$. 
By construction, the one step resolution given by the quotient map $Q(x,p) \to CFlx(x,p)$,
and the graded abelian decomposition $\Delta_Q$ of $Q(x,p)$, is a strict resolution. 

With the one step, strict (well-structured) resolution $Q(x,p) \to CFlx(x,p,a)$,
with the abelian decomposition $\Delta_Q$,
we associate its completion, according to the construction that appears in definition 1.12 of [Se2].
We denote the obtained completion, $Comp_{CFlx}$.

By the construction of the completion, $Comp_{CFlx}$, the solid limit group $(Sld(x,p),E_S)$ is naturally mapped
into it (onto the subgroup that is associated with its top level, that is isomorphic to
$(Q(x,p),E_Q)$).
Let $\tau: Sld(x,p) \to Comp_{CFlx}$ be this natural map.

\noindent
A homomorphism, $h:Sld(x,p) \to A^1*\ldots *A^{\ell}$, of the solid limit group over free products, $(Sld(x,p),E_S)$,
is called a $strictly$ $solid$ $homomorphism$ ($specialization$) 
if it is non-degenerate (definition 11.1 in [Se1]), and
for every index $i$, $1 \leq i \leq s$, it does not factor 
through the homomorphisms
$\tau_i:Sld(x,p) \to Comp_{CFlx_i}$. 
Note that a strictly solid specialization is, in particular,
a solid specialization (with respect to the given finite collection of covers).
\endproclaim

Solid limit groups over free products may be infinitely presented. For completeness, and since it will be required in the next sections,
we further associate with a solid limit group its collection of $asymptotic$ $strictly$ $solid$ sequences of homomorphisms. 

\vglue 1pc
\proclaim{Definition 1.12} Let $(Sld(x,p),E_S)$ be a solid limit group over free products,
Suppose that the rank of the solid limit group, $(Sld,E_S)$, is at most $d$. With $(L,E_L)$ we can associate a sequence of f.p.\ groups,
$F_d \to G_1 \to G_2 \ldots$, that converges into $Sld$, so that the graded abelian decomposition of $(Sld,E_S)$ lifts to
abelian decompositions of the f.p.\ groups $\{G_n\}$. With each f.p.\ group $G_n$ we can formally associate the completion,
$Comp_n$, of the identity map, $G_n \to G_n$ (note that $G_n$ is not a limit group over free products, but only its abelian decompositions
(that lifts the one of $(Sld,E_S)$ is needed for defining the completion).  

A sequence of homomorphisms $\{h_n: G_n \to A^1_n* \ldots * A^{\ell}_n \}$, 
that converges to a limit quotient $(U,E_U)$ of $(Sld,E_S)$ is called $asymptotically$ $strictly$ $solid$ for $(Sld,E_S)$, if there is no
subsequence of homomorphisms (still denoted) $\{h_n\}$, and homomorphisms $\{u_n: G_n \to A^1_n*\ldots *A^{\ell}_n\}$, so that for each index $n$
the pair of homomorphisms, $(h_n,u_n)$, extends to a homomorphism, $\tau_n : Comp_n \to A^1_n*\ldots *A^{\ell}_n$, where the restrictions of
$\tau_n$ to the two copies of $G_n$ in $Comp_n$ are the homomorphisms $h_n$ and $u_n$, and so that the sequence $\{u_n\}$
converges into a limit quotient of one of the
covers that are associated with $(Sld,E_S)$ according to theorem 1.5.
\endproclaim

Like over free groups, and in a similar way to the combinatorial boundedness of rigid homomorphisms,
in theorem 1.13 we state the existence of a
   combinatorial bound for strictly solid homomorphisms that are the 
shortest in their solid family.

\vglue 1pc
\proclaim{Theorem 1.13} Let $(Sld(x,p),E_S)$ be a  solid limit group, and let:
$(CFlx_1,E_1),\ldots,(CFlx_s,E_s)$ be a
finite collection of covers of flexible quotients of $(Sld,E_S)$ that satisfies the conclusion of proposition 1.5.
 Then there exists a
constant $R_0$ so that every  strictly solid homomorphism $h: Sld(x,p,a) \to F_k$
which is among the shortest in its solid family,
 is  
$R_0$-$P$-covered. 
\endproclaim

\nfp Identical to the proof of theorem 1.7 in [Se3].

\line{\hss$\qed$}

\medskip
So far we generalized the combinatorial boundedness of rigid and shortest strictly solid homomorphisms
from free and hyperbolic groups to the context of  rigid and solid limit groups over free products, i.e., we
showed that rigid and (shortest) strictly solid homomorphisms over free products of a given f.p.\ rigid or solid
limit group (over free products) are $R_0$-$P$-covered, for some constant $R_0$ that depends only on the 
rigid or solid limit group.

Over free and hyperbolic groups we were able to use this combinatorial boundedness to obtain a global
bound on the number of rigid and families of strictly solid homomorphisms for any given value of the
defining parameters (theorems 2.5 and 2.9 in [Se3]). Over free products such a global bound can not exist.
However, once we define appropriately families of rigid and strictly solid homomorphisms over free
products, it is possible to use the combinatorial bounds to obtain a global (uniform) bound on the
number of rigid and strictly solid families, a bound that suffices for the analysis of first order 
predicates over free products in the next sections.

\vglue 1pc
\proclaim{Theorem 1.14} Let $(Rgd(x,p),E_R)$ be a  rigid limit group over free products, generated by
$x_1,\ldots,x_s$ and so that the parameter subgroup is generated by $p_1,\ldots,p_u$. Let
$(Flx_1,E_1),\ldots,(Flx_t,E_t)$ be the set of maximal flexible quotients of $(Rgd(x,p),E_R)$ 
(see proposition 1.4).

There exist finitely many (combinatorial) systems of $fractions$ of the defining parameters (that depend only on the
rigid limit group $(Rgd(x,p),E_R)$):
$$p_1=v_1 \ldots v_{i_1}, \ p_2=v_{i_1+1} \ldots v_{i_2}, \ \ldots \ ,p_u=v_{i_{u-1}+1} \ldots v_{i_u}$$
$$x_j=x_j(v_1,\ldots,v_{i_u},a_1,\ldots,a_f) \ j=1,\ldots,s.$$
(where the indices $i_1,\ldots,i_u$ and the words $x_1,\ldots,x_s$ 
may depend on the system of fractions), so that:
\roster
\item"{(i)}" with each triple consisting of a fixed free product, $A^1*\ldots*A^{\ell}$, a  value of the defining parameters $p_0$ from the
free product $A^1*\ldots*A^{\ell}$, and
a (combinatorial) system of fractions (one out of the finitely many possible ones), 
there is either a unique or no associated values of the fractions, $v_1,\ldots,v_{i_u}$
(that are all taken from $A^1* \ldots * A^{\ell}$).

\item"{(ii)}" for each rigid homomorphism: $h: Rgd(x,p) \to A^1* \ldots * A^{\ell}$, that satisfies $h(p)=p_0$,
there exist at least one of 
the (finitely many combinatorial) systems of 
fractions, so that the associated values of the fractions, $v_1,\ldots,v_{i_u}$ (that depend only on $p_0$ and the 
combinatorial system of fractions),
together with some  elements, $a_1,\ldots,a_f \in A^1 \cup \ldots \cup A^{\ell}$, satisfy:
$$h(x_j)=x_j(v_1,\ldots,v_{i_u},a_1,\ldots,a_f) \ j=1,\ldots,s.$$
\endroster
\endproclaim

\nfp By theorem 1.10 there exists some integer, $R_0$, so that every rigid homomorphism is $R_0$-$P$-covered. This implies the existence of
finitely many combinatorial systems of fractions that have the properties that are listed in the theorem, i.e., 
that the values of a fixed generating set, $x_1,\ldots,x_s$,  under any given  rigid homomorphism into a free product, are 
specified by  fixed words that are associated with the combinatorial systems, and these words are in terms of fractions of values
of a generating set of the parameter subgroup, $p_1,\ldots,p_u$, and finite number of elements from the factors of the free
product.

To complete the proof of the theorem one is still required to show that for each combinatorial system of fractions, there is a global
bound $b_R$ (that does not depend on the parameters value), so that for each free product, $G=A^1* \ldots *A^{\ell}$, 
and each given value of the parameters,
$p_1,\ldots,p_u$, in $G$, it possible to "cut" the parameters in finitely many ways to obtain the values of (the fractions),
$v_1,\ldots,v_{i_u}$, to obtain all the rigid specializations that are associated with the value of $p_1,\ldots,p_u$, 
and the number of ways we need to "cut" the parameters to obtain the fractions, is bounded by the global bound $b_R$. 

The existence of such a bounded follows from the proof of theorem 2.5 in [Se3]. Note that the argument that is used
to prove theorem 2.5 in [Se3] proves exactly the statement of theorem 1.13 (also in the case of free products), just that in the case
of free (or hyperbolic) groups this statement implies a bound on the 
number of rigid specializations for each value of the defining parameters (theorem 2.5 in [Se3]), and this can not be deduced (in general)
over free products.

\line{\hss$\qed$}

The bound on the number of rigid families of homomorphisms of a rigid limit group over free products
that is stated in theorem 1.14 has an analogue for strictly solid families of homomorphisms of a solid
limit group over free products.

\vglue 1pc
\proclaim{Theorem 1.15} Let $Sld(x,p)$ be a  solid limit group over free products, generated by
$x_1,\ldots,x_s$ and so that the parameter subgroup is generated by $p_1,\ldots,p_u$. Let
$(CFlx_1,E_1),\ldots,(CFlx_t,E_t)$ be a finite collection of covers of flexible quotients of $(Sld(x,p),E_S)$ 
(see proposition 1.5).

There exist finitely many (combinatorial) systems of $fractions$ of the defining parameters (that depend only on the
solid limit group $(Sld(x,p),E_S)$ and the fixed collection of covers):
$$p_1=v_1 \ldots v_{i_1}, \ p_2=v_{i_1+1} \ldots v_{i_2}, \ \ldots \ ,p_u=v_{i_{u-1}+1} \ldots v_{i_u}$$
$$x_j=x_j(v_1,\ldots,v_{i_u},a_1,\ldots,a_f) \ j=1,\ldots,s.$$
(where the indices $i_1,\ldots,i_u$ and the words $x_1,\ldots,x_s$ 
may depend on the combinatorial system of fractions), so that:
\roster
\item"{(i)}" with a triple that consists of a free product, $A^1* \ldots * A^{\ell}$, a  value of the defining parameters $p_0$ from the
free product $A^1* \ldots * A^{\ell}$, and 
a combinatorial system of fractions (one out of finitely many), there is either a unique or no associated values of 
the fractions, $v_1,\ldots,v_{i_u}$
(that are all taken from $A^1*\ldots *A^{\ell}$).

\item"{(ii)}" for each (almost) shortest strictly solid  homomorphism: $h: Sld(x,p) \to A^1* \ldots *A^{\ell}$, 
that satisfies $h(p)=p_0$,
there exist at least one of 
(the finitely many combinatorial) systems of 
fractions, so that the associated values of the fractions, $v_1,\ldots,v_{i_u}$ (that depend only on $p_0$ and the 
combinatorial system of fractions),
together with some elements, $a_1,\ldots,a_f \in A^1 \cup \ldots \cup A^{\ell}$, satisfy:
$$h(x_j)=x_j(v_1,\ldots,v_{i_u},a_1,\ldots,a_f) \ j=1,\ldots,s.$$
\endroster
\endproclaim

\nfp The existence of finitely many combinatorial systems of fractions that are good for all almost shortest strictly solid specializations
follows by theorem 1.11. The existence of a global bound on the number of "cuts" or fractions that are needed in order to cover all the
almost shortest specializations follows from the proof of theorem 2.9 in [Se3].

\line{\hss$\qed$}

Graded limit groups, their graded (relative) abelian JSJ  decompositions (over free products),
their graded shortening quotients, and rigid and solid limit groups and their flexible
quotients, allow us to associate a graded Makanin-Razborov diagram over free products with any given f.p.\ group
$G(x,p)$. As in the ungraded case, the diagram is not canonical, but as over free groups, it encodes 
all the homomorphisms of the group $G(x,p)$ into free products for all the possible values of the
parameters $p$.   

To construct the graded Makanin-Razborov diagram over free products, we need two basic objects, that were 
constructed and used in the ungraded case. First given a f.p.\ group $G(x,p)$ we need to associate with
$G(x,p)$ (or a quotient of it) a collection of (well-structured) graded resolutions. 
Then given a 
well-structured graded resolution we
need to associate with it a cover that satisfies similar properties to the cover of an ungraded resolution (see
theorem 24 in [Ja-Se]). The entire collection of covers of resolutions of $G(x,p)$ should encode all
the homomorphisms of $G(x,p)$ into free products, for every possible value of the defining parameters.
Finally, to construct the graded Makanin-Razborov diagram we will show that there exist
a finite collection of covers of graded resolutions through which all the homomorphisms from the given graded 
limit group into free products do factor.

\vglue 1pc
\proclaim{Proposition 1.16} Let $G(x,p)$ be a f.g.\ group, and let 
$L(x,p)$ be a limit quotient of $G(x,p)$ (over free products).
Then there exists a finite sequence of limit groups over free products:
$$L=L_0(x,p)
  \to L_1(x,p) \to L_2(x,p) \to \ldots \to L_s(x,p)$$
for which:
\roster
\item"{(i)}" non-trivial elliptic elements in $L_i$ are mapped  to non-trivial elliptic elements in $L_{i+1}$,
$i=1,\ldots,s-1$.

\item"{(ii)}" if $L_i$
is not rigid nor solid, then $L_{i+1}$ is a shortening
quotient of $L_i$, for $i=1,\ldots,s-1$. 

\item"{(iii)}" 
 in case $L_i$ is a free product of a  solid limit group with (possibly)
elliptic factors and (possibly) a free group, for some $i<s$, then $L_{i+1}$ is a free product of a  quotient of one of the covers that are
associated with the solid factor (according to proposition 1.5) 
with (possibly) the same elliptic factors and (possibly) the same free group, for
$i=1,\ldots,s-1$.

\item"{(iv)}" The epimorphisms along the sequence are proper epimorphisms.

\item"{(v)}" $L_s$ is a free product of a rigid or a solid limit group with (possibly) some elliptic
factors and (possibly) a free group. With the rigid limit group (which is a factor in $L_s$)
we associate its finite
collection of maximal flexible quotients, and with a solid limit group (which is a factor of $L_s$), we
associate a finite collection of covers of its flexible quotients that satisfies the conclusion of
proposition 1.5.

\item"{(vi)}" the resolution:
$L_0  \to L_1 \to L_2 \to \ldots \to L_s$ is a graded strict resolution ([Se1],5), i.e., in each level
non-QH, non-virtually-abelian
vertex groups and edge groups 
in the graded  JSJ decomposition (over free products) are mapped monomorphically into the limit group in the next
level, and QH vertex groups are mapped into non-virtually-abelian, non-elliptic subgroups).

\item"{(vii)}" the 
constructed resolution is well-structured (see definition 1.11 in [Se2] for a well-structured resolution). 
As a corollary, the graded limit group (over free products) $L=L_0$ is embedded into the completion 
of the well-structured
resolution:
$$L_0  \to L_1 \to L_2 \to \ldots \to L_s$$ 
so that all the elliptic elements in $L$ are mapped into conjugates of the elliptic subgroups of the completion.
\endroster
\endproclaim

\nfp Given a sequence of homomorphisms $\{h_n\}$ of a given group $G$ into free products, $A^1_n*\ldots *A^{\ell}_n$, that converges to a
limit group $(L,E_L)$ over free products, it is shown in theorem 18 of [Ja-Se] how to associate a strict resolution with a 
subsequence of the sequence of homomorphisms. with the notions of graded limit groups over free products, their graded
abelian decompositions, rigid and solid graded limit groups over free products, their flexible quotients, and asymptotically rigid and 
strictly solid sequences, the construction of an (ungraded) resolution over free products generalizes naturally to the construction
of a graded resolution.

\line{\hss$\qed$}

Proposition 1.16 associates graded resolutions with a f.g.\ group. In the ungraded case, in order to associate a
(non-canonical) Makanin-Razborov diagram over free products with a given f.p.\ group, we have replaced each resolution of some  
limit quotient (over free products) of the given f.p.\ group with a $cover$ of that limit quotient, and  the entire resolution with 
a cover resolution, so that each cover resolution has a f.p.\ completion (and terminates in a f.p.\ group). This existence of
cover resolutions enables one to show that it is sufficient to choose finitely many cover resolutions, so that every homomorphism of
the original (given) f.p.\ group into free products factors through at least one of the finite (non-canonical) collection of cover
resolutions (see theorems 24-27 in [Ja-Se]).

Our approach to constructing a graded Makanin-Razborov diagram of a given f.p.\ group $G(x,p)$ is conceptually similar. However,
we will need to modify the notion of a cover resolution to the graded case. As our goal is to obtain a finite collection of
(graded) covers so that for each value of the defining parameters the homomorphisms of $G(x,p)$ into free factors will factor
through the given finite collection, we need to define cover resolutions in a finite way, i.e., graded covers should be 
embedded into f.p.\ completions (over free products). Also, we need to define the covers of flexible quotients that are
attached to (terminal) rigid and solid limit groups in such a way that they can be embedded into f.p.\ completions. 
To achieve these goals we need to slightly modify the definitions and the constructions that we used, so
that covers can be defined in a "finite" way, and all our previous results, especially those for rigid and solid limit
groups remain valid.

\vglue 1pc
\proclaim{Theorem 1.17} Let $(Rgd(x,p),E_R)$ be a rigid limit group over free products. Suppose that $(Rgd(x,p),E_R)$ is
a (rigid) limit quotient of a f.p.\ group, $G(x,p)$, 
and let $\{h_n: G(x,p) \to A^1_n*\ldots*A^{\ell}_n \}$ be an asymptotically rigid sequence of homomorphisms (of $Rgd$) that converges into $(Rgd(x,p),E_R)$. 

There exists an approximating rigid limit group of ($Rgd(x,p),E_R)$, that we denote $(APRgd,E_{APR})$, with the following
properties:
\roster
\item"{(1)}" $(Rgd,E_R)$ is a quotient limit of $(APRgd,E_{APR})$. $APRgd$ is a rigid limit quotient of $G(x,p)$ (over free products).

\item"{(2)}" there is an ungraded  resolution (over free products) $WRes$, $APRgd \to W_1 \to \ldots \to W_s$, so that $W_s$
is f.p.\ and 
 there exists a subsequence of homomorphisms, (still denoted) $\{h_n\}$, that factor through this resolution. In particular,
the completion of the resolution $WRes$ is f.p.\ and the approximating limit group, $(APRgd,E_{APR})$, is embedded into this
f.p.\ completion.

\item"{(3)}" with $(APRgd,E_{APR})$ we associate finitely many  limit quotients (over free products) of $G(x,p)$, $CF_1,\ldots,CF_g$, 
with ungraded resolutions (over free products), $Res_1,\ldots,Res_g$, that terminate in f.p.\ limit groups. Hence,
the completions of $Res_1,\ldots,Res_g$ are f.p.\ and the limit groups (over free products) , $CF_1,\ldots,CF_g$, 
can be embedded in
f.p.\ completions  
(that have  f.p.\ terminal limit groups). 

\item"{(4)}" the finite collection of limit groups, $CF_1,\ldots,CF_g$, are all limit quotients of $G(x,p)$, and they
dominate all the flexible quotients of
$(APRgd,E_{APR})$ (although they need not be quotients of the rigid limit group $(APRgd,E_{APR})$). Furthermore, every flexible
homomorphism of the rigid approximation, $(APRgd,E_{APR})$, factors through at least one of the resolutions, $Res_1,\ldots,Res_g$.  

\item"{(5)}" from the subsequence of homomorphisms, $\{h_n\}$, that factor through the resolution, $WRes$, it is possible to extract
a subsequence that does not factor through any of the 
limit groups over free products, $CF_1,\ldots,CF_g$, and in particular they do not factor through the resolutions, $Res_1,\ldots,Res_g$.
This implies that  the homomorphisms, $\{h_n\}$, are rigid homomorphisms of the rigid limit group (over free products)
$(APRgd,E_{APR})$.
\endroster
\endproclaim

\nfp Let $(Rgd,E_R)$ be a rigid limit group over free products, and let
$\{h_n: G(x,p) \to A^1_n* \ldots *A^{\ell}_n\}$ be a sequence of asymptotically rigid homomorphisms that that converges into the rigid
limit group over free products, $(Rgd,E_R)$.

By theorem 18 in [Ja-Se], starting with the sequence of homomorphisms, $\{h_n\}$, we can pass to a subsequence, and
associate with the rigid
limit group over free products, $(Rgd,E_R)$,  
a strict well-structured ungraded resolution, that we denote $Res$, $Rgd=W_0 \to W_1 \to \ldots \to W_t$, so that
$W_t$ is the free product of (possibly) finitely many elliptic factors and (possibly) a free group.

We set $Comp_W$ to be the completion of this ungraded  resolution. Note that $Comp_W$ is generated by the terminal limit group $W_t$,
in addition to finitely many elements and relations. However, $W_t$ may be not  finitely presented.

\noindent
We modify $Comp_W$ by replacing $W_t$ with a sequence of f.p.\ approximations of it.
Starting with $Comp_W$ we define a sequence of completions, obtained by replacing the terminal limit group 
(over free products) in $Comp_W$, $W_t$, with approximating f.p.\ groups, that we denote $W_t^m$. These groups, $W_t^m$, are obtained by replacing
each of the elliptic factors of $W_t$, with a f.p.\ group that is generated by some fixed generating set of the corresponding elliptic
factor of $W_t$, and imposing only the relations up to length $m$ in this (fixed) generating set that are equal to the identity in
the corresponding factor of $W_t$. For every index $m$, $W_t$ is a quotient of $W_t^m$, and for $m$ large enough, the group obtained
by replacing $W_t$ with $W_t^m$ in $Comp_W$, is a completion of a strict well-structured resolution (over free products), and we denote the obtained
completion, $Comp_W^m$. 

Note that the approximating completions, $Comp_W^m$, are all finitely presented. With $Comp_W^m$, we naturally associate an
ungraded resolution over free products, that we denote $Res^m$, $W_0^m \to W_1^m \to \ldots \to W_t^m$. The limit groups over
free products, $W_0^m, \ldots , W_t^m$, are all embedded in the f.p.\ completion, $Comp_W^m$. The limit groups (over free products), 
$W_0, \ldots , W_t$, are quotients of the  limit groups over free products,
$W_0^m, \ldots , W_t^m$, in correspondence, and the direct limit of the limit groups over free products, $W_0^m,\ldots,W_t^m$, are
$W_0,\ldots,W_t$, in correspondence.

By construction, with each of the graded limit groups, $W_i$, there is an associated virtually abelian decomposition,
which is the JSJ decomposition (over free products) of $W_i$. 
By the construction of $Comp_W^m$, the JSJ decompositions of the  limit groups (over free products), $W_i$, lift to 
virtually abelian decompositions of the (approximating) graded limit groups over free products, $W_i^m$. 

By the proof of theorem 27 in [Ja-Se], for large enough index $m$, these  virtually abelian decompositions of
the limit groups $W_i^m$ (that are lifted from $W_i$) 
are their JSJ decompositions over free products. Furthermore, by the same argument (that appears in detail in the proof of theorem
27 in [Ja-Se]), for large enough $m$, the graded JSJ decomposition of $W_0^m$ is trivial, which means that for large $m$,
$W_0^m$ is rigid. 

So far we have shown that for large index $m$, the resolutions $Res^m$ and their initial rigid limit group $W_0^m$,
satisfy properties (1) and (2) of the theorem. To get the other parts, we first need to show that for large index
$m$, there exists some index $n_m$, so that for every $n>n_m$, the homomorphism $h_n$ is a rigid homomorphism of
the rigid limit group, $W_0^m$.

The sequence of homomorphisms, $\{h_n\}$, is an asymptotically rigid sequence of homomorphisms with respect to the
original rigid limit group, $(Rgd,E_R)$, and it converges to $(Rgd,E_R)$. Since the completions, $Comp_W^m$, are finitely
presented, for each index $m$, there exists some index, $a_m$, so that for every $n>a_m$, the sequence of
of homomorphisms, $\{h_n\}$, factor through $Comp_W^m$, and in particular through the rigid limit group, $W_0^m$.

Suppose that there exists a subsequence of indices, still denoted $m$, so that for every index $m$, there is a subsequence
of homomorphisms (where the subsequence may depend on $m$), still denoted, $\{h_n\}$, 
that are not rigid homomorphisms of $W_0^m$, i.e., they factor through a flexible
quotient of $W_0^m$. In this case, the sequence of homomorphisms, $\{h_n\}$,  converges into a flexible quotient of $(Rgd,E_R)$,
a contradiction to $(Rgd,E_R)$ being rigid. Therefore, for large enough index $m$, there exists an index $n_m$, so that the
for every $n>n_m$, the  homomorphisms $\{h_n\}$ that do all factor through $W_0^m$, are rigid homomorphisms.

Furthermore,  
for large enough index $m$, $W_0^m$ is rigid, hence, we can associate with $W_0^m$ its finite collection of maximal flexible
quotients, 
$(Flx_1^m,E^m_{F_1}),\ldots,(Flx^m_v,E^m_{F_{v(m)}})$
 (such a collection of maximal flexible quotients exists by proposition 1.4). If there exists a subsequence of the indices $m$, so
that the original rigid limit group, $(Rgd,E_R)$, which is a limit quotient of $W_0^m$, is a limit  quotient of at least one
of the maximal flexible quotients of $W_0^m$, then $(Rgd,E_R)$, can not be a rigid limit group, a contradiction. Therefore,
for $m$ large enough, $(Rgd,E_R)$ is not a  limit quotient of any of the maximal flexible quotients of $W_0^m$.

We pick such a large index $m$. Since $(Rgd,E_R)$ is not a limit quotient of any of the maximal flexible quotients:
$(Flx_1^m,E^m_{F_1}),\ldots,(Flx^m_v,E^m_{F_{v(m)}})$, and since $G(x,p)$ maps epimorphically onto its maximal flexible
quotients and onto $(Rgd,E_R)$,  there exist elements $r_1,\ldots,r_{v(m)} \in G(x,p)$, so that for each $j$, $1 \leq j \leq v(m)$,
$r_j$ is either mapped into non-trivial element in $(Rgd,E_R)$ and to the identity element in $Flx^m_j$, or it is mapped into a 
non-elliptic element in $(Rgd,E_R)$, and to an elliptic element in $(Flx^m_j,E_{F^m_j})$.

With each of the maximal flexible quotients, $(Flx^m_j,E^m_{F_j})$, we associate its (ungraded) Makanin-Razborov diagram according to
theorem 27 in [Ja-Se]. Given one of these flexible quotients, $(Flx^m_j,E^m_{F_j})$, and a resolution that appears in its (strict) ungraded
Makanin-Razborov diagram, we take the completion of this resolution, into which $G(x,p)$ is mapped. We further take a f.p.\ 
approximation of this completion into which $G(x,p)$ is mapped as well. Clearly,
we can choose the f.p.\ approximation of the completion, so that the element $r_j \in G(x,p)$ is mapped either to a trivial or to an elliptic
element, depending whether
it is mapped to a trivial or to an elliptic element in $(Flx^m_j,E^m_{F_j})$. Finally we set the groups, $(CF_1,\ldots,CF_g)$ to be the
images of $G(x,p)$ in the corresponding f.p.\ completions of the resolution in the (strict) ungraded Makanin-Razborov diagrams of
the maximal flexible quotients: 
$(Flx_1^m,E^m_{F_1}),\ldots,(Flx^m_v,E^m_{F_{v(m)}})$. By construction, this finite collection of limit groups with f.p.\ completions
(and corresponding strict well-structured resolution), together with
the rigid limit group, $W_0^m$, and its f.p.\ completion, satisfy properties (1)-(5) of the theorem.

\line{\hss$\qed$}

Theorem 1.17 associate (non-canonically) a rigid limit group that can be embedded into a f.p.\ completion with any given rigid limit
group over free products. To construct a graded Makanin-Razborov diagram over free products we further need to associate a limit
group that can be embedded into a f.p.\ completion with any given solid limit group. In order to associate similar limit groups
with a solid limit group over free products, 
we need to weaken the notion of a solid limit group, to a $weakly$ $solid$ limit group.

\vglue 1pc
\proclaim{Definition 1.18} Let $(L(x,p),E_L)$ be a graded limit group over free products. Suppose that $(L,E_L)$ does not admit a non-trivial free
decomposition in which the parameter subgroup, $<p>$, and  all the elliptic subgroups of $L(x,p)$ are elliptic (i.e., 
elliptic subgroups are contained in conjugates of
the factors of the free decomposition). Let $\Lambda$ be the graded virtually abelian JSJ decomposition of $(L,E_L)$ over free products, 
and suppose that it is
non-trivial. 

Let $GMod(L,E_L)$ be the graded modular group of $(L,E_L)$ over free products. Even though $(L,E_L)$ need not be solid, with the graded limit
group, $(L,E_L)$, 
we can associate its collection of flexible quotients according to definition 1.2. Hence, following the proof of proposition 1.5,
with $(L,E_L)$ we can associate a finite collection of covers that dominate all the flexible quotients of $(L,E_L)$, and satisfy 
the properties of covers (that are associated with a solid limit group) that are listed in proposition 1.5.

We say that a graded limit group over free products, $(L(x,p),E_L)$, that does not admit a graded free decomposition over free
products, and that has a non-trivial graded JSJ decomposition, is $weakly$ $solid$, if there exists
a finite collection of covers of the flexible
quotients of $L$, $(CF_1,\ldots,CF_b)$, so that there are homomorphisms of $(L,E_L)$ into free products that are $weakly$ $strictly$ $solid$,
i.e.,
they can not be extended to homomorphisms of any of the completions of the one step resolutions, $L \to CF_i$ (cf. definition 1.11).
\endproclaim

\vglue 1pc
\proclaim{Proposition 1.19} Let $(L(x,p),E_L)$ be a weakly solid limit group, and let $(CF_1,\ldots,CF_b)$, be a finite collection of
covers that dominate all the flexible quotients of $(L,E_L)$. Then there exists a constant, $R_0$, so that every (almost) shortest weakly
strictly solid homomorphism of $(L,E_L)$ into free products is $R_0$-$P$-covered. Furthermore, the conclusion of theorem
1.15 (for shortest strictly solid homomorphisms of a solid limit group  over free products) hold for weakly strictly solid
homomorphisms of $(L,E_L)$.
\endproclaim

\nfp The same arguments that are used to prove  proposition 1.8, and theorems 1.13  and 1.15 for strictly solid homomorphisms of 
a solid limit group (over
free products) remain valid for weakly strictly solid homomorphisms of a weakly solid limit group.

\line{\hss$\qed$}

The notion of weakly solid limit group, and the ability to generalize the bounds on strictly solid homomorphisms of solid
limit groups to weakly strictly solid homomorphisms of weakly solid limit groups, enable us to state an analogue of theorem 1.16 for 
approximations of solid
limit groups over free products. Theorems 1.17 and 1.20 finally allow us to generalize the construction of the Makanin-Razborov diagram
over free products from the ungraded to the graded case. 

\vglue 1pc
\proclaim{Theorem 1.20} Let $(Sld(x,p),E_S)$ be a solid limit group over free products. Suppose that $(Sld(x,p),E_S)$ is
a (solid) limit quotient of a f.p.\ group, $G(x,p)$, 
and let $\{h_n: G(x,p) \to A^1_n*\ldots *A^{\ell}_n \}$ be an asymptotically strictly solid  sequence of homomorphisms (of $Sld$) that converges into $(Sld(x,p),E_S)$. 

There exists an approximating weakly solid limit group of ($Sld(x,p),E_S)$, (see definition 1.18 for weakly solid) 
that we denote $(WSld,E_{WS})$, with the following
properties:
\roster
\item"{(1)}" $(Sld,E_S)$ is a quotient limit of $(WSld,E_{WS})$. $WSld$ is a weakly solid  limit quotient of $G(x,p)$ (over free products).

\item"{(2)}" there is an ungraded  resolution (over free products) $Res$, $WSld \to V_1 \to \ldots \to V_r$, so that $V_r$
is f.p.\ and 
 there exists a subsequence of homomorphisms, (still denoted) $\{h_n\}$, that factor through this resolution. In particular,
the completion of the resolution $Res$ is f.p.\ and the approximating weakly solid limit group, $(WSld,E_{WS})$, is embedded into this
completion.

\item"{(3)}" with $(WSld,E_{WS})$ we associate finitely many  limit groups (over free products), $CF_1,\ldots,CF_g$, 
with ungraded resolutions (over free products), $Res_1,\ldots,Res_g$, that terminate in f.p.\ limit groups. Hence,
the completions of $Res_1,\ldots,Res_g$ are f.p.\ and the limit groups (over free products) , $CF_1,\ldots,CF_g$, 
can be embedded in
f.p.\ completions  
(that have  f.p.\ terminal limit groups). 

\item"{(4)}" the finite collection of limit groups, $CF_1,\ldots,CF_g$, are all limit quotients of $G(x,p)$, and they
dominate (as quotients of the f.p.\ group $G(x,p)$) all the flexible quotients of
$(WSld,E_{WS})$ (although they need not be quotients of the weakly solid limit group $(WSld,E_{WS})$). Furthermore, for each
homomorphism $f: WSld \to A^1*\ldots*A^{\ell}$ that is not weakly strictly solid, there exists a homomorphism, $u: WSld \to A^1* \ldots *A^{\ell}$, such that 
the pair, $(f,u)$, extends to a homomorphism of the completion of the identity resolution, $WSld \to WSld$, and the homomorphism
$u$ as a homomorphism of the f.p.\ group, $G(x,p)$, 
factors through at least one of the resolutions, $Res_1,\ldots,Res_g$.  

\item"{(5)}" from the sequence of homomorphisms, $\{h_n\}$, it is possible to extract a subsequence 
of weakly solid homomorphisms of
the approximating solid limit group, $(WSld,E_{WS})$. Furthermore, for each homomorphism $h_n$ in this subsequence, there does
not exist a homomorphism, $u_n : WSld \to A^1_n*\ldots*A^{\ell}_n$, so that the pair of homomorphisms, $(h_n,u_n)$, extends to a homomorphism
of the completion of the one step identity resolution, $WSld \to WSld$, and the homomorphism $u_n$, as a homomorphism from the f.p.\
group $G(x,p)$
(into $A^1_n*\ldots*A^{\ell}_n$), factors through one of the limit groups: $CF_1,\ldots,CF_g$.
\endroster
\endproclaim

\nfp The argument that we use follows the proof of theorem 1.17.
Let $(Sld,E_S)$ be a solid limit group, and  let
$\{h_n: G(x,p) \to A^1_n*\ldots*A^{\ell}_n\}$ be a sequence of asymptotically strictly solid homomorphisms  that converges into the solid
limit group over free products, $(Sld,E_S)$.

Starting with the sequence of homomorphisms, $\{h_n\}$, we can pass to a subsequence, and
associate with the solid
limit group over free products, $(Sld,E_S)$,  
a strict well-structured ungraded resolution, that we denote $Res$, $Sld=W_0 \to W_1 \to \ldots \to W_t$, so that
$W_t$ is the free product of (possibly) finitely many elliptic factors and (possibly) a free group.

We set $Comp_W$ to be the completion of this ungraded  resolution. $Comp_W$ is generated by the terminal limit group $W_t$,
in addition to finitely many elements and relations. However, $W_t$ may be not  finitely presented. 
We modify $Comp_W$, by replacing $W_t$ with a sequence of f.p.\ approximations of it precisely as we did in the rigid case, i.e., 
for each index $m$, we replace
each factor of $W_t$ by a f.p.\ group, so that the presentation of each (new) factor includes only relations up to length $m$ from the
relations of the corresponding factor of $W_t$. We denote the obtained f.p.\ group $W_t^m$, and the completion that is obtained from $Comp_W$,
by replacing $W_t$ with $W_t^m$, $Comp_W^m$.

For each $m$, $Comp_W^m$, is finitely presented. With $Comp_W^m$, we naturally associate an
ungraded resolution over free products, that we denote $Res^m$, $W_0^m \to W_1^m \to \ldots \to W_t^m$. The limit groups over
free products, $W_0^m, \ldots , W_t^m$, are all embedded in the f.p.\ completion, $Comp_W^m$. The limit groups (over free products), 
$W_0, \ldots , W_t$, are quotients of the  limit groups over free products,
$W_0^m, \ldots , W_t^m$, in correspondence, and the direct limit of the limit groups over free products, $W_0^m,\ldots,W_t^m$, are
$W_0,\ldots,W_t$, in correspondence.

By construction, with each of the graded limit groups, $W_i$, there is an associated virtually abelian decomposition,
which is the JSJ decomposition (over free products) of $W_i$. 
By the construction of $Comp_W^m$, the JSJ decompositions of the  limit groups (over free products), $W_i$, lift to 
virtually abelian decompositions of the (approximating) graded limit groups over free products, $W_i^m$. 
By the proof of theorem 27 in [Ja-Se], for large enough index $m$, these  virtually abelian decompositions of
the limit groups $W_i^m$ (that are lifted from $W_i$) 
are their JSJ decompositions over free products. 

The sequence of homomorphisms, $\{h_n\}$, is an asymptotically strictly solid sequence of homomorphisms with respect to the
original rigid limit group, $(Sld,E_S)$, and it converges to $(Sld,E_S)$. Since the completions, $Comp_W^m$, are finitely
presented, for each index $m$, there exists some index, $a_m$, so that for every $n>a_m$, the sequence of
of homomorphisms, $\{h_n\}$, factor through $Comp_W^m$, and in particular through the limit group, $W_0^m$.

Suppose that there exists a subsequence of indices, still denoted $m$, so that for every index $m$, there is a subsequence
of homomorphisms (where the subsequence may depend on $m$), still denoted, $\{h_n\}$, 
that are not weakly strictly solid  homomorphisms of $W_0^m$, i.e., they can be extended to a completion of a one step 
resolution: $W^0_m \to CFlx_m$, where $CFlx_m$ is a (cover of a) flexible quotient of the graded limit group $W_0^m$ (note that $W^0_m$
need not be solid, still one can define a flexible quotient of it as it appears in definition 1.18). 
In this case, the original sequence of homomorphisms, $\{h_n\}$, is not an asymptotically strictly solid sequence of
homomorphisms with respect to the solid limit group, $(Sld,E_S)$, a contradiction. Therefore, for large enough $m$, the approximating
limit groups, $W_0^m$ are weakly solid, and for every index $m$, there exists an index $n_m$, so that for every $n>n_m$, the
homomorphisms, $\{h_n\}$, are weakly strictly solid homomorphisms of $W_0^m$.

For each index $m$, we look at some finite collection of covers of its flexible quotients. We further look at a strict (ungraded) 
Makanin-Razborov resolution of each of the (ungraded) resolutions in this diagram, and with each completion of a resolution that
 appears in the diagram we associate a f.p.\ completion, that is obtained by replacing the terminal limit
group of such an ungraded resolution with some f.p.\ approximation of it.  By looking at the image of the f.p.\ group $G(x,p)$ in each of these
(finitely many) f.p.\ completions, we obtained finitely many limit groups (covers), that we denote $CF^m_1,\ldots,CF^m_{v(m)}$, that
are quotients of $G(x,p)$,
and can be embedded into f.p.\ completions, and cover all the flexible quotients of the weakly solid limit groups (over free products) $W_0^m$.

The limit groups, $\{W_0^m\}$, do converge into the original solid limit group, $(Sld,E_S)$. As the flexible quotients of $W_0^m$ are
(by definition) quotients of $W_0^m$, by choosing the approximations properly, i.e., choosing enough of the relations of the 
terminal limit group of each completion, to obtain the covers (that are not necessarily quotients of $W_0^m$),
$CF^m_1,\ldots,CF^m_{v(m)}$, we can guarantee that any sequence of homomorphisms, $s_m : CF^m_{i(m)} \to A^1_m*\ldots*A^{\ell}_m$, that converges to a limit
group, converges into a quotient of $(Sld,E_S)$.

Suppose that there exists a subsequence of indices, still denoted $m$, so that for every index $m$, there is a subsequence
of homomorphisms (where the subsequence may depend on $m$), still denoted, $\{h_n\}$, so that for every $h_n$ there exists a
homomorphism, $u_n: W_0^m \to  A^1_n*\ldots*A^{\ell}_n$, that as a homomorphism of the f.p.\ group, $G(x,p)$,  factors through one of the finitely many covers, 
$CF^m_{i(m)}$, and the pair of homomorphisms of $WSld$, $(h_n,u_n)$, extend to a homomorphism of the completion of the identity
resolution, $W_0^m \to W_0^m$  (with the graded  JSJ decomposition of WSld as the associated virtually abelian decomposition). 
In this case, it is possible to extract a subsequence of the homomorphisms, $\{h_n\}$, 
with homomorphisms $\{u_n: W_0^m \to A^1_n*\ldots*A^{\ell}_n\}$, so that the
pairs, $(h_n,u_n)$, extend to homomorphisms of the identity resolution, $W_0^m \to W_0^m$, and the homomorphisms, $\{u_n\}$, converge
into a flexible quotient of $(Sld,E_S)$, a contradiction to our assumption that the sequence, $\{h_n\}$, is weakly strictly solid. 

Therefore, for large enough $m$, the weakly solid limit group, $W_0^m$, that can be embedded into a f.p.\ completion, and the covers,
$CF^m_1,\ldots,CF^m_{v(m)}$, that can be embedded into f.p.\ completions as well, satisfy properties (1)-(5) of the theorem.

\line{\hss$\qed$}

Given a sequence of homomorphisms from a fixed graded group into free products, theorem 1.16 enables us to associate a graded resolution with
a subsequence of the homomorphisms. 
Theorems 1.17 and 1.20 finally allow us to generalize the notion of a cover of an ungraded resolution (theorem 24 in [Ja-Se]), to define
a cover of a resolution, so that, in particular, given a sequence of homomorphisms of a f.p.\ group into free products, there is a subsequence 
of a  sequence of homomorphisms that factor through  a cover of the resolution that was constructed in proposition 1.16.

\vglue 1pc
\proclaim{Theorem 1.21} Let $G(x,p)$ be a f.p.\ group, and let $\{h_n:G(x,p) \to A^1_n* \ldots *A^{\ell}_n\}$ be a sequence of homomorphisms
of $G(x,p)$ into free products that converges into a limit group (over free products), $(L(x,p),E_L)$.
Let $L(x,p)= L_0 \to L_1 \to \ldots L_s$ be a graded resolution, that we denote $GRes$, of the limit group $(L(x,p),E_{L})$, that is obtained from a
subsequence of the sequence of homomorphisms, $\{h_n\}$, according to proposition 1.16.

There exists a f.g.\ graded limit quotient of $G(x,p)$, $CM$, with a set of elliptics, $E_{CM}$, and
a well-structured graded resolution of $CM$, $CGRes$, $CM=CM_0 \to CM_1 \to \ldots \to CM_s$, 
that satisfies similar properties to the ones listed in proposition 1.16. The graded resolution, $CGRes$, satisfies the 
following properties:

\roster
\item"{(1)}" for each index $i$, $i=1,\ldots,s$, there is an epimorphism of limit groups over free products, $\tau_i: CM_i \to L_i$. The epimorphisms
$\tau_i$ 
commute with the quotient maps in the two graded resolutions, $GRes$ and $CGRes$. 

\item"{(2)}" non-trivial elliptic elements in $CM_i$ are mapped  to non-trivial elliptic elements in $CM_{i+1}$,
$i=1,\ldots,s-1$.

\item"{(3)}" the epimorphisms along $CGRes$ are proper epimorphisms.

\item"{(4)}" the resolution $CGRes$
 is a graded strict resolution ([Se1],5), and a well-structured
resolution. All the graded abelian decompositions that are associated with the various limit groups
(over free products), $CM_i$, $i=1,\ldots,s$, are their graded JSJ decompositions over free products. Furthermore, 
the graded JSJ decompositions (over free products) of
the limit groups $CM_i$, $i=1,\ldots,s$, have the same structure as the corresponding graded JSJ decompositions (over free products) 
of the limit groups, $L_i$, where the
difference is only in the rigid vertex groups of the JSJ decompositions  and in the elliptic factors.
 
\item"{(5)}" if no factor of $L_i$
is  rigid nor solid, then $CM_{i+1}$ is a cover of a shortening
quotient of $CM_i$, for $i=1,\ldots,s-1$. 

\item"{(6)}" if a factor of $L_i$ is rigid then the corresponding factor of $CM_i$ is rigid. This can happen only for $i=s$. 
 If a factor of $L_i$ is solid then the corresponding factor of $CM_i$ is weakly solid (definition 1.18).

\item"{(7)}" in case $L_i$ is a free product of a  solid limit group with (possibly)
elliptic factors and (possibly) a free group, for some $i<s$, then by (6) the corresponding factor of $CM_i$ is weakly solid, 
and $CM_{i+1}$ is a free product of a cover of a
flexible quotient of 
the weakly  solid factor in $CM_i$ with (possibly)  the same   elliptic factors of $CM_i$ 
and (possibly) the same free group (as in $CM_i$).

\item"{(8)}" $L_s$, the terminal limit group of the graded resolution $GRes$, 
 is a free product of a rigid or a solid limit group with (possibly) some elliptic
factors and (possibly) a free group. $CM_s$ is a free product of a rigid or a weakly solid factor, that we denote $RSF_s$, in correspondence with
the factor of $L_s$,  with (possibly) f.p.\ covers of the  elliptic factors of $L_s$, 
and a free group (of the same rank as in $L_s$).

With the rigid or the weakly solid factor of $CM_s$, $RSF_s$, we associate finitely many  limit groups (over free products), $CF_1,\ldots,CF_g$, 
with ungraded resolutions (over free products), $Res_1,\ldots,Res_g$, that terminate in f.p.\ limit groups (cf. theorems 1.17 and 1.20). Hence,
the completions of $Res_1,\ldots,Res_g$ are f.p.\ and the limit groups (over free products), $CF_1,\ldots,CF_g$, 
can be embedded in
f.p.\ completions  
(that have  f.p.\ terminal limit groups). 

\item"{(9)}" the finite collection of limit groups, $CF_1,\ldots,CF_g$, are all limit quotients of $G(x,p)$, and they
dominate all the flexible quotients of $RSF_s$
(although they need not be quotients of $RSF_s$). 

By theorem 1.16, 
from the sequence of homomorphisms, $\{h_n\}$, it is possible to extract
a subsequence, still denote $\{h_n\}$,  on  which it is possible to perform iterative shortenings, and obtain another
sequence of homomorphisms, $\{u_n\}$, which is asymptotically rigid or asymptotically 
strictly solid and converges into the rigid or solid factor of $L_s$,
the terminal graded limit group in $GRes$.   

The homomorphisms $\{u_n\}$ factor through $CM_s$, and they restrict to rigid or weakly strictly solid homomorphisms of
$RSF_s$ (with respect to that cover). Furthermore, in the rigid case, every flexible homomorphism of $RSF_s$ 
factors through at least one of the resolutions, $Res_1,\ldots,Res_g$, that are associated with the limit groups, $CF_1,\ldots,CF_g$.
In the weakly solid case, for every non-weakly strictly solid homomorphism of $RSF_s$,
$f: RSF_s \to A^1* \ldots *A^{\ell}$, there exists a homomorphism, $u: RSF_s \to A^1*\ldots *A^{\ell}$, such that 
the pair, $(f,u)$, extends to a homomorphism of the completion of the identity resolution, $RSF_s \to RSF_s$, and the homomorphism
$u$ as a homomorphism of the f.p.\ group, $G(x,p)$, 
factors through at least one of the resolutions, $Res_1,\ldots,Res_g$.

\item"{(10)}"   with the terminal graded limit group $CM_s$, it is possible to associate an ungraded strict well-structured
resolution:
$CM_s=V_0 \to V_1 \to \ldots \to V_t$, so that the terminal (ungraded) limit group (over free products) $V_t$ is a free product of
(possibly) elliptic factors and (possibly) a free group. Every non-trivial elliptic element in each of the limit groups
(over free products) $V_i$, $i=0,\ldots,t-1$, is mapped to a non-trivial element in $V_{i+1}$. 

With the combined (ungraded) resolutions, $CM_0 \to \ldots \to CM_s \to V_1 \to \ldots \to V_t$ it is possible to associate a completion,
$Comp_{CM}$, which is f.p.\ (and terminates in a f.p.\ group).
\endroster
\endproclaim

\nfp Given the sequence of homomorphisms, $\{h_n:G(x,p) \to A^1_n*\ldots *A^{\ell}_n\}$, we can extract a subsequence from it 
(still denoted $\{h_n\}$), from which it 
is possible to obtain a graded  strict, well-structured resolution over free products,  $L_0 \to L_1 \to \ldots \to L_s$ , according to theorem
1.16.
Furthermore, with
each homomorphism $h_n$ (from that subsequence) it is possible to associate a homomorphism $u_n : L \to A^1_n*\ldots*A^{\ell}_n$, so that the homomorphisms
$\{u_n\}$ restrict to an asymptotically rigid or to an asymptotically strictly solid homomorphisms with respect to the rigid or solid factor
of the terminal limit group, $L_s$.

Following theorem 18 in [Ja-Se], starting with the sequence of homomorphisms, $\{u_n\}$, we can pass to a further subsequence,
and associate with the limit
group over free products, $L_s$ ,  a strict well-structured ungraded resolution, $L_s=W_0 \to W_1 \to \ldots \to W_t$, so that
$W_t$ is the free product of (possibly) finitely many elliptic factors and (possibly) a free group.

At this point we look at the combined (ungraded) resolution:
$$L_0 \to L_1 \to \ldots \to L_s=W_0 \to W_1 \to \ldots \to W_t$$ 
We set $Comp_L$ to be the completion of this ungraded combined resolution. Note that all the limit groups over free products, $L_0,\ldots,L_s$
and $W_1,\ldots,W_t$, are embedded in the completion $Comp_L$.

Starting with $Comp_L$ we define a sequence of completions, obtained by replacing the terminal limit group 
(over free products) in $Comp_L$, $W_t$, with approximating f.p.\ groups, that are obtained from the terminal limit group
$W_t$, by imposing only relations of length at most $m$ from each factor of $W_t$ (with respect to a fixed set
of generators of these factors). We denote the approximating f.p.\ subgroups, $W_t^m$, and the group that is obtained from
$Comp_L$ be replacing $W_t$ with $W_t^m$, we denote $Comp_L^m$.

\noindent
For large enough index $m$, we do obtain a new completion, hence, for large $m$, $Comp_L^m$ is a completion of a strict well-structured
ungraded resolution over free products. We denote the resolution that is associated with $Comp_L^m$, $Res^m$, and this is the resolution: 
$$L_0^m \to L_1^m \to \ldots \to L_s^m=W_0^m \to W_1^m \to \ldots \to W_t^m$$ 
where the groups $L_i^m$ and $W_j^m$ approximate the corresponding subgroups, $L_i$ and $W_j$.

By construction, with each of the graded limit groups, $L_i$, there is an associated graded (virtually) abelian decomposition, which is the
graded virtually abelian JSJ decomposition of $L_i$ over free products. Also, with each of the ungraded limit group, $W_j$, there
is an associated (virtually) abelian decomposition which is the  virtually abelian decomposition of $W_j$ over free products. 
These graded and ungraded JSJ decompositions lift to virtually abelian decompositions of the approximating limit groups, $L_i^m$ and $W_j^m$.

Moreover, by the proof of theorem 27 in [Ja-Se], for large enough index $m$, these graded virtually abelian decompositions of
the limit groups $L_i^m$ (that are lifted from $L_i$) 
are their graded  JSJ decompositions over free products, and the virtually abelian
decompositions of the limit groups $W_j^m$ (that are lifted from $W_j$) are their  JSJ decompositions over free products.

By theorems 1.17 and 1.20, if $L_i$  is solid, then for large $m$, $L_i^m$ is weakly solid. If $L_s$ is rigid (it's the only graded limit group
along the resolution $Res$ that can be rigid), then for large $m$, $L_s^m$, is rigid. Therefore, for large enough $m$, the resolution, $Res^m$,
satisfy properties (1)-(7) and (10). Properties (8) and (9) hold for $Res^m$, for large enough $m$, by the conclusions of theorems 1.17 and
1.20.

\line{\hss$\qed$}

Theorem 1.21 associates a collection of covers of resolutions (over free products) with a given f.p.\ group.
Every homomorphism of a given f.p.\ group into free products factors through one of these cover resolutions. As in the ungraded
case (theorem  26  in [Ja-Se]), if we start with a f.p.\ group, $G(x,p)$, there exists a finite collection of cover
resolutions through which all the homomorphism from the f.p.\ group into free products do factor. Such a finite
collection is not canonical, but as in the ungraded case, we view each such finite collection as the
$graded$ $Makanin$-$Razborov$ diagram of the graded f.p.\ group, $G(x,p)$, over free products.

\vglue 1pc
\proclaim{Theorem 1.22} Let $G(x,p)$ be a f.p.\ group. There exist finitely many well-structured resolutions:  
$CM=CM_0 \to CM_1 \to \ldots \to CM_s$, where $CM$ is a graded limit quotient of $G(x,p)$ (over free
products), and $CM_s$ is a free product of  either a rigid or a weakly solid factor with (possibly) 
finitely many elliptic factors, and (possibly) a free  group,
that are all covers of graded resolutions, that satisfy properties (1)-(10) that are listed in theorem 1.21. 

With each of the finitely many cover graded resolutions we associate the collection of homomorphisms
that factor through it. Homomorphisms of the terminal graded limit group, $CM_s$, are obtained from 
rigid or weakly strictly solid homomorphisms of the rigid or weakly solid factor of $CM_s$ that are combined with homomorphisms
of the elliptic factors (into elliptic subgroups), and arbitrary homomorphisms of the (possible) free
factor. We further require that the restriction of weakly strictly solid homomorphisms to edge groups and to abelian
and QH vertex groups in the graded abelian JSJ decomposition of the weakly solid factor are not (entirely) elliptic.

Homomorphisms of the limit groups, $CM_{s-1},\ldots,CM_0$ (and hence of the given f.p.\ group
$G(x,p)$), are obtained from  homomorphisms of $CM_s$, in a similar way to homomorphisms that factor
through graded resolutions over free and hyperbolic group, i.e., in the following way:
\roster
\item"{(1)}" in case $CM_i$, $i=0,\ldots,s-1$,  does not have a factor which is  rigid or weakly solid, homomorphisms of $CM_i$
are obtained from homomorphisms of $CM_{i+1}$, by precomposing the homomorphisms of $CM_{i+1}$ with
graded modular automorphisms of $CM_i$. 

\item"{(2)}" in case $CM_i$, $i=0,\ldots,s-1$, is weakly solid, we look at the following homomorphisms.
With $CM_i$ we associate a one step (well-structured) resolution $CM_i \to CM_{i+1}$. With this one step
resolution we associate its completion, $Comp_i$. Both $CM_i$ and $CM_{i+1}$ are naturally embedded in
$Comp_i$. With $CM_i$ we associate all its homomorphisms that can be extended to homomorphisms of
the completion, $Comp_i$, so that the restrictions of these homomorphisms of $Comp_i$ to the image of
$CM_{i+1}$ in $Comp_i$ is a homomorphism that is associated with $CM_{i+1}$.

\item"{(3)}" the restriction of the homomorphisms of the various groups, $CM_i$, $i=0,\ldots,s-1$, to the
edge groups and to abelian
and QH vertex groups in the graded abelian JSJ decomposition of theses are not (entirely) elliptic.
\endroster

Every homomorphism of $G(x,p)$ into a free product factors through at least one of the
finitely many (cover) resolutions that are associated with the f.p.\ group $G(x,p)$ (i.e., the cover resolutions that appear in its
(non-canonical) Makanin-Razborov diagram).
\endproclaim

\nfp The cover graded resolutions that appear in theorem 1.21, can be defined by finitely many generators, relations, and generators of
subgroups, since these are determined by finitely many f.p.\ completions and finitely may f.g.\ subgroups of these f.p.\ completions.
The construction of cover graded resolutions in theorem 1.21 allows one to consider only homomorphisms that satisfy properties
(1)-(3) as homomorphisms that factor through a cover graded resolution.

By theorem 1.21, every homomorphism of the f.p.\ group, $G(x,p)$,  factors through at least one of the cover graded resolutions
that can be associated with $G(x,p)$. Therefore, by 
ordering the countable set of cover graded resolutions, and applying the same argument that was used to prove theorem 26 in [Ja-Se], i.e., the
argument that was used to show that the set of all homomorphisms of a f.p.\ group into free products factor through finitely many
cover ungraded resolutions, there exists a finite collection of cover graded resolutions, through which all the homomorphisms of the f.p.\ group, $G(x,p)$,
do factor. 

\line{\hss$\qed$}

The finite collection of graded resolutions that are associated with a f.p.\ group, $G(x,p)$, according
to theorem 1.22, and through which all the homomorphisms of $G(x,p)$ into free products do factor, form
a graded Makanin-Razborov diagram of $G(x,p)$ over free products. Note that the collection of graded
resolutions is not canonical, hence, the diagram we constructed is not canonical. 

Also, note that over
free and hyperbolic groups we needed to study the singular locus of the graded resolutions that appear in
the (in these cases canonical) graded Makanin-Razborov diagram. This was crucial in studying sentences and
predicates over these groups. The diagram that we constructed over free 
products, and particularly the conditions on homomorphisms that factor through each graded resolution in
the diagram, guarantee that the collection of homomorphisms that are associated with each graded resolution from
the diagram lies outside the singular locus, and these collections cover the entire set of homomorphisms from
the f.p.\ group $G(x,p)$ into free products. Hence, there is no need to study the singular locus of the
constructed graded resolutions.

\vglue 1.5pc
\centerline{\bf{\S2. Formal Solutions and Formal Limit Groups}}
\medskip
So far we have generalized the construction of the Makanin-Razborov diagram from a free 
group to a free product of groups [Ja-Se], and in section 1 we have generalized the study of
systems of equations with parameters (over free products), and associated a graded Makanin-Razborov
with such a system.
 
To analyze the structure of elementary sets over a free product, we first need to generalize $formal$
$solutions$, and their (formal) limit groups ([Se2]), to formal solutions and formal limit groups over free products.
As in the construction of the Makanin-Razborov diagrams (graded and ungraded), we study formal solutions
over the entire class of free products, and not over a particular free product. This approach will finally lead us
to the construction of uniform  proofs of a given sentence (or predicate),  uniform proofs that
lead to a (uniform) reduction of a sentence over a free product to sentences over the various factors.

 We start this section with a generalization of (a special case of) Merzlyakov
original theorem, mainly to demonstrate that basic concepts over free groups generalize to free
products.

\vglue 1pc
\proclaim{Theorem 2.1 (cf. [Me])} Let $G=A^1*\ldots *A^{\ell}$ be a non-trivial free product that is not isomorphic to
$D_{\infty}$.
Let
$y=(y_1,\ldots,y_{\ell})$ and  $x=(x_1,\ldots,x_{q})$. Let
$w_1(x,y,a)=1,\ldots,w_s(x,y,a)=1$ be a system of equations over $G$, where $a$ is a tuple of elements from the factors,
$A^1,\ldots,A^{\ell}$.
 Suppose that
the sentence:
$$\forall y \ \exists x \ \ w_1(x,y,a)=1,\ldots,w_s(x,y,a)=1$$
is a truth sentence over $G$. Then there exists a $formal$ solution
$x=x(y,a')$, where $a'$ is a tuple of elements from the factors, $A^1,\ldots,A^{\ell}$,  
so that  each of the words $w_j(x(y,a'),y,a)$ is the trivial word in
the free product $<y>*A^1*\ldots *A^{\ell}=<y>*G$.
\endproclaim

\nfp In a non-trivial free product, $G=A^1*\ldots *A^{\ell}$, which is not isomorphic to $D_{\infty}$, there exists a free,
purely hyperbolic subgroup,
so that its Cayley graph is quasi-isometrically mapped, under the embedding of the free group,
 into the Bass-Serre tree that corresponds to the free product, $G=A^1*\ldots*A^{\ell}$. We can clearly write $G$ as a non-trivial free
product, $G=C*B$, where $C$ is not isomorphic to $Z_2$.
If $c_1,c_2 \in C$ are distinct non-trivial elements, and $b \in B$ is non-trivial, then
$w_n=(c_1b)^n$ and $t_n=(c_2b)^n$, generate such a free group of rank 2, for $n$ large enough (e.g. $n>5$).

The sentence is a truth sentence over the free product, $G$, hence, we may assign arbitrary values to the universal variables $y$, and for
every chosen values (in $G$), there exist values for the existential variables $x$, so that the equalities, $w_1=1,\ldots,w_s=1$, hold
in $G$. 

In $G$ there is a copy of the free group $F_2$ that is quasi-isometrically embedded in the Bass-Serre tree that corresponds to the free
product, $G=A^1*\ldots *A^{\ell}$. We assign a sequence of values in $G$ to the universal variables $y$. We give the universal
variables $y$, a sequence of values in $G$ that we denote $\{y(n)\}$, so that these values are in (the image of) the 
(quasi-isometrically embedded) free group $F_2$ in $G$, that form a test sequence in $F_2$ 
(see the proof of theorem 1.1 and definition 1.20 in [Se2] for the notion of test
sequence). 

For every index $n$, there are values in $G$ of the existential values, $x$, that we denote $x(n)$, so that $w_i(x(n),y(n),a)=1$,
$i=1,\ldots,s$. For each $n$ we choose the shortest possible $x(n)$ (in the Bass-Serre tree that is associated with the free
product, $G=A^1*\ldots *A^{\ell}$, with its simplicial metric) that satisfies these equalities. 

By the same (geometric) argument that was used to prove Merzlyakov's theorem  over free groups (theorem 1.1 in [Se2]), by possibly 
iteratively modifying the shortest values $x(n)$ finitely many times, the sequence of elements $\{(x(n),y(n),a)\}$ has a subsequence
that converges to the limit group over free products, $L(x,y,a) \, = \, <y>*A_1*\ldots*A_{\ell}*E_1*\ldots*E_t*F_r$, where:
\roster
\item"{(1)}" $<y>$ is a free group generated by the universal variables: $y_1,\ldots,y_{\ell}$.

\item"{(2)}" $A_i$, $i=1,\ldots,\ell$, is elliptic and contains the  tuple of (fixed) elements $a$ from the factor $A^i$.

\item"{(3)}" $E_1,\ldots,E_r$ are elliptic. $F_r$ is a (possibly trivial) free group of rank $r$.

\item"{(4)}" $L$ contains the elements $x$, and the elements, $w_1(x,y,a),\ldots,w_s(x,y,a)$, represent the trivial element
in the limit group $L$.
\endroster

We replace $L$ by its (quotient) factor, $R(x,y,a) \, = \, <y>*A_1*\ldots *A_{\ell}$, which is also a (limit) quotient of $L$
 ($R$ is in fact a retract of $L$).
$L$, and hence $R$, are f.g.\ groups, but in general they need not be finitely presented, i.e., the elliptic factors $A_1,\ldots,A_{\ell}$,
are
f.g.\ but not necessarily finitely presented. Still, the words $w_1,\ldots,w_r$ represent the trivial element in $R$, as they represent
the trivial element in $L$, and $R$ is a retract of $L$. 
Hence, if we look at a free group, $F_y$,  generated by the elements $y_1,\ldots,y_{\ell}$ and finite
(fixed) generating sets for $A_1,\ldots,A_{\ell}$, then  there is a finite collection of relations, $r_1,\ldots,r_f$, 
 in the (fixed) generating sets of $A_1,\ldots,A_{\ell}$, relations that hold in these factors,
such that if these relations are imposed on $F$, then the words, $w_1,\ldots,w_s$, represent the trivial element in the obtained quotient
f.p.\ group.

$L$ is a limit group that is obtained as a limit from a sequence of elements, $\{(x(n),y(n),a)\}$. Hence, for large enough $n$, the relations
$r_1,\ldots,r_f$ that hold in $A_1,\ldots,A_{\ell}$, hence in $L$, 
hold for the specializations, $\{(x(n),y(n),a)\}$, and the fixed set of generators of $A_1,\ldots,A_{\ell}$ are elliptic, hence, 
their  specializations are in the factors $A^1,\ldots,A^{\ell}$ correspondingly. Therefore, using the specialization $\{(x(n),y(n),a)\}$,
for large $n$,
it is possible to find  elements
$x \in <y>*A^1*\ldots*A^{\ell}$, 
so that the words, $w_1(x,y,a),\ldots,w_s(x,y,a)$ are the trivial elements in the free product, $<y>*G$.

\line{\hss$\qed$}

To analyze sentences and predicates over free products, we will need
a generalization of Merzlyakov
theorem
to a truth sentence defined over an arbitrary (given) variety. 
Before generalizing the results that were proved over free group for general varieties,
we generalize Merzlyakov theorem to a coefficient free sentence 
that contains inequalities. 
For such sentences the conclusion can be stated in a uniform way for
all free products for which the sentence is a truth sentence.

\vglue 1pc
\proclaim{Theorem 2.2} 
Let
$w_1(x,y)=1,\ldots,w_s(x,y)=1$ be a system of equations (over a group),
and let $v_1(x,y),\ldots,v_r(x,y)$ be a collection of words in
the free group generated by  $\{x,y\}$.
Let:
$$\forall y \ \exists x \ \ w_1(x,y)=1,\ldots,w_s(x,y)=1 \ \wedge \
v_1(x,y) \neq 1,\ldots,v_r(x,y) \neq 1$$
be a sentence over groups. Then there exist finitely many f.p.\ limit groups over free products:
$$H_j=<y>*S_1*\ldots*S_{m_j}*F_{r_j} \ \ j=1,\ldots,t$$ 
and tuples of elements (formal solutions), $x_j \in H_j$, 
so that:
$w_1(x_j,y)=1,\ldots,w_s(x_j,y)=1$ in the limit group $H_j$. 

Furthermore,
 for every non-trivial free product,
$G=A^1*\ldots*A^{\ell}$, that is not isomorphic to $D_{\infty}$, and so that the given sentence is
a truth sentence over $G$, for at least one index $j$, $1 \leq j \leq t$, there exists a homomorphism
$\tau :H_j \to G*<y>$, that maps each factor $S_i$, $1 \leq i \leq m_j$, 
into an elliptic subgroup in $G$,  
the factor $<y>$ in $H_j$ isomorphically onto the factor
$<y>$ in $G*<y>$, and the free factor $F_{r_j}$ into $G$,  
such that 
%
the sentence:
$$\exists y \ \ v_1(\tau (x_j),y) \neq 1,\ldots,v_r(\tau (x_j),y) \neq 1$$
is a truth sentence in $G$. 
\endproclaim

\nfp The argument that we use is a combination of the proof of theorem 2.1, with the proof of Merzlyakov theorem with inequalities
(theorem 1.2 in [Se2]), and the construction of formal limit groups in section 2 of [Se2].

As we have pointed out in the proof of theorem 2.1, in  a non-trivial free product, $G=A^1*\ldots*A^{\ell}$, that is not isomorphic to
$D_{\infty}$, 
there exists a free, purely hyperbolic subgroup,
so that its Cayley graph is quasi-isometrically mapped, under the embedding of the free group,
 into the Bass-Serre tree that corresponds to the free product, $G=A^1*\ldots*A^{\ell}$.

We look at all the non-trivial free products, $G=A^1*\ldots*A^{\ell}$ (for an arbitrary finite $\ell$), that are not isomorphic
to $D_{\infty}$, over which the given sentence  is a truth
sentence. For
 each such free product $G$, we look at all the f.g.\ free groups that are embedded in $G$, are purely hyperbolic in $G$, and the Cayley
graph of $F$ is quasi-isometrically mapped into the Bass-Serre tree that is associated with the free product of $G$, $G=A^1*\ldots*A^{\ell}$. 

We look
at all such free products $G=A^1*\ldots*A^{\ell}$ (over which the given sentence is a truth sentence), 
the quasi-isometrically embedded free groups,  $F$,
all possible test sequence for the universal variables $y$ in $F$, and for each value of the universal variables $y$, the shortest possible 
values for the existential variables $x$ (in the simplicial metric of the Bass-Serre tree that is associated with
the free product $G=A^1*\ldots*A^{\ell}$), for which both the equalities, $w_1,\ldots,w_s$, and the inequalities, $v_1,\ldots,v_r$, 
hold (in $G$).

By the construction of formal limit groups that appears in section 2 of [Se2], and by the argument that was used to prove
that a f.p.\ group has only finitely many maximal limit quotients over free products (theorem 21 in [Ja-Se]), 
with the entire collection of these test sequences it is
possible to associate finitely many (maximal) formal limit groups (over free products). 
Since these limit groups were constructed from sequences that restrict
to test sequences for the universal variables, and shortest possible values for the existential variables, each of the finitely many
maximal formal limit groups has the structure: 
$$H_j=<y>*S_1*\ldots*S_{m_j}*F_{r_j} \ \ j=1,\ldots,t$$ 
 where the factors, $S_1,\ldots,S_{m_j}$, are elliptic, and $F_{r_j}$ is free. Furthermore, each such (maximal) formal limit group,
$H_j$, is finitely presented, and with each such formal limit group, $H_j$, there are
associated tuples of elements, $x_j \in H_j$, so that the words, $w_1(x_j,y),\ldots,w(x_j,y)$, are trivial in $H_j$, and the words,
$v_1(x_j,y),\ldots,v_r(x_j,y)$, are non-trivial in $H_j$. 

Suppose that $G=A^1*\ldots*A^{\ell}$ is a non-trivial free product, that is not isomorphic to $D_{\infty}$, and the given coefficient free
sentence is true over $G$. Then there exists a free group $F_2$, with an embedding into $G$ so that the Cayley
graph of $F_2$ is mapped quasi-isometrically into the Bass-Serre tree that is associated with the free product, $G=A^1*\ldots*A^{\ell}$. 
We assign (tuples of) values, $y(n)$,  from (the embedding in $G$ of) a test sequence of $F_2$ to the universal variables $y$. 
Given each $y(n)$, we set
$x(n)$ to have the shortest possible values in $G$ (with respect to the simplicial metric on the Bass-Serre tree that
is associated with $G=A^1*\ldots*A^{\ell}$).

Since the limit groups over free products, $\{H_j\}$, are finitely presented, the sequence $\{(x(n),y(n))\}$  has a subsequence 
(still denoted $\{(x(n),y(n))\}$, that factors
thorough one of the limit groups, $H_j$. Furthermore, we may assume that for every index $n$, the restriction of the specializations,
$\{(x(n),y(n)\}$, to the 
(elliptic) factors of $H_j$, $S_1,\ldots,S_{m_j}$, are elliptic subgroups in $G$, and the elements, $v_1(x(n),y(n)),\ldots,v_r(x(n),y(n))$, 
are all non-trivial in $G$.

For every index $n$,  we can define a  homomorphism $h_n: H_j \to G$, given by the specialization $\{(x(n),y(n)\}$. We modify each homomorphism
$h_n$, to a homomorphism: $u_n: H_j \to <y>*G$, by sending the subgroup $<y>$ in $H_j$ isomorphically onto the subgroup $<y>$ in $<y>*G$,
and setting $u_n$ restricted to the factor $S_1*\ldots*S_{m_j}*F_{r_j}$ to be identical to $h_n$ restricted to that factor. $u_n$ is clearly
a homomorphism, so all the elements, $w_1,\ldots,w_s$, which are the identity element in $H_j$ are mapped to the identity
element in $<y>*G$. Furthermore, $v_1,\ldots,v_r$ are mapped to non-trivial elements in $G$ by the homomorphisms $h_n$, hence, they must
be mapped to non-trivial elements by $u_n$, as $h_n = \nu_n \circ u_n$, where $\nu_n: <y>*G \to G$, $nu_n(y)=y(n)$.

Finally, we set $\tau=u_n$ for an arbitrary $n$. The equalities $w_i(\tau(x),y)=1$ hold  in $G$ (as they hold in
$<y>*G$), $i=1,\ldots,s$. 
Furthermore, the sentence:
$$\exists y \ \ v_1(\tau (x_j),y) \neq 1,\ldots,v_r(\tau (x_j),y) \neq 1$$
is a truth sentence in $G$, as the inequalities hold for $y=y(n)$.
  
\line{\hss$\qed$}

To generalize Merzlyakov's original theorem to sentences (over a free group) that hold over some varieties and
not only over the entire affine set, 
we have associated with each well-structured resolution a $completion$, and with a completion
we have associated $closures$ of it, that are obtained by adding roots to a finite collection of elements in 
abelian vertex groups that appear along the well-structured resolution (see definitions 1.12 and 1.15 in [Se2]).

Given a completion of a well-structured resolution over a free group, and a (finite) collection of
closures of that completion, we call the finite collection of closures, a $covering$ $closure$ of the completion
(definition 1.16 in [Se2]),
if every homomorphism that factors through the well-structured resolution can be extended to a completion
of at least one of the closures of the completion (from the given finite collection). 

These notions generalize to well-structured resolutions over free products, although the generalizations
require some modifications.  
The construction of a completion generalizes naturally
and canonically from well-structured resolutions over a free group (definition 1.12 in [Se2]) to
well-structured resolutions over free products. The notion of a closure of the completion over free products,
requires some changes in comparison with  the similar object over a free group.

\vglue 1pc
\proclaim{Definition 2.3} Let $Res(y)$ be a well-structured 
resolution over free products, and let
 $Comp(Res)(z,y)$ be
its completion. 

\noindent
Let $E_1,\ldots,E_r$ be the terminal elliptic subgroups (factors) in the well-structured resolution, $Res(y)$. Note
that these are also the terminal elliptic subgroups in the completion, $Comp(Res)(z,y)$.  
Let $Ab_1,\ldots,Ab_d$ be the non-conjugate, non-cyclic, maximal abelian subgroups that appear along the 
completion, $Comp(Res)(z,y)$, and are mapped onto
a non-cyclic, (non-elliptic)  abelian factor in a free decomposition associated with one
of the levels of the completion. 

\noindent
Let 
 $PAb_1,\ldots,PAb_{pd}$ be the non-conjugate, non-cyclic, (non-elliptic) maximal $pegged$ abelian groups
that appear along the completed resolution,
i.e., maximal non-cyclic abelian subgroups in
$Comp(Rlim)(z,y)$,
that are mapped onto a non-cyclic, (non-elliptic)  abelian vertex group in some  abelian decomposition
associated with some level of the completed resolution $Comp(Res)(z,y)$, and this abelian
 vertex is connected to the other vertices of the completed decomposition of that level by an edge
 with (maximal) cyclic (non-elliptic) stabilizer.
 We call the maximal cyclic subgroup of 
a pegged abelian group connecting it to the other vertices of the corresponding
completed decomposition, the $peg$ of the $pegged$ abelian group $PAb$.

Let $S_1,\ldots,S_d$ be free abelian groups so that $Ab_1<S_1,\ldots,Ab_d<S_d$ are
subgroups of finite index.
Let $PS_1,\ldots,PS_{pd}$ be free abelian groups
so that $PAb_1<PS_1,\ldots,PAb_{pd}<PS_{pd}$ are subgroups of finite 
index, and the pegs $peg_1,\ldots,peg_{pd}$  are primitive elements in the
ambient free abelian groups $PS_1,\ldots,PS_{pd}$.

A $closure$ of the completed resolution $Comp(Res)(z,y)$ is obtained by replacing
the free abelian groups $Ab_1,\ldots,Ab_d$ by the free abelian groups
$S_1,\ldots,S_d$, and 
the pegged abelian groups $PAb_1,\ldots,PAb_{pd}$ by the free
abelian groups 
$PS_1,\ldots,PS_{pd}$ in correspondence, along the entire 
 completed resolution, i.e., from the top level through the bottom level
in which a subgroup of the pegged abelian group appears along the completed
resolution.

We also associate with the closure new elliptic subgroups, $D_1,\ldots,D_{f}$, where $\ell \leq r$, and $B_1,\ldots,B_t$, with the 
following properties:
\roster
\item"{(1)}"  (up to a change in the order of the elliptic factors $E_1,\ldots,E_r$) $E_i$ is 
mapped into $D_i$, $i=1,\ldots,f$. 

\item"{(2)}" for each index $i$, $f+1 \leq i \leq r$, we add a new (free) generator,
$c_i$, and map $E_i$ into $D_{j_i}$. 
%
\endroster

The closure is obtained from the completion, $Comp(Res)$, by possibly enlarging the maximal abelian and maximal pegged abelian
subgroups (by finite index supergroups), and further replacing each of the terminal elliptic subgroups, $E_1,\ldots,E_f$,
with the corresponding subgroup, $D_1,\ldots,D_{f}$, and replacing each of the subgroups, 
$E_{f+1},\ldots,E_r$, by its image in $D_{j_i}$ conjugated by $c_i$, $i=f+1,\ldots,r$. The terminal limit group of the closure
is the free product of the elliptic subgroups, $D_1*\ldots*D_{f}$, with the free product of the (additional) elliptic subgroups,
$B_1*\ldots*B_t$  (possibly) free product with a free group and finitely many abelian (non-elliptic) factors and closed surface groups.
The
completion, $Comp(Res)(z,y)$, is mapped naturally into a closure of it.
\endproclaim

Having defined closures of a completion, we can generalize the notion of a $covering$ $closure$.

\vglue 1pc
\proclaim{Definition 2.4} Let $Res(y)$ be a well-structured 
resolution over free products, let $Comp(Res)(z,y)$ be its completion,
and let $Cl_1(Res),\ldots,Cl_v(Res)$ be a finite set of closures of $Comp(Res)$. 

Let $G=A^1*\ldots*A^{\ell}$ be a (non-trivial) free product that is not isomorphic to $D_{\infty}$.
We say that the given finite collection of closures is a $covering$ $closure$ of the completion,
$Comp(Res)$, over the free product $G$, if every tuple of specializations in $G$ of the variables $y$,
  that factors through the resolution $Res(y)$, can be extended to a specialization that factors through at least
one of the closures, $Cl_1,\ldots,Cl_v$. 
\endproclaim

The completion of a resolution $Comp(Res)(z,y)$,
 its closures $Cl(Res)(s,z,y)$,
 and the notion of a $covering$ $closure$, finally allow us to 
present $formal$
$solutions$ associated with a well-structured resolution of a limit group over free products.

\vglue 1.5pc
\proclaim{Theorem 2.5} 
Let $u_1(y),\ldots,u_m(y)$ be a collection of words in the free group, $<y>$, and let $(L,E_L)$ be
a limit group over free products that is a limit quotient (over free products) of the f.p.\ group:
$G(y)=<y \, | \, u_1(y),\ldots,u_m(y)>$.  
Let $Res(y)$ be a  well-structured 
resolution of the limit group $(L,E_L)$, 
and let $Comp(Res)(z,y)$  be the completion of the resolution
$Res(y)$. 
 
Let
$w_1(x,y)=1,\ldots,w_s(x,y)=1$ be a system of equations,
and let $v_1(x,y),\ldots,v_t(x,y)$ be a collection of words in
the alphabet $\{x,y\}$.
Let
the sentence:
$$\forall y \ (u_1(y)=1,\ldots,u_m(y)=1) \ \exists x \ \
 w_1(x,y)=1,\ldots,w_s(x,y)=1 \wedge $$
$$ \wedge \ 
 v_1(x,y) \neq 1,\ldots,v_t(x,y) \neq 1 $$
be a sentence over groups.

 There exists a finite collection of closures of the resolution $Res(y)$ over free products:
 $Cl(Res)_1(s,z,y),\ldots,Cl(Res)_q(s,z,y)$, so that for each index $i$, $1 \leq i \leq q$,
there exists a limit group over free products, $H_i=Cl_i(Res)*F_{d_i}$, with
a tuple of elements $x_i \in H_i$, for which the words,
 $w_1(x_i,y),\ldots,w_s(x_i,y)$, are the trivial elements in the limit groups $H_i$, and the words, $v_1(x_i,y),\ldots,v_t(x_i,y)$,
are non-trivial elements in the limit groups $H_i$, for $i=1,\ldots,q$.

In addition, let $G=A^1*\ldots*A^{\ell}$ be a non-trivial free product that is not isomorphic to $D_{\infty}$,
and
suppose that the given sentence is a truth sentence over $G$. Let $E_1,\ldots,E_r$ be the terminal elliptic 
subgroups in the completion, $Comp(Res)$, of the resolution, $Res(y)$. Then:
\roster
\item"{(1)}" the closures, $Cl_1,\ldots,Cl_q$, form a covering closure (definition 2.4) of the completion, $Comp(Res)$,
over the free product $G$.

\item"{(2)}" let $(z_0,y_0)$ be a tuple of specializations from $G$ that factors through the completed resolution, $Comp(Res)(z,y)$.
There exists an index $i$, $1 \leq i \leq q$, 
for which the specialization $(z_0,y_0)$ extends to a specialization
$(s_0,z_0,y_0)$ of the closure, $Cl_i(Res)$. 

 the specialization $(z_0,y_0)$ restricts to homomorphisms of the elliptic subgroups of the completion,
$Comp(Res)$, 
$E_1,\ldots,E_r$, into the free product $G$. These homomorphisms extend to a homomorphism $h: H_i \to G$
for which for every index $j$, $1 \leq j \leq r$:
$v_j(h(x_i),h(y)) \neq 1$ in $G$.
\endroster
\endproclaim

\nfp To construct the set of closures that is associated with a given sentence, we start with the collection of free products of
non-trivial groups,
$G=A^1*\ldots*A^{\ell}$, where $\ell>1$ is (an arbitrary)  positive integer, and  $G$ is not isomorphic to $D_{\infty}=Z_2*Z_2$, 
and over which the given sentence is a truth sentence. Given the collection of all
 these groups, $\{G\}$, we look at the collection of all the test sequences of the given resolution over free products, $Res(y)$. Note
that the collection of test sequences is divided into finitely many subsets, where in each subset it is specified which of the terminal
elliptic subgroups of $Res(y)$, $E_1,\ldots,E_r$, are mapped into conjugate factors in $G$. If for some test sequence of
$Res(y)$ over a free products $G=A^1*\ldots*A^{\ell}$, two terminal elliptic factors, 
$E_i$ and $E_j$, are mapped into conjugate factors of $G$,
we add new elements 
to the (generators of the) resolution, $Res(y)$, that conjugate the factors into which $E_i$ and $E_j$ are mapped. 

Since the sentence is a truth sentence over the free products $G$ that we consider, for the specializations of each test sequence
(in the free products $G$), it is 
possible to add specializations to the existential variables $x$, so that both the equalities, $w_1,\ldots,w_s$, and the inequalities,
$v_1,\ldots,v_t$, do hold in $G$. Given each specialization in a test sequence, we choose the shortest possible specialization of the
existential variables $x$ (in the simplicial metric of the Bass-Serre tree that is associated with the free product, $G=A^1*\ldots,A^{\ell}$).

Given each test sequence of $Res(y)$, and the extensions to shortest possible existential variables,
it is possible to extract a subsequence, that converges into a limit group (over free products), that has the form,
$Cl(Res)*F_d$, for some closure of the given resolution, $Res(y)$, and (possibly trivial) free group, $F_d$. 
As we did in proving theorem 21 in [Ja-Se], we can define a partial
order on the collection of these obtained limit groups. We can replace each of the obtained closures, $Cl(Res)$, with a closure that is 
obtained from the completion, $Comp(Res)$ , by adding finitely many generators and relations, and 
 so that the elements, $w_1,\ldots,w_s$ are still trivial, and the elements, $v_1,\ldots,v_t$, are 
non-trivial in the relatively finitely presented closure. By the ability to perform such a replacement,  every 
maximal element with respect to the partial order on the obtained limit groups,
$Cl(Res)*F_d$, is obtained from $Comp(Res)$ by adding finitely many generators and relations. 
By the same argument that was used to prove the finiteness of maximal limit quotients
(over free products) of a f.p.\ group (theorem 21 in [Ja-Se]), there are only finitely many maximal limit groups of the form,
$Cl(Res)*F_d$, that dominate all the limit groups that are obtained from test sequences of the resolution, $Res$, and its extension to
shortest existential variables.

We denote these maximal limit groups, $H_1,\ldots,H_q$, where for each $i$, $i=1,\ldots,q$, $H_i=Cl_i(Res)*F_{d_i}$. Since the limit groups
$\{H_i\}$ dominate all the possible test sequences over groups for which the sentence is a truth sentence  (together with their
extensions to the existential variables), the closures, $Cl_1,\ldots,Cl_q$, do form a covering closure for the resolution, $Res$, for
every non-trivial free product $G$, which is not $D_{\infty}$, and  over which the sentence is a truth sentence. 
By construction, the elements $w_1,\ldots,w_s$ represent the
trivial element in each of the groups $H_i$. The elements $v_1,\ldots,v_t$ represent non-trivial elements. Part (2) follows
since the limit groups
$H_i$, $i=1,\ldots,q$, dominate (the tails of) all the test sequences of the resolution, $Res(y)$, 
over free products for which the given sentence is a truth
sentence (and their extensions to shortest  existential variables).

\line{\hss$\qed$}

Theorem 2.5 proves, in particular, that if a given sentence is a truth sentence over a variety $V$ that
is defined over some free product $G$, then there exist  formal solutions that prove the validity of the
sentence for generic points in the variety $V$. As over free and hyperbolic groups, we will need a 
uniform way to pick these formal solutions. To get the type of uniformity that we need, we still need to
collect the entire set of formal solutions. These can be encoded using formal limit groups (over free
products) and graded formal limit groups.

As over free groups (definition 2.1 in [Se2]), given a completion of a well-structured resolution
over free products,
$Comp(Res)(z,y)$, and a system of equations, $\Sigma(x,y)=1$, we define a $formal$ $limit$ $group$ 
that is associated with the completion and the system of equations, as a limit of a sequence of homomorphisms:
$h_n:<x,z,y> \to A^1_n*\ldots*A^{\ell}_n$, where the restrictions, $h_n:<z,y> \to A^1_n*\ldots*A^{\ell}_n$, form a test sequence of the
completion, $Comp(Res)(z,y)$, and $\Sigma(h_n(x),h_n(y))=1$ in $A^1_n*\ldots*A^{\ell}_n$.

In [Ja-Se] we have associated well-structured resolutions with a limit group over free products.
If we combine the results and the techniques of [Ja-Se], with the construction of formal resolutions 
of a formal limit group (over a free group)
in the second chapter
of [Se2], then given a formal limit group, $FL(x,z,y)$, we are able to  associate with it its collection of formal resolutions. 
Furthermore, in [Ja-Se] it is shown that given a f.p\ group it is possible to associate with it finitely
many well-structured resolutions, so that every homomorphism from the given f.p.\ group into a free
product factors through at least one of these finitely many resolutions, and so that the completion of each
of the finitely many resolutions is finitely presented (such a finite collection of resolutions is taken to be 
the (non-canonical) Makanin-Razborov diagram of the given f.p.\ group over free products).

Similarly,  let $Res(y)$ be a well-structured resolution  with a f.p.\ completion,
$Comp(Res)(z,y)$ (hence, $comp(Res)(z,y)$ terminates in a free product of finitely many f.p.\ elliptic factors and (possibly) a
f.g.\ free group).   Given a (finite) system of equations, $\Sigma(x,y)=1$, it is possible to associate
with the  given system of equations, $\Sigma(x,y)=1$, and the given resolution over free products, $Res(y)$, 
a finite collection of (well-structured) formal resolutions, 
with f.p.\ completions, so that:
\roster
\item"{(1)}" each resolution terminates in a f.p.\ group, $Cl_i(Res)*F_{d_i}$, 
where $Cl_i(Res)$ is a closure of the completion, $Comp(Res)$,  and 
the closure, $Cl_i(Res)$, and its terminating elliptic  factors are all 
finitely presented.

\item"{(2)}" every formal solution of the given completion and the given system of equations, over any non-trivial free product
that is not isomorphic to $D_{\infty}$,   factors
through at least one of the finitely many formal resolutions.
\endroster

Like the Makanin-Razborov diagram (over free products) that we associated with a f.p.\ group, the finite collection
of formal resolutions that we associated with a given f.p.\ completion and a finite system of equations is
not canonical. However, it encodes all the formal solutions that can be associated with the given completion
and the given system of equations, and it suffices to analyze sentences and predicates. Hence, we call such
 a finite collection of formal resolutions (that satisfy properties (1) and (2)), a $formal$ $Makanin$-$Razborov$
diagram of the completion, $Comp(Res)$, and the finite system of equations, $\Sigma$, over free products.
 
\medskip
To analyze sentences and predicates over free groups, we needed to collect not only formal solutions that
are defined over a given completion and a given system of equations, but rather to collect all the formal
solutions that are defined over a given graded completion, $Comp(Res(y,z,p))$, and a  system
of equations with parameters, $\Sigma(x,y,p)=1$. 

\noindent
Note that a technical difficulty that exists in the graded case, and does not
appear in the ungraded case,  is that  a (graded)  completion
of a graded resolution over free products, can not be assumed to be f.p.\ even if we start with a graded f.p.\ group. f.p.\ 
is essential in our approach to constructing (graded) Makanin-Razborov diagrams over free products. However,
this technical difficulty can be overcome by applying the techniques and the results that appear in theorems 1.17, 1.20, and 1.21,
while studying graded resolutions and their covers.

We start the construction of a graded formal Makanin-Razborov diagram over free products, by following the construction of 
formal graded limit groups and their formal graded Makanin-Razborov diagrams over free groups, that appears in section 3 in [Se2].
Let $Res(y,p)$ be a graded well-structured 
resolution over free products. We assume that $Res(y,p)$ terminates in a free product of a rigid or a solid limit group 
with (possibly) finitely many elliptic factors and (possibly) a f.g.\ free group. In case the terminal graded limit
group of $Res(y,p)$ contains a solid factor we assume that with the solid factor, there is an associated finite collection of covers
of its flexible quotients (see proposition 1.5).
 Let $Comp(Res)(z,y,p)$ be the graded completion
of $Res(y,p)$.

As in the ungraded case (over free products), given the graded completion, $Comp(Res)(z,y,p)$, of the graded well-structured resolution
over free products, $Res(y,p)$,
and a (finite) system of equations, $\Sigma(x,y,p)=1$, we define a $graded$ $formal$ $limit$ $group$ 
that is associated with the completion and the system of equations, as a limit of a sequence of homomorphisms
$h_n:<x,z,y,p> \to A^1_n*\ldots*A^{\ell}_n$, where the restrictions, $h_n:<z,y,p> \to A^1_n*\ldots*A^{\ell}_n$, 
form a graded test sequence of the
graded completion, $Comp(Res)(z,y,p)$, and $\Sigma(h_n(x),h_n(y),h_n(p))=1$ in $A^1_n*\ldots*A^{\ell}_n$. Furthermore, 
by possibly passing to a subsequence 
of the homomorphisms, $\{h_n\}$, it is possible to use the subsequence of homomorphisms (still denoted $\{h_n\}$), and associate with
them a graded formal resolution, that terminates in a graded limit group of the form, $Cl(Res)*F_d$, where $Cl(Res)$ is a closure
of the given completion, $Comp(Res)$,  and $F_d$ is a possibly trivial free group. Note that 
even if the completion, $Comp(Res)$, is
f.p.\ it does not imply that the closure, $Cl(Res)$, is f.p.\ as well.

Since the closure, $Cl(Res)$, may be infinitely presented, to construct finitely many formal graded resolutions that will cover all the
graded formal resolutions that are associated with $Comp(Res)$, we need to replace each graded formal resolution by a $cover$, in
a similar way to the construction of covers of graded resolutions over free products that appear in theorem 1.21.

\vglue 1pc
\proclaim{Theorem 2.6} Let $G(y,p)$ be a f.p.\ group,   let $Res(y,p)$ be a well-structured resolution over free products of a (graded)
limit quotient of $G(y,p)$ that terminates in a free product of a rigid or solid factor and possibly finitely many elliptic factors and 
a (possibly trivial) free group. Let $Comp(Res)(z,y,p)$ be the completion of $Res(y,p)$, and let $\Sigma(x,y,p)=1$ be a (finite) 
system of equations. 

Let $\{h_n: <x,z,y,p> \to A^1_n*\ldots*A^{\ell}_n\}$ (where $A^1_n*\ldots*A^{\ell}_n$  are  non-trivial free products that are not
isomorphic to $D_{\infty}$) be a sequence of homomorphisms, 
where the restrictions, $h_n:<z,y,p> \to A^1_n*\ldots*A^{\ell}_n$, form a graded test sequence of the
graded completion, $Comp(Res)(z,y,p)$, and $\Sigma(h_n(x),h_n(y),h_n(p))=1$ in $A^1_n*\ldots*A^{\ell}_n$. 
Furthermore, we assume that the sequence 
$\{h_n\}$ converges into a formal limit group over free products, $FL(x,z,y,p)$, 
and that from the sequence, $\{h_n\}$, it is possible to construct
a formal resolution (over free products), $FGRes$:
$$FL(x,z,y,p)= FL_0 \to FL_1 \to \ldots \to FL_s=Cl(Res)*F_d$$ 
where $Cl(Res)$ is a graded closure of the completion, $Comp(Res)(z,y,p)$, and  $F_d$ is a possibly
trivial free group.

There exists a f.g.\ graded formal limit quotient of $G(y,p)$, $CF$,  and
a well-structured graded formal resolution of $CF$, $CFGRes$, $CF=CF_0 \to CF_1 \to \ldots \to CF_s$, that covers the graded
formal resolution that is constructed from the homomorphisms, $\{h_n\}$, i.e., the formal graded resolution,
$CFGRes$,  has the following properties (cf. theorem 1.21):
\roster
\item"{(1)}" the (cover) graded formal resolution, $CFGRes$,  is a strict well-structured resolution over free products. In particular
it maps non-trivial elliptic elements  to non-trivial elliptic elements. Furthermore, 
 for each index $i$, $i=1,\ldots,s$, there is an epimorphism of limit groups over free products, $\tau_i: CF_i \to FL_i$. The epimorphisms
$\tau_i$ 
commute with the quotient maps in the two graded resolutions, $FGRes$ and $CFGRes$. 

\item"{(2)}" the epimorphisms along $CFGRes$ are proper epimorphisms.

\item"{(3)}" 
all the graded formal abelian decompositions that are associated with the various formal limit groups
(over free products), $CF_i$, $i=1,\ldots,s$, are their graded formal JSJ decompositions over free products. Furthermore, 
the graded formal JSJ decompositions (over free products) of
the formal limit groups, $CF_i$, $i=1,\ldots,s$, have the same structure as the corresponding graded formal JSJ decompositions (over free products) 
of the formal limit groups, $FL_i$, where the
difference is only in the rigid vertex groups of the JSJ decompositions  and in the elliptic factors.
 
\item"{(4)}" if no factor of $FL_i$
is (formal)  rigid nor solid, then $CF_{i+1}$ is a cover of a (formal) shortening
quotient of $CF_i$, for $i=1,\ldots,s-1$.

\item"{(5)}" for $i<s$  no factor of $FL_i$ can be  (formal) rigid. 
If for $i<s$ a factor of $FL_i$ is  (formal) solid then the corresponding factor of  $CF_i$ is 
 formal weakly solid (see definition 1.18 for a weakly solid group).
In this (solid) case, when  $i<s$,  $CF_{i+1}$ is a free product of a cover of a
(formal) flexible quotient of 
the (formal) weakly  solid factor in $CF_i$ with (possibly) the same elliptic factors and (possibly) the same free group as
in $CF_i$.

\item"{(6)}" $FL_s$, the terminal formal graded limit group of the formal graded resolution, $FGRes$, 
 is of the form $Cl(Res)*F_d$, where $Cl(Res)$ is a closure of the completion, $Comp(Res)$, and $F_d$ is a possibly
trivial free group. $Cl(Res)$ terminates in a free product of either a (non-formal) rigid or solid factor with
(possibly) finitely many elliptic factors and (possibly) a free factor.

$CF_s$, the terminal limit group of the cover formal graded resolution, $CFGRes$, is of the form $CF_s=CCl(Res)*F_d$, where
$CCl(Res)$ is a cover of the closure, $Cl(Res)$, which is a factor of $FL_s$. The cover closure, $CCl(Res)$, has the same
structure as the closure, $Cl(Res)$, and it differs from the closure $Cl(Res)$ only in the terminating limit group.
The cover closure $CCl(Res)$ terminates in a free product of a (non-formal) rigid or a weakly solid factor, 
that we denote $RSF_s$, in correspondence with
the factor of the terminal limit group of $Cl(Res)$,  with (possibly) f.p.\ covers of the elliptic factors in
the terminal limit group of $Cl(Res)$, and a free group (of the same rank as in the terminal limit group of $Cl(Res)$).

With the rigid or the weakly solid factor of the terminal limit group of $CCl(Res)$, $RSF_s$, 
we associate finitely many  limit groups (over free products), $CFl_1,\ldots,CFl_g$, 
with ungraded resolutions (over free products), $Res_1,\ldots,Res_g$, that terminate in f.p.\ limit groups (cf. theorems 1.17 and 1.20). Hence,
the completions of $Res_1,\ldots,Res_g$ are f.p.\ and the limit groups (over free products), $CFl_1,\ldots,CFl_g$, 
can be embedded in
f.p.\ completions  
(that have  f.p.\ terminal limit groups). 

\item"{(7)}" the finite collection of limit groups, $CFl_1,\ldots,CFl_g$, are all limit quotients of the f.p.\ group 
$G(y,p)*<x>$, and they
dominate all the flexible quotients of $RSF_s$
(although they need not be quotients of $RSF_s$). 

From the sequence of homomorphisms, $\{h_n\}$, 
it is possible to perform iterative shortenings, and obtain another
sequence of homomorphisms, $\{u_n\}$, that restrict to  asymptotically rigid or asymptotically 
strictly solid homomorphisms that converge into the rigid or solid factor of the terminal limit group of the closure, 
$Cl(Res)$, that appears
as a factor in a free product, $FL_s=Cl(Res)*F_d$, where $FL_s$ is
the terminal graded limit group in $FGRes$.   

The homomorphisms $\{u_n\}$ factor through $CFL_s$, and they restrict to rigid or weakly strictly solid homomorphisms of
$RSF_s$ (with respect to the cover, $CFl_1,\ldots,CFl_g$). Furthermore, in the rigid case, every flexible homomorphism of $RSF_s$ 
factors through at least one of the resolutions, $Res_1,\ldots,Res_g$, that are associated with the limit groups, 
$CFl_1,\ldots,CFl_g$.
In the weakly solid case, for every non-weakly strictly solid homomorphism of $RSF_s$,
$f: RSF_s \to A^1*\ldots*A^{\ell}$, there exists a homomorphism, $u: RSF_s \to A^1*\ldots*A^{\ell}$, such that 
the pair, $(f,u)$, extends to a homomorphism of the completion of the identity resolution, $RSF_s \to RSF_s$, and the homomorphism
$u$ as a homomorphism of the f.p.\ group, $G(y,p)*<x>$, 
factors through at least one of the resolutions, $Res_1,\ldots,Res_g$.

\item"{(8)}"   with the terminal graded limit group of the cover closure, $CCl(Res)$, that we denote, $TCCl(Res)$, 
it is possible to associate an ungraded strict well-structured
resolution:
$TCCl=V_0 \to V_1 \to \ldots \to V_t$, so that the terminal (ungraded) limit group (over free products) $V_t$ is a free product of
(possibly) f.p.\ elliptic factors and (possibly) a free group. Every non-trivial elliptic element in each of the limit groups
(over free products) $V_i$, $i=0,\ldots,t-1$, is mapped to a non-trivial element in $V_{i+1}$.

It is clearly possible to combine the cover formal graded resolution, $CFGRes$, with the ungraded resolution,  
$TCl=V_0 \to V_1 \to \ldots \to V_t$, and obtain a combined ungraded resolution over free products. With this resolution we can
naturally associate an (ungraded) completion, that we denote, $Comp_{CF}$. Then the completion, $Comp_{CF}$, and its terminal limit
group, $V_t$, are finitely presented. This implies that all the formal limit groups that appear along
the cover formal graded resolution, $CFGRes$, and in particular the terminal cover closure, $CCl(Res)$, can all be
embedded in a f.p.\ completion.
\endroster
\endproclaim

\nfp With the techniques for the construction of graded formal limit groups (over free groups) that appear in section 3
of [Se2], the proof is identical to the proof of theorem 1.21.

\line{\hss$\qed$}

Given a graded well-structured resolution, $Res(y,p)$, of a graded limit quotient of a f.p.\ group, $G(y,p)$, 
its completion, $Comp(Res)(z,y,p)$, a system
of equations, $\Sigma(x,y,p)=1$, and a sequence of homomorphisms, $\{h_n: <x,z,y,p> \to A^1_n*\ldots*A^{\ell}_n\}$ 
that restricts to
a test sequence of the completion, $Comp(Res)$, we first passed to a subsequence of the homomorphisms, $\{h_n\}$, from
which we constructed a formal graded resolution over a closure of the completion, $Comp(Res)$, and then by theorem 2.6 we
have associated (non-canonically) a cover formal graded resolution with the constructed formal graded resolution. By the
properties of the cover formal resolution, all the graded limit groups that are involved in its construction can be embedded
in f.p.\ graded completions (see theorem 2.6).

Therefore, with the graded resolution, $Res(y,p)$, and the system of equations, $\Sigma(x,y,p)=1$, theorem 2.6 associates
a collection of cover formal graded resolutions, so that every formal solution, that is associated with the resolution and the 
system of equations, factors through at least one of these cover formal graded resolutions.  
As in the ungraded
case (theorem  26  in [Ja-Se]), and as in the graded case (theorem 1.22), from this collection
of cover formal graded resolutions it is possible to find a finite subcollection, so that every formal solution that is
associated with the graded resolution, $Res(y,p)$, and the system of equations, $\Sigma(x,y,p)=1$, factors through
at least one of the cover formal graded resolutions that belong to the finite subcollection. Hence, we can
view this (non-canonical) finite subcollection as a $formal$ $graded$ $Makanin$-$Razborov$ $diagram$ of the resolution,
$Res(y,p)$, and the system of equations, $\Sigma(x,y,p)=1$.

\vglue 1pc
\proclaim{Theorem 2.7} Let $G(y,p)$ be a f.p.\ group,
let $Res(y,p)$ be a well-structured resolution over free products of a (graded)
limit quotient of $G(y,p)$ that terminates in a free product of a rigid or solid factor, 
and possibly finitely many elliptic factors, and 
a (possibly trivial) free group. Let $Comp(z,y,p)$ be the completion of $Res(y,p)$, and let $\Sigma(x,y,p)=1$ be a (finite) 
system of equations.

There exist finitely many well-structured formal graded resolutions that satisfy properties (1)-(8) of theorem 2.6:
$CF=CF_0 \to CF_1 \to \ldots \to CF_s$,
where $CF$ is a graded formal limit group over free products that is associated with the f.p.\ group $G(y,p)$ and the system of equations,
$\Sigma(x,y,p)=1$, 
and $CF_S$ is a free product $CF_s=CCl(Res)*F_d$, i.e., a free product of a cover of a closure of $Comp(Res)(z,y,p)$ 
with a (possibly trivial)
free group.

With each of the finitely many cover formal graded resolutions we associate the collection of formal solutions
that factor through it. These formal solutions, are obtained from homomorphisms of the terminal graded limit group, $TCCl(Res)$ of
the cover closure, $CCl(Res)$, which is a factor of $CF_s$. These homomorphisms of $TCCl(Res)$ 
are  obtained from 
rigid or weakly strictly solid homomorphisms of the rigid or weakly solid factor of $TCCl(Res)$,
 that are combined with homomorphisms
of the f.p.\ elliptic factors (into elliptic subgroups), and arbitrary homomorphisms of the (possible) free
factor (of $TCCl(Res)$). We further require that the restriction of weakly strictly solid homomorphisms to edge groups and to abelian
and QH vertex groups in the  graded abelian JSJ decomposition of the weakly solid factor are not (entirely) elliptic.

Formal solutions that are associated with the given cover formal graded resolution, and with the
 formal limit groups, $CF_{s-1},\ldots,CF_0=CF$, 
are obtained from  homomorphisms of $CF_s$, in a similar way to homomorphisms that factor
through graded resolutions over free and hyperbolic group, i.e., according to parts (1)-(3) that are listed in theorem 1.22.

Every formal solution that is associated with the given resolution, $Res(y,p)$, and the system of equations, 
$\Sigma(x,y,p)=1$,  factors through at least one of the
finitely many (cover) formal graded resolutions that are associated with the f.p.\ group $G(y,p)$,
the resolution, $Res(y,p)$, and  the system of equations, $\Sigma(x,y,p)=1$.
\endproclaim

\nfp Given the construction of a cover formal graded resolution that appears in theorem 2.6, the proof of theorem 2.7 is identical 
to the proof of theorem 1.22 (and to theorem 26 in [Ja-Se]).

\line{\hss$\qed$}

The finite collection of formal graded resolutions that are associated (non-canonically)
 with a f.p.\ group, $G(y,p)$,  a well-structured resolution
of a limit quotient (over free products) of it, $Res(y,p)$, and a system of equations, $\Sigma(x,y,p)=1$, according  
to theorem 2.7, and through which all the formal solutions that are associated with this triple  do factor, form
a $formal$ $graded$ $Makanin$-$Razborov$ $diagram$  over free products. Note that the collection of formal graded
resolutions is not canonical, hence, the diagram we constructed is not canonical. 

As we have already indicated in constructing the graded Makanin-Razborov diagram, the construction of the formal
graded Makanin-Razborov diagram does not require a separate study of the singular locus of the constructed
formal graded resolutions, as by construction, this singular locus is covered by other formal graded resolutions from
the constructed finite collection.

\vglue 1.5pc
\centerline{\bf{\S3. 
AE Sentences}}
\medskip

In [Ja-Se] and the first 2 sections in this paper we have generalized the results and notions that were 
presented in 
[Se1], [Se2]
and [Se3], for studying varieties, sentences and predicates defined over a free group,
to free products. In this section we show how to use these notions and constructions to reduce an AE sentence
over free products to a sentence in the Boolean algebra of AE sentences in the factors. Like all our 
constructions in the previous sections and in [Ja-Se], the reduction of an AE sentence over free products
to a sentence in the factors is uniform, and does not depend on a specific free product (although it does depend on the
number of factors in the free product).

In [Se4] we have associated a $complexity$ with each well-structured resolution, or an induced 
resolution (see section 3 of [Se4] for an induced resolution). This complexity is slightly modified in studying AE sentences over hyperbolic groups
(definition 4.2 in [Se7]). Over free products we use the same complexity of well-structured resolutions
(and their completions), as the one that is used over hyperbolic groups.

\vglue 1pc
\proclaim{Definition 3.1 ([Se4],4.2)} Let $Comp(Res)(t,y)$ be a completion of a  well-structured resolution,
$Res(y)$, over free products,
 with
 (possibly) reduced
modular groups associated with each of its various $QH$ subgroups.
Let $Q_1,\ldots,Q_m$ be the $QH$ subgroups that appear
in the completion, $Comp(Res)(t,y)$, and let $S_1,\ldots,S_m$ be the
 (punctured) surfaces associated with the reduced modular group associated with
each of the $QH$ vertex group.
To each (punctured) surface $S_j$ we may associate an ordered couple
  $(genus(S_j),|\chi(S_j)|)$.
We will assume that the $QH$ subgroups  $Q_1,\ldots,Q_m$ are ordered according to the lexicographical
(decreasing) order of the ordered couples associated with their corresponding surfaces.
Let $rk(Res(y))$ be the rank of the  free group that is dropped along the
resolution $Res(y)$, let $fact(Res(y))$ be the number of 
elliptic terminal  factors of the resolution $Res(y)$,  
and 
let $Abrk(Res(y))$ be the sum of the ranks of the kernels of the mappings of (free) abelian groups 
that appear as vertex groups along
the  resolution $Res(y)$ (see definition
1.15 in [Se4]).

\noindent
We set the complexity of the  resolution, $Res(y)$, and the completion, $Comp(Res)$, denoted
 $Cmplx(Res(y))$, to be:
$$ Cmplx(Res(y)) \, = \, (fact(Res(y))+rk(Res(y)), \
 (genus(S_1),|\chi(S_1)|),\ldots$$
$$\ldots,(genus(S_m),|\chi(S_m)|), \ 
Abrk(Res(y))).$$
On the set of complexities of completed resolutions with (possibly) reduced modular groups 
we can define a linear order.
Let $Res_1(y)$ and $Res_2(y)$ be two completed resolutions with (possibly)
reduced modular groups. 
 We say that 
$Cmplx(Res_1(y)) \, = \,  Cmplx(Res_2(y))$ if the tuples defining the two complexities
are identical. We say that
$Cmplx(Res_1)(y)) \, < \,  Cmplx(Res_2(y))$ if:
\roster
\item"{(1)}"  the "Kurosh" rank, $fact(Res_1(y))+rk(Res_1(y))$ is smaller than
the Kurosh rank
 $fact(Res_2(y))+rk(Res_2(y))$.

\item"{(2)}" the above ranks are equal and 
the tuple: $$((genus(S^1_1),|\chi(S^1_1)|),\ldots,(genus(S^1_{m_1}),|\chi(S^1_{m_1}|))$$
is smaller in the lexicographical order than the tuple:
 $$((genus(S^2_1),|\chi(S^2_1)|),\ldots,(genus(S^2_{m_2}),|\chi(S^2_{m_2}|)).$$

\item"{(3)}"  the above ranks and tuples are equal and $Abrk(Res_1(y)) \, < \, 
Abrk(Res_2(y))$. 
\endroster
\endproclaim

Given a complete, well-separated resolution, $Res(t,y)$,
over a free group, and a subgroup $<y>$ of its associated limit group, we have constructed
in section 3 of [Se4], the induced resolution, $Ind(Res(t,y))(u,y)$. The construction 
of the induced resolution generalizes directly to well-separated resolutions over free products,
hence, we omit its detailed description.

In section 4 of [Se4] a procedure for a validation of an AE sentence over a free group is presented. This
procedure is generalized to torsion-free hyperbolic groups in [Se8]. In this section our goal is to analyze
uniformly an AE sentence over all non-trivial free products. Hence, we are not aiming at validating
an AE sentence, but rather our goal is to find a uniform way to reduce an AE sentence over a free product to
a sentence over its  factors.

In the procedure that is presented in [Se4] we used a single formal solution (at times  particular families of formal
solutions), at each step of the procedure. As we look for a uniform construction, we will need to
collect all the formal
solutions that are defined over a completion (or rather on  closures of it) in each step of the procedure. Also, as we
will see in the sequel, to obtain the reduction of an AE sentence from a free product to its factors,
it is easier to analyze the places in which the iterative
procedure of the type that is presented in [Se4] fails, rather
than the places it succeeds (these possible "failures" are provided by the iterative procedure for validation of
an AE sentence over free and hyperbolic groups 
(see section 4 in [Se4]), and they are later used
in the quantifier elimination procedure over these groups).

\smallskip
Let:
$$  \forall y \ \exists x \
 \Sigma(x,y)=1 \, \wedge \, \Psi(x,y) \neq 1$$ 
be a sentence over groups. Let $F_y=<y>$, be a free group that is freely generated by (copies of) the (universal)
variables $y$. Let $G=A^1* \ldots *A^{\ell}$ be a non-trivial free product  that is not isomorphic to $D_{\infty}$. If the
given sentence is a false sentence over $G$, then there is a homomorphism $h:F_y \to G$, so that for the
corresponding values of the universal variables $y$ (i.e., the image of the variables $y$ under the homomorphism $h$), 
there exist no values for the existential variables
$x$, for which both the equalities $\Sigma(x,y)=1$ and the inequalities $\Psi(x,y) \neq 1$ hold.

We look at all the possible sequences of homomorphisms, $\{h_n:F_y \to A^1_n*\ldots*A^{\ell}_n\}$, where $\ell \geq 2$ is an arbitrary positive
integer, and the free products, $A^1_n*\ldots*A^{\ell}_n$, are non-trivial and not isomorphic to $D_{\infty}$. We further assume that each
of the homomorphisms $h_n$
fails the given AE sentence for $A^1_n*\ldots*A^{\ell}_n$, i.e., for each $h_n$
there are no values for the existential variables $x$ (in $A^1_n*\ldots*A^{\ell}_n$), so that for the tuple, $(x,h_n(y))$,
both the equalities 
$\Sigma$ and the inequalities $\Psi$ hold (in $A^1_n*\ldots*A^{\ell}_n$).

By theorem 18 in [Ja-Se] given  such a sequence, $\{h_n\}$, we can pass to a subsequence that converges into a 
well-structured  (even well-separated) resolution: $L_0 \to L_1 \to \ldots \to L_s$, where $L_s$ is
 a free product of (possibly) a free group and (possibly) finitely many elliptic factors. We denote this resolution
$BRes$. Note that the terminal elliptic subgroups in $BRes$ are f.g.\ but they may be infinitely presented.

\noindent
As the the terminal elliptic factors of the resolution, $BRes$, and its completion, may be infinitely presented, we start by iteratively
approximating it by resolutions with the same structure that have f.p.\ completions and terminal elliptic factors, approximations (or covers)
that
we used in constructing the ungraded Makanin-Razborov diagram over free products (see theorem 25 in [Ja-Se]).

Let $E_1,\ldots,E_r$ be the elliptic factors in the free decomposition of the terminal limit group of the 
resolution, $BRes$, $L_s$. $E_1,\ldots,E_r$ are all f.g. but they may be infinitely presented. Hence,
we fix a system of f.p.\ 
approximations of $E_1,\ldots,E_r$, that we denote, $E^m_1,\ldots,E^m_r$, that are obtained
from $E_1,\ldots,E_r$ by fixing a generating set of $E_1,\ldots,E_r$, and keeping only the relations of length
up to $m$ in each of the elliptic factors, $E_1,\ldots,E_r$. 

For sufficiently large index $m$, we set $Res^m$ to be the resolution (over free products) that is obtained from $BRes$,
by replacing the elliptic factors, $E_1,\ldots,E_r$, by the f.p.\ factors, $E^m_1,\ldots,E^m_r$ (note that
for sufficiently large $m$ it is guaranteed that the all the retractions in the resolution, $BRes$, lift to
corresponding retractions in the resolutions, $Res^m$). Since
the resolution $BRes$ is well-structured (and well-separated), the resolutions, $Res^m$, are well-structured
and well-separated for large $m$. Since the elliptic factors, $E^m_1,\ldots,E^m_r$, are all finitely presented,
 for each
index $m$, there exists some index $n_m$, so that for all $n>n_m$, the homomorphisms $\{h_n\}$ factor through 
the resolution $Res^m$.

In section 2 we have shown that given a resolution over free products, $Res(y)$, with a f.p.\ completion, and a
finite system of equations, $\Sigma(x,y)=1$, it is possible to  associate with them (non-canonically) a formal
Makanin-Razborov diagram that encodes all the formal solutions that are defined over (a closure of) $Res(y)$, and over every
non-trivial free product, $G=A^1*\ldots*A^{\ell}$, that is not isomorphic to $D_{\infty}$ 
(see theorem 2.5). Furthermore, each formal resolution
in the formal Makanin-Razborov diagram terminates in a closure of $Res(y)$, and this closure, as well as its terminating elliptic 
factors are all finitely presented.

\noindent
Therefore, with each approximating resolution of the constructed resolution $BRes$, $Res^m$, 
we associate (non-canonically) a formal Makanin-Razborov diagram (over free products). Note that
the completions of
the formal resolutions in these formal Makanin-Razborov diagrams are all finitely presented.

For each index $m$, there exists an index $n_m$, so that for all $n>n_m$, the homomorphisms
$\{h_n:A^1_n*\ldots*A^{\ell}_n\}$  factor through the resolution, $Res^m$. Hence, with each such homomorphism, $h_n$,
we can associate values (specializations) with each of the elliptic subgroups, $E^m_1,\ldots,E^m_r$.  

If there exists an index $m$, for which there is an infinite subsequence of homomorphisms (still denoted),
$\{h_n\}$, so that for each of the specializations of the elliptic subgroups, $E^m_1,\ldots,E^m_r$, that are
associated with the homomorphisms, $\{h_n\}$, there exists a test sequence of $Res^m$ that extends the values of the elliptic factors,
that does not extend to formal solutions over (a closure of) $Res^m$, or 
it does extend to formal solutions over (a closure of) $Res^m$, but for each such formal solution at least one of the inequalities in the
system, $\Psi(x,y) \neq 1$, does not hold (i.e., it is an equality and not an inequality),
we reached a terminal point of the iterative procedure. In this case, the associated output is the (approximating)
resolution, $Res^m$, and its associated formal Makanin-Razborov diagram.

Suppose that for every index $m$, there is no subsequence of homomorphisms $\{h_n\}$, for which for the values of the elliptic
factors of $Res^m$ that are associated with each homomorphism from the subsequence, there exists a test sequence that extends these values,
and either
there exists no
formal solution (over a closure of $Res^m$) that extends the specializations of the test sequence, or
there exist such formal solutions but for all of them the system of inequalities, $\Psi(x,y) \neq 1$ does not hold. 
In this case for each index $m$, there exists
an index $k_m > n_m > m$, so that $h_{k_m}$ factors through $Res^m$, and with the specializations of the elliptic subgroups, 
$E^m_1,\ldots,E^m_r$, that are associated with $h_{k_m}$, it is possible to associate a formal solution $x_m$ that
does  satisfy $\Sigma(x_m,y)=1$ and $\Psi(x_m,y) \neq 1$ for generic  $y$ (generic for the fiber that is associated with
the resolution, $Res(y)$, and the corresponding specializations of the elliptic factors, $E^m_1,\ldots,E^m_r$).

By construction, the elliptic factors in the resolutions, $\{Res^m\}$, $E^m_1,\ldots,E^m_r$, converge into the
elliptic factors in the resolution, $BRes$, $E_1,\ldots,E_r$ (as they are f.p.\ approximations of $E_1,\ldots,E_r$). 
The formal solutions, $\{x_m\}$, are defined
over closures of the resolutions, $\{Res^m\}$.
Using the techniques to construct formal limit groups over free products, that were presented in the previous 
section, from the sequence of formal solutions, $\{x_m\}$, it is possible to extract a subsequence (still
denoted) $\{x_m\}$, that converges into a formal limit group over a closure of the (limit) resolution $BRes$, $FL(x,z,y)$. 
By the 
construction of the formal limit group $FL(x,z,y)$, the equations from the system, $\Sigma(x,y)=1$, represent the 
trivial word in $FL(x,z,y)$, whereas each of the inequations, $\Psi(x,y) \neq 1$, represent a non-trivial 
element in $FL(x,z,y)$ (as these elements are non-trivial for the formal solutions, $\{x_m\}$, and generic values of
the variables $y$).

At this point we look at the sequences of specializations, $\{(x_m,y_m,z_1(m),\ldots,z_t(m))\}$, of the formal
solution $x_m$, the universal variables $y$, and its successive shortenings $(z_1,\ldots,z_t)$, that take their
values in the free products, $\{A^1_{k_m}*\ldots*A^{\ell}_{k_m}\}$, and each of the specializations,
$(y_m,z_1(m),\ldots,z_t(m))$, factors through the resolution, $Res^m$. The elements, $\{y_m\}$, are precisely the subsequence
of values of the universal variables for which the sentence fail to hold for the free products, $\{A^1_{k_m}*\ldots*A^{\ell}_{k_m}\}$. 

Given this sequence of specializations, we apply the first step of the procedure for validation of an
AE sentence, that is presented in section 4 of [Se4], and extract a subsequence, that converges into a
quotient resolution of the one that is associated with the formal limit group, $FL(x,z,y)$ (where the last one is a closure of the
original resolution, $BRes$). Since the formal
solutions, $\{x_m\}$, were assumed to satisfy both the equalities, $\Sigma(x,y)=1$, and the inequalities,
$\Psi(x,y) \neq 1$, the obtained quotient resolution is not a closure of the resolution, $BRes$, that we have started
with, but rather a resolution of "reduced complexity" (in the sense of the iterative procedure that is presented in section
4 in [Se4]).

We continue iteratively. At each step we start with a quotient resolution, $QRes$, that was constructed in the
previous step of the procedure, 
using the general step of the iterative procedure that is presented in section 4 in
[Se4], and a sequence of homomorphisms into free products, $\{h_n\}$, that converges into the (completion of the)
quotient resolution, $QRes$. The homomorphisms, $\{h_n\}$, that are associated with $QRes$, are constructed from a subsequence of the 
sequence of homomorphisms that are associated with the quotient resolution that was constructed in the previous step, in addition to 
elements that are associated with the formal solution that is imposed on that quotient resolution.
Furthermore, the sequence $\{h_n\}$, restricts to specializations of the universal variables $y$, that demonstrate that
the given AE sentence is a false sentence for the corresponding free
products (i.e., for these values of the universal variables there exist no specializations of the existential variables for
which both the equalities, $\Sigma$, and the inequalities, $\Psi$, do hold).  

\noindent
We start the current (general) step, by fixing a sequence of finitely presented approximations, 
$QRes_m$, of the quotient resolution,
$QRes$. Note that for each index $m$, there exists an index $n_m$, so that for every $n>n_m$, the homomorphisms
$\{h_n\}$ factor through $QRes_m$, and for each homomorphism $h_n$ there are associated specializations of the 
elliptic factors of $QRes_m$ in the free product that is associated with $h_n$. 

With each f.p.\ completion (of $QRes_m$) we associate (not in a canonical way)
a formal Makanin-Razborov diagram. Note that each formal resolution in such a formal diagram has a f.p.\
completion. If there exists an index $m$, and a subsequence of the given sequence of homomorphisms, $\{h_n\}$,
that we still denote $\{h_n\}$,
so that over the specializations of the elliptic factors of $QRes_m$ that are associated
with the subsequence of homomorphisms, $\{h_n\}$,  there exists a test sequence that does not extend to any formal solutions over $QRes_m$,
 so that 
both the equalities, $\Sigma(x,y)=1$, and the inequalities, $\Psi(x,y) \neq 1$, hold for generic $y$'s (over the corresponding sequence
of free products, $\{A^1_n*\ldots*A^{\ell}_n\}$), we
reached a terminal point of our iterative procedure. In this case the final output of the iterative
procedure is the quotient resolution, $QRes_m$, (that has a f.p.\ completion), and its (non-canonical)
formal Makanin-Razborov diagram, where each  formal resolution in this diagram has a f.p.\ completion as well,
and this formal Makanin-Razborov diagram encodes all the formal solutions over $QRes_m$ that satisfy the
equalities, $\Sigma(x,y)=1$. 

Suppose that there is no such index $m$, and no such subsequence of the homomorphisms, $\{h_n\}$. In this case
for each index $m$, there exists some index $k_m>m$, so that there is a formal solution $x_m$ that
is defined over the resolution, $QRes_m$, with  values of its elliptic factors that are associated
with the homomorphism, $h_{k_m}$, so that both the equalities, $\Sigma(x_m,y)=1$, and the inequalities,
$\Psi(x_m,y) \neq 1$, hold for generic values of the universal variables $y$ (generic in the corresponding fibers of the resolution $QRes_m$,
and the values of its terminal elliptic factors that are associated with the homomorphisms, $\{h_{k_m}\}$). 
In this case we look at the sequence,
$\{(x_m,h_{n_m})\}$, that has a subsequence that converges into a quotient resolution, that is 
obtained using the general step of the iterative procedure for the analysis of an AE sentence that
is presented in section 4 of [Se4]. Since the formal solutions are guaranteed to satisfy both the
equalities, $\Sigma(x_m,y)=1$, and the inequalities, $\Psi(x_m,y) \neq 1$, the complexity of the
obtained quotient resolution is strictly smaller than the complexity of the resolution, $QRes$, that we have 
started the current step with.

By theorem 4.12 in [Se4], this iterative procedure terminates after finitely many steps. If we reached a
terminal quotient resolution along the process, we found a f.p.\ resolution, $QRes_m$,
and a subsequence
of the original sequence of homomorphisms, for which for the specializations of the elliptic factors of
$QRes_m$ that are associated with the subsequence of homomorphisms (these elliptic factors take their values in  the free products,
$\{A^1_n*\ldots*A^{\ell}_n\}$, that are associated with the homomorphisms, $\{h_n\}$), 
there exist test sequences that do not extend to  formal solutions that
are defined over (closures of)  $QRes_m$, and for which both the equalities and the inequalities hold for generic value
of the universal variables $y$ (generic in the fibers that are associated with $QRes_m$ and the corresponding values of the elliptic factors
of $QRes_m$ in $A^1_n*\ldots*A^{\ell}_n$).

If we didn't find such a f.p.\ resolution along the iterative procedure, it continued until we reached 
a quotient resolution which is a free product of elliptic factors (as by theorem 4.12 in [Se4] the iterative procedure terminates
after finitely many steps). Once again, we look at a sequence
of f.p.\ approximations of these elliptic factors. These f.p.\ approximations are free products of elliptic
factors as well. 

Since the given sequence of homomorphisms $\{h_n\}$, testifies that the given sentence is not valid
for the free products, $\{A^1_n*\ldots*A^{\ell}_n\}$, there must exist a f.p.\ approximation of the terminal elliptic factors,
and a subsequence of the given sequence of homomorphisms, $\{h_n\}$,
that we still denote $\{h_n\}$,
so that over the specializations of the elliptic factors of the f.p.\ approximation that are associated 
with the subsequence of homomorphisms, $\{h_n\}$,  no formal solution over the f.p.\ approximating free 
product
can be constructed so that 
both the equalities, $\Sigma(x,y)=1$, and the inequalities, $\Psi(x,y) \neq 1$, hold for the corresponding value of
the universal variables $y$.

Therefore, whatever terminal point of the iterative procedure we have reached, we found a resolution (over free
products) with f.p.\ completion, through which a subsequence of the given sequence of homomorphisms do factor,
and for which for some test sequences of this resolution that are associated with 
those values of the elliptic factors of the  terminal limit group of the resolution,
that are associated with the given subsequence of homomorphisms, the specializations in these test sequences can not be extended to 
formal solutions (over  closures  of the resolution) that satisfy both the equalities, $\Sigma(x,y)=1$, and the 
inequalities, $\Psi(x,y) \neq 1$,
for generic values of the variables $y$.

This procedure that starts with a given AE sentence, and with a sequence of values of the universal variables $y$ 
of that sentence, that testifies that the
given AE sentence fails (over the  free products in which the universal variables take their values), and extracts 
a subsequence that factors through
a resolution (over free products) with a f.p.\ completion,    is what we need in order to reduce a sentence uniformly from
free products to their factors. It 
enables one to associate with a given AE sentence, a finite collection of resolutions with f.p.\
completions, so that if the given sentence is false for a non-trivial free product, $A^1*\ldots*A^{\ell}$ (which is not $D_{\infty}$), 
it must be false for a generic 
point of one of the constructed resolutions, where the elliptic factors take some value in $A^1*\ldots*A^{\ell}$. This is the key
for reducing an AE sentence over a free product to a sentence over the factors of the free product, and a 
key for our approach to the analysis of general sentences and predicates over free products.

\vglue 1pc
\proclaim{Theorem 3.2} Let:
$$  \forall y \ \exists x \
 \Sigma(x,y)=1 \, \wedge \, \Psi(x,y) \neq 1$$ 
be a sentence over groups. Then there exist finitely many resolutions over free products:
 $Res_1(z,y),\ldots,Res_d(z,y)$ with the following properties:

\roster
\item"{(1)}" the completion of each of the resolutions, $Res_i(z,y)$, is finitely presented.

\item"{(2)}" with each resolution, $Res_i(z,y)$, we associate (non-canonically) its formal Makanin-Razborov
diagram over free products with respect to the system of equations, $\Sigma(x,y)=1$. Every resolution in these formal Makanin-Razborov diagrams
has a f.p.\ completion.

\item"{(3)}" for every non-trivial free product, $A^1*\ldots*A^{\ell}$, which is not $D_{\infty}$, and for which the given
sentence is false over $A^1*\ldots*A^{\ell}$, there exists
an index $i$, $1 \leq i \leq d$, and a specialization of each of the terminal elliptic  factors of $Res_i(z,y)$ in elliptic factors
in the free products, $A^1*\ldots*A^{\ell}$,
so that for these specializations of the elliptic factors of $Res_i(z,y)$, there exists a  test sequence that is associated with these 
values of the elliptic factors, that does not extend to
formal solutions over (closures of) $Res_i(z,y)$, for which both the equalities, $\Sigma(x,y)=1$, and the inequalities,
$\Psi(x,y) \neq 1$, hold for generic values of $y$. 
\endroster

\noindent
In other words,
there exists a finite collection of resolutions (with f.p.\ completion) over free products,
so that the failure of an AE sentence over general free products can be demonstrated by the lack of the existence
of a formal solution over a (generic point in a) fiber of at least one of these resolutions.
\endproclaim

\nfp Let $A^1_n,\ldots,A^{\ell}_n$, for some $\ell>1$,  be a sequence of non-trivial free products, 
that are not isomorphic to $D_{\infty}$, over which the given AE sentence is false.  Let
$\{y_n\}$ be a sequence of specializations 
of the universal variables $y$, in the free products, $A^1_n,\ldots,A^{\ell}_n$,  that fail the given AE sentence over 
these free products, i.e., the existential sentences (with coefficients) over the free products, $A^1_n*\ldots*A^{\ell}_n$:
$$   \exists x \
 \Sigma(x,y_n)=1 \, \wedge \, \Psi(x,y_n) \neq 1$$ 
are false. 

Starting with the sequence, $\{y_n\}$, the terminating iterative procedure that we have presented
constructs a resolution, $VRes$, with the following
properties:
\roster
\item"{(i)}" the resolution $VRes$ has a f.p.\ completion and terminal limit group.

\item"{(ii)}" there exists a subsequence of the sequence of specializations, $\{y_n\}$, that extend to specializations 
that factor through the resolution, $VRes$. Hence, with each specialization from this subsequence of the sequence, $\{y_n\}$, 
specializations of the f.p.\ terminal limit group of the resolution, $VRes$, can be associated. For each such specialization of the terminal
limit group of $VRes$, there exists an associated subsequence, that does not extend to any
formal solutions (over closures of $VRes$)  that satisfy both the equalities, $\Sigma(x,y)=1$, and the inequalities, $\Psi(x,y) \neq 1$,
for generic values of the variables $y$.

\item"{(iii)}" with the resolution $VRes$, and the system of equations, $\Sigma(x,y)=1$, we associate (non-canonically) a
formal Makanin-Razborov diagram (see section 2). Every (formal) resolution in this diagram has a f.p.\
completion and terminal limit groups as well.
\endroster

At this point we are able to apply the argument that we used in constructing the ungraded and graded Makanin-Razborov diagrams
(theorems 26 in [Ja-Se] and 1.22 in this paper). We look at all the non-trivial free products, $A^1*\ldots*A^{\ell}$
(possibly for varying $\ell>1$), that are not isomorphic to $D_{\infty}$, over which the given AE sentence is false.
We further look at sequences of specializations of the universal variables $y$,  $\{y_n\}$, that testifies that the given
AE sentence fail over the corresponding sequence of non-trivial free products (that are not $D_{\infty}$),
$A^1_n*\ldots*A^{\ell}_n$. From every such  sequence we use our terminating iterative procedure, and extract a subsequence
of the specializations, $\{y_n\}$, and a resolution,  
$VRes$, that has the properties (i)-(iii), and in particular,   the subsequence of specializations, (still denoted) $\{y_n\}$, 
extend to specializations that factor through the resolution $VRes$,
and no formal solutions that satisfy both the equalities, $\Sigma(x,y)=1$, and the inequalities, $\Psi(x,y) \neq 1$,
for generic values (i.e., some test sequences) of the variables $y$ can be constructed, 
for those values of the elliptic factors in the terminal 
limit group of $VRes$, that are associated with the subsequence of specializations of the sequence, $\{y_n\}$.

The completion of each of the constructed resolutions, $VRes$, is finitely presented, and so are the resolutions in its 
associated formal Makanin-Razborov diagram. Hence, we can define a linear  order on this (countable) collection of
resolutions ($VRes$), and their (non-canonically) associated formal Makanin-Razborov diagrams. By the same argument that
was used in constructing the Makanin-Razborov diagram (theorem 26 in [Ja-Se]), there exists a finite subcollection
of these resolutions that satisfy properties (1)-(3) of the theorem.

\line{\hss$\qed$}

Theorem 3.2 constructs a finite collection of resolutions over free products, with f.p.\ completions, that
demonstrate the failure of an AE sentence over general free products. This is precisely what is required
in order to reduce an AE sentence over  a free product to a sentence which is in the Boolean algebra of AE sentences over 
its factors, which is a special
case of our general goal. Before stating this reduction of a sentence from a free product to its factors, we need to study the
$singular$ $locus$ of a resolution over free products.
  
Recall that in constructing the graded and ungraded Makanin-Razborov diagrams over free products, and the formal
Makanin-Razborov diagram, we have considered the set of specializations that 
factor through a resolution in one of these diagrams, as those homomorphisms that factor through the resolutions (in the sense that
they are  obtained from specializations of
the terminal limit group of the resolution by iteratively precomposing these homomorphisms with automorphisms from the modular groups
(over free products) that are associated with the various levels), and in addition we required that the associated specializations of
the various limit groups along the resolution restrict to non-elliptic specializations of all the abelian edge (and vertex) groups, and all
the QH vertex groups along the resolution. 

This assumption on the non-ellipticity of abelian edge groups, and QH vertex groups, 
allows us to ignore the singular locus in studying (graded) resolutions and
the formal resolutions that are defined over them. However, it is essential to determine the singular locus of the resolutions,
$Res_1,\ldots,Res_d$, that were constructed in theorem 3.2, 
as specializations that do factor through the singular locus are not being considered in further analyzing a given
resolution, in order to reduce the AE sentence over the free product $G$ to an AE sentence over the factors of $G$,
and this singular locus needs to be defined uniformly (for all free products that are not isomorphic to $D_{\infty}$).

\vglue 1pc
\proclaim{Definition 3.3} Let $Res(y)$ be a well-structured, coefficient free resolution over free products. Let $G=A^1*\ldots*A^{\ell}$
be a non-trivial free product that is not isomorphic to $D_{\infty}$. A specialization of the terminal limit group of $Res(y)$ in $G$
is said to be in the $singular$ $locus$ of $Res(y)$ over the free product $G$ if it does not extend to a test sequence over $G$, 
or equivalently,
if for every  specialization of  $Res(y)$ in $G$ that is obtained from the   given specialization of the terminal limit group
of $Res(y)$,  by a finite sequence of automorphisms in the modular groups that are associated with the various levels
of $Res(y)$, the associated specialization of at least one of the edge groups in the virtually abelian decompositions that
are associated with the various levels of $Res(y)$ is elliptic or trivial, or the specialization of
at least one of the QH vertex groups in these abelian decompositions is elliptic or trivial or virtually abelian. 
\endproclaim 

To obtain a (uniform) reduction of AE sentences, we need a uniform description of the singular locus of a given resolution (over 
all free products).

\vglue 1pc
\proclaim{Proposition  3.4} Let $Res(y)$ be a well-structured, coefficient free
resolution over free products with a f.p.\ completion (and terminal limit 
group). Then there exist finitely many graded resolutions with f.p.\ completions and terminal limit groups, 
$SLRes_1(y),\ldots,SLRes_u(y)$, so that for every non-trivial free product,
$G=A^1*\ldots*A^{\ell}$, that is not isomorphic to $D_{\infty}$, a specialization of the terminal limit group of $Res(y)$ is in the singular
locus of $Res(y)$, if and only if it extends to a specialization of (at least) one of the terminal  limit groups of the resolutions,
$SLRes_1(y),\ldots,SLRes_u(y)$.
\endproclaim

\nfp Let the resolution, $Res(y)$, be given by the sequence of epimorphisms: $L_0 \to L_1 \to \ldots \to L_s$.
Given an edge or a QH vertex in the virtually abelian decomposition that is associated with level $i$ of the resolution $Res(y)$, 
$0 \leq i \leq s-1$,
we look at the partial resolution, $Res_{i+1}$, that is given by: $L_{i+1} \to L_{i+2} \to \ldots \to L_s$. Since $Res(y)$ is
well-structured and coefficient-free, so is $Res_{i+1}$. Since $L_s$, the terminal limit group of $Res(y)$, is finitely presented,
so is the completion of $Res_{i+1}$. 
 
An edge group or a QH vertex group in the virtually abelian decomposition that is associated with $L_i$, is mapped into $L_{i+1}$.
Given the edge or a QH vertex group in the abelian decomposition that is associated with $L_i$, $0 \leq i \leq s-1$,
we look at all the test sequences of the completion of the resolution, $Res_{i+1}$, 
over arbitrary non-trivial free products that are not $D_{\infty}$,
for which the image of the  given edge group is elliptic or trivial or the image of the given QH vertex group in $L_{i+1}$ is either 
elliptic or trivial or virtually abelian.  

Every such test sequence subconverges into a limit quotient $U$ (over free products)
of the completion, $Comp(Res_{i+1})$, in which the image of the
given edge group is elliptic or trivial, or the image  of the given QH vertex group in $L_{i+1}$ is elliptic or trivial or
virtually abelian. As the limit quotient $U$ is obtained as a limit of
a test sequence of specializations of $Comp(Res_{i+1})$, $U$ has the same structure as the completion, $Comp(Res_{i+1})$, 
but its terminal limit group
is a quotient of the terminal limit group of $Res(y)$, $L_s$. Note that the terminal limit group of $L_s$ is f.g.\ but it
may be infinitely presented.

On the set of limit quotients that are obtained from convergent sequences of test sequences of $Res_{i+1}$ in which the image of
the  given edge group is trivial or elliptic, or the the image of the given QH vertex group in $L_{i+1}$ is either elliptic or trivial
or virtually abelian, we can naturally define a partial 
order, that is identical to the partial order that is defined on limit quotients (over free products) of a given f.g.\ group
in definition 12 in [Ja-Se]. Since the completion, $Comp(Res_{i+1})$, is finitely presented, and we require that finitely many more elements 
(the image of a generating set of an edge group or a QH vertex group in $L_{i+1}$) are elliptic or trivial, or that the image of a QH
vertex group is virtually abelian (which means that it has an abelian subgroup of index 2 in our case), every limit quotient $U$ that
is obtained as a limit of a test sequence of $Res_{i+1}$ in which the image of
the  given edge group is trivial or elliptic, or the image of the given QH vertex group in $L_{i+1}$ is either elliptic or trivial
or virtually abelian, is dominated (under the natural
partial order) by such a limit quotient
which is finitely presented.

Therefore, by the same argument that was used to prove the finiteness of maximal limit quotients (over free products) of
a f.p.\ group (theorem 21 in [Ja-Se]), with a given edge or a QH vertex in the abelian decomposition that is associated
with $L_i$ in $Res(y)$, $0 \leq i \leq s-1$, we can associate finitely many limit quotients $U_1,\ldots, U_d$, that are all
f.p.\ and do all have the same structure as $Res_{i+1}$, such that every test sequence of $Res_{i+1}$ for which
the image of the  given edge group is trivial or elliptic, or the image of the given QH vertex group in $L_{i+1}$ is either 
elliptic or trivial or virtually abelian, has a subsequence that
factors through one of the limit groups, $U_1,\ldots,U_d$. 

Hence, by taking the finite union of these limit quotients for all the (finitely many) edge and QH vertex groups in the abelian
decompositions that are associated with all the limit groups, $L_i$, $0 \leq i \leq s-1$, we constructed a finite set of limit quotients,
that are all completions, and together they satisfy the conclusion of the proposition.

\line{\hss$\qed$}

To reduce an AE sentence from the ambient free group to its factors, we also need to construct auxiliary resolutions, that will enable
one to decide uniformly if over a given free product $G$, the set of specializations that do factor through a given
resolution are covered by a given finite collection of closures. Following the construction of the $Root$ resolutions (over free groups)
in sections 1 and 3
of [Se5], we associate a finite collection of $Root$ $resolutions$ with a given resolution and a finite collection of its closures (over
free products).     

\vglue 1pc
\proclaim{Proposition  3.5} Let $Res(y)$ be a well-structured, coefficient free
resolution over free products with a f.p.\ completion (and terminal limit 
group), and let $Cl_1(Res),\ldots,Cl_f(Res)$ be a given finite set of closures of $Res(y)$ (see definition 2.3). 

Each  abelian vertex group that appears in the abelian decompositions that are associated with the various levels of
$Res(y)$ is contained as a finite index subgroup in an abelian vertex group  in each of the closures,
$Cl_1(Res),\ldots,Cl_f(Res)$. 
We set $Pind$ to be the product of all  the indices of the abelian vertex groups of $Res(y)$ in the corresponding abelian vertex
groups in the closures, $Cl_1(Res),\ldots,Cl_f(Res)$. 

Let $pg_1,\ldots,pg_t$ be the pegs of pegged abelian groups in $Res(y)$. With each collection of positive integers, $d_1,\ldots,d_t$,
that do all divide $Pind$, we associate a finite (possibly empty) collection of closures of $Res(y)$ with f.p.\ completions and
terminal limit groups. We denote the entire
finite collection of such closures, $\{RootRes_r\}$, and call them $Root$ $resolutions$.

Let $G=A^1*\ldots*A^{\ell}$ be a non-trivial free product that is not isomorphic to $D_{\infty}$. We fix a specialization of the terminal
limit group of $Res(y)$ in $G$. The pegs of abelian vertex groups in $Res(y)$, $pg_1,\ldots,pg_t$, have roots of orders
$d_1,\ldots,d_t$, for every test sequence (in $G$) over the given specialization of the terminal limit group of $Res(y)$, if and
only if the specialization of the terminal limit group of $G$ can be extended to at least one of the Root resolutions, $RootRes_r$,
that are associated with the sequence, $d_1,\ldots,d_t$.   
\endproclaim

\nfp Let $d_1,\ldots,d_t$ be a tuple of integers that divide $Pind$. We look at the collection of all the test sequences of $Res(y)$,
over all possible non-trivial 
free products that are not isomorphic to $D_{\infty}$, for which the pegs, $pg_1,\ldots,pg_t$, have roots of orders
$d_1,\ldots,d_t$. By the construction of the formal Makanin-Razborov diagram over free products in section 2, with this collection
of test sequences, it is possible to associate (non-canonically) a finite collection of closures of $Res(y)$, that do all have f.p.\
completions. This finite collection of closures has the property, that every test sequence of $Res(y)$ (over an arbitrary 
non-trivial free product that is not isomorphic to $D_{\infty}$), for which the pegs in $Res(y)$, $pg_1,\ldots,pg_t$, have
roots of orders, $d_1,\ldots,d_t$, in correspondence,  has a subsequence that factors through one of the closures
from the finite collection. 

By taking the finite union of the finite collections of closures for every tuple of integers, $d_1,\ldots,d_t$, for which each $d_i$ 
divides $Pind$, we get a (non-canonical) finite collection that satisfies the conclusion of the theorem.  

\line{\hss$\qed$}

Proposition 3.5 associates Root resolutions with a given coefficient-free 
resolution over free products, and a finite collection of closures of this resolution.
Using these Root resolutions, it is possible to write a sentence, which is a disjunction of conjunctions of sentences over the factors
of a free product, that determines (uniformly) if for a given free product the finite collection of closures forms a covering
closure of the given resolution. 

\vglue 1pc
\proclaim{Corollary  3.6} With the notation of proposition 3.5, and for every positive integer $\ell>1$,
the set of closures, $Cl_1(Res),\ldots,Cl_f(Res)$, is a covering 
closure of the resolution, $Res(y)$, 
over a non-trivial free product, $G=A^1*\ldots*A^{\ell}$, that is not isomorphic to $D_{\infty}$  (see definition 2.4 for
a covering closure), if and only if a (finite) conjunction of disjunctions of 
coefficient free AE sentences over the factors, $A^1,\ldots,A^{\ell}$ is a truth sentence. 
\endproclaim

\nfp The given closures, $Cl_1(Res),\ldots,Cl_f(Res)$, do not form a covering closure of the resolution, $Res(y)$, over a
free product,  
$G=A^1*\ldots*A^{\ell}$, that is not isomorphic to $D_{\infty}$, if and only if there exist  specializations of the elliptic
factors in the terminal
limit group of $Res(y)$, for which:
\roster
\item"{(i)}" the specializations  extend to specializations of the elliptic factors in the terminal limit groups of 
(one of finitely many possible) prescribed subset of  the associated Root 
resolutions (see proposition 3.5) and to any of the terminal limit groups of the Root resolutions in the complement of the 
prescribed set.

\item"{(ii)}" the specializations 
 do not extend to specializations of the elliptic factors in the terminal limit groups of 
(one of finitely many possible) prescribed subset of  the closures, $Cl_1(Res),\ldots,Cl_f(Res)$, that given condition (i) on the
Root resolutions, cover the corresponding  fiber of $Res(y)$ that is associated with the specializations of the elliptic  factors. 
\endroster

By counting all the possible finite subsets of Root resolutions, and all the possible finite subsets of closures,
$Cl_1(Res),\ldots,Cl_f(Res)$, that satisfy part (ii), given the chosen finite subset of Root resolutions, the set
of closures, $Cl_1(Res),\ldots,Cl_f(Res)$,
do not form a covering closure of $Res(y)$ over $G$, if and only if a (finite) disjunction of conjunctions of $EA$ sentences over the 
factors $A^1,\ldots,A^{\ell}$ is a truth sentence. Hence, the converse is a finite conjunction of disjunctions  of AE sentences
over the factors, $A^1,\ldots,A^{\ell}$.

\line{\hss$\qed$}

Proposition 3.4 on the structure of the singular locus of a resolution, the construction of Root resolutions in proposition
3.5, together with theorem 3.2 that associates with a given AE sentence
over free products, finitely many resolutions,
and their formal Makanin-Razborov diagrams, enable us to reduce (uniformly) a given AE sentence from free products to their factors.

\vglue 1pc
\proclaim{Theorem 3.7} Let:
$$  \forall y \ \exists x \
 \Sigma(x,y)=1 \, \wedge \, \Psi(x,y) \neq 1$$ 
be a sentence over groups. 

Then for every positive integer, $\ell>1$, there exists a coefficient-free 
sentence over free products, 
which is a disjunction of conjunctions of AE sentences, where each of these last AE sentences involves elliptic 
elements from the same factor, $A^1,\ldots,A^{\ell}$,
such that for every non-trivial free product, $G=A^1*\ldots*A^{\ell}$,  which is not $D_{\infty}$, 
the original
AE sentence over the free product $G=A^1*\ldots*A^{\ell}$ is a truth sentence, if and only if the sentence which is a 
disjunction of conjunctions of AE sentences over the factors $A^1,\ldots,A^{\ell}$  is a truth sentence.
\endproclaim 

\nfp By theorem 3.2, with the given AE sentence, it is possible to associate finitely many well-structured resolutions,
$Res_1,\ldots,Res_d$, that do all have f.p.\ completion, and terminal limit groups.
With each resolution, $Res_i$, and the system of equations, $\Sigma(x,y)=1$,
it is possible to associate (non-canonically) its formal Makanin-Razborov diagram. We denote each resolution in the formal
Makanin-Razborov diagram of $Res_i$, $FRes_i^j$. Since the completion of $Res_i$ is f.p.\ each formal resolution, $FRes_i^j$, 
has a f.p.\ completion and terminal
limit group.

By the construction of the formal resolutions, the elements that correspond to the equations in the system, $\Sigma(x,y)=1$, represent
the identity element in the completion, $Comp(FRes_i^j)$. We further look at all the test sequences of each of the formal
resolutions, $FRes_i^j$, for which at least one of the words in the system of inequalities,
$\Psi(x,y) \neq 1$, holds as equality. By the construction of the formal Makanin-Razborov diagram,
this associates a finite collection of  resolutions  with each of the formal resolutions, $FRes_i^j$, that do all have  f.p.\ completion and
terminal limit groups. In each of these resolutions, at least one of the words in the system of inequalities,
$\Psi(x,y) \neq 1$ is the trivial element. We denote each of these resolutions, $\Psi FRes_i^{j,k}$.

By proposition 3.4, with the resolutions, $Res_1,\ldots,Res_d$, that are associated with the given AE sentence, it is possible to associate
finitely many resolutions, $SLRes_i^v$, that have  f.p.\ completions and terminal limit groups. A specialization of the 
terminal limit group of $Res_i$
is in the singular locus, if and if only this specialization extends to a specialization of the terminal limit group of one of the
resolutions, $SLRes_i^v$. Similarly we associate with each formal resolution, $FRes_i^j$, finitely many formal resolutions that are
associated with its singular locus, and we denote of these resolutions, $SLFRes_i^{j,v}$.

Each of the formal resolutions, $FRes_i^j$, terminates with a closure of the resolution $Res_i$. Hence with the resolution, $Res_i$, we can
associate finitely many closures, $Cl(Res_i)^j$. By proposition 3.5, with the resolution, $Res_i$, and its closures, $Cl(Res_i)^j$, we can 
associate finitely many Root resolutions, that we denote $RootRes_i^r$, that collect all the specializations of the terminal
limit group of $Res_i$, for which the pegs in the resolution, $Res_i$, have roots of prescribed orders for every test sequence that extends the 
given specialization of the terminal limit group of $Res_i$.

Therefore, by theorem 3.2 and propositions 3.4 and 3.5, the given AE sentence is false over a non-trivial free product, $G=A^1*\ldots*A^{\ell}$,
that is not isomorphic to $D_{\infty}$, if and only if
there exists a specialization in $G$ of  the terminal limit group of a resolution
$Res_i$ (one of the resolutions, $Res_1,\ldots,Res_d$, that were associated with the sentence by theorem 3.2), i.e.,
specializations of the elliptic factors of the terminal limit group of $Res_i$ in the factors, $A^1,\ldots,A^{\ell}$, for which:
\roster
\item"{(1)}" the specialization does not extend to a specialization of the terminal limit group of one of the resolutions,
$SLRes_i^v$ (the singular locus of $Res_i$).

\item"{(2)}" the specialization extends to the terminal limit group of some prescribed Root resolutions (possibly only the trivial roots), 
$RootRes_i^r$, 
and not to other ones. 

\item"{(3)}" the specialization does not extend to 
a (finite) collection of specializations of  the formal resolutions, $FRes_i^j$, that are associated
with the resolution $Res_i$, that satisfy:

\itemitem{(i)} these specializations are not  in the singular locus of $FRes_i$, i.e., they don't extend to specializations of the terminal
limit group of one of the resolutions, $SLRes_i^{j,v}$.

\itemitem{(ii)} these specializations do not  extend to specializations of the terminal limit group of  any of the
resolutions, $\Psi FRes_i^{j,k}$. 

\itemitem{(iii)} the fibers that are associated with these specializations and the corresponding formal resolutions, $FRes_i^j$, form a 
covering 
closure of the fiber that is associated with $Res_i$ and the given specialization of its terminal limit group. 
\endroster

Finally, the existence of such a specialization of the terminal limit group of one of the resolutions, $Res_i$, that satisfies
properties (1)-(3), is clearly a finite disjunction of finite conjunctions of EA sentences over the factors of $G$, $A^1,\ldots,A^{\ell}$. Hence,
the AE sentence, which is exactly the negation of this sentence is a finite disjunction of finite conjunctions of AE sentences over the
factors, $A^1,\ldots,A^{\ell}$.

\line{\hss$\qed$}

\vglue 1.5pc
\centerline{\bf{\S4. 
EAE sentences and predicates}}
\medskip

In  the previous section we used the iterative procedure that is presented in section 4 of [Se4]
to reduce AE sentences from the ambient free product to its factors. In this section we use the
procedures that were used to prove quantifier elimination over free groups, to analyze EAE predicates
and sentences over free products. 

Let:
$$  EAE(p) \ = \ \exists w \forall y \ \exists x \
 \Sigma(x,y,w,p)=1 \, \wedge \, \Psi(x,y,w,p) \neq 1$$ 
be a predicate over groups. 
Let $G=A^1*\ldots*A^{\ell}$ be a non-trivial free product (for some $\ell>1$) that is not isomorphic to $D_{\infty}$.
If $p_0 \in EAE(p)$ over the free product $G$, then 
there exists some specializations $w_0$ of the existential variables $w$, so that the AE sentence:
$$   \forall y \ \exists x \
 \Sigma(x,y,w_0,p_0)=1 \, \wedge \, \Psi(x,y,w_0,p_0) \neq 1$$ 
is a truth sentence over $G$.

In the previous section we studied AE sentence over free products, and we have finally shown that an AE sentence over free
products is a truth sentence over a non-trivial free product $G=A^1*\ldots *A^{\ell}$ (that is not isomorphic to $D_{\infty}$) if and only if 
a disjunction of conjunctions of sentences over the factors $A^1,\ldots,A^{\ell}$ does hold. 

\noindent
To prove a similar reduction from a given EAE sentence over free products to a (finite) disjunction of (finite) conjunctions
of sentences over the factors of the free product,  
we use the procedure for quantifier elimination over free groups that was presented in [Se5] and [Se6]. Recall that the
analysis of an EAE set over a free group is divided into two steps. In the first step one analyzes uniformly the AE 
sentences:
$$   \forall y \ \exists x \
 \Sigma(x,y,w,p)=1 \, \wedge \, \Psi(x,y,w,p) \neq 1$$ 
(uniformly in the pair $(w,p)$, which is equivalent to the analysis of the corresponding AE set, $AE(w,p)$). Then one
uses this uniform analysis of AE sentences, and the (iterative) sieve procedure that is presented in [Se6], to analyze an
EAE set over a free group. 

We start with the uniform analysis of AE sentences, which is equivalent to the analysis of an AE set, that combines
what we did in the previous section with the iterative procedure that is presented in section 2 in [Se5]. Let:
$$AE(w,p) \ = \   \forall y \ \exists x \
 \Sigma(x,y,w,p)=1 \, \wedge \, \Psi(x,y,w,p) \neq 1$$ 
be an AE set. Let $G=A^1*\ldots*A^{\ell}$ be a non-trivial free product that is not isomorphic to $D_{\infty}$, and suppose that
$(w_0,p_0)$ is a specialization of the free variables, $(w,p)$, in the free product $G$, 
that is not in the definable set $AE(w,p)$ over the free
product $G$. Then there exists a specialization $y_0$ of the universal variables $y$ in $G$, for which for all the possible 
values of the existential variables $x$ in $G$ either the equalities, $\Sigma(x,y_0,w_0,p_0)=1$, do not hold, or at least one
of the inequalities, $\Psi(x,y_0,w_0,p_0) \neq 1$, does not hold in $G$.

We continue by combining the construction of the diagram that is associated with an AE sentence over free products, that
was constructed in the previous section, with the procedure for the analysis of AE sets (over free groups) that is presented in section
2 of [Se5].
We start by looking at all the sequences of specializations of tuples, $\{(p_n,w_n)\}$, that take their values in non-trivial free
products, $G_n=A^1_n*\ldots*A^{\ell}_n$, that are not isomorphic to $D_{\infty}$, specializations of the tuple $(p,w)$ 
that are not in the definable set, $AE(w,p)$. This implies that for every
pair, $(w_n,p_n)$, in the sequence, there exists a specialization of the universal variables, $y_n$, in the free product,
$G_n$, 
for which there are no values for the existential variables $x$ (in $G_n$) so that both the equalities 
$\Sigma(x,y,w,p)=1$ and the inequalities $\Psi(x,y,w,p) \neq 1$ hold.

By proposition 1.16,
 given  such a sequence of triples, $\{(p_n,w_n,y_n)\}$, we can pass to a subsequence that converges into a 
well-structured  (even well-separated) graded resolution with respect to the parameter subgroup $<w,p>$, that we denote $GRes$: 
$$L_0 \to L_1 \to \ldots \to L_s$$ where $L_s$ is
 a free product of a rigid or a solid limit group over free products (that contains the subgroup, $<w,p>$), 
and (possibly) a free group and (possibly) finitely many elliptic factors. 
Note that the terminal limit group $L_s$   of $GRes$ is f.g.\ but it may be infinitely presented.

As the terminal limit group $L_s$ of $GRes$ is f.g.\ but perhaps not f.p.\ we  
fix a sequence of  approximating covers of the resolution $GRes$, that we denote, $\{CGRes_m\}$.
These approximating covers  are constructed according to the construction of a cover approximating resolution that appears in theorem 1.21, 
hence, the
cover approximating resolutions, $\{CGRes_m\}$, satisfy all the properties 
that are listed in theorem 1.21. In particular, the graded completions of the covers, $CGRes_m$, can be embedded into
f.p.\ completions, and if the terminal limit group $L_s$ of $GRes$ is a free product of a rigid (solid) factor with elliptic
and free factors, then the terminal limit group of the covers, $CGRes_m$, is a free product of a rigid (weakly solid) factor
with (f.p.\ approximating) elliptic and free factors. With each of the graded covers, $CGRes_m$, 
there is an associated cover of its flexible quotients, and each cover (of flexible quotients)
from this finite collection can be embedded into a f.p.\ completion. Also, for each cover, $CGRes_m$, there exists a subsequence
of the triples, $\{(p_n,w_n,y_n)\}$, that factor through it, and restrict to rigid or weakly strictly solid specializations
of the rigid or weakly solid factor of the terminal limit group of $CGRes_m$. By the construction of approximating
cover resolutions,  we  further assume that the sequence of 
cover approximations, $\{CGRes_m\}$, converges into the original graded resolution, $GRes$, and in particular the terminal limit groups
of the covers, $\{CGRes_m\}$, that we denote, $L_s^m$, converge into the terminal limit group of $GRes$,  $L_s$.

With each of the approximating covers, $CGRes_m$, and the system of equations, $\Sigma(x,y,w,p)=1$,
 we associate (non-canonically) a graded
formal Makanin-Razborov diagram. Note that by theorem 2.7, the (formal graded) completion of each of the resolutions in 
such a  formal graded Makanin-Razborov
diagram can be embedded into a f.p.\ (formal graded) completion, and that with the terminal limit group of the cover closure of
$CGRes_m$, that is associated with a resolution in the graded formal Makanin-Razborov diagram, there is an associated finite
collection of covers of the flexible quotients of the rigid or weakly solid factor of the terminal limit group, and these covers
are all embedded in f.p.\ completions (see theorem 2.7).

By the construction of the approximating covers, $\{CGRes_m\}$, for each index $m$, there is a subsequence of the sequence of tuples,
$\{(p_n,w_n,y_n)\}$, that factor through $CGRes_m$. In particular, for each index $m$, with the subsequence of tuples, 
(still denoted) $\{(p_n,w_n,y_n)\}$,
that factor through $CGRes_m$, we can associate a sequence of specializations of the terminal limit group of $CGRes_m$, $L_s^m$,
that restricts to rigid or weakly strictly solid specializations of the rigid or weakly solid factor of $L_s^m$, and to 
specializations of the elliptic factors of $L_s^m$.

If there exists an index $m$, for which there is an infinite subsequence of tuples (still denoted),
$\{(p_n,w_n,y_n)\}$, that factor through $CGRes_m$, and 
so that the rigid or weakly strictly solid specializations of the rigid or weakly solid factor in $L_s^m$,
and the specializations of the elliptic factors of $L_s^m$, that are
associated with the tuples, $\{(p_n,w_n,y_n)\}$, extend to test sequences (of $CGRes_m$) that do not extend to formal solutions,
or they
do extend to formal solutions (over  closures of $CGRes_m$), but for each such formal solution, at least one of the inequalities in the
system, $\Psi(x,y,w,p) \neq 1$, does not hold (for generic value of the universal variables $y$),
we reached a terminal point of the iterative procedure. In this case, the associated output is the approximating cover, $CGRes_m$,
its associated (finite) collection of covers of the flexible quotients of the rigid or weakly solid factor of $L_s^m$, and its formal
graded Makanin-Razborov diagram.

Suppose that for every index $m$, there is no such subsequence of tuples $\{(p_n,w_n,y_n)\}$. 
In this case for each index $m$, there exists
an index $n_m > m$, so that the tuple $(p_{n_m},w_{n_m},y_{n_m})$ factors through $CGRes_m$, and with the associated rigid or 
weakly strictly solid
specialization  of the rigid or weakly solid factor of $L_s^m$, and the associated specializations of the elliptic factors of $L^m_s$, 
it is possible to associate a test sequence and a formal solution $x_{n_m}$ that
does  satisfy $\Sigma(x_{n_m},y,w_{n_m},p_{n_m})=1$ and $\Psi(x_{n_m},y,w_{n_m},p_{n_m}) \neq 1$ for generic  $y$.

By construction, the factors of the terminal limit groups of the graded resolutions, $\{CGRes_m\}$, $L^m_s$, converge into the
corresponding factors of the terminal limit group of the graded resolution, $GRes$, $L_s$. The formal solutions, $\{x_{n_m}\}$, are defined
over the approximating cover resolutions, $\{CGRes_m\}$.
Using the techniques to construct graded formal limit groups over free products, that were presented in section 2,
from the sequence of formal solutions, $\{x_{n_m}\}$, and the specializations of the terminal limit groups, $L^m_s$, that are
associated with the tuples, $\{(p_{n_m},w_{n_m},y_{n_m})\}$, it is possible to extract a subsequence 
that converges into a graded formal limit group over the original (limit) graded resolution $Res$, $FL(x,z,y,w,p)$. 
By the 
construction of the graded formal limit group $FL(x,z,y,w,p)$, the equations from the system, $\Sigma(x,y,w,p)=1$, represent the 
trivial word in $FL(x,z,y,w,p)$, whereas each of the inequations in the system, $\Psi(x,y,w,p) \neq 1$, represent a non-trivial 
element in $FL(x,z,y,w,p)$.

At this point we look at the sequences of specializations: $$\{(x_{n_m},y_{n_m},z_1(n_m),\ldots,z_t(n_m),w_{n_m},p_{n_m})\}$$ 
of the formal
solution $x_{n_m}$, the universal variables $y_{n_m}$,
and its successive shortenings $(z_1,\ldots,z_t)$, and the parameters, $(p_{n_m},w_{n_m})$, that take their
values in the free products, $\{A^1_{n_m}*\ldots*A^{\ell}_{n_m}\}$, and each of the tuples,
$(p_{n_m},w_{n_m},y_{n_m})$, factors through the cover graded resolution, $CGRes_m$. The elements, $\{y_m\}$, are precisely the subsequence
of values of the universal variables for which the sentence fail to hold for the free products, $\{A^1_{n_m}*\ldots*A^{\ell}_{n_m}\}$. 

Given this sequence of specializations, we apply the first step of the procedure for the analysis of an AE set, 
that is presented in section 2 of [Se5], and extract a subsequence, that converges into a
quotient graded resolution of the one that is associated with the formal limit group, $FL(x,z,y,w,p)$. Since the formal
solutions, $\{x_{n_m}\}$, were assumed to satisfy both the equalities, $\Sigma(x,y,w,p)=1$, and the inequalities,
$\Psi(x,y,w,p) \neq 1$, the obtained quotient resolution is not a graded closure of the original graded resolution, $GRes$, that we
 have started
with, but rather a resolution of "reduced complexity" (in the sense of the iterative procedure that is presented in section
2 in [Se5]).

We continue iteratively. At each step we start with a quotient graded resolution (with respect to the parameter
subgroup), $<w,p>$, $QRes$, that was constructed in the
previous step of the procedure, 
using the general step of the iterative procedure that is presented in section 2 in
[Se5], and a sequence of homomorphisms into free products, that extends a subsequence of the original sequence
of tuples, $\{(p_n,w_n,y_n)\}$, that converges into the (completion of the)
quotient resolution, $QRes$. Note that the completion of such a quotient resolution need not be finitely presented. 

We start the current (general) step, by fixing a sequence of approximating cover graded resolutions, $QRes_m$, of the quotient resolution, $QRes$, 
that satisfy the properties that
are listed in theorem 1.21. By the construction of cover resolutions (theorem 1.21), 
we may further assume that the sequence of cover resolutions, $QRes_m$, converges into the original 
quotient resolution, $QRes$.
Note that for each index $m$, there exists a subsequence of the original sequence of specializations (that converges into the completion of the
quotient resolution, $QRes$) that factor through the cover (graded) resolution, $QRes_m$. In particular, this subsequence restricts to rigid 
or weakly strictly solid specializations of the rigid or weakly solid factor of the terminal limit group of the cover resolution, 
$QRes_m$. 

With each cover (approximating) graded resolution, $QRes_m$,  we associate (not in a canonical way)
its graded formal Makanin-Razborov diagram (see theorem 2.7). Note that the completion of each resolution in the formal
Makanin-Razborov diagram can be embedded in a f.p.\ completion, and that with each such formal resolution there is an associated finite
collection of covers of the flexible quotients of the rigid or weakly solid factor of its terminal limit group, and each cover from this finite
collection can be embedded in a f.p.\ completion (see theorem 2.7). Also, recall that a formal Makanin-Razborov diagram of $QRes_m$ encodes all
the formal solutions that can be defined over it for all possible free products, $A^1*\ldots*A^{\ell}$.

If there exists an index $m$, and a subsequence of the original sequence of specializations that factor through the cover approximating resolution,
$QRes_m$,
so that over the rigid or the weakly strictly solid specializations of the rigid or weakly solid factor of the terminal limit group of
$QRes_m$, and the specializations of the elliptic factors of the terminal limit group of $QRes_m$, that are associated
with the subsequence,  there exist test sequences that can not be extended to  formal solutions over  closures of $QRes_m$,  so that 
both the equalities, $\Sigma(x,y,w,p)=1$, and the inequalities, $\Psi(x,y,w,p) \neq 1$, hold for generic $y$'s
(over the corresponding subsequence
of free products, $\{A^1_n*\ldots*A^{\ell}_n\}$), we
reached a terminal point of our iterative procedure. In this case the final output of the iterative
procedure is the graded quotient resolution, $QRes_m$,  and its (non-canonical)
graded formal Makanin-Razborov diagram.

Suppose that there is no such index $m$, and no such subsequence of the specializations that factor through $QRes_m$.
In this case
for each index $m$, there exists a specialization from the sequence that converges to $QRes$ (of index bigger than $m$ in that
sequence) that factors through $QRes_m$,  and a formal solution $x_{m}$ that
is defined over (a closure of) the graded resolution, $QRes_m$, so that for the specialization of the terminal limit group of $QRes_m$ that is
associated with this specialization, and for the formal solution, $x_m$, both the equalities, $\Sigma(x,y,w,p)=1$, and the inequalities,
$\Psi(x,y,w,p) \neq 1$, hold for generic values of the universal variables $y$ (generic values that are associated with the specialization
of the terminal limit group of $QRes_m$, and the specialization that factors through it from the sequence that converges to $QRes$).

In this case we look at the sequence of these specializations (that factor through $QRes_m$), 
and their associated formal solutions, $\{x_m\}$, that are defined over $QRes_m$.
We analyze the  combined sequence (of pairs of a specialization that factors through $QRes_m$ and the formal solution $x_m$), 
according to the analysis of quotient graded resolutions, which is part of the general step of the iterative procedure for the
analysis of AE sets (over free groups), that is presented in section 2 of [Se5], and    
finally extract a subsequence that converges into 
a quotient graded resolution of the resolutions that we have this step with, $QRes$ (see section 2 of [Se5]).
Since the formal solutions are guaranteed to satisfy both the
equalities, $\Sigma(x,y,w,p)=1$, and the inequalities, $\Psi(x,y,w,p) \neq 1$, for generic values of $y$, the complexity of the
obtained quotient graded resolution is strictly smaller than the complexity of the graded resolution, $QRes$, that we have 
started the current step with.

By theorem 2.10 in [Se5], this iterative procedure terminates after finitely many steps (the proof of the termination over free products is identical
to the one over free groups). If we reached a
terminal quotient graded resolution along the process, we found a cover graded  resolution, $QRes_m$, that satisfies the properties of theorem 1.21.
Furthermore, there exists a subsequence of the original sequence of tuples, $\{(p_n,w_n,y_n)\}$, that extend to specializations that factor through 
this cover resolution, $QRes_m$, and these extended specializations restrict to rigid or weakly strictly solid specializations of the rigid
or weakly solid factor of the terminal limit group of $QRes_m$. Each extended specialization from this subsequence restricts to a 
specialization of the terminal graded limit group of $QRes_m$, and there  exists a test sequence that is associated with it, that
can not be extended to a formal solution that is defined over (a closure of) $QRes_m$, for which both
the equalities, $\Sigma(x,y,w,p)=1$, and the inequalities, $\Psi(x,y,w,p) \neq 1$, hold for generic values of the universal variables $y$, that
are associated with $QRes_m$, and the associated  specialization of the terminal limit group of $QRes_m$.

If we didn't find such an  approximating cover  resolution, $QRes_m$,  along the iterative procedure, it continued until we reached 
a quotient graded  resolution which is a free product of a rigid  limit group over free products (with respect to the parameter subgroup
$<w,p>$) and elliptic factors. Once again we look at a sequence
of  approximating covers of the rigid limit group (over free products) and the elliptic factors, and we construct these covers according to the
construction that appears in theorem 1.21. With each such approximating cover,
and the system of equations, $\Sigma(x,y,w,p)=1$,  we
associate (non-canonically) its graded formal Makanin-Razborov diagram, according to theorem 2.7.

Hence, whatever terminal point of the iterative procedure we have reached, we found an approximating cover  (graded) resolution (over free
products) that satisfies the properties of approximating covers that appear in theorem 1.21. 
A subsequence of the given sequence of specializations, $\{(p_n,w_n,y_n)\}$, extend to specializations that do factor through this approximating
cover. Each extended specialization restricts to a specialization of the terminal limit group of the approximating cover resolution.
For each associated specialization of the terminal limit group of the cover resolution, there exists a test sequence for which
there is no formal solution  that is defined over the approximating cover resolution, 
and for which the equalities, $\Sigma(x,y,w,p)=1$, and the inequalities, $\Psi(x,y,w,p) \neq 1$, do hold for generic values of the universal
variables $y$, in the the fiber that is associated with the approximating cover resolution, and the  specialization of its terminal
graded limit group.

As in the previous section, the existence of such a cover resolution, that satisfies the properties that are listed in theorem 1.21, is 
sufficient to enable one to associate with a given AE set, a finite collection of (cover) graded resolutions that satisfy the properties
that are listed in theorem 1.21, and their (non-canonical) formal graded Makanin-Razborov diagrams, 
 so that if a specialization of the tuple $(w,p)$ in some  free product, $A^1*\ldots*A^{\ell}$ (for some $\ell >1$), is not in the definable set $AE(w,p)$,
then the corresponding sentence  must be false for a generic 
point of one of the constructed (cover) graded resolutions, where the terminal  factors of the graded resolutions take some value in the free
product, $A^1*\ldots*A^{\ell}$. This is the key
for analyzing an AE set over free products, i.e., the key for obtaining quantifier elimination from a predicate over the ambient free product
to a predicate over the factors.

\vglue 1pc
\proclaim{Theorem 4.1} Let:
$$AE(w,p) \ = \   \forall y \ \exists x \
 \Sigma(x,y,w,p)=1 \, \wedge \, \Psi(x,y,w,p) \neq 1$$ 
be an AE set  over groups. Then there exist finitely many (cover) graded resolutions over free products (with respect to the
parameter subgroup $<w,p>$) that satisfy the properties that are listed in theorem 1.21:
 $GRes_1(z,y,w,p),\ldots,GRes_d(z,y,w,p)$, with the following properties:

\roster
\item"{(1)}" with each graded resolution, $GRes_i(z,y)$, we associate (non-canonically) its graded formal Makanin-Razborov
diagram over free products with respect to the system of equations, $\Sigma(x,y,w,p)=1$.

\item"{(2)}" each of the graded resolutions, $GRes_i$, can be extended to an ungraded resolution with a f.p.\ completion,
$Comp_i$ (see theorem 1.21).

\item"{(3)}" let $G=A^1*\ldots*A^{\ell}$ be a non-trivial free product (for some $\ell >1$) that is not isomorphic to $D_{\infty}$.
Let $(w_0,p_0)$ be a specialization
of the parameters $(w,p)$ in the free product $A^1*\ldots*A^{\ell}$, for which $(w_0,p_0) \notin AE(w,p)$ over the free product $G=A^1*\ldots*A^{\ell}$.
Then there exists
an index $i$, $1 \leq i \leq d$,  a rigid or a weakly strictly solid specialization of the rigid or weakly solid factor of the terminal limit group
of $GRes_i$ in $G=A^1*\ldots*A^{\ell}$, that extends $(w_0,p_0)$, 
and specializations of the elliptic factors of that terminal limit group in elliptic subgroups in
$G=A^1*\ldots*A^{\ell}$,
 so that the combined  specialization of the terminal limit group of
the graded resolution, $GRes_i$, can be extended to a specialization of the f.p.\ completion, $Comp_i$, and   with this combined 
specialization of the terminal limit group of $GRes_i$ there exists an associated test sequence for which there
is no formal solution over (a closure of) $GRes_i(z,y,w,p)$ for which both the equalities, $\Sigma(x,y,w,p)=1$, and the inequalities,
$\Psi(x,y,w,p) \neq 1$, do hold for generic values of the universal variables $y$.
\endroster

\noindent
In other words, given an AE set,
there exists a finite collection of cover graded resolutions  over free products (with the properties that are listed in theorem 1.21),
so that the failure of a specialization of the parameters (free variables) to be in the AE set over any given free product,
can be demonstrated by the lack of the existence
of a formal solution over a test sequence of at least one of these (finitely many cover) graded resolutions.
\endproclaim

\nfp Let $G_n=A^1_n*\ldots*A^{\ell}_n$ (possibly for varying $\ell>1$), be a sequence of non-trivial free products 
that are not isomorphic to $D_{\infty}$. Let
$\{(p_n,w_n,y_n)\}$ be a sequence of tuples in $G_n$, so that $(p_n,w_n) \notin AE(p,w)$
over $G_n$, and $y_n$ testifies for the failure of  the pair $(p_n,w_n)$ to be in $AE(p,w)$, i.e.
the following existential sentence (with coefficients):
$$ \exists x \
 \Sigma(x,y_n,w_n,p_n)=1 \, \wedge \, \Psi(x,y_n,w_n,p_n) \neq 1$$ 
is false over $G_n$.

Starting with the sequence, $\{(p_n,w_n,y_n)\}$, the terminating iterative procedure that we have presented
constructs a graded resolution, $GRes$, with the following
properties:
\roster
\item"{(i)}" the graded resolution $GRes$ satisfies the properties of a cover graded resolution that are 
listed in theorem 1.21. In particular, its completion can be extended to (an ungraded) f.p.\ completion. With the rigid
or weakly solid factor of the terminal  limit group of $GRes$, there is a finite collection of covers of its
flexible quotients that can all be embedded into f.p.\ completions.

\item"{(ii)}" there exists a subsequence of the sequence of specializations, $\{(p_n,w_n,y_n)\}$, that extend to 
specializations 
that factor through the resolution, $GRes$, and to specializations of the f.p.\ completion of the ungraded resolution that
extends the graded resolution, $GRes$ (see theorem 1.21). 
Hence, with each specialization from this subsequence of the tuples, 
$\{(p_n,w_n,y_n)\}$, 
specializations of the  terminal limit group of the graded resolution, $GRes$, can be associated. Furthermore, 
these specializations of the terminal limit group of $GRes$
 can be extended to specializations of the f.p.\ completion of some ungraded resolution of the terminal limit group of
$GRes$. With this subsequence of specializations of the terminal limit group of the graded resolution, $GRes$, there are associated
test sequences for which no
 formal solutions that satisfy both the equalities, $\Sigma(x,y,w,p)=1$, and the inequalities, $\Psi(x,y,w,p) \neq 1$, for generic values
of the universal variables $y$,
can be constructed.

\item"{(iii)}" with the resolution $GRes$, and the system of equations, $\Sigma(x,y,w,p)=1$, we associate (non-canonically) a
formal graded Makanin-Razborov diagram (see theorem 2.7). The completion of every (formal) resolution in this 
diagram can be embedded into a f.p.\
completion, and the finite collection of covers that is associated with the rigid or weakly solid
factor of the terminal limit group of each resolution in these formal graded Makanin-Razborov diagrams,
can be embedded in a f.p.\ completion as well.
\endroster

At this point we are able to apply the argument that was used to prove theorem 3.2 (for AE sentences).
We look at all the sequences of non-trivial free products, $G_n=A^1_n*\ldots*A^{\ell}_n$
(possibly for varying $\ell>1$), that are not isomorphic to $D_{\infty}$, an associated sequence of 
 the tuples, $(p_n,w_n) \notin AE(p,w)$,
over $G_n$,
and specializations $y_n$ (in $G$) of the universal variables $y$ that testify for that. 
Given every such  sequence we use our terminating iterative procedure, and extract a subsequence
of the tuples, (still denoted) $\{(p_n,w_n,y_n)\}$, and a graded resolution,  
$GRes$, that has the properties (i)-(iii), and in particular,   the subsequence of tuples, (still denoted) $\{(p_n,w_n,y_n)\}$, 
extend to specializations that factor through the resolution $GRes$, and to specializations of the f.p.\ completion of
the ungraded resolutions that extends the graded resolution $GRes$.
Furthermore, no formal solution that satisfies both the equalities, $\Sigma(x,y,w,p)=1$, 
and the inequalities, $\Psi(x,y,w,p) \neq 1$,
for generic values of the variables $y$ can be constructed, for some test sequences that are associated with those values of the  terminal 
limit group of $GRes$, that are associated with the subsequence of tuples, $\{(p_n,w_n,y_n)\}$.

The completion of each of the constructed resolutions, $GRes$, can be extended to an ungraded f.p.\ completion, 
and so are the resolutions in its 
associated formal graded Makanin-Razborov diagram. Furthermore, the (finite collection of) covers of flexible quotients
of the rigid or weakly solid factors of the terminal limit groups of the resolutions, $GRes$, and of the formal graded
resolutions that are associated with the resolutions $GRes$, can all be embedded into (ungraded) resolutions with f.p.\
completions. Hence, we can define a linear  order on this (countable) collection of
resolutions ($GRes$), and their (non-canonically) associated formal Makanin-Razborov diagrams, and (finite collections of)
covers of flexible quotients.  By the same argument that
was used in constructing the Makanin-Razborov diagram (theorem 26 in [Ja-Se]), there exists a finite subcollection
of these graded resolutions that satisfy properties (1)-(3) of the theorem.

\line{\hss$\qed$}

\medskip
Given an AE set, theorem 4.1 constructs a finite collection of (cover) graded resolutions over free products, where the completions of all
these graded resolutions can be embedded into (ungraded) f.p.\ completions (cf. theorem 1.21), that
demonstrates the failure of specializations of the parameters (free variables) to be in the given  AE set over general free products. 
This is an  analogue of the uniformization of proofs (the construction of the tree of stratified sets) that was proved in section
2 of [Se5] for AE sets over free groups. 

\noindent
To use these (finitely many) constructed resolutions for the analysis of EAE sets, we still
need to associate with the graded resolutions that are constructed in theorem 4.1, a (finite) collection of resolutions that contain
their singular locus, and a finite collection of resolutions that classifies those values of their terminal limit groups for which
the pegs of these graded resolutions have roots of prescribed orders (cf. propositions 3.4 and 3.5).

As in the ungraded case, in constructing the graded and formal graded Makanin-Razborov diagrams over free products, 
we have considered the set of specializations that 
factor through a resolution in one of these diagrams, as those homomorphisms that factor through the graded resolutions,
and in addition we required that the associated specializations of
the various limit groups along the resolution restrict to non-elliptic specializations of all the abelian edge (and vertex) groups, and all
the QH vertex groups along the resolution. 

In applying theorem 4.1 to  analyze EAE sets, and general definable sets, over free products,
 it is essential to determine the singular locus of the graded resolutions that are associated with an AE set by the theorem, as we need
to ignore 
 specializations that do factor through the singular locus. 

\noindent
The singular locus of a graded resolution, $GRes(y,p)$, is defined precisely in the same way as in the ungraded case (definition 3.3).
To analyze the singular locus of a graded resolution we prove a graded version
of proposition 3.4. Given the techniques of the first two sections, especially the construction of the formal graded Makanin-Razborov
diagram (theorem 2.7), the generalization of proposition 3.4 to the graded case is rather straightforward.

\vglue 1pc
\proclaim{Proposition  4.2} Let $GRes(y,p)$ be a well-structured, coefficient free, graded
resolution over free products.  
There exist finitely many graded approximating cover resolutions, that satisfy the properties of the resolutions in the formal
graded Makanin-Razborov diagram (theorems 2.6 and 2.7), 
$SLRes_1(y,p),\ldots,SLRes_u(y,p)$, and in particular the (graded) completion of each of these graded resolutions embeds into an (ungraded)
f.p.\ completion.   

For every non-trivial free product,
$G=A^1*\ldots*A^{\ell}$, that is not isomorphic to $D_{\infty}$, a specialization of the terminal limit group of $GRes(y)$ is in the singular
locus of $GRes(y,p)$, if and only if it extends to a specialization of (at least) one of the terminal  limit groups of the resolutions,
$SLGRes_1(y,p),\ldots,SLGRes_u(y,p)$, and further extends to a specialization of the f.p.\ (ungraded) completion that contains the (graded) completion
of the corresponding resolution, $SLGRes_i(y,p)$.
\endproclaim

\nfp The statement of the theorem follows by exactly the same argument that was used to prove the analogous statement in the
ungraded case (proposition 3.4), where instead of using the construction of the (ungraded) formal Makanin-Razborov diagram, we apply
the construction of the graded formal Makanin-Razborov diagram (theorem 2.7).

\line{\hss$\qed$}

As in the ungraded case, to analyze an EAE and general definable sets over free products,  we also need to construct auxiliary resolutions, 
that will enable
one to decide uniformly if over a given free product $G$, the set of specializations that do factor through a given
resolution are covered by a given finite collection of closures. The construction of these $graded$ $Root$ resolutions essentially
follows the construction of ungraded Root resolutions in proposition 3.5.

\vglue 1pc
\proclaim{Proposition  4.3} Let $GRes(y,p)$ be a well-structured, coefficient free
resolution over free products with a f.p.\ completion (and terminal limit 
group), and let: $$Cl_1(Res),\ldots,Cl_f(Res)$$ be a given finite set of graded closures of $GRes(y,p)$ (see definition 2.3). 

Each of the abelian vertex groups that appear in the abelian decompositions that are associated with the various levels of
$GRes(y,p)$ is contained as a finite index subgroup in an abelian vertex group in each of the closures,
$Cl_1(Res),\ldots,Cl_f(Res)$. We set $Pind$ to be the product of these indices (where the product is over  all the abelian 
vertex groups in $GRes(y,p)$, and
all the closures, $Cl_1(Res),\ldots,Cl_f(Res)$). 

Let $pg_1,\ldots,pg_t$ be a generating set for the pegs of pegged abelian groups in $GRes(y,p)$. Note that in the graded case,
the peg of a pegged abelian vertex group is in general a f.g.\ free abelian group, hence, its generating set may consist of more than a single
element). With each collection of positive integers, $d_1,\ldots,d_t$,
that do all divide $Pind$, we associate a finite (possibly empty) collection of graded closures of $GRes(y,p)$ with f.p.\ completions and
terminal limit groups. We denote the entire
finite collection of such graded closures, $\{GRootRes_r\}$, and call them (graded) $Root$ $resolutions$.

Let $G=A^1*\ldots*A^{\ell}$ be a non-trivial free product that is not isomorphic to $D_{\infty}$. We fix a specialization of the terminal
limit group of $GRes(y,p)$ in $G$. The pegs of abelian vertex groups in $GRes(y,p)$, $pg_1,\ldots,pg_t$, have roots of orders
$d_1,\ldots,d_t$, for every test sequence (in $G$) over the given specialization of the terminal limit group of $GRes(y,p)$, if and
only if the specialization of the terminal limit group of $G$ can be extended to at least one of the graded
Root resolutions, $GRootRes_r$,
that are associated with the sequence, $d_1,\ldots,d_t$.   
\endproclaim

\nfp Identical to the proof of the analogous statement in the ungraded case (proposition 3.5), were as in the construction of the 
graded resolutions that are associated with the singular locus of a graded resolution (proposition 4.2), 
we apply the construction of the graded
formal Makanin-Razborov diagram (theorem 2.7), instead of the construction of the (ungraded) formal Makanin-Razborov diagram  that
was used in proving proposition 3.5.

\line{\hss$\qed$}

In a similar way to the ungraded case, the graded Root resolutions that are constructed in proposition 4.3, enable one to  associate
a predicate with the set of values of the parameters $p$, for which a given set of (graded) closures is a covering closure. 

\noindent
The finite collection of graded resolutions that are associated with an AE set in theorem 4.1, their formal graded Makanin-Razborov diagrams,
the singular loci of all these resolutions, and their associated graded Root resolutions, enable us to analyze EAE sets over free products. For such
analysis 
we need
a variation  of the sieve procedure over free groups that is presented in [Se6] for free products.

Let:
$$  EAE(p) \ = \ \exists w \ \forall y \ \exists x \
 \Sigma(x,y,w,p)=1 \, \wedge \, \Psi(x,y,w,p) \neq 1$$ 
be a predicate over groups. With the given EAE set we naturally look at the AE set:
$$  AE(p,w) \ = \  \forall y \ \exists x \
 \Sigma(x,y,w,p)=1 \, \wedge \, \Psi(x,y,w,p) \neq 1.$$ 
By theorem 4.1 with this AE set it is possible to associate finitely many (cover) graded resolutions that satisfy the properties that are listed
in theorem 1.21. With each of these (finitely many) graded resolutions, we have  associated 
(non-canonically) a graded formal Makanin-Razborov diagram, a finite collection of graded 
resolutions that are associated with the singular
loci of all these resolutions, and a finite collection of graded Root resolutions, that are associated with each of the graded
resolutions that are associated with the AE set by theorem 4.1, and the closures that are associated with the resolutions in
the formal Makanin-Razborov diagrams of these graded resolutions.  These  finite collection of graded (and graded formal) resolutions, 
demonstrate the failure of a specialization of a pair $(p,w)$ to be in the
definable set $AE(p,w)$ over an arbitrary non-trivial free product, $G=A^1*\ldots*A^{\ell}$ (for an arbitrary $\ell>1$), where
 $G$ is
not isomorphic to $D_{\infty}$.

Let: $GRes_1(z,y,w,p),\ldots,GRes_d(z,y,w,p)$ be the (cover) graded resolutions that we have associated with the set $AE(p,w)$. Each of these graded
resolutions terminates in a graded limit group which is a free product of a rigid or a weakly solid limit group with (possibly) a finite collection
of f.p.\ elliptic factors, and the completion of each of these cover graded resolutions can be embedded into a f.p.\  (ungraded)  completion, that we denote $Comp_1,\ldots,Comp_d$.

With each of the graded resolutions, $GRes_i$, and the system of equations, $\Sigma(x,y,w,p)=1$,  
we have associated (non-canonically) its graded formal Makanin-Razborov diagram, 
and the completion of each formal resolution in one
of these diagrams can be embedded into a f.p.\ completion as well. We denote each of the (finitely many) resolutions
in the formal graded Makanin-Razborov diagram of a graded resolution, $GRes_i$, $FGRes_i^j$. 

With each of the formal graded resolutions, $FGRes_i^j$, we further associate a finite collection of graded resolutions.
For each formal graded resolution, $FGRes_i^j$, we look at all of its test sequences, for which at least one of the 
inequalities in the system, $\Psi(x,y,w,p) \neq 1$, fails to be an inequality and is in fact an equality for the entire
test sequence. With the collection of all such test sequences it is possible to associate a finite collection of graded
resolutions, that do all have the same structure as the formal graded resolution, $FGRes_i^j$, just that each factor in the 
terminal limit group of $FGRes_i^j$ is replaced by a quotient. We denote each of the (finitely many) obtained resolutions,
that are associated with $FGRes_i^j$,
$CollFGRes_i^{j,k}$. 

With each of the graded resolutions, $GRes_i$, we have associated a finite collection of graded resolutions that are associated
with its singular locus. We denote these graded resolutions, $SLGRes_i^v$. We also associate finitely many resolutions with each
graded formal resolution, $FGRes_i^j$, that we denote, $SLFGRes_i^{j,v}$. 

\noindent
We further associate graded Root resolutions with each
of the graded resolutions, $GRes_i$, and the collection of closures that are associated with the graded resolutions in its
formal graded Makanin-Razborov diagram. We denote these graded Root resolutions, $GRootRes_i^r$.
 
We start the analysis of the set $EAE(p)$ with all the sequences of specializations of the tuple, $(p,w)$, $\{(p_n,w_n)\}$, that take values
in non-trivial free products, $G_n=A^1_n*\ldots*A^{\ell}_n$ (for some $\ell>1$),  which is not isomorphic to $D_{\infty}$. 
We further assume that for every index $n$,
$p_n \in EAE(p)$ over the
free product, $G_n$, and that $w_n$ is a witness for $p_n$, i.e., a specialization of the existential variables $w$ in $G_n$, so
that the AE sentence:
$$   \forall y \ \exists x \
 \Sigma(x,y,w_n,p_n)=1 \, \wedge \, \Psi(x,y,w_n,p_n) \neq 1$$ 
is a truth sentence over $G_n$.

By theorem 1.16, 
given such a sequence of pairs, 
 $\{(p_n,w_n)\}$, we can pass to a subsequence that converges into a 
well-structured  (even well-separated) graded resolution with respect to the parameter subgroup $<p>$: 
$M_0 \to M_1 \to \ldots \to M_s$, where $M_s$ is
 a free product of a rigid or a solid group over free products (that contains the subgroup, $<p>$), 
and (possibly) a free group and (possibly) finitely many elliptic factors. 
We denote this graded resolution
$WRes$. Note that the terminal limit group $M_s$   in $WRes$ is f.g.\ but it may be infinitely presented.

As we did in analyzing the set, $AE(p,w)$, we  
fix a sequence of approximating cover resolutions of $WRes$, that we denote, $\{WRes_m\}$, that are constructed following the construction that
appears in theorem 1.21, and hence each of them satisfies the properties of cover graded resolutions that are listed in
theorem 1.21. By construction we may further assume that the sequence of approximating cover resolutions, $\{WRes_m\}$, converge into
the original graded resolution, $WRes$. 

By the construction of the approximating cover graded resolutions that
appears in theorem 1.21, for every index $m$, there is a subsequence of the sequence of specializations, $\{(p_n,w_n)\}$, that
factor through the approximating cover resolution, $WRes_m$, and they further extend to specializations of
the f.p.\ completion of the ungraded resolution that extends the approximating cover, $WRes_m$ (see theorem 1.21).
In particular, with this subsequence it is
possible to associate specializations of the terminal limit group of $WRes_m$, that restrict to rigid or weakly
strictly solid specializations of the rigid or weakly solid factor of the terminal (graded) limit group
of $WRes_m$.

First, suppose that there exists an approximating resolution, $WRes_m$, for which there exists a subsequence of pairs, (still
denoted) $\{(p_n,w_n)\}$, so that for each specialization of the terminating limit group of $WRes_m$ that is associated with  
a pair $(p_n,w_n)$ from the subsequence, there exists a test sequence  in the fiber that is associated with
such a specialization  of the terminal graded limit group (which takes its values in $G_n$) that satisfies $p_n \in EAE(p)$
(over $G_n$), and the values of the variables  $w$ from the test sequence are witnesses for $p_n$ (over $G_n$). 

In this case with the approximating graded resolution, $WRes_m$, we
associate finitely many graded resolutions, that have similar properties to the resolutions in
 the formal graded Makanin-Razborov
diagram (see theorem 2.7).

First, we look at
all the test sequences of $WRes_m$, that can be extended to specializations of one of the  terminal limit
groups of the graded resolutions, $GRes_1,\ldots,GRes_d$,  that restrict to rigid or weakly strictly solid 
specializations of the rigid or weakly solid factor of that terminal limit group, and so that the specializations
of the terminal limit group, $GRes_i$,  extends to a specializations of the f.p.\ completion, $Comp_i$, into which the completion of
$GRes_i$ is embedded.

By the construction the formal graded Makanin-Razborov diagram (theorem 2.7),  with the approximating (cover) graded resolution,
$WRes_m$, and the collection of test sequences of it, that can be extended to specializations of the terminal
limit groups of the graded resolutions, $GRes_1,\ldots,GRes_d$, it is possible to associate (non-canonically) a finite collection
of (formal like) graded resolutions, that do all terminate in  graded cover closures of the graded resolution, $WRes_m$, $CCl(WRes_m)$, possibly
free product with a free group. Each resolution from this finite collection has the same properties as the resolutions
in the formal Makanin-Razborov diagram (see theorems 2.6 and 2.7). 
We denote each of the constructed (finitely many, graded, cover) resolutions, $ExWRes_m$.

In a similar way to the formal Makanin-Razborov diagram (cf. theorem 2.7), given
a test sequence of $WRes_m$ that can be extended to specializations of the terminal limit group of one of the
graded resolutions, $GRes_1,\ldots,GRes_d$,  and its associated f.p.\ completion $Comp_i$, it is possible to
extract a subsequence of the extended specializations that factor through one of the resolution, $ExWRes_m$.

Given this finite collection of closures, $ExWRes_m$, that collect all the test sequences of the graded resolution, 
$WRes_m$, that can be extended to specializations of the terminal limit group of one of the graded resolutions,
$GRes_1,\ldots,GRes_d$, we further associate with each of the graded cover  resolutions, $ExWRes_m$, finitely many graded resolutions,
that terminate in (cover) closures of the graded resolutions, $ExWRes_m$.
These include 6 types of cover graded resolutions, each for collecting test sequences of $ExWRes_m$ with different properties. With
each collection of test sequences it is possible to associate a finite collection of cover graded resolutions, in a similar way,
and with similar properties to those of the resolutions in the formal Makanin-Razborov diagram (see theorems 2.6 and 2.7):  
\roster
\item"{(1)}" first we look at all the test sequences of $ExWRes_m$ for which the restrictions to specializations of the rigid or weakly solid
factor of the terminal limit group of one of the resolutions, $GRes_1,\ldots,GRes_d$, is not rigid or not
weakly strictly solid.

\item"{(2)}" second, we look at all the test sequences of $ExWRes_m$ for which the restrictions to specializations of 
 the terminal limit group of one of the resolutions, $GRes_1,\ldots,GRes_d$, can be extended to specializations of 
the terminal limit group of one the graded resolutions, $SLGRes_i^v$, that are associated with the singular locus of the
corresponding graded resolution, $GRes_i$.

\item"{(3)}" third, we look at all the test sequences of $ExWRes_m$ for which the restrictions to specializations of  
 the terminal limit group of one of the  resolutions, $GRes_1,\ldots,GRes_d$, can be extended to specializations of 
the terminal limit group of one the graded Root resolutions, $GRootRes_i^r$, that are associated with the 
corresponding graded resolution, $GRes_i$.

\item"{(4)}" we look at all the test sequences of $ExWRes_m$, for which to the restriction to specializations
of the terminal limit group of one of the graded resolutions, $GRes_1,\ldots,GRes_d$, it is possible to add a specialization of
the terminal limit group of one of the associated graded formal (cover) resolution, $FGRes_i^j$.

\item"{(5)}"  given the finite collection of graded cover resolutions that are associated with the test sequences that are
described in part (4), we look at those test sequences of these graded cover resolutions ,for which the restrictions to
the specializations of the rigid or weakly solid factor of the terminal limit group of $FGRes_i$ are not rigid or
not weakly strictly solid.

\item"{(6)}"  given the finite collection of graded cover resolutions that are associated with part (4), we look at those 
test sequences of these (finitely many) graded resolutions, for which the restrictions to the terminal limit groups of one
of the resolutions, $FGRes_i^j$, either factor through the terminal limit group of one of the associated collapse formal graded
resolutions, $CollFGRes_i^{j,k}$, or the can be extended to specializations of the terminal limit group of one of the resolutions,
$SLFGRes_i^{j,v}$, that are associated with the singular locus of the corresponding resolution, $FGRes_i^j$.
\endroster

Finally, with the resolution, $WRes_m$,  itself, we associate finitely many graded resolutions, that are associated with its singular
locus, according to proposition 4.2. We further associate with $WRes_m$,   and the finite collection of closures of it, that
include the resolutions, $ExWRes_m$, and the resolutions that are associated with it according to parts (1)-(6), a finite
collection of graded Root resolutions, according to proposition 4.3.

\smallskip
Suppose that there is no 
approximating resolution, $WRes_m$, with a subsequence of (the original sequence of) tuples $\{(p_n,w_n)\}$ 
that factor through it, and that can be extended to specializations of the f.p.\ completion that extends $WRes_m$, for
which for the fiber of $WRes_m$ that contains the tuple, $(p_n,w_n)$, the restrictions of a test sequence in the fiber
to the existential variables $w$ 
testify that $p_n \in EAE(p)$ over $G_n$, i.e., for these $w$'s, the tuple $(p_n,w)$ is not in the associated set
$AE(p,w)$ (over $G_n$).

In this case we can associate with the original graded resolution, $WRes$, a graded closure that is constructed from a limit of
test sequences of a sequence of approximating resolutions, $WRes_m$, that can be extended to specializations of the terminal
limit group of one 
of the graded resolutions, $GRes_1,\ldots,GRes_d$, that restrict to rigid or weakly strictly solid specializations
of the rigid or weakly solid factor of that terminal group, where the specializations of the terminal limit group of one
of the graded resolutions, $GRes_1,\ldots,GRes_d$, can not be extended to
formal solutions that are defined over these graded resolutions, and these formal solutions form a covering closure (definition 2.4),
and
satisfy both the equalities,
$\Sigma(x,y,w,p)=1$, and the inequalities, $\Psi(x,y,w,p) \neq 1$. We denote such a closure, $ExWRes$.

Now, we continue in one of two possible ways, depending on the sequence of specializations, $\{(p_n,w_n\}$, and the fibers
that contain them in the approximations, $WRes_m$.
The first possibility is to look at a sequence of cover approximations of the closure, $ExWRes$, (see theorem
1.21 for the construction of such approximations), and force one of 3 possible  collapse conditions on it:
\roster
\item"{(i)}" we require that
the additional rigid or weakly strictly solid specialization of the rigid or weakly solid factor of the
terminal limit group of the associated graded resolution, $GRes_1,\ldots,GRes_d$, that was added to $ExWRes$, 
will be either non-rigid or non weakly strictly solid.

\item"{(ii)}" we require that the specialization of the terminal limit group of the associated graded resolution,
$GRes_i$, extends to a specialization of the terminal limit group of one of the graded resolutions, $SLGRes_i^v$,  that are
associated with the graded resolution, $GRes_i$.

\item"{(iii)}" we require that the specialization of the terminal limit group of the associated graded resolution,
$GRes_i$, extends to a specialization of the terminal limit group of one of the graded Root resolutions, $GRootRes_i^r$,  that are
associated with the graded resolution, $GRes_i$ (and the restrictions of the test sequences from which the resolution, $ExWRes$,
was constructed, do not extend to specializations of this graded Root resolution).
\endroster
By passing to a further subsequence, and applying the first step of the sieve
procedure, that is presented in [Se6], we obtain a  quotient resolution of the closure, $ExWRes$.

The second possibility is to look at a sequence of cover approximations of the closure, $ExWRes$, and force the existence 
of a formal solution that is
associated with the specialization of the terminal limit group of one of the resolutions, $GRes_1,\ldots,GRes_d$, that
was added to specializations of $ExWRes$. i.e., to each specialization in a test sequence of a cover approximation
of $ExWRes$ that we consider, we add specializations of the terminal limit group of one of the formal resolutions,
$FGRes_i^j$, that are associated with the graded resolutions, $GRes_1,\ldots,GRes_d$. From the collection of test sequences
that factor through cover approximations of $ExWRes$, and the additional specializations of the terminal limit groups of 
one of the formal resolutions that are associated with one of the resolutions, $GRes_1,\ldots,GRes_d$, we extract
 a subsequence that converges into a quotient resolution of the closure, $ExWRes$. The structure of such a
quotient resolution is determined by the first step of the sieve procedure [Se6], in the same way that it was
used in proving the equationality of Diophantine sets in section 2 of [Se9].

We continue iteratively, at each step we first look at a sequence of approximating (cover) resolutions of a quotient resolution
that was constructed in the previous step, where these approximating (cover) resolutions satisfy the properties that
are listed in theorem 1.21, and in particular they can be extended to ungraded resolutions with f.p.\ completions.
If there exists an approximating cover resolution for which there exists a subsequence
of pairs, (still
denoted) $\{(p_n,w_n)\}$, so that for each specialization of the terminating limit group of the approximating resolution
 that is associated with  
a pair $(p_n,w_n)$ from the subsequence, there exists a test sequence in the fiber that is associated with
such a specialization  (which takes its values in $G_n$) that satisfies $p_n \in EAE(p)$
(over $G_n$), and the restrictions of the test sequence to the existential variables $w$ are witnesses for $p_n$, 
we do what we did in this case in the first step. 

This means that we associate
with the approximate resolution finitely many auxiliary resolutions. First, we look at all the test sequences of
the approximating cover resolution that can be extended to specializations of the terminal limit group of one of the
graded resolutions, $GRes_1,\ldots,GRes_d$, and use a construction similar to the construction of the formal graded
Makanin-Razborov diagram to associate with the cover approximating resolution finitely many resolutions that are similar
to the resolutions, $ExWRes_m$, that were constructed in the first step. Then we look at all the test sequences of
these constructed resolutions that 
satisfy properties (1)-(6), where for each such collection of test sequences, we use a construction similar
to the one that was used to construct a formal Makanin-Razborov diagram (theorems 2.6 and 2.7). Finally ,with the approximating resolution and
the auxiliary resolutions that are associated with it, we associate finitely many graded resolutions that are associated
with the singular locus of the approximating resolution, and with it and the finite set of closures that are associated with its auxiliary
resolutions, according to propositions 4.2 and 432. Note that each of the constructed (auxiliary)
resolutions has the same properties of a resolution in the formal Makanin-Razborov diagram (see theorem 2.6), so in particular it terminates 
in a cover closure of the approximating cover resolution, and each of the constructed resolution can
be extended to an ungraded resolution with a f.p.\ completion.

If there is no such approximating resolutions for a constructed quotient resolution, we associate with it a graded resolution
that is constructed from test sequences that can be extended to specializations of the terminal
limit groups of one of the graded resolutions: $GRes_1,\ldots,GRes_d$ (that can not be extended to formal
solutions that are defined over these resolutions and satisfy both the equalities, $\Sigma(x,y,w,p)=1$, and the
inequalities, $\Psi(x,y,w,p) \neq 1$). We further force one of the 3 possible  collapse conditions over the constructed 
closure of the quotient resolution, or the existence of a formal solution that is associated with the
additional specialization of the terminal limit group of one of the graded resolutions,
$GRes_1,\ldots,GRes_d$, and apply the general step of the sieve procedure [Se6], and associate a quotient resolution of
the quotient resolution that we have started this (general) step with.

\noindent
By an argument which is similar to the argument that guarantees the termination of the sieve procedure in [Se6], and the
argument that was used to prove the equationality of Diophantine sets over  free and hyperbolic groups (proposition 2.2 in [Se9]),  we
obtain a termination of the this procedure.

\vglue 1pc
\proclaim{Theorem 4.4} The procedure for the analysis of an EAE set over free products terminates after finitely
many steps.
\endproclaim

\nfp
At each step of the procedure, on  the quotient resolution, $WRes$, that is analyzed in that step,  we impose one of two types of
restrictions.
The first adds to  generic points of fibers (test sequences) of the resolution, $WRes$, specializations of the terminal limit group
of one of the graded resolutions, $GRes_1,\ldots,GRes_d$, that restrict to rigid or weakly strictly solid specializations of the rigid
or weakly solid factor of that terminal limit group, and associate with $WRes$ a (cover) closure.
 Then we impose one of 3 possible  Diophantine conditions on the obtained (cover) closure, that either forces
the additional rigid or weakly strictly solid specialization to be either non-rigid or non weakly strictly solid, or forces the
specialization of the terminal limit group of one of the graded resolution, $GRes_1,\ldots,GRes_d$, to be extended to a 
specialization of  the terminal limit group of one of the resolutions that are associated with the singular locus,
$SLGRes_i^v$, or to the terminal limit group of one of the graded Root resolutions,
$GRootRes_i^r$. 

The second type adds specializations of the terminal limit groups of one of the associated formal graded (cover) resolutions, $FGRes_i^j$, that
restrict to rigid or weakly strictly solid specializations of the rigid or weakly solid factor in the terminal limit group of that formal
graded resolution, and these specializations do not extend to a closure of the quotient resolution from the previous step. 

These (3 possible) Diophantine conditions are similar in nature to the Diophantine condition 
that is imposed on quotient resolutions in the general step of the sieve procedure (over free groups)
in [Se6]. The additional specializations of the second type (that do not extend to a closure of the previous
quotient resolution), are similar to the rigid or strictly solid specializations that are added in each step in the iterative procedure that
is used to prove the equationality of Diophantine
sets over free and hyperbolic groups (section 2 in [Se9]). 
However, we need to modify the argument for the termination of the sieve procedure over free groups, to guarantee
the termination of the procedure for the analysis of an EAE set over free products.

The sieve procedure over free groups [Se6], analyzes a sequence of quotient resolutions, where each quotient resolution is obtained
from the previous one by imposing one of finitely many Diophantine conditions on a closure of the quotient resolution that was
obtained in the previous step. The complexity of such a quotient  resolution is measured by some finite collection of associated 
graded limit groups
(over free groups), a finite collection of core resolutions, their associated induced resolutions, and possibly a finite collection
of sculpted resolutions and carriers. Since the analysis and the definitions of these objects are rather involved, we won't elaborate
on them, and the interested reader is refered to [Se6] for the complete details.

The analysis of quotient resolutions over free products, along the steps of the iterative procedure for the analysis of an
EAE set, is done according to the general step of the sieve procedure. Hence, with each quotient resolution we also
associate (finitely many) graded limit groups (over free products), core resolutions, induced resolutions, and possibly carriers
and sculpted resolutions, precisely as they appear in the sieve procedure [Se6].

As the graded limit groups, the core resolutions, and the induced resolutions, that are constructed along the iterative procedure for the
analysis of an EAE set over free products, are similar to the ones that are constructed along the sieve procedure over free groups
in [Se6], 
the proof of theorem 22 in [Se6] that guarantees the termination of the sieve procedure over free groups, implies that 
the limit groups, the core resolutions,
and the induced resolutions, that appear along the procedure for the analysis of an EAE set over free products, can change only finitely many times
along the procedure. Hence, by the same proof of theorem 22 in [Se6], 
if the procedure does not terminate, the number of sculpted resolutions that are associated with the
quotient resolutions that are constructed in the various steps can not be bounded.

Suppose that the procedure for the analysis of an EAE set over free products does not terminate. By the argument that is used in proving
theorem 22 in [Se6], the larger and larger numbers of sculpted resolutions that are constructed along the iterative procedure,
in combination with the argument that was used to prove the combinatorial bounds on rigid and (weakly) strictly solid specializations of rigid
and (weakly) solid limit groups over free products (theorems 1.10 and 1.13), together with the argument that was
used to prove the existence of
finitely many systems of fractions for rigid and (weakly) strictly solid specializations in theorems 1.14 and 1.15,
guarantee the existence of a subsequence of steps along the iterative procedure, with an unbounded number of associated sculpted resolutions,
for which the additional specializations of the rigid or weakly solid factors of the terminal limit groups of one of the graded
resolutions,$GRes_1,\ldots,GRes_d$, satisfy the conclusions of theorems 1.10, 1.13, 1.14 and 1.15 (even though these specializations 
are not rigid nor weakly strictly solid in general).

Over free groups, such a subsequence of steps gives an immediate contradiction, as theorems 1.14 and 1.15 in the case of a free
group imply a global bound on the number of (families of) the specializations of the rigid or weakly solid factors of the terminal limit groups
of $GRes_1,\ldots,GRes_d$, that contradicts the lack of a bound on the number of sculpted resolutions along the
subsequence of steps of the procedure. Over free products, theorems 1.14 and 1.15 do not imply similar bounds on the number of
specializations of the rigid or weakly solid factors, hence, we need to slightly modify the argument.

The quotient resolutions that are constructed along the procedure, are obtained by either imposing one of 3 possible Diophantine
conditions on a closure of the quotient resolution that was constructed in the previous step,
or by requiring the existence of an additional specialization of the terminal limit group of one of the formal graded
resolutions, $FGRes_i^j$, that are associated with the graded resolutions, $GRes_1,\ldots,GRes_d$, a specialization that does not extend to a 
closure of the previously constructed quotient resolution. The second type of requirement, the existence of an additional specialization of
one of the terminal limit groups of the associated graded formal resolutions, can occur for only boundedly many sculpted 
resolutions, by the same argument that was used in the proof of theorem 22 in [Se6] (over free groups), or alternatively by the argument
 that was applied in proving the equationality of Diophantine sets over free and hyperbolic groups (theorem 2.2 in [Se9]).

Hence, to prove the termination of the iterative procedure for the analysis of an EAE set, we need to prove that the requirement
that the additional specializations of the terminal limit group of one of the graded resolutions, $GRes_1,\ldots,GRes_d$, satisfy
one of 3 possible Diophantine conditions, can  occur for boundedly many sculpted resolutions as well.  

Each of the 3 possible Diophantine requirements, is expressed by an extension of the specialization of the additional terminal
limit group of one of the graded resolutions, $GRes_1,\ldots,GRes_d$, to a specialization of one of finitely many possible 
graded limit groups.
With each such graded limit group we associate a (non-canonical) graded Makanin-Razborov diagram (see theorem 1.22) 
with respect to the parameter 
subgroup that is generated by the fixed generating set of the  (original)  terminal limit group of one of the
graded resolutions, $GRes_1,\ldots,GRes_d$. Each such graded resolution in the graded Makanin-Razborov diagram terminates in a free
product of a rigid or
a weakly solid limit group  with finitely many elliptic factors and possibly a free group.

Therefore, the Diophantine conditions that are imposed in each step of the iterative procedure 
on the additional specializations of the terminal limit group of one
of the graded resolutions, $GRes_1,\ldots,GRes_d$, can be expressed by the possibility to extend this specialization to a
specialization of one of the terminal limit groups  in the finitely many graded Makanin-Razborov diagrams that are associated with
the 3 possible Diophantine conditions, that restrict to  
a rigid or a weakly strictly solid specialization of the rigid or weakly solid factor of the corresponding terminal graded
limit group.

We have already pointed out that in case the iterative procedure doesn't terminate,
there exists a subsequence of steps for which the additional specializations of the terminal
limit groups of the graded resolutions, $GRes_1,\ldots,GRes_d$, that are associated with unbounded collections of sculpted
resolutions, satisfy the conclusions of theorems 1.10,1.13,1.14 and 1.15. Precisely the same argument guarantees that we can pass
to a further subsequence for which the extended specializations to terminal limit groups of the resolutions in the graded
Makanin-Razborov diagrams that are associated with the 3 possible Diophantine conditions, satisfy the conclusions of these 
(1.10,1.13.1.14 and 1.15) theorems as well.

Now we can apply the argument that was used to prove the termination of the sieve procedure over free groups (theorem  22 in [Se6]), 
and conclude that by the conclusions of theorems 1.14 and 1.15 for the extended specializations of the terminal limit groups in the
graded Makanin-Razborov diagrams that are associated with the Diophantine conditions, if there is no bound on the number of sculpted resolutions
that are associated with the quotient resolutions along the iterative procedure, then there must exist a specialization of the terminal
limit group of one of the graded resolutions, $GRes_1,\ldots,GRes_d$, that  already satisfies one of the associated Diophantine conditions 
within the
 sculpted (or developing) resolution  in which it was constructed, which is a contradiction to the way sculpted and developing resolutions
are constructed in the general step of the iterative procedure.

Therefore, there exists a global bound on the number of possible sculpted resolutions that can be associated with a 
quotient resolution along the iterative procedure.  This global bound implies the termination of the iterative procedure for the analysis of an
EAE set over free products, since as we have already pointed out, 
 in case it doesn't terminate, there is no bound on the number of associated sculpted resolutions. Therefore, the iterative procedure terminates 
after finitely many steps. 

\line{\hss$\qed$}

By theorem 4.4 the procedure for the analysis of an EAE set terminates after finitely many steps. When it terminates we are left
with a quotient resolution, that has all the properties of a cover approximating (graded) resolution that are listed in theorem
1.21,  and a finite collection of auxiliary resolutions, that do all satisfy the properties of the resolutions in the
formal Makanin-Razborov diagram (theorems 2.6 and 2.7)

This terminating iterative procedure that allows one to associate an approximating cover   graded resolution,  and its finite
collection of auxiliary resolutions (that can all be extended to
ungraded resolutions with f.p.\ completions),
with a subsequence of
any given sequence of pairs, $\{(p_n,w_n)\}$, of specializations in free products $G_n$ for which $p_n \in EAE(p)$ over $G_n$, and $w_n$
are witnesses for $p_n$, enables one to construct finitely many such approximating cover  resolutions with their (finitely many)
auxiliary resolutions, that can be associated with the set
$EAE(p)$. In the case of an EAE sentence, the construction of these  (finitely many) approximating and auxiliary  resolutions,
allows one to reduce (uniformly) the given EAE sentence over free products to a (finite)  disjunction of conjunctions of sentences over the factors
of the free product.

\vglue 1pc
\proclaim{Theorem 4.5} Let:
$$EAE(p) \ = \ \exists w \  \forall y \ \exists x \
 \Sigma(x,y,w,p)=1 \, \wedge \, \Psi(x,y,w,p) \neq 1$$ 
be an EAE set  over groups. Then there exist finitely many (coefficient-free) 
graded resolutions over free products (with respect to the
parameter subgroup $<p>$):
 $WRes_1(z,y,w,p),\ldots,WRes_e(z,y,w,p)$, with the following properties:

\roster
\item"{(1)}"each of  the graded resolutions, $WRes_i$, is a graded cover resolution that satisfies the properties of cover
graded resolutions that are listed in theorem 1.21. In particular, each of the graded 
resolutions, $WRes_i$, can be continued to an ungraded
resolution with a f.p.\ completion. With the rigid or weakly solid factor of the terminal limit group of $WRes_i$ there are
finitely many associated covers of the flexible quotients of that rigid or weakly solid factor, and each of these  covers can be embedded in an 
ungraded resolution with a f.p.\ completion (see theorem 1.21).

\item"{(2)}" with each graded resolution, $WRes_i(z,y,w,p)$, we associate (non-canonically) finitely many graded 
(auxiliary) resolutions, precisely in the same way that we have associated such auxiliary resolutions with an approximating
cover resolution, $WRes_m$, in the first step of the iterative procedure, in case the iterative procedure terminates in
the first step.  These auxiliary graded resolutions have the same properties as the resolutions in a formal graded
Makanin-Razborov diagram (theorems 2.6 and 2.7), and they are constructed from test sequences of each of the graded
resolutions, $WRes_i$, $i=1,\ldots,e$, that can be extended to specializations of the terminal limit group of one of
the graded resolutions, $GRes_1,\ldots,GRes_d$, and test sequences that can be extended to specializations that satisfy
one of the properties (1)-(6), that are listed in constructing the auxiliary resolutions that are associated with
the resolution, $WRes_m$, in the first step of the iterative procedure (in case this procedure terminates in the first 
step). We further associate auxiliary resolutions that are associated with the singular loci of the resolutions, $WRes_i$, (see proposition
4.2) and 
auxiliary resolutions that are graded Root resolutions that are associated with the graded resolutions, $WRes_i$, and the (finite) set
of closures that are associated with all the previously constructed auxiliary resolutions (see proposition 4.3).

\item"{(3)}" let $G=A^1*\ldots*A^{\ell}$ be a non-trivial free product, which is not isomorphic to
$D_{\infty}$. Let $p_0$ be a specialization
of the parameters $p$ in the free product $G$, for which $p_0 \in EAE(p)$ over $G$.
Then there exists
an index $i$, $1 \leq i \leq e$,  a rigid or a weakly strictly solid specialization of the rigid or weakly solid factor 
of the terminal limit group
of $WRes_i$ in $G$ (that restricts to $p_0$), and specializations of the elliptic factors of that terminal limit group 
in elliptic factors in $G$, 
so that for the combined  specialization of the terminal limit group of
the graded resolution, $WRes_i$,   there  is a test sequence in the fiber that is associated with
this specialization of the terminal limit group (of $WRes_i$), that either can  not be extended to 
specializations of the terminal
limit group of one of the graded resolutions, $GRes_1,\ldots,GRes_d$, or if it can be extended to such specializations,
then either their restrictions to the rigid or weakly solid factor in the terminal limit group of
$GRes_1,\ldots,GRes_d$, is not rigid or not weakly solid, or the specializations of the terminal limit group of
$GRes_1,\ldots,GRes_d$, can be extended to specializations of the terminal limit group of 
one of the graded resolutions, $SLGRes_i^v$, that are associated with
the singular locus of the graded resolutions, $GRes_1,\ldots,GRes_d$, or
the specializations of the terminal limit group
of $GRes_1,\ldots,GRes_d$, can be further extended to specializations of the terminal limit groups of
 one or more formal graded resolutions,
$FGRes_i^j$, that are associated with, $GRes_1,\ldots,GRes_d$, so that the associated fibers of the formal resolutions, $FGRes_i^j$, 
form a covering closure of the corresponding fiber of the
graded resolution $GRes_i$, and these specializations of the terminal limit groups of the formal resolutions, $FGRes_i$, 
can not be further extended to
the terminal limit groups of the collapse resolutions, $CollFRes_i^{j,k}$, that are associated with the formal
graded resolutions, $FRes_i^j$, and not to specializations of the terminal limit group of one of the graded resolutions,
$SLFGRes_i^{j,v}$, that are associated with the singular locus of the formal resolutions $FGRes_i^j$  
(see the construction of the auxiliary resolutions that are associated with $WRes_m$ in
the first step of the iterative procedure for the analysis of an EAE set, and in particular properties (1)-(6)).
\endroster

In other words, given an EAE set,
there exists a finite collection of graded resolutions (that can be extended to ungraded resolutions with f.p.\ completions) 
over free products, and finitely many associated auxiliary resolutions of these graded resolutions,
so that the inclusion of  a specialization of the parameters (free variables)  in the EAE set over any given free product,
can be demonstrated by  a generic point in at least one of these resolutions, and (possibly) their associated 
auxiliary resolutions.
\endproclaim

\nfp The argument that we use is similar to the one that was used to prove theorems 4.1  and 3.2, that is based on the
arguments that were used in constructing the ungraded and graded Makanin-Razborov diagrams (theorem 26 in [Ja-Se] and theorem
1.22).

\noindent
Let $G_n=A^1_n*\ldots*A^{\ell}_n$ (possibly for varying $\ell>1$), be a sequence of non-trivial free products 
that are not isomorphic to $D_{\infty}$. Let
$\{(p_n,w_n)\}$ be a sequence of tuples in $G_n$, so that $p_n \in EAE(p)$ over $G_n$, and $w_n$ is a witness for $p_n$,
i.e., $(p_n,w_n)$ is in the associated set $AE(p,w)$ over $G_n$, which means that the following sentence:
$$ \forall y \ \exists x \
 \Sigma(x,y,w_n,p_n)=1 \, \wedge \, \Psi(x,y,w_n,p_n) \neq 1$$ 
is true over $G_n$.

Starting with the sequence, $\{(p_n,w_n)\}$, the terminating iterative procedure for the analysis of an EAE set 
that we have presented
constructs a graded resolution, $WRes$, with the following
properties:
\roster
\item"{(i)}" the graded resolution $WRes$ satisfies the properties of a cover graded resolution that are 
listed in theorem 1.21. In particular, its completion can be extended to (an ungraded) f.p.\ completion. With the rigid
or weakly solid factor of the terminal  limit group of $WRes$, there is a finite collection of covers of its
flexible quotients that can all be embedded into f.p.\ completions.

\item"{(ii)}" there exists a subsequence of the sequence of specializations, $\{(p_n,w_n)\}$, that extend to 
specializations 
that factor through the resolution, $WRes$, and to specializations of the f.p.\ completion of the ungraded resolution that
extends the graded resolution, $WRes$ (see theorem 1.21). 
Hence, with each specialization from this subsequence of the tuples, 
$\{(p_n,w_n)\}$, 
specializations of the  terminal limit group of the graded resolution, $WRes$, can be associated. Furthermore, 
these specializations of the terminal limit group of $WRes$
 can be extended to specializations of the f.p.\ completion of some ungraded resolution of the terminal limit group of
$WRes$. Restrictions of generic points in the fibers that are associated with these specializations of the terminal limit group
of $WRes$, to the variables $(p,w)$, that we denote, $(p_n,w)$, are in the set $AE(p,w)$ that is associated with the given 
EAE set, $EAE(p)$. i.e., for every index $n$, and for generic $w$ in the fiber that is associated with $(p_n,w_n)$ (i.e., for the
restriction of some test sequence
in this fiber to the existential variables $w$),
 the sentence:
$$ \forall y \ \exists x \
 \Sigma(x,y,w,p_n)=1 \, \wedge \, \Psi(x,y,w,p_n) \neq 1$$ 
is true over $G_n$.

\item"{(iii)}" with the resolution $WRes$,  we associate (non-canonically) a finite collection of graded
auxiliary resolutions, that are part of the output of the terminating iterative procedure for the analysis of an
EAE set, i.e., auxiliary resolutions of the same type as those that were associated with an approximating cover resolution,
$WRes_m$, in the first step of the procedure, in case it terminates in the first step. These auxiliary resolutions have the
same properties of the resolutions in a formal graded Makanin-Razborov diagram (theorems 2.6 and 2.7), and in particular
they can be extended to ungraded resolutions with f.p.\ completions.
\endroster

Now we can apply the argument that was used to prove theorems 3.2 and 4.1.
We look at all the sequences of non-trivial free products, $G_n=A^1_n*\ldots*A^{\ell}_n$
(possibly for varying $\ell>1$), that are not isomorphic to $D_{\infty}$, an associated sequence of 
 tuples, $p_n \in EAE(p)$ over $G_n$, and witnesses $w_n$ for $p_n$, i.e., a sequence of pairs $(p_n,w_n) \in AE(p,w)$
over $G_n$.
From every such  sequence we use our terminating iterative procedure, and extract a subsequence
of the tuples, (still denoted) $\{(p_n,w_n)\}$, and a graded resolution,  
$WRes$, that has the properties (i)-(iii), and in particular,   the subsequence of tuples, (still denoted) $\{(p_n,w_n)\}$, 
extend to specializations that factor through the resolution $WRes$, and to specializations of the f.p.\ completion of
the ungraded resolutions that extends the graded resolution $WRes$.
Furthermore, the restrictions of generic points in the fibers that are associated with the pairs, $\{(p_n,w_n)\}$,
are in the set $AE(p,w)$ (i.e., generic values of $w$ in these fibers are witnesses for $p_n \in EAE(p)$
over $G_n$).

The completion of each of the constructed resolutions, $WRes$, can be extended to an ungraded f.p.\ completion, 
and so are its associated (finitely many) auxiliary resolutions.
Furthermore, the (finite collection of) covers of flexible quotients
of the rigid or weakly solid factors of the terminal limit groups of the resolutions, $WRes$, and of the auxiliary resolutions
that are associated with $WRes$, 
can all be embedded into (ungraded) resolutions with f.p.\
completions. Hence, we can define a linear  order on this (countable) collection of
graded resolutions ($WRes$), and their (non-canonically) associated auxiliary resolutions, and (finite collections of)
covers of flexible quotients.  By the same argument that
was used in constructing the Makanin-Razborov diagram (theorem 26 in [Ja-Se]), there exists a finite subcollection
of these graded resolutions that satisfy properties (1)-(3) of the theorem.

\line{\hss$\qed$}

As we have already pointed out, the existence of finitely many graded resolutions that satisfy the  properties
that are listed in theorem 4.5, and their (finitely many) auxiliary resolutions, 
allows one to reduce an EAE sentence from free products to a sentence over its factors in a uniform way.

\vglue 1pc
\proclaim{Theorem 4.6} Let:
$$  \exists w \ \forall y \ \exists x \
 \Sigma(x,y,w)=1 \, \wedge \, \Psi(x,y,w) \neq 1$$ 
be an EAE  sentence over groups. 

Then there exists a coefficient-free 
sentence over free products, 
which is a (finite) disjunction of conjunctions of EAE sentences, where each of these last EAE sentences involves elliptic 
elements from the same factor,
such that for every non-trivial free product, $G=A^1*\ldots*A^{\ell}$, that is not isomorphic to $D_{\infty}$,
the original
EAE sentence over the free product $G=A^1*\ldots*A^{\ell}$ is a truth sentence, if and only if the sentence which is a 
(finite) disjunction of conjunctions of EAE sentences over the factors $A^1,\ldots,A^{\ell}$  is a truth sentence.
\endproclaim 

\nfp
 By theorem 4.5, with a given EAE set, it is possible to associate finitely many graded resolutions,
$WRes_1,\ldots,WRes_e$, that do all satisfy the properties that are listed in  theorem 1.21, and in particular they can all be embedded into
f.p.\ completions.  
With each graded resolution, $WRes_i$, we have associated 
finitely many auxiliary resolutions, according to parts (1)-(6) that appear in the first step of the iterative procedure for the analysis of
an EAE set, as well as finitely many graded resolutions that are associated with the singular locus of $WRes_i$ (according to proposition 
4.2), and finitely many graded Root resolutions that are associated with $WRes_i$ and its collection of auxiliary resolutions according
to proposition 4.3. Note that all the (finitely many) graded resolutions that are associated with $WRes_i$ have the properties of
the resolutions in the formal graded Makanin-Razborov diagram as listed in theorems 2.6 and 2.7, and in particular they can all be
embedded into f.p.\ completions.

In case of an EAE sentence, the same constructions that enable one to associate  graded  resolutions and their
auxiliary resolutions with an EAE set (i.e., the iterative procedure for the analysis of an EAE set, and the proof of theorem 4.5),
enable one to associate with a given EAE sentence a (non-canonical) finite collection of (ungraded) resolutions 
(over free products) with f.p.\ completions, 
and with each resolution finitely many (ungraded) auxiliary resolutions that do all have f.p.\ completions, and these auxiliary resolutions
have the same properties and they are constructed in the same way as the (graded) auxiliary resolutions that are constructed in the case of
an EAE set. 

We (still) denote the (ungraded) resolutions that are associated with the given EAE sentence, $WRes_1,\ldots,WRes_e$.
By theorems 4.5 and 3.2, and propositions 4.2 and 4.3, the given EAE sentence is true over a non-trivial free product, 
$G=A^1*\ldots*A^{\ell}$,
that is not isomorphic to $D_{\infty}$, if and only if
there exists a specialization in $G$ of  the terminal limit group of a resolution
$WRes_i$ (one of the resolutions, $WRes_1,\ldots,WRes_e$, i.e.,
specializations of the elliptic factors of the terminal limit group of $WRes_i$ in the factors, $A^1,\ldots,A^{\ell}$, for which:
\roster
\item"{(1)}" the specialization does not extend to a specialization of the terminal limit group of one of the auxiliary resolutions,
that are associated with the singular locus of $WRes_i$ according to proposition 4.2.
The specialization extends to the terminal limit group of some prescribed Root resolutions (possibly only the trivial roots), 
that are associated with $WRes_i$ according to theorem 4.3, and does not extend to the other Root resolutions.

\item"{(2)}" the specialization does not extend to 
a (finite) collection of specializations of  the  terminal limit groups of the auxiliary resolutions that are constructed from test
sequences of $WRes_i$ 
that extend to specializations of the terminal limit groups of the graded resolutions, $GRes_1,\ldots,GRes_d$, that are associated
with the AE set, $AE(w)$:
$$  \forall y \ \exists x \
 \Sigma(x,y,w)=1 \, \wedge \, \Psi(x,y,w) \neq 1$$ 
that satisfy:

\itemitem{(i)} these specializations extend to specializations of the terminal limit groups of a prescribed finite collection (possibly empty) 
of auxiliary resolutions 
(that are associated with $WRes_i$) that are constructed  from test sequences for which  the additional specializations of $GRes_1,\ldots,GRes_d$,
are either collapsed, i.e., the restriction of the additional specializations to the rigid or weakly solid factor are not rigid nor weakly
strictly solid, 
or they can be extended to either the the terminal limit groups of resolutions that are associated with the singular locus of the corresponding resolution,
$GRes_1,\ldots,GRes_d$, or to the terminal limit groups of resolutions that are associated with the graded Root resolutions that are
associated with the resolutions, $GRes_1,\ldots,GRes_d$. The specializations extend only to specializations of the terminal limit groups of 
such auxiliary resolutions from the prescribed set, and not to specializations of terminal limit groups of auxiliary resolutions that do
not belong to the prescribed set. 

\itemitem{(ii)} these specializations extend to specializations of the terminal limit groups of a prescribed finite collection (possibly empty) 
of auxiliary resolutions 
(that are associated with $WRes_i$) that are constructed  from test sequences for which  the additional specializations of $GRes_1,\ldots,GRes_d$,
extend to specializations of the terminal limit groups of formal graded resolution, $FGRes_i^j$, that are associated with
the graded resolutions, $GRes_i$.
The specializations extend only to specializations of the terminal limit groups of 
such auxiliary resolutions from the prescribed set, and not to specializations of terminal limit groups of auxiliary resolutions that do
not belong to the prescribed set. 

\itemitem{(iii)} the extended specializations of the terminal limit groups of the auxiliary resolutions that are
associated with the formal graded resolutions, $FGRes_i^j$, extends further  to specializations of the terminal limit groups of a 
prescribed finite collection (possibly empty) 
of auxiliary resolutions 
(that are associated with $WRes_i$) that are constructed  from test sequences for which  the additional specializations of the terminal limit
groups of the formal graded resolutions, $FGRes_i$,
are either collapsed, i.e., the restriction of the additional specializations to the rigid or weakly solid factor are not rigid nor weakly
strictly solid, 
or they can be extended to either the the terminal limit groups of resolutions that are associated with the singular locus of the 
corresponding graded formal resolution, $FGRes_i^j$, 
or to the terminal limit groups of one of the graded resolutions, $\Psi FGRes_i^{j,k}$.
Again, the specializations extend only to specializations of the terminal limit groups of 
such auxiliary resolutions from the prescribed set, and not to specializations of terminal limit groups of auxiliary resolutions that do
not belong to the prescribed set.

\itemitem{(iv)} the fibers that are associated with the specializations of the terminal limit groups of auxiliary resolutions from
the prescribed (finite) sets, form a covering closure of the fiber that is associated with the specialization of the terminal limit
group of $WRes_i$.
\endroster

Finally, by going over the finitely many possibilities for prescribed sets of auxiliary resolutions that satisfy part (iv),
the existence of  a specialization of the terminal limit group of one of the resolutions, $WRes_i$, that satisfies
properties (1)-(2), is clearly a finite disjunction of finite conjunctions of EAE sentences over the factors of $G$, $A^1,\ldots,A^{\ell}$,
i.e., conditions (i)-(iv) can be easily written as a disjunction of conjunctions of such sentences in the factors by using all the
constructed auxiliary resolutions.

\line{\hss$\qed$}

\vglue 1.5pc
\centerline{\bf{\S5. 
AEAE sentences and predicates}}
\medskip

In  the previous two sections we used the iterative procedure for the analysis of an AE sentence over free groups [Se4], and the sieve
procedure that was used for quantifier elimination over free groups [Se6], to analyze AE and EAE predicates and sentences. In
particular we showed that AE and EAE sentences over free products can be reduced to (finite) disjunction of conjunctions of sentences
over the factors of the free product. In this section we further continue towards similar statements for general predicates and
sentences. We use what we proved for AE and EAE predicates, to analyze AEAE predicates and sentences. We expand the structure
that we are working with, by adding new quantifiers over elliptic factors, and prove that an AEAE set can be defined using a
predicate that requires only 3 quantifiers over the ambient free product, and additional quantifiers over elements that are in
elliptic factors. The tools that are used in this reduction enable us to reduce an AEAE sentence over the ambient free product
to a disjunction of conjunctions of sentences over the elliptic factors. These results and the tools and techniques that
are used in their proof, are the basis for our general results for sentences and predicates (with an arbitrary (finite) number of quantifiers)
over free products, and therefore they are the key for all our results on the first order theory of free products of groups.

Let:
$$  AEAE(p) \ = \ \forall t \ \exists w \ \forall y \ \exists x \
 \Sigma(x,y,w,t,p)=1 \, \wedge \, \Psi(x,y,w,t,p) \neq 1$$ 
be a predicate over groups. 
With the set $AEAE(p)$ we naturally associate an EAE set:
$$  EAE(p,t) \ = \  \exists w \ \forall y \ \exists x \
 \Sigma(x,y,w,t,p)=1 \, \wedge \, \Psi(x,y,w,t,p) \neq 1$$ 
and with the sets, $AEAE(p)$ and $EAE(p,t)$, we naturally associate an AE set:  
$$  AE(p,t,w) \ = \   \forall y \ \exists x \
 \Sigma(x,y,w,t,p)=1 \, \wedge \, \Psi(x,y,w,t,p) \neq 1.$$

Recall that with the set $AE(p,t,w)$ we have associated (in theorem 4.1) finitely many graded resolutions,
$GRes_1,\ldots.GRes_d$, with respect to the parameter subgroup $<p,t,w>$, that do all satisfy the properties of graded
cover resolutions that are listed in theorem 1.21. With these graded resolutions we have associated graded resolutions that are
associated with their singular locus (proposition 4.2), graded Root resolutions (proposition 4.3), and graded
formal resolutions, $FGRes_i^j$, with which we associated graded resolutions that are associated with their singular locus,
and collapse graded formal resolutions,
$CollFGRes_i^{j,k}$.  A specialization $(p_0,t_0,w_0) \notin AE(p,t,w)$,
over some non-trivial free product $G=A^1*\ldots *A^{\ell}$, which is not isomorphic to
$D_{\infty}$, 
if and only if there is a specialization of the terminal graded limit group of one of the graded resolutions,
$GRes_1,\ldots,GRes_d$, that extends the specialization $(p_0,t_0,w_0)$, a specialization 
that is composed from a rigid or a weakly strictly solid specialization of the rigid or solid factor of that 
graded limit group,
and specializations of the f.p.\ elliptic factors of the terminal graded limit group, 
so that the restrictions to the universal variables $y$ 
of generic elements (i.e., a test sequence) in the fiber that is associated with the specialization of
the terminal limit group, are witnesses that $(p_0,t_0,w_0) \notin AE(p,t,w)$.

With the EAE set, $EAE(p,t)$, we associated finitely many graded resolution (in proposition 4.3) that we denoted, $WRes_1,\ldots,WRes_e$,
that satisfy the properties of theorem 1.21.
With each of these graded resolutions we have associated finitely
many auxiliary resolutions that satisfy the properties of the resolutions in the formal Makanin-Razborov diagram (theorems
2.6 and 2.7). Both the graded resolutions, $GRes_i$, and their associated auxiliary resolutions can be
extended to ungraded resolutions with f.p.\ completions. 
By part (3) of theorem 4.5,  a specialization $(p_0,t_0) \in EAE(p,t)$
over some non-trivial free product $G=A^1*\ldots *A^{\ell}$ which is not isomorphic to $D_{\infty}$,
if and only if there is a specialization of the terminal graded limit group of one of the graded resolutions, 
$WRes_1,\ldots,WRes_e$, that extends $(p_0,t_0)$, 
so that a generic point (i.e., a test sequence) in the fiber that is associated with this specialization of the terminal limit group, 
restricts to specializations of the existential variables $w$ that are witnesses that  $(p_0,t_0) \in EAE(p,t)$.

Both in the analysis of AE sets, and in the
analysis of EAE sets, the existence of finitely many graded resolutions with the properties listed above enabled us to reduce a
2 quantifier and a 3 quantifier sentence to a finite disjunction of conjunctions of sentences over the factors of 
the free product (theorems 3.7 and 4.6).

In this section we use our results on AE and EAE sets over free products to analyze AEAE sets over free products, and in particular 
we associate finitely many graded resolutions, that satisfy the properties of theorem 1.21, and finitely many 
auxiliary resolutions with each of these graded resolutions (that satisfy the properties of formal graded resolutions
that are listed in theorem 2.6), with any given
AEAE set. The construction of these graded and auxiliary resolutions uses the structural results that were proved 
in sections 1 and 2, and is once
again a variation of the sieve procedure that was presented in [Se6].

We start the analysis of the set $AEAE(p)$ with all the sequences of specializations of the tuple, $(p,t)$, $\{(p_n,t_n)\}$, that take values
in non-trivial free products, $G_n=A^1_n*\ldots *A^{\ell}_n$ (for possibly varying $\ell>1$),  which are not
isomorphic to $D_{\infty}$. We further assume that for every index $n$,
$p_n \notin AEAE(p)$ over the
free product, $G_n$, and that $t_n$ is a witness for $p_n$, i.e., a specialization of the universal variables $t$ in $G_n$, so
that the EAE sentence:
$$   \exists w \ \forall y \ \exists x \
 \Sigma(x,y,w,t_n,p_n)=1 \, \wedge \, \Psi(x,y,w,t_n,p_n) \neq 1$$ 
is a false sentence over $G_n$.

By theorem 1.16, 
given such a sequence of pairs, 
 $\{(p_n,t_n)\}$, we can pass to a subsequence that converges into a 
well-structured  (even well-separated) graded resolution with respect to the parameter subgroup $<p>$: 
$T_0 \to T_1 \to \ldots \to T_s$, where $T_s$ is
 a free product of a rigid or a solid group over free products (that contains the subgroup, $<p>$), 
and (possibly) a free group and (possibly) finitely many elliptic factors. 
We denote this graded resolution
$TRes$. As we noted in similar constructions in the previous sections, the terminal limit group $T_s$   in $TRes$ is f.g.\ 
but it may be not finitely presented.

We apply the techniques of theorem 1.21, and fix a sequence of cover approximating graded 
resolutions of the resolution, $TRes$,
that we denote, $\{TRes_m\}$, that satisfy the properties of cover approximating resolutions that are listed in theorem 1.21,
and converge into the graded resolution, $TRes$. For each index $m$, there is a subsequence of the sequence of pairs,
$\{(p_n,t_n)\}$, that factors through the cover approximating (graded) resolution, $TRes_m$.

For each index $m$, and  each  pair, $(p_n,t_n)$, that factors through $TRes_m$,  
we can associate a specialization
of the terminating limit groups of the (approximating) resolution, $TRes_m$ (in the free product $G_n$). 
Such specialization is composed from a rigid or  a weakly strictly
solid specialization of the rigid or weakly solid factor of the terminal limit group of $TRes_m$, and specializations 
of the f.p.\ elliptic factors of that
terminal graded limit group.

As in the analysis of an EAE set in section 4, we first  assume that there exists an approximating resolution, $TRes_m$, 
for which there exists a subsequence of pairs, (still
denoted) $\{(p_n,t_n)\}$, so that for each specialization of the terminating limit group of $TRes_m$ that is associated with  
a pair $(p_n,t_n)$ from the subsequence, a generic  pair in the fiber (i.e., a test sequence) that is associated with
such a specialization  of the terminal graded limit group of $TRes_m$ (which  values in $G_n$) 
satisfies $p_n \notin AEAE(p)$
(over $G_n$), and a generic $t$ in the fiber is a witness for that. i.e., for generic $t$ in the fiber, the pairs
$(p_n,t) \notin EAE(p,t)$ over $G_n$.

In this case with the approximating (cover) graded resolution, $TRes_m$, that has the properties that are listed in theorem 1.21,
we associate finitely many auxiliary resolutions, that have similar properties to the resolutions in the formal graded
Makanin-Razborov resolution (theorems 2.6 and 2.7). These  auxiliary resolutions are similar in nature to the
auxiliary resolutions that were associated in section 4 with an EAE set, and with a cover approximating resolution, $WRes_m$
(see the first step of the iterative procedure for the analysis of an EAE set in section 4). The auxiliary resolutions that are
associated with the approximating cover, $TRes_m$, are constructed from test sequences of $TRes_m$ that have the following
properties: 
\roster
\item"{(1)}" we start by associating (non-canonically) finitely many graded resolutions with the singular locus of $TRes_m$ 
(according to proposition 4.2), that we denote $SLTRes$, and finitely many graded Root resolutions (according to proposition
4.3), that we denote $GRootTRes$.
 
\item"{(2)}" we  look at all the test sequences of $TRes_m$, that can be extended to specializations of 
the terminal limit group of one of the cover graded resolutions, $WRes_i$, that are associated with the EAE set, $EAE(p,t)$, and
to the f.p.\ completion into which this terminal limit group embeds. By the construction of the formal graded
Makanin-Razborov diagram, with this collection of test sequences it is possible to associate finitely many graded resolutions
that have the same properties as the resolutions in the formal Makanin-Razborov diagram (theorems 2.6 and 2.7). In particular,
the completion of each of the constructed resolutions embeds into a f.p.\ completion. We denote the constructed 
resolutions, $WTRes$.

\item"{(3)}" we look at test sequences of the graded resolutions that are constructed in part (2), for which 
the restriction of the additional specializations of the terminal limit group of the resolution, $WRes_i$, to a specialization of the
rigid or weakly solid factor of that terminal limit group, is non-rigid or  non weakly solid. By the construction of
the formal Makanin-Razborov diagram, it is possible to associate finitely many graded (cover) resolutions with these test
sequences that have the properties of the resolutions in the formal Makanin-Razborov diagram. We denote the constructed
resolutions, $CollWTRes$.

Similarly, we look at test sequences of the graded resolutions that are constructed in part (2), for which the restrictions of the
additional 
specializations of the terminal limit group of the resolution, $WRes_i$, can be extended to specializations of the terminal limit
group of one the auxiliary resolutions that are associated with the singular locus of $WRes_i$, or of one of the graded Root resolutions that
are associated with $WRes_i$. We denote the obtained auxiliary resolutions, $SLWTRes$ and $GRootWTRes$, in correspondence.

\item"{(4)}"  we look at test sequences that are constructed in part (2),  that can be further extended to specializations 
of the terminal limit group of one of the auxiliary resolutions that are associated with the graded resolution, $WRes_i$, and was constructed
from test sequences that can be extended to the terminal limit group of one of the graded resolutions, $GRes_1,\ldots,GRes_d$,
that are associated with the AE set, $AE(p,t,w)$. We further require that the specialization of the terminal limit
group of this auxiliary resolution can be extended to the f.p.\ completion that is associated with the auxiliary resolution,
and the terminal limit group (of the auxiliary resolution) is embedded into it. Once again with this collection of sequences
it is possible to associate finitely many graded resolutions that have the properties of the resolutions in a formal
graded Makanin-Razborov diagram (theorems 2.6 and 2.7). We denote the constructed resolutions, $YWTRes$.

\item"{(5)}" we look at test sequences that were considered in part (4) for which the additional specialization of the
terminal limit group of the associated auxiliary resolution, restricts to non-rigid non-weakly strictly solid 
specialization of the rigid or weakly solid factor of the terminal limit group of the auxiliary resolution. With
this collection of sequences we associate finitely many graded resolutions following the construction of the formal 
graded Makanin-Razborov diagram. We denote the constructed resolutions, $CollYWTRes$.

Similarly, we look at test sequences of the graded resolutions that are constructed in part (4), for which the restrictions of the
additional 
specializations of the associated auxiliary  resolution  of $WRes_i$, can be extended to specializations of the terminal limit
group of one the auxiliary  resolutions that are associated with $WRes_i$ and with the singular locus of $GRes_i$, or with
 one of the graded Root resolutions that
are associated with $GRes_i$. We denote the obtained auxiliary resolutions, $SLYWTRes$ and $GRootYWTRes$, in correspondence.

\item"{(6)}"  we look at test sequences that are constructed in part (4),  that can be further extended to specializations 
of the terminal limit group of one of the auxiliary resolutions that are associated with the graded resolution, $WRes_i$, and were constructed
from test sequences that can be extended to the terminal limit group of one of the graded formal resolutions, 
$FGRes_i^j$, that are associated with the graded resolution $GRes_i$, and both 
are associated with the AE set, $AE(p,t,w)$. We further require that the specialization of the terminal limit
group of this auxiliary resolution can be extended to the f.p.\ completion that is associated with the auxiliary resolution,
and the terminal limit group (of the auxiliary resolution) is embedded into it. With this collection of
sequences we, once again, associate finitely many resolutions, in a construction that follows the construction of the
formal graded Makanin-Razborov diagram (theorem 2.7). We denote the constructed resolutions, $XYWTRes$.

\item"{(7)}" we look at test sequences that were considered in part (6) for which the additional specialization of the
terminal limit group of the associated auxiliary resolution (that is associated with $FGRes_i^j$), 
restricts to non-rigid non-weakly strictly solid 
specialization of the rigid or weakly solid factor of the terminal limit group of the auxiliary resolution. With
this collection of sequences we associate finitely many graded resolutions, following the construction of the formal
graded Makanin-Razborov diagram. We denote the constructed resolutions, $CollXYWTRes$.

Similarly, we look at test sequences of the graded resolutions that are constructed in part (6), for which the restrictions of the
additional 
specializations of the terminal limit group of one the auxiliary resolutions that are associated with $WRes_i$ and one of the formal
graded resolution, $FGRes_i^j$, can be extended to specializations of the terminal limit
group of one the auxiliary resolutions that are associated with $WRes_i$ and with the singular locus of $FGRes_i^j$. 
We denote the obtained auxiliary resolutions $SLXYWTRes$.

\item"{(8)}"  we look at test sequences that are constructed in part (6),  that can be further extended to specializations 
of the terminal limit group of one of the auxiliary resolutions that are associated with the graded resolution, $WRes_i$, and were constructed
from test sequences that can be extended to the terminal limit group of one of the collapse graded formal resolutions, 
$CollFGRes_i^{j,k}$, that are associated with the graded formal resolution, $FGRes_i^j$, that is associated with the graded
resolution, $GRes_i$, and all the three 
are associated with the AE set, $AE(p,t,w)$. We further require that the specialization of the terminal limit
group of this auxiliary resolution can be extended to the f.p.\ completion that is associated with the auxiliary resolution,
and the terminal limit group (of the auxiliary resolution) is embedded into it. With this collection of
sequences we associate finitely many resolutions, in a construction that follows the construction of the
formal graded Makanin-Razborov diagram (theorem 2.7). We denote the constructed resolutions, $\Psi XYWTRes$, as well. 

\item"{(9)}" we look at test sequences that were considered in part (8) for which the additional specialization of the
terminal limit group of the associated auxiliary resolution (that is associated with $CollFGRes_i^{j,k}$), 
restricts to non-rigid non-weakly strictly solid 
specialization of the rigid or weakly solid factor of the terminal limit group of the auxiliary resolution. With
this collection of sequences we associate finitely many graded resolutions, following the construction of the formal
graded Makanin-Razborov diagram. We denote the constructed resolutions, $Coll \Psi XYWTRes$.
\endroster
Hence, with the cover approximating graded resolution, $TRes_m$, that satisfies the properties of theorem 1.21, we associated
finitely many auxiliary resolutions that do all satisfy the properties of resolutions in the formal
Makanin-Razborov diagram (theorems 2.6 and 2.7).

\medskip
Suppose that there is no 
approximating cover resolution, $TRes_m$, with a subsequence of the original sequence, $\{(p_n,t_n)\}$ (for which
$p_n \notin AEAE(p)$ (over $G_n$), and $t_n$ (a specialization of the universal variables $t$ in $G_n$) is a witness for that
(in $G_n$)), that factor through $TRes_m$, so that
for each specialization of the terminal limit group of $TRes_m$ that is associated with  
a pair $(p_n,t_n)$ from the subsequence, a generic  pair (i.e., a test sequence) in the fiber that is associated with
such a specialization  of the terminal graded limit group (which takes its values in $G_n$), satisfies $p_n \notin AEAE(p)$
(over $G_n$), and $t$ (a generic value of $t$ in the associated fiber) is a witness for that.

As in analyzing EAE sets in sections 4,  in this case (in which there is no approximating
cover $TRes_m$ with the desired properties) 
we can associate with the original 
graded resolution, $TRes$, another graded resolution, which is obtained by applying the first step of the sieve procedure
for the analysis of quotient resolutions [Se6]. The new graded resolution has a smaller complexity than the original graded
resolution $TRes$, and the original sequence of specializations, 
$\{(p_n,t_n)\}$, is guaranteed to
have a subsequence that extends to specializations that converge into the obtained quotient resolution (which is of smaller complexity).
This enables us to continue iteratively, in a similar way to the sieve procedure [Se6], and to the iterative procedure that 
was used to prove the 
equationality of Diophantine sets over free and hyperbolic groups in section 2 of [Se8].

With the resolution, $TRes$, we have associated a sequence of approximating cover resolutions, $TRes_m$, that
satisfy the properties that are listed in theorem 1.21. In particular, for each index $m$, there is a subsequence
of the original sequence of pairs, $\{(p_n,t_n)\}$, that factor through $TRes_m$. Recall that by our assumptions,
for every index $m$, there is no subsequence of the  pairs, $\{(p_n,t_n)\}$, that do
factor through the approximating cover, $TRes_m$, for which  for generic pair $(p_n,t)$ (i.e., for a restriction of a test sequence
to the universal variables $t$) in the fiber that
contains the pair $(p_n,t_n)$ (in the graded variety that is associated with $TRes_m$), 
$p_n \notin AEAE(p)$ over $G_n$, and the generic 
$t$ (in the fiber) is a witness for that, i.e., the generic pairs, $(p_n,t) \notin EAE(p,t)$ over $G_n$.

We construct a sequence of specializations over $G_n$ as follows.
We go over the indices $m$, and given an index  $m$ we pick an index $n_m > m$, and a tuple, 
$(p_{n_m}, \hat t_{n_m}, \tilde w_{n_m}, \tilde t_{n_m})$, that takes its values in $G_{n_m}$,
with the following properties:
\roster
\item"{(1)}" $(p_{n_m},t_{n_m})$ factors through the graded resolution, $TRes_m$, and $\hat t_{n_m}$ is the 
specialization of the terminal limit group of the graded resolution, $TRes_m$, that contains the pair, $(p_{n_m},t_{n_m})$.

\item"{(2)}" $\tilde t_{n_m}$ is a specialization of the variables $t$, from the fiber that is associated with
$\hat t_{n_m}$, i.e., the fiber that contains the pair $(p_{n_m},t_{n_m})$ in the graded variety that is associated with
$TRes_m$. Furthermore, for every index $m$, $(p_{n_m}, \tilde t_{n_m}) \in EAE(p,t)$ over $G_n$.

\item"{(3)}" the sequence $\{(p_{n_m},\tilde t_{n_m}) \}$ is a graded test sequence that converges into the graded resolution, 
$TRes$. 

\item"{(4)}" the tuples $(p_{n_m},\tilde t_{n_m}, \tilde w_{n_m})$, are specializations of the terminal limit group of
one of the  finitely many graded resolutions, $WRes$ (that are associated with the set $EAE(p,t)$), that do not extend to  specializations
of the terminal limit groups of the graded resolutions that are associated with the singular locus of $WRes_m$, and these
specializations restrict to  rigid  or  weakly 
strictly solid specializations of the rigid or weakly solid factor of that terminal limit group 
(with respect to the parameter group $<p,t>$). 
Furthermore, a generic value of the existential variables $w$ in the fiber of $WRes$ that is associated with the tuple,
$(p_{n_m},\tilde t_{n_m}, \tilde w_{n_m})$,  
is a witness that:
$(p_{n_m}, \tilde t_{n_m}) \in EAE(p,t)$ over $G_{n_m}$.
\endroster

From the sequence of tuples,
$(p_{n_m}, \hat t_{n_m}, \tilde w_{n_m}, \tilde t_{n_m})$, we can extract a subsequence that converges into a closure of
the graded resolution, $TRes$, that we have started with. We denote this closure, $ExWTRes$. 
Recall that we have assumed that the original sequence of pairs,
$\{(p_n,t_n)\}$, satisfies $p_n \notin AEAE(p)$, and $t_n$ is a witness for that (i.e., $(p_n,t_n) \notin EAE(p,t)$).
On the other hand, in the sequence of tuples that
we chose (in the fibers that are associated with a subsequence of the original sequence
of specializations), a generic value of the existential variables $w$ in the fiber that is associated with the tuple,
$(p_{n_m},\tilde t_{n_m}, \tilde w_{n_m})$,  
is a witness that the pair $(p_{n_m}, \tilde t_{n_m}) \in EAE(p,t)$.

Given the closure, $ExWTRes$, we continue in a similar way to what we did in analyzing EAE sets in section 4, a way that is
adapted to sets with 4 quantifiers. The closure, $ExWTRes$, contains an (additional) rigid or weakly strictly solid specialization of the 
terminal 
existential variables $w$, 
that for generic
value of the universal variables $t$ (i.e., for restrictions of a test sequence to the variables $t$) are witnesses that 
the pairs, $(p_{n_m}, t) \in EAE(p,t)$. 
Hence, we apply the construction that was used in theorem 1.21, and associate  a sequence of approximating covers 
of $ExWTRes$ that
satisfy the properties that are listed in theorem 1.21. We denote this sequence, $\{ExWTRes_r\}$. 

By construction, for each index $r$, it is possible to extract a subsequence of the original sequence of tuples:
$(p_{n_m}, \hat t_{n_m}, \tilde w_{n_m}, \tilde t_{n_m})$, that factor through $ExWTRes_r$. The tuples $(p_{n_m},t_{n_m})
\notin EAE(p,t)$, whereas the tuples $(p_{n_m}, \tilde t_{n_m}) \in EAE(p,t)$, and generic values (a test sequence) 
of the existential variables
$w$ in the fiber that is associated with the tuple,
$(p_{n_m},\tilde t_{n_m}, \tilde w_{n_m})$,  
are witnesses 
for this last
inclusion. Hence, on the tuples of  variables $(p,t,w)$, that are associated with the additional tuples,
$(p_{n_m},\tilde t_{n_m}, \tilde w_{n_m})$, in the approximating cover $ExWTRes_r$, 
 at least one of the following 
collapse forms can be imposed (cf. the list of  test sequences that are associated with the approximating cover $TRes_m$  in the first step):
\roster
\item"{(1)}" the additional specialization of the variables $(p,t,w)$ in $ExWTRes_r$, that is assumed to restrict to rigid or weakly strictly solid
specialization of the rigid or weakly solid factor of the terminal limit group of a resolution, $WRes$ (that is associated with $EAE(p,t)$)
restricts to  non-rigid or non weakly 
strictly solid specialization of that factor. This forces a Diophantine condition on the specializations of $ExWTRes_r$, similar to the
ones that are forced along the sieve procedure in [Se6].

\item"{(2)}" the additional specialization of the variables $(p,t,w)$ in $ExWTRes_r$, extends to a specialization of 
the terminal limit group of one of the graded
resolutions that are associated with the singular locus of the corresponding resolution, $WRes_i$.
This  forces a Diophantine condition on the specializations of $ExWTRes_r$. 

\item"{(3)}" we look at all the test sequences of $ExTWRes_r$, for which  the specialization of the variables 
$(p,t,w)$ in  $ExWTRes_r$, can be extended to  a specialization of the terminal limit group of an 
auxiliary resolution that is associated with graded Root resolution of $WRes$ (and with the set $EAE(p,t)$), and not to a specializations
of an auxiliary resolutions that is associated with graded Root resolutions of $WRes$ of higher order roots, 
and for which the specializations of the variables, $(p,t,w)$, can be extended to specializations of 
the terminal limit groups of a prescribed set of auxiliary resolutions 
that are associated with the  resolution
$WRes$ (and with the set $EAE(p,t)$),  auxiliary resolutions
that were constructed from test sequences of $WRes$  that can be extended to specializations of the terminal limit group of one of
the graded resolutions, $GRes_1,\ldots,GRes_d$, or with the formal graded resolutions, $FGRes^i_j$, or with their associated graded Root
resolutions, that are associated with $AE(p,t,w)$. 
It is further required that the prescribed set of 
auxiliary resolutions that are associated with the graded resolutions, $GRes_1,\ldots,GRes_d$, form a covering closure of the 
graded resolution $WRes$, but that they do not form a covering closure if we take out the auxiliary resolutions that are associated with
the graded formal resolutions, $FGRes_i^j$. 

By the techniques for the construction of the
graded formal Makanin-Razborov diagram (theorem 2.7), with this collection of test sequences and their extended specializations
(for all possible graded root resolutions and prescribed sets of auxiliary resolutions),
it is possible to associate finitely many graded resolutions that have the properties of resolutions in the formal
Makanin-Razborov diagrams (theorems 2.6 and 2.7), and they are all closures of $ExWTRes_r$. We denote each of these finitely many
graded resolutions, $XYExWTRes_r^s$.

On each of the closures, $XYExWTRes_r^s$, we further impose Diophantine conditions (in parallel). 
We require that the  additional specializations
of the terminal limit group of one of the  auxiliary resolutions of $WRes$ that is associated with $FGRes_i^j$ in $XYExWTRes_r^s$, 
extend to  specializations of the terminal limit group of either an auxiliary resolution that is associated with the singular
locus of the corresponding resolution, $FGRes_i^j$, or to the terminal limit group of an auxiliary resolution that is associated with
a resolution, $\Psi FGRes_i^{j,k}$, or to an auxiliary resolution that is associated with specializations of the terminal
limit group of a corresponding graded resolution, $FGRes^i_j$, that do not restrict to rigid or weakly strictly solid specializations of
the rigid or weakly solid factor of the terminal limit group of $FGRes_i^j$, or that the specializations of the terminal limit group of 
such an auxiliary resolution of $WRes$ do
not restrict to rigid or weakly strictly solid specializations  of the rigid or weakly solid factor of their terminal limit group.

Another possible Diophantine restrictions on $XYExWTRes$ is that the specializations of the terminal limit group of the graded resolution
$WRes$, extend to specializations of the terminal limit group of an auxiliary resolution of $WRes$, that is associated with
a  graded Root resolution of $WRes$ or with a graded Root resolution of one of the graded resolutions, $GRes_1,\ldots,GRes_d$,
from the prescribed set,  with roots of higher orders than the prescribed ones.

\item"{(4)}" we further look at all the specializations of the graded resolution, $ExWTRes_r$, for which to the specialization of
the the variables
$(p,t,w)$ in  $ExWTRes_r$, it is possible to add  additional specialization of
the terminal limit group of one of the finitely many  auxiliary resolutions that are associated with $WRes$ (and with $EAE(p,t)$),  
auxiliary resolutions
that were constructed from test sequences of $WRes$  that can be extended to specializations of the terminal limit group of one of
the graded resolutions, $GRes_1,\ldots,GRes_d$, that are associated with $AE(p,t,w)$, and the combined specialization does not
extend to a specialization that factors through (or in the same weakly strictly solid family of a specialization in) 
one of the closures, $XYExWTRes_r^s$, that  were constructed in part (3). We further assume that these additional 
specializations
of the terminal limit group of an auxiliary resolution of $WRes$ (that is associated with $GRes_i$), restrict to rigid or weakly
strictly solid specializations of the rigid or weakly solid factor of the terminal limit group of the auxiliary resolution, and they 
can not be extended to specializations of the terminal limit groups of auxiliary resolutions that are associated with the singular locus of 
$GRes_i$, or of auxiliary resolutions that are associated with specializations of the terminal limit group of $GRes_i$ that do
not restrict to rigid or weakly strictly solid specializations of the rigid or weakly solid factor of the terminal limit group of
$GRes_i$.
\endroster

For each index $r$,  we can find an index $n_r > r$, so that the specialization, $(p_{n_r},t_{n_r})$ (in $G_{n_r}$),
extends to a specialization that satisfies the  conditions that are specified by one the parts, (1)-(4) (for
$ExWTRes_r$).
We further apply the 
analysis of quotient resolutions that appear in the first step of the sieve procedure [Se6], and
pass to a subsequence of the extended specializations (of the sequence  $\{(p_{n_r},t_{n_r})\}$, that converges into a
quotient resolution of smaller complexity (in terms of the sieve procedure) than the complexity of the original resolution, $TRes$ (note that the
obtained quotient resolution may be a proper closure of $TRes$, but there is a global bound (depending on $TRes$) on the number of steps
for which the obtained quotient resolution is a closure of $TRes$. After a number of steps that is bigger than the bound, 
the obtained quotient resolution
is no longer a closure, and hence, there is a reduction in other parts of the complexity of the obtained quotient resolution).

We continue iteratively. At each step we first look at a sequence of approximating resolutions of a quotient resolution
that was constructed in the previous step, where these approximating resolutions satisfy the properties that are
listed in theorem 1.21.
If there exists an approximating resolution for which there exists a subsequence
of pairs, (still
denoted) $\{(p_n,t_n)\}$, so that for each specialization of the terminating limit group of the approximating resolution
 that is associated with  
a pair $(p_n,t_n)$ from the subsequence, a generic  pair (i.e. a restriction of a test sequence to the pair $(p,t)$) 
in the fiber that is associated with
such a specialization  of the terminal graded limit group (which takes its values in $G_n$) forms a witness for  $p_n \notin AEAE(p)$
(over $G_n$),  we do what we did in this case in the first step (i.e., we look at all the test sequences of the approximating
resolutions that satisfy properties (1)-(9) and associate finitely many auxiliary resolutions with it). 

If there is no such approximating resolutions for a constructed quotient resolution, we associate with it a closure
that is  constructed from test sequences that can be extended to additional (new) specialization 
  of the terminal limit group of one of
the resolutions, $WRes$ (that are associated with $EAE(p,t)$), in a similar way to the construction of the graded closure, $ExWTRes$,
in the first step of the procedure.
We further force a collapse condition over the constructed 
closure of the quotient resolution, according to parts (1)-(4), and apply the general step of the sieve procedure [Se6], 
extract a subsequence of the original sequence of tuples, $\{(p_n,t_n)\}$, and construct a quotient resolution
that has smaller complexity than the quotient resolution that was constructed in the previous step (smaller complexity in light
of the sieve procedure).

By an argument which is similar to the argument that guarantees the termination of the sieve procedure in [Se6], and to the termination of the
procedure for the analysis of EAE sets,  we
obtain a termination of the this procedure.

\vglue 1pc
\proclaim{Theorem 5.1} The procedure for the analysis of an AEAE set over free products terminates after finitely
many steps.
\endproclaim

\nfp The argument that we use is similar to the proof of theorem 4.4 (termination of the procure for the analysis of an EAE set).
At each step of the procedure, on  the quotient resolution, $TRes$, that is analyzed in that step  we impose one of finitely
many  restrictions. The first ones impose a non-trivial Diophantine condition on specializations of a closure of the resolution, $TRes$. 
This is similar
to the Diophantine conditions that are imposed on quotient resolutions in the general step of the iterative procedure for
the analysis of an EAE set. Hence, by the termination of this iterative procedure (theorem 4.4),
these type of restrictions can occur only at finitely 
many steps along the iterative procedure for the analysis of an AEAE set.

The second type of restrictions, adds specializations of the terminal limit groups of one of the associated auxiliary resolutions that are
associated with the graded resolution, $WRes$, where these specializations restrict to rigid or weakly strictly
solid specializations of the rigid or weakly solid factor of the terminal limit group. We further require that the additional
specializations do not factor through (or are not in the same (weak) strictly solid families as specializations of) closures that are
built from test sequences of $TRes$ that can be extended to specializations of such terminal limit groups of auxiliary resolutions.

Therefore, the argument that was used to prove the termination of the iterative procedure for the analysis of an EAE set
(theorem 4.4), proves that restrictions of the
second type can also occur in only finitely many steps (this type of restrictions are  precisely the ones that are analyzed 
in proving the equationality of Diophantine sets
(over free and hyperbolic) groups in [Se9]). Since the restrictions of both the first and the second type can occur at only finitely many
steps, the iterative procedure for the analysis of an AEAE set terminates after finitely many steps.

\line{\hss$\qed$}

As we argued for AE and EAE sets, the termination of the procedure for the analysis of AEAE sets, enables one to associate
with the given AEAE set a finite collection of graded resolutions, that do all satisfy the properties that are listed in theorem 1.21,
and for each of these graded resolutions an associated finite collection of auxiliary resolutions, that do all satisfy the 
properties of resolutions in the formal Makanin-Razborov diagram (theorems 2.6 and 2.7).

This finite collection of graded resolutions, enable us to obtain a form of quantifier elimination
for AEAE sets over free products, i.e., it enables one to show that an AEAE set can be defined by a predicate that uses only 3 quantifiers
over the free product, and additional quantifiers over the various factors.  
In the case of an AEAE sentence, the construction of these  (finitely many) graded resolutions (and their auxiliary resolutions),
allows one to reduce the given AEAE sentence over free products to a (finite) 
disjunction of conjunctions of sentences over the factors
of the free product. These statements are generalized to arbitrary definable sets and predicates in the next section, 
and they are the 
key for all the results in this paper.

\vglue 1pc
\proclaim{Theorem 5.2} Let:
$$AEAE(p) \ = \ \forall t \ \exists w \  \forall y \ \exists x \
 \Sigma(x,y,w,t,p)=1 \, \wedge \, \Psi(x,y,w,t,p) \neq 1$$ 
be an AEAE set  over groups. Then there exist finitely many graded resolutions over free products (with respect to the
parameter subgroup $<p>$):
 $$TRes_1(z,y,w,t,p),\ldots,TRes_g(z,y,w,t,p)$$ with the following properties:

\roster
\item"{(1)}" The graded resolutions, $TRes_i(z,y,w,t,p)$, satisfy the properties that are listed in theorem 1.21. In particular, with
each of them there is an associated f.p.\ completion, into which the completion of them  embeds.

\item"{(2)}" with each graded resolution, $TRes_i(z,y,w,t,p)$, we associate (non-canonically) a finite collection of
auxiliary resolutions, that do all have the properties of resolutions in the formal graded Makanin-Razborov diagram 
(theorems 2.6 and 2.7). These auxiliary resolutions are constructed from test sequences of $TRes_i$, 
precisely in the way we associated
auxiliary resolutions with an approximating cover, $TRes_m$, in the first part of the procedure, i.e., according to parts
(1)-(9) in the first step. 
We denote these auxiliary resolutions: $WTRes$, $CollWTRes$, $YWTRes$, $CollYWTRes$, $XYWTRes$, $CollXYWTRes$, $\Psi XYWTRes$,
and
$Coll \Psi XYWTRes$, and their associated graded Root resolutions and resolutions that are associated with the corresponding
singular locus (see parts (1)-(9) in the first step of the procedure).

\item"{(3)}" let $G=A^1*\ldots*A^{\ell}$ be a non-trivial free product, which is not isomorphic to $D_{\infty}$.
Let $p_0$ be a specialization
of the parameters $p$ in the free product $G$, for which $p_0 \notin AEAE(p)$ over $G$.
Then there exists
an index $i$, $1 \leq i \leq g$,  a specialization of the terminal limit group of $TRes_i$ (in $G$) that restricts to $p_0$, 
that extends to a specialization
of the f.p.\ completion into which the completion of $TRes_i$ embeds, and restricts to a 
rigid or a weakly strictly solid  specialization of the rigid or weakly solid factor of the terminal limit group
of $TRes_i$, such that one of the following holds: 
\itemitem{(i)}  the specialization of the terminal limit group of $TRes_i$, does not extend to a specialization of the
 terminal limit group of any graded resolutions, $SLTRes$, that are associated with the singular locus of $TRes$. It may extend
to specializations of the terminal limit groups of some of the auxiliary resolutions, $WTRes$, but the corresponding fibers do not
cover the fiber of $TRes$ that is associated with the given specialization.

\itemitem{(ii)} the specialization of the terminal limit group of $TRes_i$ does not extend to the terminal limit group of any of
the resolutions, $SLTRes$, but it extends to specializations of the terminal limit groups of some graded resolutions, $WTRes$,
that do not extend to specializations of the terminal limit groups of the auxiliary resolutions, $SLWTRes$ and $CollWTRes$, 
and the corresponding fibers do cover the fiber of $TRes$ that is associated with the given specialization.

In that case, for each extension of the specialization of the terminal limit group of $TRes_i$, to a specialization of the terminal
limit group of one of the auxiliary resolutions, $WTRes$, there can be  further extensions to specializations of the terminal limit
groups of auxiliary resolutions, $YWTRes$. Furthermore, these specializations of the terminal limit group of $YWTRes$ may 
extend to specializations of the terminal limit groups of auxiliary resolutions, $XYWTRes$. For every finite collection of
extensions of the specialization
of the terminal limit group of $TRes_i$ to the terminal limit groups of (finitely many)  auxiliary resolutions $WTRes$, 
the fibers that are associated with these extensions (and with the auxiliary resolutions $WTRes$), minus the fibers that are associated 
with further extensions to specializations of the terminal limit groups of auxiliary resolutions $YWTRes$, that are not covered by fibers that are
associated with further extensions to the terminal limit groups of auxiliary resolutions $XYWTRes$ (that can not be extended
to fibers of auxiliary resolutions $CollXYWTRes$, $\Psi XYWTRes$, and $SLXYWTRes$), do not cover the fiber that is
associated with the original specialization of the terminal limit group of $TRes_i$. 
\endroster

In other words, given an AEAE set,
there exists a finite collection of graded resolutions  over free products (that satisfy the properties that are listed
in theorem 1.21), and a finite collection of auxiliary resolutions (that satisfy the properties that are listed in theorem 2.6),
so that the inclusion of  a specialization of the parameters (free variables) in the complement of the given  AEAE set over 
any given free product, which is not $D_{\infty}$,
can be demonstrated by  a generic point (test sequences) in a  disjunction of conjunctions of (fibers of)  these resolutions.
\endproclaim

\nfp The argument that we use is similar to the one that was used to prove theorems 4.5,4.1  and 3.2, that is based on the
arguments that were used in constructing the ungraded and graded Makanin-Razborov diagrams (theorem 26 in [Ja-Se] and theorem
1.22).

\noindent
Let $G_n=A^1_n*\ldots*A^{\ell}_n$ (possibly for varying $\ell>1$), be a sequence of non-trivial free products 
that are not isomorphic to $D_{\infty}$. Let
$\{(p_n,t_n)\}$ be a sequence of tuples in $G_n$, so that $p_n \notin AEAE(p)$ over $G_n$, and $t_n$ is a witness for $p_n$,
i.e., $(p_n,t_n) \notin EAE(p,t)$  over $G_n$.

Starting with the sequence, $\{(p_n,t_n)\}$, the terminating iterative procedure for the analysis of an AEAE set 
that we have presented,
constructs a graded resolution, $TRes$, with the following
properties:
\roster
\item"{(i)}" the graded resolution $TRes$ satisfies the properties of a cover graded resolution that are 
listed in theorem 1.21. 
In particular, its completion can be extended to (an ungraded) f.p.\ completion. 
With the rigid
or weakly solid factor of the terminal  limit group of $WRes$, there is a finite collection of covers of its
flexible quotients that can all be embedded into f.p.\ completions.

\item"{(ii)}" there exists a subsequence of the sequence of specializations, $\{(p_n,t_n)\}$, that extend to 
specializations 
that factor through the resolution, $TRes$, and to specializations of the f.p.\ completion of the ungraded resolution that
extends the graded resolution, $TRes$ (see theorem 1.21). 
Hence, with each specialization from this subsequence of the tuples, 
$\{(p_n,t_n)\}$, 
specializations of the  terminal limit group of the graded resolution, $TRes$, can be associated. Furthermore, 
these specializations of the terminal limit group of $TRes$
 can be extended to specializations of the f.p.\ completion of some ungraded resolution of the terminal limit group of
$TRes$. Restrictions of generic points in the fibers that are associated with these specializations of the terminal limit group
of $TRes$, to the variables $(p,t)$, that we denote, $(p_n,t)$, are not in the set $EAE(p,t)$ that is associated with the given 
AEAE set, $AEAE(p)$.

\item"{(iii)}" with the resolution $TRes$,  we associate (non-canonically) a finite collection of graded
auxiliary resolutions, according to parts (1)-(9) of the first step of the iterative procedure for the analysis of an AEAE set. These
auxiliary resolutions are part of the output of this terminating procedure.
The auxiliary resolutions have the
same properties of the resolutions in a formal graded Makanin-Razborov diagram (theorems 2.6 and 2.7), and in particular
they can be extended to ungraded resolutions with f.p.\ completions.
\endroster

Now we can apply the argument that was used to prove theorems 3.2, 4.1 and 4.5.
We look at all the sequences of non-trivial free products, $G_n=A^1_n*\ldots*A^{\ell}_n$
(possibly for varying $\ell>1$), that are not isomorphic to $D_{\infty}$, an associated sequence of 
 tuples, $p_n \notin AEAE(p)$ over $G_n$, and witnesses $t_n$ for $p_n$, i.e., a sequence of pairs $(p_n,t_n) \notin EAE(p,t)$
over $G_n$.
From every such  sequence we use our terminating iterative procedure, and extract a subsequence
of the tuples, (still denoted) $\{(p_n,t_n)\}$, and a graded resolution,  
$TRes$, that has the properties (i)-(iii), and in particular,   the subsequence of tuples, (still denoted) $\{(p_n,t_n)\}$, 
extend to specializations that factor through the resolution $TRes$, and to specializations of the f.p.\ completion of
the ungraded resolutions that extends the graded resolution $TRes$.
Furthermore, the restrictions of generic points in the fibers that are associated with the pairs, $\{(p_n,t_n)\}$,
are not in the set $EAE(p,t)$ (i.e., generic values of $t$ in these fibers are witnesses for $p_n \notin AEAE(p)$
over $G_n$).

The completion of each of the constructed resolutions, $TRes$, can be extended to an ungraded f.p.\ completion, 
and so are its associated (finitely many) auxiliary resolutions.
Furthermore, the (finite collection of) covers of flexible quotients
of the rigid or weakly solid factors of the terminal limit groups of the resolutions, $TRes$, and of the auxiliary resolutions
that are associated with $TRes$, 
can all be embedded into (ungraded) resolutions with f.p.\
completions. Hence, we can define a linear  order on this (countable) collection of
graded resolutions ($TRes$), and their (non-canonically) associated auxiliary resolutions, and (finite collections of)
covers of flexible quotients.  By the same argument that
was used in constructing the Makanin-Razborov diagram (theorem 26 in [Ja-Se]), there exists a finite subcollection
of these graded resolutions that satisfy properties (1)-(3) of the theorem.

\line{\hss$\qed$}

As we argued for AE and EAE sets, the existence of a finite collection of graded resolutions and their auxiliary resolutions,
with the properties that are listed in
theorem 5.2, 
allows one to reduce an AEAE sentence from free products to a sentence over its factors in a uniform way.

\vglue 1pc
\proclaim{Theorem 5.3} Let:
$$  \forall t \ \exists w \ \forall y \ \exists x \
 \Sigma(x,y,w,t)=1 \, \wedge \, \Psi(x,y,w,t) \neq 1$$ 
be an AEAE  sentence over groups. 

Then there exists a coefficient-free 
sentence over free products, 
which is a (finite) disjunction of conjunctions of AEAE sentences, where each of these last AEAE sentences involves elliptic 
elements from the same factor,
such that for every non-trivial free product, $G=A^1*\ldots*A^{\ell}$, that is not isomorphic to $D_{\infty}$, 
the original
AEAE sentence over the free product $G$ is a truth sentence, if and only if the sentence which is a (finite)
disjunction of conjunctions of AEAE sentences over the factors,  $A^1,\ldots,A^{\ell}$  is a truth sentence.
\endproclaim 

\nfp The argument that we use, that is based on theorem 5.2, is similar to the proof of theorem 4.6 (which is based on
theorem 4.5). 
 By theorem 5.2, with a given AEAE set, it is possible to associate finitely many graded resolutions,
$TRes_1,\ldots,TRes_g$, that do all satisfy the properties that are listed in  theorem 1.21, and in particular they can all be embedded into
f.p.\ completions.  
With each graded resolution, $TRes_i$, we have associated 
finitely many auxiliary resolutions, according to parts (1)-(9) that appear in the first step of the iterative procedure for the analysis of
an AEAE set.
Note that all the (finitely many) graded resolutions that are associated with $TRes_i$ have the properties of
the resolutions in the formal graded Makanin-Razborov diagram as listed in theorems 2.6 and 2.7, and in particular they can all be
embedded into f.p.\ completions.

In case of an AEAE sentence, the same constructions that enable one to associate  graded  resolutions and their
auxiliary resolutions with an AEAE set (i.e., the iterative procedure for the analysis of an AEAE set, and the proof of theorem 5.2),
enables one to associate with a given AEAE sentence a (non-canonical) finite collection of (ungraded) resolutions 
(over free products) with f.p.\ completions, 
and with each resolution finitely many (ungraded) auxiliary resolutions that do all have f.p.\ completions, and these auxiliary resolutions
have the same properties and they are constructed in the same way as the (graded) auxiliary resolutions that are constructed in the case of
an AEAE set.

We (still) denote the (ungraded) resolutions that are associated with the given AEAE sentence, $TRes_1,\ldots,TRes_g$.
By theorems 5.2, 4.5 and 3.2, and propositions 4.2 and 4.3, the given AEAE sentence is false over a non-trivial free product, 
$G=A^1*\ldots*A^{\ell}$,
that is not isomorphic to $D_{\infty}$, if and only if
there exists a specialization in $G$ of  the terminal limit group of a resolution
$TRes_i$ (one of the resolutions, $TRes_1,\ldots,TRes_g$), i.e.,
specializations of the elliptic factors of the terminal limit group of $TRes_i$ in the factors, $A^1,\ldots,A^{\ell}$, so that
 the specialization does not extend to a specialization of the terminal limit group of one of the auxiliary resolutions
that are associated with the singular locus of $TRes_i$ (see proposition 4.2). 
This specialization of the terminal limit
group of $TRes_i$ also satisfies: 
\roster
\item"{(1)}" it extends to the terminal limit group of an auxiliary resolution that is associated with a  Root resolution 
of $TRes_i$ according to proposition 4.3
(possibly only the trivial roots), 
and does not extend to specializations of the terminal limit group of
other Root resolutions, $RootTRes$, that are associated with higher order roots.

\item"{(2)}" it may extend to  specializations of the terminal limit groups of some of the auxiliary resolutions
$WTRes$, that do not extend to specializations of the terminal limit groups that are associated with the singular locus of the
associated graded resolution, $WRes$. Either the collection of the fibers that are associated with these specializations of
the terminal limit groups of $WTRes$ do not form a covering closure of the fiber that is associated with the specialization
of the terminal limit group of $TRes_i$, or (some of) these specializations can be extended to the specializations of the terminal
limit groups of the following auxiliary resolutions (see pars (1)-(9) in the first step for the construction and the properties of the 
auxiliary resolutions that we refer to):
\itemitem{(i)} each of the specializations of the terminal limit groups of $WTRes$ extends to a specialization of the terminal
limit group of one of the auxiliary resolutions, $GRootWTRes$, and not to 
similar auxiliary resolutions that are associated with roots
of higher orders.

\itemitem{(ii)} some of these specializations extend to specializations of the terminal limit groups of auxiliary resolutions,
$YWTRes$, that do not extend to specializations of the terminal limit groups of auxiliary resolutions,
$CollYWTRes$ and $SLYWTRes$.

Each such specialization extends to a specialization of the terminal limit group of an auxiliary resolution, $GRootYWTRes$, and not
to the terminal limit groups of such resolutions that are associated with higher order roots.

\itemitem{(iii)} some of the specializations in (ii) can further extend to specializations of the terminal limit groups of auxiliary
resolutions, $XYWTRes$, that do not further extend to specializations of the terminal limit groups of auxiliary resolutions,
$SLXWTRes$, $CollXYWTRes$, or $\Psi XYWTRes$.

\itemitem{(iv)} the fibers that are associated with the specializations of the terminal limit groups of auxiliary resolutions $WTRes$,
minus the fibers that are associated with their extensions to  specializations of the terminal limit groups of auxiliary resolutions $YWTRes$,
from which we take out fibers that are covered by fibers that are associated with further extensions to
 specializations of the terminal limit groups of auxiliary resolutions $XYWTRes$, that do not extend to terminal limit groups of auxiliary
resolutions: 
$SLXWTRes$, $CollXYWTRes$, or $\Psi XYWTRes$, do not form a cover of the fiber that is associated with the original specialization of the 
graded resolution $TRes_i$.
\endroster

Finally, by going over the finitely many possibilities for prescribed sets of auxiliary resolutions that satisfy part (iv),
the existence of  a specialization of the terminal limit group of one of the resolutions, $TRes_i$, that satisfies
properties (1)-(2), is clearly a finite disjunction of finite conjunctions of EAEA sentences over the factors of $G$, $A^1,\ldots,A^{\ell}$,
i.e., conditions (i)-(iv) can be easily written as a disjunction of conjunctions of such sentences in the factors by using all the
constructed auxiliary resolutions. Hence, the given AEAE set is a disjunction of conjunctions of AEAE sentences over the factors 
$A^1,\ldots,A^{\ell}$.

\line{\hss$\qed$}

In analyzing AE and EAE sets, we have associated finitely many graded resolutions with these sets, and using them we were
able to reduce an AE or an EAE sentence over free products to disjunctions of conjunctions of sentences over the factors.
The graded resolutions that we constructed for analyzing AEAE sets, enable us to reduce not only AEAE (and EAEA)
sentences over free products, but also AEAE (and EAEA) predicates. As we will see in the next section, this kind of reduction,
or quantifier elimination, hold for arbitrary predicates (or definable sets) over free products, and can be viewed as
a form of quantifier elimination over free products.

\vglue 1pc
\proclaim{Theorem 5.4} Let:
$$AEAE(p) \ = \   \forall t \ \exists w \ \forall y \ \exists x \
 \Sigma(x,y,w,t,p)=1 \, \wedge \, \Psi(x,y,w,t,p) \neq 1$$ 
be an AEAE  sentence over groups. 

Then there exists a coefficient-free predicate over groups that are free products, that is composed from only 3 quantifiers over 
variables that take values in the ambient
free product, and additional quantifiers over variables that take their values in the various factors of the free product, 
so that for every 
non-trivial free product, $G=A^1*\ldots*A^{\ell}$, that is not isomorphic to $D_{\infty}$, the set $AEAE(p)$ over $G$ can be defined by the following
predicate over free products:
$$AEAE(p) \ = \   \exists u \  (\forall t \ \exists w \ \forall y \ \exists x)  \ \forall v \ \exists s (\ \exists e)$$
$$ (\Sigma_1(x,y,w,t,u,v,s,p)=1 \, \wedge \, \Psi_1(x,y,w,t,u,v,s,p) \neq 1) \, \vee \, \ldots \, \vee$$
$$ (\Sigma_{k}(x,y,w,t,u,v,s,p)=1 \, \wedge \, \Psi_{k}(x,y,w,t,u,v,s,p) \neq 1)$$
where the variables $u,v,s$ take values in the ambient free product $G$, and the variables $t,w,y,x,e$ take values in the factors 
$A^1,\ldots,A^{\ell}$.  

Furthermore, an EAEA set, $EAEA(p)$, can be defined over every non-trivial  free product $G=A^1*\ldots*A^{\ell}$, that is not isomorphic to $D_{\infty}$,
by the following predicate:
$$EAEA(p) \ = \   \exists u \  (\exists t \ \forall w \ \exists y \forall x)  \ \forall v \ \exists s (\ \exists e)$$
$$ (\Sigma_1(x,y,w,t,u,v,s,p)=1 \, \wedge \, \Psi_1(x,y,w,t,u,v,s,p) \neq 1) \, \vee \, \ldots \, \vee$$
$$ (\Sigma_{k}(x,y,w,t,u,v,s,p)=1 \, \wedge \, \Psi_{k}(x,y,w,t,u,v,s,p) \neq 1)$$
where (as in the AEAE case) the variables $u,v,s$ take values in the ambient free product $G$, and the variables $t,w,y,x,e$ 
take values in the factors $A^1,\ldots,A^{\ell}$. 
\endproclaim 

\nfp Our argument for the analysis of AEAE and EAEA predicates combines theorem 5.3, that proves the reduction of sentences from the ambient
free product to the factors, with theorem 5.2 that associates with an AEAE or EAEA predicates a finite collection of graded resolutions
and their auxiliary resolutions that enable one to perform the reduction of sentences uniformly.

 By theorem 5.2, with a given AEAE set, it is possible to associate finitely many graded resolutions,
$TRes_1,\ldots,TRes_g$, that do all satisfy the properties that are listed in  theorem 1.21, and in particular they can all be embedded into
(ungraded) f.p.\ completions. Furthermore, with the rigid or weakly strictly solid factor of the terminal limit group of each of the
graded resolutions, $TRes_i$, we have associated finitely many (graded limit group) covers of its flexible quotient, 
and each of these covers embeds into an (ungraded)  f.p.\ completion as well (see theorem 1.21).
 
With each graded resolution, $TRes_i$, we have associated 
finitely many auxiliary resolutions, according to parts (1)-(9) that appear in the first step of the iterative procedure for the analysis of
an AEAE set.
Note that all the (finitely many) graded resolutions that are associated with $TRes_i$ have the properties of
the resolutions in the formal graded Makanin-Razborov diagram as listed in theorems 2.6 and 2.7, and in particular they can all be
embedded into f.p.\ completions, and with the rigid or weakly solid factor of their terminal limit groups, there are associated
finitely many covers of its flexible quotients, and these covers are embedded into (ungraded) f.p.\ completions. 

Theorem 5.3 analyzes logically an AEAE sentence over free products, and shows how the graded resolutions, $TRes_i$, can be used to
reduce (uniformly) such an AEAE sentence from the ambient free product to a (finite) disjunction of conjunctions of AEAE sentences
over the factors of the free product.

For the analysis of sentences, we have used ungraded resolution over free products, that terminate in a free product of elliptic
factors and possibly a free group. To analyze predicates we have to use graded resolutions, that terminate in a free product of
a rigid or a weakly solid factor with (possibly) finitely many elliptic factors and possibly a free factor. Hence, to apply the same 
analysis
that was used in the analysis of AEAE sentences to analyze AEAE predicates, we need to further find a way to encode all the rigid
or all the weakly strictly solid specializations of the rigid or weakly solid factors,  and to encode the fact that two specializations
belong to the same rigid or weakly strictly solid family of a rigid or a weakly solid limit group (over free products).

By theorems 1.14, for any given rigid  limit group, $Rgd(x,p)$, and any given (finite) covers of its flexible quotients,
there exist finitely many combinatorial systems of fractions (that depend only on the rigid limit group and its generating
set),
so that for every value $p_0$ of the parameter subgroup $<p>$, we can associate at most one value of the fractions for each combinatorial system (these
values depend on the value of the defining parameters $p_0$), so that every rigid specialization that is associated with $p_0$ can be expressed as
a fixed word (that depends only on the combinatorial system, i.e., depends only on the group and its generating set), in the fractions (that
take their values in the ambient free product), and in finitely many elliptic elements (see theorem 1.14  for the exact statements).
By theorem 1.15 the same holds for almost shortest weakly strictly solid homomorphisms that are associated with the parameters value $p_0$.

These finite collection of combinatorial systems that describe the structure of all the rigid or almost shortest weakly strictly solid 
homomorphisms of rigid and weakly solid limit groups, together with theorems 5.2 and 5.3, allow us to replace an AEAE (or EAEA) predicate
with a predicate that uses only 3 quantifiers on elements from the ambient free product, as described in the statement of
the theorem.

Given an AEAE set, $AEAE(p)$, we can write a predicate of the form: 
$$AEAE(p) \ = \   \exists u \  (\forall t \ \exists w \ \forall y \ \exists x)  \ \forall v \ \exists s (\ \exists e)$$
$$ (\Sigma_1(x,y,w,t,u,v,s,p)=1 \, \wedge \, \Psi_1(x,y,w,t,u,v,s,p) \neq 1) \, \vee \, \ldots \, \vee$$
$$ (\Sigma_{k}(x,y,w,t,u,v,s,p)=1 \, \wedge \, \Psi_{k}(x,y,w,t,u,v,s,p) \neq 1)$$
where the existential variables $u$ represent the fractions in all the possible system of fractions that are associated 
(by theorems 1.14 and 1.15) with the
rigid or weakly solid factors of the terminal limit groups of the graded resolutions, $TRes_1,\ldots,TRes_g$, 
and their auxiliary resolutions, that are associated with
the set $AEAE(p)$ by theorem 5.2. The universal variables $v$ and the existential variables $s$ and $e$ (the variables $e$ are contained in a
factor $A^1,\ldots,A^{\ell}$ of the free product), enable us to guarantee that the
values of the fractions $u$ satisfy the conclusions of theorems 1.14 and 1.15, i.e., enable the covering of  all the families of rigid or
weakly strictly solid families that are associated with any given value of the defining parameters.

The variables $t,w,y,x$ take their values in the various elliptic factors, $A^1, \ldots,A^{\ell}$, of the given free product $G$, and
this part of the predicate is given by the description that appears in theorem 5.3, only that the elements $t,w,y,x$ are either 
specializations of the  elliptic factors of the terminal limit groups of the graded resolution, $TRes_1,\ldots,TRes_g$, and of
their auxiliary resolutions, or they are elliptic elements that appear in the combinatorial description of rigid and almost
shortest weakly strictly solid specializations (given by theorems 1.14 and 1.15) of the rigid or weakly solid factors of
these terminal limit groups. The universal variables $v$ (that get values in the ambient free product) enable us to guarantee that
the existential variables $w$ and $x$ represent rigid or weakly strictly solid specializations. The existential variables  
 $s$ (that get values in the ambient free product) enable us to guarantee that
the universal variables $t$ and $y$ represent rigid or weakly strictly solid specializations, and that two given specializations
of a weakly solid limit group belong to the same family.

Once we have used the variables $u,v,s$ (that take their values in the ambient free product) and the elliptic variables $e$, to
go over all the rigid and weakly strictly solid specializations, the theorem for AEAE sets follows by the proof of theorem 5.3. 
The proof in the EAEA case, is naturally identical.

\line{\hss$\qed$}

\vglue 1.5pc
\centerline{\bf{\S6. Definable sets and sentences over free products}} 
\medskip

In  the previous 3 sections we used the iterative procedure for the analysis of an AE sentence over free groups [Se4], and the sieve
procedure that was used for quantifier elimination over free groups [Se6], to analyze AE, EAE, and AEAE predicates and sentences. In
particular, we showed that AEAE and EAEA sentences over free products can be reduced to a (finite) disjunction of conjunctions of sentences
over the factors of the free product (theorem 5.3), and that every AEAE or EAEA set over free products can be defined by a predicate that 
contains only 3 quantifiers over the ambient free product, and additional quantifiers over elements in the various factors (theorem 5.4). These
were conclusions of the existence of finitely many graded resolutions that satisfy the properties that are listed in theorem 1.21, and
their (finitely many) auxiliary resolutions (that satisfy the properties of resolutions in the formal graded Makanin-Razborov diagram),
so that a specialization $p_0$ of the defining parameters (free variables) is in the AEAE (EAEA) set if and only if a generic
element in a fiber of one of these resolutions is a witness for that, and the existence of such a generic witness can be 
reduced to the terminal limit groups of these graded resolutions, and their auxiliary resolutions, 
and the condition on the terminal limit groups can be expressed by a predicate of
the indicated form.

In this section we use the same techniques that were used to analyze AEAE and EAEA sets in  section 5, to analyze
general definable sets and sentences over free products. First we apply a finite induction process and associate (non-canonically)
a finite collection
of graded resolutions with any given definable set over free products, graded resolutions that have similar properties to the ones
that are associated with an AEAE and EAEA sets (see theorem 5.2). The existence of finitely many graded resolutions with such properties allows
one to deduce that every definable set over free products can be defined by a predicate that contains only 3 quantifiers over the
ambient free product and finitely many quantifiers over elements in the factors of the free product. These resolutions also
enable one to deduce that every sentence over free products can be reduced to a finite disjunction of finite conjunctions of
sentences over the factors over the free product. As we show in the next section, these reductions have somewhat surprising 
(uniformity) corollaries for sentences over free products, and we are sure that they will find quite a few generalizations and further
applications in the 
near future.

Let:
$$  E{(AE)}^k(p) \ = \ \exists t \ \forall y_1 \ \exists x_1 \ \ldots \ \forall y_k \ \exists x_k \
 \Sigma(t,y_1,x_1,\ldots,y_k,x_k,p)=1 \, \wedge $$
$$ \wedge \, \Psi(t,y_1,x_1,y_k,x_k,p) \neq 1$$ 
be a predicate (with 2k+1 quantifiers) over groups. Our goal is to associate with this set a finite collection of graded
resolutions that have the properties of the graded resolutions that are associated with an AEAE set according to theorem 
5.2. 

\vglue 1pc
\proclaim{Theorem 6.1} Let $E{(AE)}^k$ be a definable set over free products.
Then there exist finitely many graded resolutions over free products (with respect to the
parameter subgroup $<p>$):
 $$DRes_1(z,x_k,y_k,\ldots,x_1,y_1,t,p),\ldots,DRes_g(z,x_k,y_k,\ldots,x_1,y_1,t,p)$$ with the following properties:

\roster
\item"{(1)}" the resolutions, $DRes_i$, $i=1,\ldots,g$, satisfy the properties that are listed in theorem 1.21.
In particular, with
each of them there is an associated f.p.\ completion, into which the completion of them  embeds.

\item"{(2)}" with each of the graded resolutions, $DRes_i$, $i=1,\ldots,g$, we associate (non-canonically) a finite collection of
auxiliary resolutions, that do all have the properties of resolutions in the formal graded Makanin-Razborov diagram 
(theorems 2.6 and 2.7). These auxiliary resolutions are constructed from test sequences of $DRes_i$, 
precisely in the way we associated
auxiliary resolutions with the graded resolutions, $TRes_i$, that are associated with an AEAE set (part (2) of theorem 5.2).
These include auxiliary resolutions that are associated with the singular locus of $DRes_i$ according to proposition 4.2, that we denote $SLDRes$, 
and graded resolutions that are associated with
 possible roots of pegs in $DRes_i$, according to proposition 4.3, that we denote $GRootDRes$. Besides these resolutions, the auxiliary 
resolutions that are associated with $DRes_i$ are similar to the ones that are constructed in that case of an AEAE set
in part (2) of theorem 5.2.

We denote these auxiliary resolutions, $Y_1DRes$, $X_1Y_1DRes$, $\ldots$, $X_1Y_1 \ldots X_kY_kDRes$,
$\Psi X_1Y_1 \ldots X_kY_kDRes$, and their associated
collapsed resolutions, 
$CollY_1DRes$, $CollX_1Y_1DRes$, $\ldots$, $CollX_1Y_1 \ldots X_kY_kDRes$, 
$Coll \Psi X_1Y_1 \ldots X_kY_kDRes$, 
that indicate that certain  specializations that are assumed to be rigid or weakly strictly solid, are 
non-rigid or non weakly strictly solid, and similar auxiliary resolutions that are associated with the singular locus of the 
corresponding resolutions,  $SLY_1DRes, \ldots$, and auxiliary resolutions that are associated with possible roots of pegs
in the associated graded resolutions, $GRootY_1DRes,\ldots$. 

\item"{(3)}" let $G=A^1*\ldots*A^{\ell}$ be a non-trivial free product, which is not isomorphic to $D_{\infty}$.
Let $p_0$ be a specialization
of the parameters $p$ in the free product $G$, for which $p_0 \in E{(AE)}^k(p)$ over the free product $G$.
Then there exists
an index $i$, $1 \leq i \leq g$,  a specialization of the terminal limit group of $DRes_i$ (in $G$) that restricts to $p_0$, 
that extends to a specialization
of the f.p.\ completion into which the completion of $DRes_i$ embeds, does not extend to a specialization of an auxiliary resolution
that is associated with the singular locus, and restricts to a 
rigid or a weakly strictly solid  specialization of the rigid or weakly solid factor of the terminal limit group
of $DRes_i$, such that one of the following holds: 

\itemitem{(i)}  the specialization of the terminal limit group of $DRes_i$ does not  extend to specializations of the terminal limit group
of  auxiliary resolutions, $Y_1DRes$, that are associated with $DRes_i$.

\itemitem{(ii)} suppose that the specialization of the terminal limit group of $DRes_i$ extends to specializations of the terminal 
limit groups of some
of the auxiliary resolutions, $Y_1DRes$. Then each such specialization of the terminal limit group of an auxiliary resolution, $Y_1DRes$,
 either extends to  specializations of the terminal limit groups of  auxiliary resolutions, $CollY_1DRes$ and $SLY_1DRes$, so that the corresponding
fibers of these resolutions cover the fiber of $Y_1DRes$, or the specialization of the terminal limit group of $Y_1DRes$ extends to a 
specialization of the terminal limit group of
 an auxiliary resolution, $X_1Y_1DRes$, that does not extend to specializations of the terminal limit groups of any auxiliary resolutions,
$CollX_1Y_1DRes$, $SLX_1Y_1DRes$, and $Y_2X_1Y_1DRes$.

In this case, for each extension of the original specialization (of the terminal limit group of $DRes$) to the terminal limit group of
$Y_1DRes$, the fiber that is associated with such extended specialization (and with $Y_1DRes$), is covered by the fibers that are associated with
extended specializations to the terminal limit groups of auxiliary resolutions: $CollY_1DRes$, $SLY_1DRes$, and $X_1Y_1DRes$ (that do not extend 
to specializations of the terminal limit groups of auxiliary resolutions:
$CollX_1Y_1DRes$, $SLX_1Y_1DRes$, and $Y_2X_1Y_1DRes$).

\itemitem{(iii)} suppose that the specialization of the terminal limit group of $DRes_i$ extends to specializations of the terminal 
limit groups of some
of the auxiliary resolutions, $Y_1DRes$, and this extension does not satisfy part (ii) (i.e., the corresponding fibers are
 not covered by fibers of the auxiliary
resolutions:
$CollY_1DRes$, $SLY_1DRes$, and fibers of an auxiliary resolution $X_1Y_1DRes$ that satisfy the properties that are listed in part (ii)).

In this case we continue the conditions on extensions of specializations iteratively. If we got to an extension of the
specialization of the terminal limit group of $DRes_i$, to a specialization of the terminal limit group of the terminal limit
group of an auxiliary resolution, $Y_jX_{j-1}Y_{j-1} \ldots X_1Y_1DRes$, $j < k$, then  every such specialization extends to 
specializations of the terminal limit groups of 
$CollY_jX_{j-1}Y_{j-1} \ldots X_1Y_1DRes$ ,or of
$SLY_jX_{j-1}Y_{j-1} \ldots X_1Y_1DRes$, or to  specializations of the terminal limit groups of  auxiliary resolutions, 
$X_jY_jX_{j-1}Y_{j-1} \ldots X_1Y_1DRes$, that do not extend to a specialization of the terminal limit group of 
$CollX_jY_jX_{j-1}Y_{j-1} \ldots X_1Y_1DRes$ nor of 
$$SLX_jY_jX_{j-1}Y_{j-1} \ldots X_1Y_1DRes$$ nor of
$Y_{j+1}X_jY_jX_{j-1}Y_{j-1} \ldots X_1Y_1DRes$. Furthermore, the fibers that are associated with these extended specializations (and their
associated auxiliary resolutions) cover the fiber that is associated with the specialization of  the terminal limit group of
$Y_jX_{j-1}Y_{j-1} \ldots X_1Y_1DRes$.

In case we got to an extension of the
specialization of the terminal limit group of $DRes_i$, to a specialization of the terminal limit group of the terminal limit
group of an auxiliary resolution, $Y_kX_{k-1}Y_{k-1} \ldots X_1Y_1DRes$, then  every such specialization extends to a
specialization of the terminal limit group of 
$CollY_kX_{k-1}Y_{k-1} \ldots X_1Y_1DRes$ or 
$SLY_kX_{k-1}Y_{k-1} \ldots X_1Y_1DRes$, 
or to a specialization of the terminal limit group of an auxiliary resolution, 
$X_kY_kX_{k-1}Y_{k-1} \ldots X_1Y_1DRes$, that does not extend to a specialization of the terminal limit group of 
$CollX_kY_kX_{k-1}Y_{k-1} \ldots X_1Y_1DRes$ nor of 
$SL X_kY_kX_{k-1}Y_{k-1} \ldots X_1Y_1DRes$, 
nor to a specialization of the terminal limit group of an auxiliary resolution,
$ \Psi X_kY_kX_{k-1}Y_{k-1} \ldots X_1Y_1DRes$, that does not extend to a specialization of the terminal limit group of
$$Coll \Psi X_kY_kX_{k-1}Y_{k-1} \ldots X_1Y_1DRes$$ 
Furthermore, the fibers that are associated with these extended specializations (and their
associated auxiliary resolutions) cover the fiber that is associated with the specialization of  the terminal limit group of
$Y_kX_{k-1}Y_{k-1} \ldots X_1Y_1DRes$. 
\endroster

In other words, given an $E{(AE)}^k$ set,
there exists a finite collection of graded resolutions (with completions that can be embedded into f.p.\ completions) over free products,
with finitely many auxiliary resolutions (that have the properties of resolutions in the formal Makanin-Razborov diagram),
such that the inclusion of  a specialization of the parameters (free variables) in the $E{(AE)}^k$ set over any given non-trivial
free product (which is not $D_{\infty}$),
can be demonstrated by  a generic point in a disjunction of conjunctions of (fibers of)  these resolutions and their auxiliary resolutions.
\endproclaim

\nfp
We prove theorem 6.1 by induction on the index $k$, that determines the number of quantifiers that are used in defining the set,
$E{(AE)}^k(p)$. Our induction hypothesis is therefore that theorem 6.1 holds for ${(AE)}^k$ sets. Note that
in section 3 and 5 we proved theorem 6.1 for AE and AEAE sets, hence, the induction hypothesis holds for $k=1,2$.

Having the induction hypothesis, we prove theorem 6.1 following our analysis of EAE and AEAE sets.
We start with all the sequences of specializations of the tuple, $(p,t)$, $\{(p_n,t_n)\}$, that take values
in non-trivial free products, $G_n=A^1*\ldots*A^{\ell}$ (possibly for varying $\ell$), that are not isomorphic to $D_{\infty}$. 
We further assume that for every index $n$,
$p_n \in E{(AE)}^k(p)$ over the
free product, $G_n$, and that $t_n$ is a witness for $p_n$, i.e., $(p_n,t_n) \in {(AE)}^k(p,t)$ over $G_n$.

Given such a sequence of pairs, 
 $\{(p_n,t_n)\}$, we can pass to a subsequence that converges into a 
well-structured  (even well-separated) graded resolution with respect to the parameter subgroup $<p>$: 
$T_0 \to T_1 \to \ldots \to T_s$, where $T_s$ is
 a free product of a rigid or a solid group over free products (that contains the subgroup, $<p>$), 
and (possibly) a free group and (possibly) finitely many elliptic factors. 
We denote this graded resolution
$DRes$. The terminal limit group $T_s$   in $DRes$ is f.g.\ 
but it may be  infinitely presented.

Following our arguments from previous sections, we apply the construction that appears in theorem 1.21, and fix a sequence of  
approximating covers of $T_s$, that we denote, $DRes_m$. The approximating covers, $DRes_m$, satisfy the properties that are
listed in theorem 1.21, and they converge to  the resolution $T_s$, precisely like the approximating covers that we chose
in analyzing EAE and AEAE sets in sections 4 and 5.

By the construction of approximating covers, with each approximating cover, $DRes_m$, 
there exists a subsequence of the original sequence of tuples, $\{(p_n,t_n)\}$, that factor through it.
As in the analysis of an AEAE set in section 5, we first  assume that there exists an approximating cover, $DRes_m$, 
for which there exists a subsequence of tuples, (still
denoted) $\{(p_n,t_n)\}$, so that for each specialization of the terminal limit group of $DRes_m$ (in $G_n$), that is associated with  
a tuple $(p_n,t_n)$ from the subsequence, a generic  pair (i.e., a test sequence)  in the fiber that is associated with
such a specialization  of the terminal graded limit group  satisfies $p_n \in E{(AE)}^k(p)$
(over $G_n$), and (a generic) $t$ is a witness for that. 

In this case with the approximating graded resolution, $DRes_m$, we associate finitely many auxiliary resolutions, that do
all have the properties of resolutions in the formal Makanin-Razborov diagram (theorems 2.6 and 2.7), in a similar way to what
we did in the procedure for the analysis of an AEAE set in section 5.
The auxiliary resolutions that are
associated with the approximating cover, $DRes_m$, are constructed from test sequences of $DRes_m$ that have the following
properties: 
\roster
\item"{(1)}" With $DRes_m$ we associate finitely many graded resolutions that are associated with its singular locus
(according to proposition 4.2), that we denote $SLDRes$, and given all the auxiliary resolutions that are constructed below, a finite collection
of graded Root resolutions (proposition 4.3), that we denote $GRootDRes$.

\item"{(2)}" we start by looking at all the test sequences of $DRes_m$, that can be extended to specializations of 
the terminal limit group of one of the cover graded resolutions, $Y_1Res$, that are associated with the  set, ${(AE)}^k(p,t)$, and
to the f.p.\ completion into which this terminal limit group embeds. By the construction of the formal graded
Makanin-Razborov diagram, with this collection of test sequences it is possible to associate finitely many graded resolutions
that have the same properties as the resolutions in the formal Makanin-Razborov diagram (theorems 2.6 and 2.7). In particular,
the completion of each of the constructed resolutions embeds into a f.p.\ completion. We denote the constructed 
resolutions, $Y_1DRes$. 

We further associate with each of these resolutions, (collapse) auxiliary  resolutions, $CollY_1DRes$, 
auxiliary resolutions that are associated with the singular
locus of the resolutions $Y_1Res$, that we denote $SLY_1DRes$, and auxiliary resolutions that are associated with the graded Root resolutions
of $Y_1Res$, that we denote $GRootY_1DRes$. 

\item"{(3)}"  we look at test sequences that are constructed in part (2),  that can be further extended to specializations 
of the terminal limit group of one of the auxiliary resolutions that are associated with the graded resolution, $Y_1Res$
(that are associated with the set ${(AE)}^k(p,t)$), and was constructed
from test sequences that can be extended to the terminal limit group of one of the graded resolutions, $X_1Y_1Res$,
that are associated with the set $E{(AE)}^{k-1}(p,t,y_1)$. 
We further require that the specialization of the terminal limit
group of this auxiliary resolution can be extended to the f.p.\ completion that is associated with the auxiliary resolution,
and the terminal limit group (of the auxiliary resolution) is embedded into it. Once again with this collection of sequences
it is possible to associate finitely many graded resolutions that have the properties of the resolutions in a formal
graded Makanin-Razborov diagram (theorems 2.6 and 2.7). We denote the constructed resolutions, $X_1Y_1DRes$.

We further associate with each of these resolutions, (collapse) auxiliary  resolutions, $CollX_1Y_1DRes$, 
auxiliary resolutions that are associated with the singular
locus of the resolutions $X_1Y_1Res$, that we denote $SLX_1Y_1DRes$, and auxiliary resolutions that are associated with the graded Root resolutions
of $X_1Y_1Res$, that we denote $GRootX_1Y_1DRes$. 

\item"{(4)}"  we continue iteratively. For each index $j \leq k$, we look at all the test sequences that were looked at
while  constructing the auxiliary resolutions, $X_{j-1}Y_{j-1} \ldots X_1Y_1DRes$, that can be extended to the terminal
limit group of one of the auxiliary resolutions that are associated with the graded resolution, $Y_1Res$, and one of the 
resolutions, 
$Y_jX_{j-1}Y_{j-1} \ldots X_1Y_1Res$, that are associated with the set ${(AE)}^{k+1-j}$. 
We further require that the specialization of the terminal limit
group of this auxiliary resolution can be extended to the f.p.\ completion that is associated with the auxiliary resolution,
and the terminal limit group (of the auxiliary resolution) is embedded into it. With this collection of
sequences we, once again, associate finitely many resolutions, in a construction that follows the construction of the
formal graded Makanin-Razborov diagram (theorem 2.7). We denote the constructed resolutions, 
$Y_jX_{j-1}Y_{j-1} \ldots X_1Y_1DRes$. With each such resolution we also associate finitely many collapse auxiliary
resolutions,
$CollY_jX_{j-1}Y_{j-1} \ldots X_1Y_1DRes$, finitely many auxiliary resolutions that are associated with the singular locus of the
corresponding resolutions,  
$SLY_jX_{j-1}Y_{j-1} \ldots X_1Y_1DRes$, and finitely many auxiliary resolutions that are associated with all the corresponding graded
Root resolutions, 
$GRootY_jX_{j-1}Y_{j-1} \ldots X_1Y_1DRes$. 

Similarly, we look at all the test sequences that were looked while constructing these last auxiliary resolutions  
that can be extended to the terminal 
limit group of one of the auxiliary resolutions that are associated with the graded resolution, $Y_1Res$, and one of the 
auxiliary resolutions, 
$X_jY_j \ldots X_1Y_1Res$, that are associated with the set $E{(AE)}^{k-j}$. 
We denote the constructed resolutions, 
$X_jY_j \ldots X_1Y_1DRes$. With such a resolution we associate finitely many collapse auxiliary resolutions, 
$CollX_jY_j \ldots X_1Y_1DRes$, auxiliary resolutions that are associated with the singular loci of the constructed resolutions, 
$SLX_jY_j \ldots X_1Y_1DRes$, and auxiliary resolutions that are associated with graded Root resolutions,
$GRootX_jY_j \ldots X_1Y_1DRes$.

\item"{(5)}"  Finally, we look at all the test sequences that were looked at
while  constructing the auxiliary resolutions, $X_{k}Y_{k} \ldots X_1Y_1DRes$, that can be extended to the terminal
limit group of one of the auxiliary resolutions that are associated with the graded resolution, $Y_1Res$, and one of the 
resolutions, 
$\Psi X_{k}Y_{k} \ldots X_1Y_1Res$, that are associated with the set ${(AE)}$. 
We denote the constructed resolutions, 
$\Psi X_kY_k \ldots X_1Y_1DRes$. With each such resolution we also associate finitely many collapse auxiliary
resolutions,
$Coll \Psi X_kY_k \ldots X_1Y_1DRes$. 
\endroster

\medskip
Suppose that there is no 
approximating cover resolution, $DRes_m$, with a subsequence of the original sequence, $\{(p_n,t_n)\}$ (for which
$p_n \in E{(AE)}^k(p)$ (over $G_n$), and $t_n$ (a specialization of the universal variables $t$ in $G_n$) is a witness for that
(in $G_n$)), that factor through $DRes_m$, so that
for each specialization of the terminal limit group of $DRes_m$ that is associated with  
a pair $(p_n,t_n)$ from the subsequence, a generic  pair in the fiber that is associated with
such a specialization  of the terminal graded limit group (which takes its values in $G_n$), satisfies $p_n \in E{(AE)}^k(p)$
(over $G_n$), and $t$ (a generic value of $t$ in the associated fiber) is a witness for that.

As in analyzing EAE and AEAE  sets in sections 4 and 5,  in this case 
we can associate with the original 
graded resolution, $DRes$, another graded resolution, of smaller complexity,
which is obtained by applying the first step of the sieve procedure
for the analysis of quotient resolutions [Se6]. Furthermore, a subsequence of the original sequence of specializations,
$\{(p_n,t_n)\}$, is guaranteed to
have a subsequence that extends to specializations that converge into the obtained quotient resolution.
This enables us to continue iteratively, in a similar way to the sieve procedure [Se6], and to the iterative procedure that 
was used to prove the 
equationality of Diophantine sets over free and hyperbolic groups in section 2 of [Se8].

With the resolution, $DRes$, we have associated a sequence of approximating cover resolutions, $DRes_m$, that satisfy the
properties that are listed in theorem 1.21.
In particular, for each index $m$, there is a subsequence
of the original sequence of pairs, $\{(p_n,t_n)\}$, that factor through $DRes_m$. By our assumptions,
for every index $m$, there is no subsequence of the  pairs, $\{(p_n,t_n)\}$, that do
factor through the approximating cover, $DRes_m$, for which  for generic pair $(p_n,t)$ in the fiber that
contains the pair $(p_n,t_n)$ (in the graded variety that is associated with $DRes_m$), 
$p_n \in E{(AE)}^k(p)$ over $G_n$, and the generic 
$t$ (a test sequence in in the fiber) is a witness for that, i.e., the generic pairs, $(p_n,t) \in {(AE)}^k(p,t)$ over $G_n$.

We construct a sequence of specializations over $G_n$, in the same way as what we did in analyzing an AEAE set in section 5
 as follows.
We go over the indices $m$, and given an index  $m$ we pick an index $n_m > m$, and a tuple, 
$(p_{n_m}, \hat t_{n_m}, \tilde {(y_1)}_{n_m}, \tilde t_{n_m})$, that takes its values in $G_{n_m}$,
with the following properties:
\roster
\item"{(1)}" $(p_{n_m},t_{n_m})$ factors through the graded resolution, $DRes_m$, and $\hat t_{n_m}$ is the 
specialization of the terminal limit group of the graded resolution, $DRes_m$, that contains the pair, $(p_{n_m},t_{n_m})$.

\item"{(2)}" $\tilde t_{n_m}$ is a specialization of the variables $t$, from the fiber that is associated with
$\hat t_{n_m}$, i.e., the fiber that contains the pair $(p_{n_m},t_{n_m})$ in the graded variety that is associated with
$DRes_m$. Furthermore, for every index $m$, $(p_{n_m}, \tilde t_{n_m}) \notin {(AE)}^k(p,t)$ over $G_n$.

\item"{(3)}" the sequence $\{(p_{n_m},\tilde t_{n_m}) \}$ is a graded test sequence that converges into the graded resolution, 
$DRes$. 

\item"{(4)}" the tuples $(p_{n_m},\tilde t_{n_m}, \tilde {(y_1)}_{n_m})$, are specializations of the terminal limit group of
one of the  finitely many graded resolutions, $Y_1Res$ (that are associated with the set ${(AE)}^k(p,t)$), 
that restrict to  rigid  or  weakly 
strictly solid specializations of the rigid or weakly solid factor of that terminal limit group of $Y_1Res$,
with respect to the parameter group $<p,t>$). 
Furthermore, a generic value of the existential variables $y_1$ in the fiber of $Y_1Res$ that is associated with the tuple,
$(p_{n_m},\tilde t_{n_m}, \tilde {(y_1)}_{n_m})$,  
is a witness that:
$(p_{n_m}, \tilde t_{n_m}) \notin {(AE)}^k(p,t)$ over $G_{n_m}$.
\endroster

From the sequence of tuples,
$(p_{n_m}, \hat t_{n_m}, \tilde {(y_1)}_{n_m}, \tilde t_{n_m})$, we can extract a subsequence that converges into a closure of
the graded resolution, $DRes$, that we have started with. We denote this closure, $ExY_1DRes$. 
Recall that we have assumed that the original sequence of pairs,
$\{(p_n,t_n)\}$, satisfies $p_n \in E{(AE)}^k(p)$, and $t_n$ is a witness for that (i.e., $(p_n,t_n) \in {(AE)}^k(p,t)$).
On the other hand, in the sequence of tuples that
we chose (in the fibers that are associated with a subsequence of the original sequence
of specializations), a generic value of the existential variables $y_1$ in the fiber that is associated with the tuple,
$(p_{n_m},\tilde t_{n_m}, \tilde {(y_1)}_{n_m})$,  
is a witness that the pair $(p_{n_m}, \tilde t_{n_m}) \notin {(AE)}^k(p,t)$.

Given the closure, $ExY_1DRes$, we continue in a similar way to what we did in analyzing AEAE sets in section 4, a way that is
adapted to sets with arbitrarily many  quantifiers. 
The closure, $ExY_1DRes$, contains an (additional) specialization of the variables $y_1$, so that for a generic 
value of the variables $t$ and $y_1$, the value of the existential variables $y_1$ are witnesses that the pairs, 
$(p_n, t) \notin {(AE)}^k(p,t)$. Hence, we look at a sequence of approximating covers of 
the closure, $ExY_1DRes$, that we denote $\{ExY_1DRes_r\}$,  that satisfy the properties that are listed
in theorem 1.21. On the additional (existential) variables $y_1$ in the approximating covers, $ExY_1DRes$, we impose 
one of the (finitely many) possible
collapse forms that are defined iteratively. These possible collapse forms are:
\roster
\item"{(1)}" the additional specialization of the terminal limit group of the graded resolutions, $Y_1Res$, 
 in $ExY_1DRes_r$, that is assumed to restrict to a rigid or a
weakly strictly solid specialization of the rigid or weakly solid factor of the terminal limit group of the graded resolution $Y_1Res$,
is non-rigid or non weakly strictly solid. This forces a Diophantine condition on the specializations of $ExY_1DRes_r$, similar to the
ones that are forced along the sieve procedure in [Se6].

Similarly, we impose the Diophantine condition that forces  
the additional specialization of the terminal limit group of $Y_1Res$ to extend to a specialization of the terminal limit group
of one of the resolutions,
$SLY_1Res$, that are associated with the singular locus of $Y_1Res$.

\item"{(2)}" to the specialization of the closure, $ExY_1DRes_r$, it is possible to add an additional specialization of the
terminal graded limit group of one of the graded resolutions, $X_1Y_1Res$, that is associated by our induction hypothesis with the set,
${(AE)}^k(p,t)$, that restricts 
to a rigid or a weakly strictly solid
specialization of the rigid or solid factor of that terminal 
graded limit group. If the extended specializations can not be  extend to a specialization of a closure of $ExY_1DRes_r$, 
then we collect these extended specializations
(for a subsequence of the original sequence of specializations, $\{(p_n,t_n)\}$), in a graded resolution that has smaller
complexity than the graded resolution, $DRes$, according to the sieve procedure [Se6], or rather according to the iterative
procedure that was used to prove the equationality of Diophantine sets in [Se9].

\item"{(3)}" suppose that there is a test sequence of the graded resolution, $ExY_1DRes$,
for which for each additional specialization of the terminal limit group of the resolution
$Y_1Res$ (that is associated with the set ${(AE)}^h(p,t)$), there is an extension to 
specializations of the terminal limit group of some auxiliary
resolutions $X_1Y_1Res$ (that are associated with the resolution $Y_1Res$ and the set
${(AE)}^k(p,t)$), so that the associated fibers covers the corresponding fiber of $Y_1Res$.

In this case, we first look at all the test sequences of $ExY_1DRes_r$, for which  the specialization of the variables 
$(p,t,y_1)$ in  $ExY_1DRes_r$, can be extended to  a specialization of the terminal limit group of an 
auxiliary resolution that is associated with graded Root resolution of $Y_1Res$ (and with the set ${(AE)}^k(p,t)$), and not to a specialization
of an auxiliary resolution that is associated with graded Root resolutions of $Y_1Res$ of higher order roots, 
and for which the specializations of the variables, $(p,t,y_1)$, can be extended to specializations of 
the terminal limit groups of a prescribed set of auxiliary resolutions, $X_1Y_1Res$  and 
$CollY_1Res$, 
that are associated with the  resolution
$Y_1Res$ (and with the set ${(AE)}^k(p,t)$),  so that the corresponding fibers of these auxiliary
resolutions cover the fiber of the resolution $Y_1Res$.

We further repeat the 
construction of the graded resolution, $ExY_1DRes_r$, and construct a closure of, $ExY_1DRes_r$, from test sequences that can be 
further extended to
specializations of the terminal limit group of one of the auxiliary resolutions, $Y_2X_1Y_1Res$, that are associated with
the set, ${(AE)}^k(p,t,y_1,x_1)$, that restrict to rigid or weakly strictly solid specializations of the rigid or solid
factor in that graded limit group. We denote this constructed graded resolution, which is a closure of $ExY_1DRes_r$, 
$ExY_2X_1Y_1DRes_r$.

First, it may be that for a subsequence of the original sequence of specializations, the additional specialization
(that is associated with the variables $Y_2$) does not restrict to a rigid or a 
weakly strictly solid specialization, or that this specialization extends to a specialization
of the terminal limit group of one of the auxiliary resolutions that are associated
with the singular locus of the associated auxiliary resolution $Y_2X_1Y_1Res$ 
(in a similar to what we did in part (1)), or that the specialization of the terminal
limit group of the graded  resolution $Y_1Res$ extends to an associated graded Root 
resolution with higher order roots. In this case
we impose a Diophantine condition similar to the ones that are imposed in part (1) on the constructed resolution, 
$ExY_2X_1Y_1DRes_r$, and as in part (1) we get a graded resolution with a smaller complexity according to the sieve procedure [Se6].

If we can not impose these Diophantine conditions,
we further collect specializations that extend to specializations of the terminal
limit groups of auxiliary resolutions, $X_2Y_2X_1Y_1Res$ (that restrict to rigid or weakly strictly solid specializations, and do not extend to the associated singular locus).
 If the extended
 specializations do not factor through a closure of $ExY_1DRes_r$, we obtained a graded resolution with smaller complexity than $DRes$. 

If
for a subsequence of the original specializations, the extended specializations do factor through a closure of $ExY_2X_1Y_1DRes_r$,
we continue as in part (3), and collect test sequences that can be extended to specializations that are associated with the 
variables $Y_3$. 

\item"{(4)}" we continue iteratively. The original sequence of specializations, $\{(p_n,t_n)\}$, are in the set,
${(AE)}^k(p,t)$, but generic points in the fibers of $ExY_1DRes_r$ that are associated with these specializations  are not in this set.
Hence, after at most $k$ steps we must end up with a graded resolution, that has a smaller complexity than the original graded
resolution, $DRes$, and this graded resolution is obtained by either forcing an additional specialization that is associated
with variables $Y_j$ to be non-rigid or non weakly strictly solid or extendable to the specialization of the terminal limit group of an auxiliary
resolution that is associated with the corresponding singular locus (a Diophantine condition), or it is obtained by collecting additional
specializations that are associated with variables $X_j$ that restrict to rigid or weakly strictly solid specializations, and do not
factor through a closure of $ExY_1DRes_r$.
\endroster

In all these cases we constructed a quotient resolution according to the first step of the sieve procedure [Se6], and
pass to a subsequence of the specializations, $\{(p_n,t_n)\}$, that extend to specializations that converge into that quotient
resolution.
Note that by the construction of these quotient resolutions that appear in [Se6], the complexity of this quotient resolution,
in light of the sieve procedure, is smaller than the complexity of the graded resolution that we have started the first step with, $DRes$ (note that the obtained quotient resolution may be 
a proper closure of the resolution that we have started the first step with, i.e. the resolution
$DRes$. However, there
is a global bound on the orders of roots that are needed to be added to such closures, hence, 
the obtained quotient resolution can be a closure of the quotient resolution that we have started the first step with only 
finitely (boundedly) many steps. After these finitely many steps, the complexity of the
obtained quotient resolution is strictly smaller).

We continue iteratively. At each step we first look at a sequence of approximating (cover) resolutions of the quotient resolution
that was constructed in the previous step, where these approximating cover resolutions satisfy the properties that are
listed in theorem 1.21, and in particular their completion can be embedded in a f.p.\ completion.
If there exists an approximating cover resolution for which there exists a subsequence
of pairs, (still
denoted) $\{(p_n,t_n)\}$, so that for each specialization of the terminating limit group of the approximating cover resolution,
 that is associated with  
a pair $(p_n,t_n)$ from the subsequence, a generic  pair (i.e., a test sequence)
in the fiber that is associated with
such a specialization  of the terminal graded limit group (which takes its values in $G_n$) forms a witness for  $p_n \in E{(AE)}^k(p)$
(over $G_n$),  we do what we did in this case in the first step (and associate with this approximating cover finitely many
auxiliary resolutions).

If there is no such approximating cover resolution for the constructed quotient resolution, we associate with it a closure
that is constructed from test sequences that can be extended to additional (new)  rigid or weakly 
strictly solid specializations of the terminal
limit groups of one of the graded resolutions $Y_1Res$, that are associated with the set ${(AE)}^k(p,t)$.
We further associate a sequence of approximating covers with that closure, and force a collapse condition over the constructed 
 approximating covers, according to parts (1)-(4). Finally, we apply the general step of the sieve procedure [Se6], and obtain a
new quotient resolution,  that has smaller complexity than the quotient resolution that we have started the current step with.

By an argument which combines the argument that guarantees the termination of the sieve procedure in [Se6], with the argument that implies
the equationality of
Diophantine sets over free groups,
we obtain a termination of the this procedure.

\vglue 1pc
\proclaim{Theorem 6.2} The procedure for the analysis of an $E{(AE)}^k$ set over free products terminates after finitely
many steps.
\endproclaim

\nfp The argument that we use is similar to the proofs of theorems 4.4 and 5.1 (termination of the procures for the analysis of AEAE and
 EAE sets).
At each step of the procedure, on  the quotient resolution, $DRes$, that is analyzed in that step  we impose one of finitely
many  restrictions. The first ones impose a non-trivial Diophantine condition on specializations of a closure of the resolution, $DRes$. 
The second type of restrictions, adds specializations of the terminal limit groups of one of the associated auxiliary resolutions that are
associated with the graded resolution, $Y_1Res$ (and with the set ${(AE)}^k(p,t)$), where these specializations restrict to rigid or weakly strictly
solid specializations of the rigid or weakly solid factor of the terminal limit group. We further require that the additional
specializations do not factor through (or are not in the same (weak) strictly solid families as specializations of) closures that are
built from test sequences of $DRes$ that can be extended to specializations of such terminal limit groups of auxiliary resolutions.

The argument that was used to prove the termination of the sieve procedure 
(theorem 22 in [Se6]), that was modified to work over free products in the proof of theorem 4.4, 
proves that restrictions of both
types can  occur in only finitely many steps. Hence,
the iterative procedure for the analysis of an $E{(AE)}^k$  set terminates after finitely many steps.

\line{\hss$\qed$}

\smallskip
So far we have shown that given a sequence,
$\{(p_n,t_n)\}$, for which  $p_n \in E{(AE)}^k(p)$ over $G_n$, and $t_n$ is a witness for $p_n$,
i.e., $(p_n,t_n) \in {(AE)}^k(p,t)$ over $G_n$, it is possible to extract a subsequence, (still denoted) $\{(p_n,t_n)\}$,
 and to construct a graded resolution, $DRes$, 
with the following
properties:
\roster
\item"{(i)}" the graded resolution $DRes$ satisfies the properties of a cover graded resolution that are 
listed in theorem 1.21. 

\item"{(ii)}" the specializations of the subsequence, $\{(p_n,t_n)\}$,  extend to 
specializations 
that factor through the resolution, $DRes$, and to specializations of the f.p.\ completion of the ungraded resolution that
extends the graded resolution, $DRes$ (see theorem 1.21). 
Hence, with each specialization from this subsequence of the tuples, 
$\{(p_n,t_n)\}$, 
specializations of the  terminal limit group of the graded resolution, $DRes$, can be associated. Furthermore, 
these specializations of the terminal limit group of $DRes$
 can be extended to specializations of the f.p.\ completion of some ungraded resolution of the terminal limit group of
$DRes$. Restrictions of generic points (i.e., test sequences) 
in the fibers that are associated with these specializations of the terminal limit group
of $DRes$, to the variables $(p,t)$, that we denote, $(p_n,t)$, are in the set ${(AE)}^k(p,t)$ that is associated with the given 
set, $E{(AE)}^k(p)$.

\item"{(iii)}" with the resolution $DRes$,  we associate (non-canonically) a finite collection of graded
auxiliary resolutions, that are part of the output of the terminating iterative procedure for the analysis of a definable set,
i.e., auxiliary resolutions of the same type as those that were associated with an approximating cover resolution,
$DRes_m$, in the first step of the procedure, in case it terminates in the first step. These auxiliary resolutions have the
same properties of the resolutions in a formal graded Makanin-Razborov diagram (theorems 2.6 and 2.7), and in particular
they can be extended to ungraded resolutions with f.p.\ completions.
\endroster

As in the proof of theorems 3.2, 4.1, 4.5 and 5.1,
we look at all the sequences of non-trivial free products, $G_n=A^1_n*\ldots*A^{\ell}_n$
(possibly for varying $\ell>1$), that are not isomorphic to $D_{\infty}$, an associated sequence of 
 tuples, $p_n \in E{(AE)}^k(p)$ over $G_n$, and witnesses $t_n$ for $p_n$, i.e., a sequence of pairs $(p_n,t_n) \in {(AE)}^k(p,t)$
over $G_n$.
From every such  sequence we use our terminating iterative procedure, and extract a subsequence
of the tuples, (still denoted) $\{(p_n,t_n)\}$, and a graded resolution,  
$DRes$, that has the properties (i)-(iii), and in particular,   the subsequence of tuples, (still denoted) $\{(p_n,t_n)\}$, 
extend to specializations that factor through the resolution $DRes$, and to specializations of the f.p.\ completion of
the ungraded resolutions that extends the graded resolution $DRes$.
Furthermore, the restrictions of generic points (test sequences) 
in the fibers that are associated with the pairs, $\{(p_n,t_n)\}$,
are in the set ${(AE)}^k(p,t)$ (i.e., generic values of $t$ in these fibers are witnesses for $p_n \in E{(AE)}^k(p)$
over $G_n$).

The completion of each of the constructed resolutions, $DRes$, can be extended to an ungraded f.p.\ completion, 
and so are its associated (finitely many) auxiliary resolutions.
Furthermore, the (finite collection of) covers of flexible quotients
of the rigid or weakly solid factors of the terminal limit groups of the resolutions, $DRes$, and of the auxiliary resolutions
that are associated with $DRes$, 
can all be embedded into (ungraded) resolutions with f.p.\
completions. Hence, we can define a linear  order on this (countable) collection of
graded resolutions ($DRes$), and their (non-canonically) associated auxiliary resolutions, and (finite collections of)
covers of flexible quotients.  By the same argument that
was used in constructing the Makanin-Razborov diagram (theorem 26 in [Ja-Se]), there exists a finite subcollection
of these graded resolutions that satisfy properties (1)-(3) of the theorem.

So far we have argued that given our induction hypothesis (the conclusion of theorem 6.1) on sets of the form, ${(AE)}^k$, the 
conclusion of theorem 6.1 follows for sets of the form $E{(AE)}^k$. Clearly, exactly the same argument proves the
conclusion of theorem 6.1 for sets of the form ${(AE)}^{k+1}$, assuming what we already proved, i.e., the conclusion of
theorem 6.1 for sets of the form $E{(AE)}^k$. Hence, we have completed the proof by induction, and theorem 6.1 follows for
every definable set over free products.

\line{\hss$\qed$}

Theorem 6.1 associates finitely many graded resolutions with any given definable set over free products.
As we argued for AE, EAE and AEAE sets, the existence of a finite collection of graded resolutions with the properties that are listed in
theorem 6.1, 
allows one to reduce a general sentence in the language of groups from free products to a sentence over its factors in a uniform way.

\vglue 1pc
\proclaim{Theorem 6.3} Let $\Phi$ be a coefficient-free sentence over groups.
Then there exists a coefficient-free 
sentence over free products, 
which is a (finite) disjunction of conjunctions of  sentences over the factors of the free product, 
such that for every non-trivial free product, $G=A^1*\ldots*A^{\ell}$, that is not isomorphic to $D_{\infty}$, 
the original
sentence over the free product $G$ is a truth sentence, if and only if the sentence which is a 
(finite) disjunction of conjunctions of sentences over the factors $A^1,\ldots,A^{\ell}$ is a truth sentence.
\endproclaim 

\nfp The argument that we use, that is based on theorem 6.1, is similar to the proof of theorem 5.3 (which is based on
theorem 5.2). 
 By theorem 6.1, with a given coefficient-free definable  set, it is possible to associate finitely many graded resolutions,
$DRes_1,\ldots,DRes_g$, that do all satisfy the properties that are listed in  theorem 1.21, and with them we have associated finitely
many auxiliary resolutions, 
according to parts (1)-(5) that appear in the first step of the iterative procedure for the analysis of
a definable set. Note that all these auxiliary resolutions have the properties of
the resolutions in the formal graded Makanin-Razborov diagram as listed in theorems 2.6 and 2.7, and in particular they can all be
embedded into f.p.\ completions.

As we remarked in the EAE and AEAE cases in sections 4 and 5, in analyzing general sentences (rather than predicates), 
the same constructions 
that enable one to associate  graded  resolutions and their
auxiliary resolutions with a definable set (i.e., the iterative procedure for the analysis of a definable set, and the proof of 
theorem 6.1),
enable one to associate with a given  sentence a (non-canonical) finite collection of (ungraded) resolutions 
(over free products) with f.p.\ completions, 
and with each resolution finitely many (ungraded) auxiliary resolutions that do all have f.p.\ completions, and these auxiliary resolutions
have the same properties and they are constructed in the same way as the (graded) auxiliary resolutions that are constructed in the case of
a definable set.

We (still) denote the (ungraded) resolutions that are associated with the given AEAE sentence, $DRes_1,\ldots,DRes_g$.
By theorem 6.2, the given sentence is true over a non-trivial free product, 
$G=A^1*\ldots*A^{\ell}$,
that is not isomorphic to $D_{\infty}$, if and only if
there exists a specialization in $G$ of  the terminal limit group of a resolution
$DRes_i$ (one of the resolutions, $DRes_1,\ldots,DRes_g$), i.e.,
specializations of the elliptic factors of the terminal limit group of $DRes_i$ in the factors, $A^1,\ldots,A^{\ell}$, so that
 the specialization does not extend to a specialization of the terminal limit group of one of the auxiliary resolutions
that are associated with the singular locus of $DRes_i$ (see proposition 4.2). 
This specialization of the terminal limit
group of $DRes_i$ also satisfies: 
\roster
\item"{(1)}" it extends to the terminal limit group of an auxiliary resolution that is associated with a  Root resolution 
of $DRes_i$ according to proposition 4.3
(possibly only the trivial roots), 
and does not extend to specializations of the terminal limit group of
other Root resolutions, $RootDRes$, that are associated with higher order roots.

\item"{(2)}" it may extend to  specializations of the terminal limit groups of some of the auxiliary resolutions
$Y_1DRes$, that do not extend to specializations of the terminal limit groups of auxiliary resolutions, $SLY_1DRes$, 
that are associated with the singular locus of the
associated graded resolution, $Y_1Res$, nor to the terminal of a collapse auxiliary resolution, $CollY_1DRes$. Either the collection 
of the fibers that are associated with these specializations of
the terminal limit groups of $Y_1DRes$ do not form a covering closure of the fiber that is associated with the specialization
of the terminal limit group of $DRes_i$, or (some of) these specializations can be extended to the specializations of the terminal
limit groups of the following auxiliary resolutions (see pars (1)-(5) in the first step of the iterative procedure for the analysis
of a definable set for the construction and the properties of the 
auxiliary resolutions that we refer to).

\item"{(3)}" each of the specializations of the terminal limit groups of $Y_1DRes$ extends to a specialization of the terminal
limit group of one of the auxiliary resolutions, $RootY_1DRes$, and not to similar auxiliary resolutions that are associated with roots
of higher orders.

\item"{(4)}" some of the specializations that extend to specializations of the terminal limit groups of the auxiliary resolutions,
$Y_1DRes$, extend further to specializations of the terminal limit groups of auxiliary resolutions,
$X_1Y_1DRes$, that do not extend to specializations of the terminal limit groups of auxiliary resolutions,
$CollX_1Y_1DRes$ and $SLX_1Y_1DRes$. 

The fibers that are associated with the specializations of the terminal limit groups of the auxiliary resolutions,
$Y_1DRes$ (and not to $SLY_1DRes$ nor $CollY_1DRes$), 
minus the fibers that are associated with specializations that can be further extended to specializations of the
terminal limit groups of auxiliary resolutions, $X_1Y_1DRes$ (and not to $SLX_1Y_1DRes$ nor $CollX_1Y_1DRes$), do not
form a covering closure of the fiber that is associated with the specialization of the terminal limit group of $DRes$.

\item"{(5)}" if none of the specializations of the terminal limit groups of the auxiliary resolutions, $X_1Y_1DRes$, that can not
be extended to specializations of the terminal limit groups of $SLX_1Y_1DRes$ nor $CollX_1Y_1DRes$, can be further
 extended to specializations of the terminal limit groups of auxiliary resolutions, $Y_2X_1Y_1DRes$, that can not be further extended
to $SLY_2X_1Y_1DRes$ nor $CollY_2X_1Y_1DRes$, we are done. 

If some of these specializations of $X_1Y_1DRes$ do extend to $Y_2X_1Y_1DRes$ we do the following. 
Each such specialization of the terminal limit group of $X_1Y_1DRes$ extends to a specialization of the terminal limit group of 
an auxiliary resolution, $GRootX_1Y_1DRes$, and not
to the terminal limit groups of such resolutions that are associated with higher order roots.

If the fibers that are associated with specializations of (the terminal limit groups of) the auxiliary resolutions, $Y_1DRes$,
minus the fibers that are associated with specializations of $X_1Y_1DRes$ from which we take out fibers that
are associated with specializations of $Y_2X_1Y_1DRes$, do not form
a covering closure of the fiber that is associated with the specialization of the terminal limit group of $DRes$ we are done.

\item"{(6)}" otherwise, we continue iteratively. At each step, we first require that the specialization of the corresponding terminal
limit groups can be extended to some associated graded Root  auxiliary resolutions, and to such auxiliary resolutions that are
associated with higher order roots.

Then we check if the inclusion exclusion combination of fibers that are associated with the auxiliary resolutions that are
associated with existential and universal variables do cover the fiber that is associated with the (original) 
specialization of the terminal
limit group of $DRes$. If they do not cover the fiber of $DRes$ we are done. Otherwise we continue to the next auxiliary resolutions.
After finitely many such steps (that depend only on the number of quantifiers in the original sentence) we have exhausted all
the possibilities for the sentence to be a true sentence.
\endroster

Finally, the given sentence is a true sentence over a non-trivial free group, $G=A^1*\ldots*A^{\ell}$, that is not
isomorphic to $D_{\infty}$, if and only if
there exists a specialization of the terminal limit group of $DRes$, that satisfies one of the possibilities that are
described in parts (1)-(6). The union of these finitely many possibilities can be expressed as a finite disjunction of
finite conjunctions of sentences over the factors,
$A^1,\ldots,A^{\ell}$.

\line{\hss$\qed$}

In  analyzing AEAE sets in section 5, we were able to reduce not only AEAE (and EAEA)
sentences over free products, but also AEAE (and EAEA) predicates. The existence of graded resolutions with the properties that are
listed in theorem 6.1, enables us to get a similar reduction for arbitrary predicates over free products.

\vglue 1pc
\proclaim{Theorem 6.4} Let $\Phi$ be an arbitrary coefficient-free predicate over groups. 
Then there exists a coefficient-free predicate over groups that are free products, that is composed from only 3 quantifiers over 
variables that take values in the ambient
free product, and additional quantifiers over variables that take their values in the various factors of the free product, 
so that for every 
non-trivial free product, $G=A^1*\ldots*A^{\ell}$, that is not isomorphic to $D_{\infty}$, 
the set $\Phi(p)$ over $G$ can be defined by the following
predicate over free products:
$$\Phi(p) \ = \   \exists u \  (\Theta(e_1)) \   \ \forall v \ \exists s \ (\exists e_2)$$
$$ (\Sigma_1(e_1,e_2,u,v,s,p)=1 \, \wedge \, \Psi_1(e_1,e_2,u,v,s,p) \neq 1) \, \vee \, \ldots \, \vee$$
 $$(\Sigma_{\ell}(e_1,e_2,u,v,s,p)=1 \, \wedge \, \Psi_{\ell}(e_1,e_2,u,v,s,p) \neq 1)$$
where the variables $u,v,s$ take values in the ambient free product $G$, and the variables $e_1,e_2$ take values in the factors 
$A^1,\ldots,A^{\ell}$,
and $\Theta$ is a predicate over elements from the factors $A^1,\ldots,A^{\ell}$.
\endproclaim 

\nfp Given theorems 6.1 and 6.3, the statement of the theorem follows by the same argument that was used to prove theorem 
5.4 (that is based on theorems 5.2 and 5.3).

The existential variables $u$ represent the fractions in all the possible systems of fractions that are associated 
(by theorems 1.14 and 1.15) with the
rigid or weakly solid factors of the terminal limit groups of the graded resolutions, $DRes_1,\ldots,DRes_g$, 
and their auxiliary resolutions, that are associated with
the given definable set  by theorem 6.1. The universal variables $v$ and the existential variables $s$ and $e_2$ 
(the variables $e_2$ are contained in the various 
factor s$A^1,\ldots,A^{\ell}$ of the free product), enable us to guarantee that the
values of the fractions $u$ satisfy the conclusions of theorems 1.14 and 1.15, i.e., enable  the covering of  all the families of rigid or
weakly strictly solid families that are associated with any given value of the defining parameters.

The variables $e_1$ take their values in the various elliptic factors, $A^1, \ldots,A^{\ell}$, of the given free product $G$, and
this part of the predicate is given by the description that appears in theorem 6.3, only that the elements $e_1$ are either 
specializations of the  elliptic factors of the terminal limit groups of the graded resolution, $DRes_1,\ldots,DRes_g$, and of
their auxiliary resolutions, or they are elliptic elements that appear in the combinatorial description of rigid and almost
shortest weakly strictly solid specializations (given by theorems 1.14 and 1.15) of the rigid or weakly solid factors of
these terminal limit groups. The universal variables $v$ (that get values in the ambient free product) enable us to guarantee that
the existential variables among the variables $e_1$ that are associated with  rigid or weakly strictly solid specializations
are indeed part of such specializations. The existential variables  
 $s$ (that get values in the ambient free product) enable us to guarantee that
the universal variables among the elliptic variables $e_1$ that  represent rigid or weakly strictly solid specializations are indeed
part of such specializations, and that two given specializations
of a weakly solid limit group belong to the same family.

Once we have used the variables $u,v,s$ (that take their values in the ambient free product) and the elliptic variables $e_1,e_2$, to
go over all the rigid and weakly strictly solid specializations, the theorem follows by the proof of theorem 6.3. 

\line{\hss$\qed$}

\vglue 1.5pc
\centerline{\bf{\S7. Basic properties of the first order theory of free products}} 
\medskip

In the previous section we associated a finite collection of graded resolutions with any given predicate (definable set)
over free products (theorem 6.1). We used these graded resolutions to reduce a general sentence over free products
to a (finite) disjunction of conjunctions of sentences over the factors of the free product (theorem 6.3), and to get
 a form of "quantifier elimination", i.e. a reduction of a predicate over free products to a 
predicate that contains only 3 quantifiers over the ambient free product and finitely many additional
quantifiers over the factors of the free product (theorem 6.4).

The (uniform) reduction of sentences from a free product to its factors has fairly immediate corollaries for
the first order theory of free product of groups. We start with a positive answer to a well known problem of 
R. Vaught (cf. [Fe-Va]). 

\vglue 1pc
\proclaim{Theorem 7.1} Let $A_1,B_1,A_2,B_2$ be groups. Suppose that $A_1$ is elementarily equivalent to $A_2$,
and $B_1$ is elementarily equivalent to $B_2$. Then $A_1*B_1$ is elementarily equivalent to $A_2*B_2$.
\endproclaim

\nfp If $A_1*B_1$ is isomorphic to $D_{\infty}$, and the free product is non-trivial, then both $A_1$ and $B_1$ are
isomorphic to $Z_2$. Hence, $A_2$ and $B_2$ are isomorphic to $Z_2$ as well, so if $A_1*B_1$ is a non-trivial
free product that is isomorphic to $D_{\infty}$ then $A_2*B_2$ is isomorphic to $D_{\infty}$ as well, and the
theorem follows in this case.    

Suppose that $A_1*B_1$ is a non-trivial free product that is not isomorphic to $D_{\infty}$. Let $\Phi$ be a coefficient
 free sentence over groups. By theorem 6.3, the sentence $\Phi$ is a truth sentence over  a non-trivial free product
$U*V$ that is not isomorphic to $D_{\infty}$, if and only if a (finite) disjunction of conjunctions of sentences over the
factors $U$ and $V$, that we denote $\alpha$,  is a truth sentence. 

Since $A_1$ is elementarily equivalent to $A_2$, and $B_1$ is elementarily equivalent to $B_2$, the sentence $\alpha$ is a truth
sentence over the factors, $A_1,B_1$, if and only if it is a truth sentence over the factors, $A_2,B_2$. Therefore, $\Phi$, is
truth over $A_1*B_1$, if and only if it is truth over $A_2*B_2$. 

\line{\hss$\qed$}

The existence of graded resolutions that are associated with a given  sentence over free products enables one to prove the following
theorem, that implies Tarski's problem for free groups.
 
\vglue 1pc
\proclaim{Theorem 7.2} Let $A,B$ be non-trivial  groups, and suppose that either $A$ or $B$ is not $Z_2$. Let $F$ be a (possibly cyclic) free group.
Then $A*B$ is elementarily equivalent to $A*B*F$.
\endproclaim

\nfp Since by Tarski's problem [Se7], all the non-abelian free groups are elementarily equivalent, by theorem 7.1 we may
assume that $F$ is a f.g.\ free group. With $A*B$ we have associated a Bass-Serre tree, $T_1$, that corresponds to the graph
of groups that contains two vertex groups, $A$ and $B$, and an edge (with trivial stabilizer) that connects them. 
With $A*B*F$ we associate a
Bass-Serre tree, $T_2$, that corresponds to a graph of groups that contains 3 vertices, one with trivial stabilizer 
and two with stabilizers $A$ and $B$ that are
connected by edges with trivial stabilizer to the vertex group with trivial stabilizer. On the vertex with trivial stabilizer
we further place $m \geq 1$ loops, if the free group $F$ is isomorphic to $F_m$. In particular, in both free products,
 $A*B$ and $A*B*F$, the elliptic elements are (only) those that can be conjugated into $A$ or $B$.

Let $\Sigma(x)=1$ be a (finite, coefficient-free) system of equations. Then, by construction, every (non-canonical)
Makanin-Razborov diagram of $\Sigma$ over  the collection of all free products, $A*B$, is a Makanin-Razborov diagram
of $\Sigma$ over the collection of all the free products, $A*B*F$. If $\Theta(x,p)=1$ is 
a coefficient-free system of equations with
parameters, $p$, then by construction, every (non-canonical) graded Makanin-Razborov diagram of $\Theta$ (that satisfy the
properties that are listed in theorem 1.22) over the collection of free products, $A*B$,  is a graded Makanin-Razborov 
diagram of $\Theta$ over the entire collection of free products $A*B*F$. This means, in particular, that a graded limit
group over the collection, $A*B$, is rigid or weakly solid, if and only if it is rigid or weakly solid over the
collection, $A*B*F$, and the same holds for flexible quotients of a rigid or a (weakly) solid limit group.

 $A*B$ naturally embeds (as a factor) in $A*B*F$, and hence each element $y_0 \in A*B$ can be naturally viewed as an element in
$A*B*F$.
Let: $$\forall y \ \exists x \ \Sigma(x,y)=1 \, \wedge \, \Psi(x,y) \neq 1$$ be a sentence over groups, and let $A*B$ be
a non-trivial free product that is not isomorphic to $D_{\infty}$. 
Suppose that $y_0 \in A*B$ is a specialization 
of the universal variables $y$,
for which there is no specialization $x_0 \in A*B$ of the existential variables $x$, so that $\Sigma(x_0,y_0)=1$ and 
$\Psi(x_0,y_0) \neq 1$ (in $A*B$). Then for $y_0$ viewed as an element of $A*B*F$, there is no specialization $x_1 \in A*B*F$ 
of the 
existential variables $x$, so that 
$\Sigma(x_1,y_0)=1$ and 
$\Psi(x_1,y_0) \neq 1$ (in $A*B*F$).  

Let $y_1 \in A*B*F$ be a specialization of the universal variables $y$, for which there is no $x_1 \in A*B*F$, such that
$\Sigma(x_1,y_1)=1$ and 
$\Psi(x_1,y_1) \neq 1$ (in $A*B*F$). Then there exists a sequence of retractions: $\tau_n : A*B*F \to A*B$, that maps the subgroup $A*B$
of $A*B*F$ identically onto $A*B$, and maps a fixed generators of $F$ into tuples of elements that form a  test sequence in $A*B$
(when $n$ grows to
infinity), so that for every index $n$, there is no element $x_n \in A*B$,
for which:
$\Sigma(x_n,\tau_n(y_1))=1$ and 
$\Psi(x_n,\tau_n(y_1)) \neq 1$ (in $A*B$).  

Hence, given a coefficient free AE sentence,  a specialization of the universal variables $y$ in $A*B$ that is a witness 
for the failure of the sentence over $A*B$, is also a witness for the failure of the sentence over $A*B*F$,
and with a specialization of the universal variables $y$ in $A*B*F$ that is
a witness for the failure of the sentence over $A*B*F$, it is possible to associate a (test) sequence of specializations in
$A*B$,  that are all witnesses for the failure of the sentence over $A*B$, and these specializations (in $A*B$)
"converge" into the specialization of the variables $y$ in $A*B*F$.

This fact, and the identity between ungraded, graded, and (graded) formal Makanin-Razborov diagrams over the collection
of non-trivial free products, $A*B$ and $A*B*F$, imply that the (non-canonical) finite
collection of (ungraded) resolutions, $Res_1(z,y),\ldots,Res_d(z,y)$, and
their finite collection of auxiliary resolutions, that
are associated in theorem 3.1 
with a given coefficient-free sentence, and the entire collection of non-trivial free products $A*B$ (that are
not isomorphic to $D_{\infty}$), can be taken to be the collection of ungraded resolutions and their  auxiliary resolutions
that is associated (by theorem 3.1) with the same coefficient free sentence, and all the non-trivial free products  $A*B*F$.
This implies that a  coefficient-free AE  sentence is a truth sentence over a non-trivial free product, $A*B$, that is not
isomorphic to $D_{\infty}$, if and only if it is a truth sentence over $A*B*F$.

By the same arguments, given a coefficient-free AE set, $AE(p)$, the (non-canonical) 
finite collection of graded resolutions and their auxiliary
resolutions, that is associated with the set $AE(p)$ and the entire collection of free products $A*B$, in theorem 4.1, can be
taken to be the finite collection of graded resolutions and their auxiliary resolutions that is associated (by theorem 4.1)
with the given AE set $AE(p)$
and the collection of free products, $A*B*F$.

We continue by induction, similar to the induction that was used in proving theorem 6.1 for general definable sets over free
products. We argued that the finite collection of graded resolutions and their auxiliary resolutions that are associated
with a given coefficient-free  AE set can be to taken to be the same over the collection of non-trivial free products $A*B$
(that are not isomorphic to $D_{\infty}$), and over the collection of non-trivial free products, $A*B*F$. Hence, 
given a coefficient-free  EAE set, the 
finite collection of graded resolutions and their auxiliary resolutions that are associated with it (by theorem 4.5) over
the entire collection of non-trivial free products, $A*B$, that are not isomorphic to $D_{\infty}$, can be taken to be the 
finite collection of such graded resolutions and auxiliary resolutions that are associated by the same theorem with the
collection of non-trivial free products, $A*B*F$. Continuing inductively, we obtain the same conclusion for the
finite collection of graded resolutions and their associated auxiliary resolutions that is associated with
a given coefficient-free definable set (by theorem 6.1) over the collection of non-trivial free products, $A*B$, and
over the collection of non-trivial free products, $A*B*F$.

Theorem 6.3 uses the finite collection of graded resolutions  and their auxiliary resolutions that is associated with a 
coefficient-free definable
set, to reduce a given coefficient free sentence from a non-trivial free product to a disjunction of conjunctions of sentences
over its factors.

The finite collection of graded resolutions and their auxiliary resolutions that is associated with
a given coefficient-free definable set, can be taken to be identical over the entire collection of non-trivial free products, 
$A*B$, and over the entire collection of free products, $A*B*F$. Hence, we may apply the proof of theorem 6.3, and deduce that
a given coefficient-free sentence over the entire collection of non-trivial free products, $A*B$, that are not
isomorphic to $D_{\infty}$, and over
the entire collection of non-trivial free products, $A*B*F$, reduces to the same (finite) disjunction of conjunctions of 
sentences over the factors $A$ and $B$. Hence, a given coefficient-free sentence is a truth sentence over a given non-trivial
free product, $A*B$, that is not isomorphic to $D_{\infty}$, if and only if it is a truth sentence over $A*B*F$. Therefore,
every such non-trivial free product $A*B$ is elementarily equivalent to $A*B*F$.  

\line{\hss$\qed$}

Note that Tarski's problem for free groups follows if we take $A$ and $B$ to be isomorphic to $Z$ in the statement of theorem 7.2. 
Also, note that exactly the same argument proves that a non-trivial free product, $A*B$, that is not isomorphic to 
$D_{\infty}$, is elementarily equivalent to a tower over $A*B$, i.e., a completion that has in its bottom level a free
product of the form, $A*B*F$, where $F$ is a (possibly trivial) free group 
(cf. theorem 7 in [Se7] in the free group case, and  theorems 7.6 and 7.10 in [Se8] in the torsion-free hyperbolic
analogue).

Theorem 7.1 implies that for every group $G$ that is a non-trivial free product, and is not the infinite dihedral group, $G$ is
elementarily equivalent to $G*F$ (where $F$ is a free group). In [Se8] it is proved that this holds for every non-elementary (torsion-free) hyperbolic 
group. By the
combination of theorems 7.1 and 7.2, the collection of groups $G$ for which $G$ is elementarily equivalent to $G*F$ is an 
elementary class, i.e., 
if $G$ is elementarily equivalent to $G*F$, and $H$ is elementarily equivalent to $G$, then $H$ is elementarily equivalent to $H*F$.
It is then natural to ask what are the properties of groups in this elementary class.

\vglue 1pc
\proclaim{Question} What are the (algebraic, first order) properties of groups $G$ for which $G$ is elementarily equivalent to $G*F$? 
\endproclaim

Other rather straightforward corollaries of theorems 6.1 and 6.3, are uniform properties of sentences over free products.
 
\vglue 1pc
\proclaim{Theorem 7.3} Let $\Phi$ be a coefficient free sentence over groups. There  exists an integer, $k(\Phi)$, so
that for every group, $H$, $\Phi$ is a truth sentence over $H_1*\ldots*H_{k(\Phi)}$, $H_i \eqsim H$, if and only if
$\Phi$ is a truth sentence over $H_1*\ldots*H_n$, $H_i \eqsim H$, for every $n \geq k(\Phi)$.
\endproclaim

\nfp  By theorem 6.1, with the coefficient free sentence, $\Phi$, it is possible to associate finitely many
(ungraded) resolutions (over free products), $DRes_1,\ldots,DRes_g$, and with each of them finitely many auxiliary resolutions,
so that the sentence $\Phi$ is a truth sentence over a non-trivial free product, $G=A^1*\ldots *A^{\ell}$, that is not
$D_{\infty}$, if and only if there exist (or there do not exist) 
specializations of the terminal limit groups of the resolutions, $DRes_1,\ldots,DRes_g$,
with properties that are listed in the statement of the theorem. 

Note that the finite collection of  resolutions, $DRes_1,\ldots,DRes_{\ell}$, is universal (although it is
not canonical), i.e., it does not depend on
the free product $G$, nor on the number of factors in the free product $G$. In theorem 6.3, we used the existence of this finite 
collection of resolutions, to reduce the coefficient-free sentence $\Phi$ to a disjunction of conjunctions of sentences in the
factors of the free product $G$, $A^1,\ldots,A^{\ell}$. This disjunction of conjunctions of sentences does depend
on the number of factors, even though the  resolutions and the auxiliary resolutions 
that are associated with $\Phi$ by theorem 6.1,
do not (depend on the number of factors).

For each subset $I$ that contains a resolution, $DRes_i$, and a subset of its (finitely many) associated auxiliary resolutions, we set 
$t_I$ to be the sum of the numbers of elliptic factors in the terminal limit groups of all the resolutions in the subset $I$ ($DRes_i$
and the subset of its associated auxiliary resolutions). Clearly, there are finitely many such subsets $I$. 
We set $k(\Phi)$ to be the maximum between the sum of the numbers $t_I$ over all the possible subsets $I$, and the number 3 (to guarantee
that the corresponding free product is not $D_{\infty}$, in case $G$ is non-trivial).  

Since the free products we are looking at are free products of the same group $G$,  theorem 6.1 implies that $\Phi$ is a truth sentence
over an iterated free product of $G$ with itself $k(\Phi)$ times, if and only if it is a truth sentence over an iterated free 
product of $G$ with itself a larger number of times.

\line{\hss$\qed$}

Note that the integer $k(\Phi)$ depends on the coefficient free sentence, $\Phi$, but it does not depend on the group, $G$.
It is easy to see that $k(\Phi)$ can not be chosen to be a universal constant, e.g.,  we can take $\Phi_m$ to be a sentence that 
specifies if the number of conjugacy classes of involutions in the group is at least $m$. For such a sentence, 
$\Phi_m$, $k(\Phi_m)=m$.

Theorem 7.3 can be further strengthened for sequences of groups.
 Let $\Phi$ be a coefficient free sentence over groups. Given any sequence of groups, $G_1,G_2,\ldots$,
we set $M_1=G_1$, $M_2=G_1*G_2$, $M_3=G_1*G_2*G_3$, and so on. The sentence $\Phi$ may be truth or false on any of the groups (free
products) $M_i$, $i=1,\ldots$. Here one can (clearly) not guarantee that the sentence $\Phi$ is constantly truth or 
constantly false 
staring at a bounded index (of the $M_i$'s). However, one can prove the following.

\vglue 1pc
\proclaim{Theorem 7.4} There exists an integer $c(\Phi)$, so that for every sequence of groups,
$G_1,G_2,\ldots$, the sentence $\Phi$ over the sequence of groups, $M_1=G_1,M_2=G_1*G_2,\ldots$ may change signs (from truth to false
or vice versa) at most $c(\Phi)$ times.
\endproclaim

\nfp The argument that we use is similar to the proof of theorem 7.3.
By theorem 6.1, with the coefficient free sentence, $\Phi$, it is possible to associate finitely many
(ungraded) resolutions, and with each of them finitely many auxiliary resolutions,
so that the sentence $\Phi$ is a truth sentence over a non-trivial free product, $G=A^1*\ldots *A^{\ell}$, that is not
$D_{\infty}$, if and only if there exist (or there do not exist) specializations of the terminal limit groups of the 
given list of resolutions, 
with properties that are listed in the statement of the theorem. 

Once again, we note that the collection of resolutions that is associated (by theorem 6.1) 
with the given coefficient free sentence $\Phi$,
and their associated finite collections of auxiliary resolutions, are universal, which means that they are good
for all the non-trivial free factors that are not isomorphic to $D_{\infty}$, and they do not depend on the number of factors
in such a free product.
In theorem 6.3, this finite collection of resolutions is  used  
to reduce the coefficient-free sentence $\Phi$ to a disjunction of conjunctions of sentences in the
factors of the free product. This disjunction of conjunctions of sentences does depend
on the number of factors, even though the  resolutions and the auxiliary resolutions 
that are associated with $\Phi$ by theorem 6.1,
do not (depend on the number of factors).

The (finitely many) resolutions and their auxiliary resolutions, that are associated with the coefficient-free
sentence $\Phi$ by theorem 6.1, allows one to reduce the sentence $\Phi$ into boundedly many sentences on the various factors 
of a given non-trivial free product, that is not isomorphic to $D_{\infty}$, 
where the given sentence $\Phi$ is true or false (over the ambient free product) if and only
if a (finite) disjunction of conjunctions of these sentences over the factors is true or false (see theorem 6.3). 
The length of this disjunction of conjunctions 
depend on the number of factors in the given free product, but the number of distinct (coefficient-free) sentences over the various factors is
uniformly bounded, and depend only on the resolutions that are associated with $\Phi$ by theorem 6.1, and on their auxiliary resolutions. 

We set $selp$ to be the sum of the numbers of elliptic factors in the terminal limit groups of all the 
(finitely many) resolutions (and auxiliary resolutions) that are associated with the coefficient free sentence $\Phi$ by theorem 6.1.
By the proof of theorem 6.3, the number of distinct coefficient-free sentences that are defined over the factors of a given non-trivial 
free product, such that the given sentence $\Phi$ is equivalent to a disjunction of conjunctions of these sentences over the various factors
of the given non-trivial free product, is bounded by the number of subsets of a set of size $selp$, i.e, it is bounded by 
$2^{selp}$.

We define the $state$ of a given factor of a free product, to be the subcollection of distinct coefficient-free sentences,  
that are defined over the factors of a free product,
and are associated with $\Phi$ by theorems 6.1 and 6.3, that are truth sentences over this given factor. Clearly, there are at most
$2^{2^{slp}}$ possible states for factors of a free product. 

Let $G_1,G_2,\ldots$ be a sequence of groups, with the corresponding  free products: $M_1=G_1,M_2=G_1*G_2,\ldots$. By the proof of theorem 
6.3, the  question whether the given 
coefficient-free sentence $\Phi$ is true or false depends only on  the  states of the factors, $G_1,\ldots,G_n$. Furthermore,
to determine if $\Phi$ is true or false, it suffices to know  for every possible factor how many factors do have this state, and for every
possible state it suffices to know the number of factors that have this state only up to $selp$, as a larger number will make no difference
for the question of $\Phi$ being true or false (i.e., for each state it suffices to know the minimum between the actual number of factors
that have this state and the number $selp$).

 Hence, with each free product $M_n$ from the sequence, we can associate a tuple of at most $2^{2^{selp}}$ non-negative integers that are all
bounded by $selp$. The sequence of tuples that are associated with the free products, $\{M_n\}$, can only increase (in the lexicographical
order), and $\Phi$ may change from true to false or vice versa, only when the tuple changes (increases). Hence,  
$\Phi$ may change from true to false or vice versa over
the sequence of free factors, $\{M_n\}$, at most at $selp \cdot 2^{2^{selp}}$ indices, and the theorem follows. 

\line{\hss$\qed$}

\vglue 1.5pc
\centerline{\bf{\S8. Stability}} 
\medskip

In the previous section we deduced several basic first order properties of free products from the reduction of sentences from a
free product to its factors (theorem 6.3), and from the association of finitely many  resolutions and their auxiliary resolutions with any
given sentence, that enables this reduction (theorem 6.1). This reduction, and the associated resolutions,
enable one to prove that a free product inherits certain first order properties from its
factors. In this section we use the scheme of argument that was used to prove
the stability of free and hyperbolic groups [Se10], to prove that stability is
inherited by a free product from its factors. This question about the possibility to lift stability from the factors to a free
product was brought to our attention by Eric Jaligot, and was the motivation for our entire work on free products.

\vglue 1pc
\proclaim{Theorem 8.1} Let $A,B$ be stable groups. Then $A*B$ is stable.
\endproclaim

\nfp We may assume that the free product $A*B$ is non-trivial and not isomorphic to $D_{\infty}$. 
As in the proof of stability of
free and (torsion-free) hyperbolic groups, we prove the stability of a free product gradually. We start with the stability of
coefficient-free varieties (that can
be deduced from complex algebraic geometry in the case of free groups), (coefficient-free) Diophantine sets, (coefficient-free) 
rigid and weakly solid sets (i.e., the set of specializations of the defining parameters for which there exists a rigid
or a strictly solid specialization of a given rigid or a solid limit group),  
coefficient-free definable sets, and finally the stability of every definable set. 

\vglue 1pc
\proclaim{Theorem 8.2} Let $A,B$ be stable groups. Then every coefficient-free variety over $A*B$ is stable.
\endproclaim

\nfp The approach that we used for proving stability of free and hyperbolic groups, associates certain finite
diagrams with some families (Diophantine, rigid and solid) of definable sets, and using these diagrams we further associate
(finite) canonical collections of $Duo$ $limit$ $groups$ with these definable sets. These finite diagrams and their associated
Duo limit groups, together with the bounds on the number of rigid and strictly solid families of specializations of
rigid and solid limit groups over free and (torsion-free) hyperbolic groups, are the main tools  that enable one to prove stability of free
and hyperbolic groups.   

Unfortunately, we do not know how to imitate this approach in the free product case, and we won't construct canonical
and universal diagrams as in the free group case. Instead we argue by contradiction. We suppose that there exists a coefficient-free
system of equations, $\Sigma(p,q)=1$, such that the corresponding variety, $V(p,q)$, is not stable over a non-trivial free
product $G=A*B$, where both $A$ and $B$ are stable. With such unstable variety $V(p,q)$, we associate (non-canonically) two
diagrams that are similar to the (canonical and universal) diagram that was constructed in the case of free and hyperbolic groups,
and from the existence of the constructed diagrams, and the stability of the factors, $A$ and $B$, we deduce a 
contradiction, that finally proves theorem 8.2.

Let $\Sigma(p,q)=1$  be a coefficient-free system of equations. Let $G=A*B$ be a non-trivial free product in which
both
$A$ and $B$ are stable, and suppose that the variety $V(p,q)$, that corresponds to the system of equations $\Sigma(p,q)=1$, 
is unstable over $G$. Since $V(p,q)$ is unstable over $G$, for every positive integer $m$, there exists two sequences of tuples
with elements in $G$,
$\{p^m_i\}_{i=1}^m$ and $\{q^m_j\}_{j=1}^m$, such that $(p^m_i,q^m_j) \in V(p,q)$ if and only if $j \leq i$.

Given the triangle of pairs, $\{(p^m_i,q^m_i) \ 1\leq i \leq m\}$, we can pass to a subtriangle (still denoted with the same 
indices), such that every sequence of pairs, $\{(p^m_{i_m},q^m_1)\}_{m=1}^{\infty}$, converges into the same limit group
(over free products), $L(p,q)$. By passing to a further subtriangle, we may assume that every such sequence,
$\{(p^m_{i_m},q^m_1)\}_{m=1}^{\infty}$, converges into the same graded resolution, 
with respect to the parameter subgroup $<q>$: 
$L=L_0 \to L_1 \to \ldots \to L_s$, where 
$L_s$ is a free product of a rigid or a solid factor with (possibly) finitely many elliptic factors, and (possibly) a free
group. We denote this graded resolution, $GRes_1$.

At this stage we pass to a further subtriangle (still denoted with the same indices). 
In this subtriangle, we may further assume that every sequence of triples, $\{(p^m_{i_m},q^m_1,q^m_2)\}_{m=2}^{\infty}$, where for
every index $m$, $i_m \geq 2$, converges into the same limit group, $U(p,q_1,q_2)$. We apply the construction of quotient
resolutions, as it appears in the first step of the sieve procedure [Se6], and pass to a further subtriangle, so that
we may assume that every sequence 
of triples, $\{(p^m_{i_m},q^m_1,q^m_2)\}_{m=2}^{\infty}$, where for
every index $m$, $i_m \geq 2$, converges into the same quotient resolution with respect to the parameter subgroup $<q_1,q_2>$
(as constructed according to the first step
of the sieve procedure). 

We continue iteratively. At each step we pass to a subtriangle, and apply the general step
of the construction of a quotient resolution, as it appears in the sieve procedure [Se6].
Note that the resolution that is constructed in step $n$ of the procedure, is graded
with respect to the parameter subgroup, $<q_1,\ldots,q_n>$. Also note that 
each quotient resolution is either a closure of the quotient resolution that was constructed
in the previous step, or it is a proper quotient resolution, which means that its complexity
(as it appears in the sieve procedure) is strictly smaller than the complexity of the resolution
that was constructed in the previous step.
 
 To conclude the first part of the argument, 
and continue with it to prove stability of varieties, 
we need to prove that this (first) iterative procedure that is associated with an unstable variety
terminates after finitely many steps.

\vglue 1pc
\proclaim{Proposition 8.3} There are only finitely many steps along the  iterative procedure that is associated with an
unstable variety, $V(p,q)$, over a given non-trivial free product, $G=A*B$, in which the constructed graded resolution
is not a closure of the graded resolution that was constructed in the previous step of the procedure.   
\endproclaim

\nfp Suppose that there are infinitely many steps in which the constructed graded resolution is not a closure of the graded
resolution that was constructed in the previous step of the procedure. This contradicts the termination of the sieve 
procedure after finitely many steps (theorem 22 in [Se6], see also the proof of theorem 4.4).

\line{\hss$\qed$}

By proposition 8.3 there exists some step $n_0$ of the iterative procedure, such that starting at this step, the constructed 
quotient resolutions along the iterative procedure, are closures of the quotient resolution that was constructed in step $n_0$.
First, we replace the quotient resolution that was constructed in step $n_0$, by a cover resolution, according to theorem 1.21.
We denote this cover resolution, $CGRes$.
By the construction of the cover resolution, $CGRes$, we can pass to a further
subtriangle of specializations, $\{p^m_i\}_{i=1}^m$, and $\{q^m_j\}_{j=1}^m$, such that:
\roster
\item"{(1)}" $(p^m_i,q^m_j) \in V(p,q)$ if and only if $j \leq i$.

\item"{(2)}" the pairs $(p^m_i,q^m_j)$ extend to  specializations that factor through the (cover) graded resolutions, $CGRes$,
for $j \leq i$, and $1 \leq j \leq n_0$. 
\endroster

The (cover) graded resolution, $CGRes$,
terminates in a free product of a rigid or weakly solid factor with (possibly) finitely
many elliptic factors. By theorems 1.14 and 1.15, with the rigid or weakly solid factor, one can associate finitely many 
(combinatorial) configurations, so that  each configuration contains finitely many fractions, and finitely
many elliptic elements, and a rigid or almost shortest (weakly) strictly solid specialization is given by
fixed words in the fractions and the elliptic elements. The value of these 
fractions depend only on the specialization of the defining parameters, and not on the (specific) rigid or 
almost shortest (weakly) strictly solid
specialization, whereas the elliptic elements do depend on the specific rigid or weakly strictly solid specialization
(and not only on the specialization of the defining parameters). See theorems 1.14 and 1.15 for the precise details.

Once again, we can pass to a further subtriangle of specializations, $\{p^m_i\}_{i=1}^m$ and $\{q^m_j\}_{j=1}^m$,  
that satisfy properties (1) and (2),
and assume that the specializations of the rigid or weakly solid factor of the terminal limit group of $CGRes$ that are
associated with the specializations, $\{(p^m_i,q^m_j)\}$, from the subtriangle, are all associated with  one (fixed) combinatorial
configuration (out of  the finitely many 
combinatorial configurations) that is presented in theorems 1.14 and 1.15,  and is associated with the rigid or
weakly solid factor of the terminal limit group of $CGRes$. 

We denote the fractions that appear in the (fixed)  combinatorial configuration (that get the same values for all 
the rigid or almost
shortest strictly solid specializations that are associated with the same specializations of the parameters ($q^m_1,\ldots,q^m_{n_0}$),
and with the fixed combinatorial configuration), $v^m_1,\ldots,v^m_f$. We denote the elliptic elements that are associated with 
the values of the parameters ($q^m_1,\ldots,q^m_{n_0}$), $eq^m_1,\ldots,eq^m_d$, and with the rigid or almost shortest (weakly) strictly 
solid specializations: 
$ep^m_1,\ldots,ep^m_g$ (see the statements of theorems 1.14 and 1.15 for these notions). For brevity we denote
the tuple $v^m_1,\ldots,v^m_f$ by $v^m$, the tuple $eq^m_1,\ldots,eq^m_d$ by $eq^m$, and the tuple $ep^m_1,\ldots,ep^m_g$ by $ep^m$.

With the variety $V$, the (cover) graded resolution, $CGRes$, the triangle of specializations that factor through $CGRes$,
the terminal limit group of $CGRes$,  
and the (fixed) combinatorial configuration that is associated with the
rigid or weakly solid factor of its terminal limit group, we associate a triangle of specializations, and a second iterative
procedure, in which each quotient resolution is ungraded, hence, terminates in a free product of f.g.\ elliptic factors.

First, with the specializations, $q^m_1,\ldots,q^m_{n_0}$, 
and the fixed combinatorial configuration  that is associated with
the rigid or weakly solid factor of the terminal limit group of $CGRes$, we associate the fixed specializations of the
elements, $v^m_1,\ldots,v^m_f$, and of the elliptic elements, $eq^m$. With each element, $p^m_i$, $n_0 < i \leq m$,
from the triangle
of specializations that is associated with $CGRes$, we associate a tuple of elliptic  specializations, $ep^m_i$, that is
associated with it and with the fixed combinatorial configuration that is associated with the rigid or weakly solid factor
of the terminal limit group of $CGRes$, and with the fixed values of the elements, $v^m$. 

We continue with the
triangle of specializations, $\{ep^m_i\}$ and $\{q^m_j\}$, $n_0 < i,j \leq m$. We swap each line of the triangle, i.e.,
we replace a pair, $(ep^m_i,q^m_i)$, with the pair, $(ep^m_{m-i+n_0+1},q^m_{m-i+n_0+1})$, for $n_0 < i \leq m$. After this swap,
the corresponding pair, $(p^m_i,q^m_j) \in V(p,q)$, $n_0 < i,j \leq m$,  if and only if $i \leq j$. 

Given the triangle of pairs, and the specialization, $v^m$ and $eq^m$, we can pass to a subtriangle (still denoted with the same 
indices), such that every sequence of pairs, $\{(ep^m_{n_0+1},q^m_{j_m},v^m,eq^m)\}_{m=1}^{\infty}$, converges into the same limit group
(over free products), $E$. By passing to a further subtriangle, we may assume that every such sequence,
converges into the same ungraded resolution (that can be viewed as multi-graded with respect to the elliptic tuples, $ep^m$ and $eq^m$): 
$E=E_0 \to E_1 \to \ldots \to E_s$, where 
$E_s$ is a free product of finitely many f.g.\ elliptic factors and possibly a free group. We denote this ungraded resolution, $ERes_1$.

As we did in the first iterative procedure, at this stage we pass to a further subtriangle (still denoted with the same indices), and
apply the construction of quotient
resolutions, as it appears in the first step of the sieve procedure [Se6], and pass to a further subtriangle, so that
we may assume that every sequence 
of tuples, 
 $\{(ep^m_{n_0+1},ep^m_{n_0+2},q^m_{j_m},v^m,eq^m)\}_{m=1}^{\infty}$, $n_0+2 < j_m \leq m$,
 converges into the same quotient resolution (which is an ungraded resolution),
as it is constructed according to the first step
of the sieve procedure. 

We continue iteratively. At each step we pass to a subtriangle, and apply the general step
of the construction of a quotient resolution, as it appears in the sieve procedure [Se6].
As in the first iterative procedure that we associated with $V(p,q)$, 
each quotient resolution is either a closure of the quotient resolution that was constructed
in the previous step, or it is a proper quotient resolution, which means that its complexity
(as it appears in the sieve procedure) is strictly smaller than the complexity of the resolution
that was constructed in the previous step. The second iterative procedure, satisfies a similar termination property
as the first iterative procedure.
 
\vglue 1pc
\proclaim{Proposition 8.4} There are only finitely many steps along the (second)  iterative procedure that is associated with an
unstable variety, $V(p,q)$, over a given non-trivial free product, $G=A*B$, in which the constructed ungraded resolution
is not a closure of the ungraded resolution that was constructed in the previous step of the procedure.   
\endproclaim

\nfp Identical to the proof of proposition 8.3.

\line{\hss$\qed$}

By proposition 8.4 starting at some step, $n_1>n_0$, all the quotient resolutions that are constructed along
the various steps of the (second) iterative procedure, are closures of the quotient (ungraded) resolution that was 
constructed in step $n_1$. 

\noindent
We denote the ungraded resolution that is constructed in step $n_1$, $ERes$. With $ERes$ we associate a cover (ungraded) resolution,
$CERes$, according to theorem 1.21. As (by theorem 1.21) the completion and the terminal limit group of the cover
resolution, $CERes$, are finitely presented, we may pass to a further subtriangle of the tuples specializations that were used to
construct the resolution $ERes$, such that all the specializations in this subtriangle factor through $CERes$.

The cover resolution, $CERes$, the subtriangle of specializations that are associated with it, together with
the stability of the factors $A$ and $B$ of the (non-trivial) free factor, $G=A*B$, contradict the instability
of the variety, $V(p,q)$, over $G$. 

\noindent
Indeed, since the constructed quotient resolutions are only replaced by closures
along the second iterative (starting with step $n_1$), after passing to a further subtriangle, with each specialization, $q^m_j$,
$n_1 < j \leq m$, we can associate a tuple of elliptic specializations $eq^m_j$. By the construction of the (ungraded)
Makanin-Razborov 
diagram of a f.p.\ group (theorem 26 in [Ja-Se]), the pairs of specializations, $(p^m_i,q^m_j) \in V(p,q)$, for $n_1 \leq i,j < m$,
if and only if the tuples of elliptic elements, $ep^m_i,eq^m_j$, satisfy one out of finitely many system of equations. Since,
$(p^m_i,q^m_j) \in V(p,q)$ if and only if $i \leq j$, from the given subtriangle of specializations, we can extract larger
and larger sequences of specializations, $\{ep^m_i\}$ and $\{eq^m_j\}$, $n_1 \leq i,j \leq m$, such that the (elliptic) tuples, 
$(ep^m_i,eq^m_j)$ satisfy a fixed coefficient-free system of equations (independent of the indices, $i,j,m$), if and only
if $i \leq j$. This clearly contradicts the stability of varieties over the factors $A$ and $B$. Hence, every coefficient-free
variety over an arbitrary free product of stable groups is stable, and we get theorem 8.2.

\line{\hss$\qed$}

Theorem 8.2 proves the stability of varieties over non-trivial free products, $G=A*B$, in case both $A$ and $B$ are stable. A similar
argument implies the stability of Diophantine sets over such non-trivial free products.

\vglue 1pc
\proclaim{Theorem 8.5} Let $A,B$ be stable groups that are not both isomorphic to $Z_2$. Then every coefficient-free 
Diophantine set over $A*B$ is stable.
\endproclaim

\nfp The argument that we used to prove stability of varieties, essentially applies to prove stability of Diophantine sets as well.
We assume that there exists an unstable  coefficient-free Diophantine set over a non-trivial free product, $G=A*B$, that is
not isomorphic to $D_{\infty}$, in which both $A$ and $B$ are stable. We apply the same iterative procedures that we applied in the
case of varieties. After finitely many steps the quotient resolutions that are produced in each of the two iterative procedures 
are closures of the quotient resolution that was constructed in the previous step. This follows 
 by the argument that proves  termination  of the iterative procedure
for the analysis of an EAE set (theorem 4.4), that is based on the termination of the sieve procedure (theorem 22 in [Se6]), 
and the proof of
equationality of Diophantine sets over free and (torsion-free) hyperbolic groups in theorem 2.2 in [Se9].

Once one proves that after finitely many steps,
quotient resolutions that are constructed along the second iterative procedure are closures of the quotient resolutions that were
constructed in the previous step, we get a contradiction to the instability of the original Diophantine set, in the same
way a contradiction was extracted for unstable varieties in the proof of theorem 8.2.

\line{\hss$\qed$}

Theorem 8.5 proves that coefficient free Diophantine sets are stable over free products of stable groups. As over
free and hyperbolic groups, our next step towards stability of the theory of free products of stable groups, is 
the analysis of rigid specializations of rigid limit groups, and (weakly) strictly solid specializations of solid limit groups
(cf. section 4 in [Se9] for the analysis of the corresponding sets over free and hyperbolic groups). 

Over free products, we looked at the set of parameters for which there are at least $s$ rigid or $s$ families of strictly
solid specializations. Unfortunately, over free products there are infinitely many such (rigid and strictly solid) families.
In theorems 1.14 and 1.15 we proved combinatorial bounds on the collection of  rigid and (weakly) strictly solid families,
that give bounds on the number of rigid and (weakly) strictly solid families over free products (that are different than the
families over free groups). However, these families and their number is not definable over free products. They are definable
if we enrich the language and allow  quantifiers over elements in the factors of the free product. Hence, instead of looking
at  sets (of parameters) with at least $s$ families of rigid or weakly strictly solid families, we look at  sets of
the defining parameters for which there exists a rigid or a weakly strictly solid specialization.

As in the case of a free group, the stability of Diophantine sets
(theorem 8.5) plays an essential role in proving the stability of these sets that are naturally associated with rigid and 
weakly solid limit groups. However, the iterative procedures that served us in proving stability of varieties and Diophantine sets
need to be further refined. 

\vglue 1pc
\proclaim{Theorem 8.6} Let $Sld(x,p,q)$ be a coefficient-free (weakly) solid  limit group over free products. Suppose that
$Sld$ embeds into a f.p.\ (ungraded) completion, and that with $Sld$ there is an associated (finite) cover of its flexible
quotients, so that every graded limit group in this cover embeds into an (ungraded) f.p.\ completion (see theorems 1.20 and 
 1.21 for
the precise details on such weakly  solid limit groups). Note that
we say that a homomorphism of $Sld(x,p,q)$ into a free product is (weakly) strictly solid, if it can be extended to
a specialization of the f.p.\ completion that is associated with $Sld$, and it does not extend to a specialization of
any of the f.p.\ completions that are associated with the (finite) cover of the flexible quotients of
$Sld$.

We set $ES(p,q)$ to be the (coefficient-free definable) set of specializations of the (free variables) parameters $p,q$, 
for which there exists
a (weakly) strictly solid  homomorphism of $Sld(x,p,q)$.  If $A$ and $B$ are
  non-trivial  stable groups that are not both isomorphic to $Z_2$, then $ES(p,q)$ is stable.
\endproclaim

\nfp The argument that we use is a modification of the argument  that was used to prove stability of varieties and Diophantine sets (theorems 8.2
and 8.5). 
Let $Sld(p,q)$ be a weakly solid limit group, and let $ES(p,q)$ be the (coefficient-free) definable set that contains those
values of the parameters (free variables) for which $Sld$ has a weakly strictly solid specialization.  Let $G=A*B$ be a non-trivial 
free product, that is not isomorphic to $D_{\infty}$, in which
both
$A$ and $B$ are stable, and suppose that the definable set, $ES(p,q)$, is unstable over $G$.
Since $ES(p,q)$ is unstable over $G$, for every positive integer $m$, there exists two sequences of tuples
with elements in $G$,
$\{p^m_i\}_{i=1}^m$ and $\{q^m_j\}_{j=1}^m$, such that $(p_i,q_j) \in ES(p,q)$ if and only if $j \leq i$.

Given the triangle of pairs, $\{(p^m_i,q^m_i) \ 1\leq i \leq m\}$, we can pass to a subtriangle (still denoted with the same 
indices), and with each pair, $(p^m_i,q^m_1)$, we associate a weakly strictly solid specialization of the weakly solid
limit group $Sld(x,p,q)$, $(x^m_{i,1},p^m_i,q^m_1)$,  such 
that every sequence of weakly strictly solid specializations, $\{(x^m_{i_m,1},p^m_{i_m},q^m_1)\}_{m=1}^{\infty}$, converges into the same limit group
(over free products), $L(x,p,q)$. By passing to a further subtriangle, we may assume that every such sequence,
$\{(x^m_{i_m,1},p^m_{i_m},q^m_1)\}_{m=1}^{\infty}$, 
converges into the same graded resolution, 
with respect to the parameter subgroup $<q>$: 
$L=L_0 \to L_1 \to \ldots \to L_s$, where 
$L_s$ is a free product of a rigid or a solid factor with (possibly) finitely many elliptic factors and (possibly) a 
free group. We denote this graded resolution, $GRes_1$.

At this stage we look for a  further subtriangle (still denoted with the same indices), for which: 
\roster
\item"{(1)}" for each pair $(p^m_i,q^m_2)$, there exists a weakly strictly solid specialization of $Sld(x,p,q)$, $(x^m_{i,2},p^m_i,q^m_2)$, 
$2 \leq i \leq m$.

\item"{(2)}" 
 every sequence of weakly strictly solid specializations (of $Sld(x,p,q)$): $$\{(x^m_{i_m,2},p^m_{i_m},q^m_2)\}_{m=1}^{\infty}$$ 
$2 \leq i_m \leq m$, converges into the same limit group, $U(x,p,q)$.

\item"{(3)}"
 every sequence of pairs of weakly strictly solid specializations:
$$\{(x^m_{i_m,1},x^m_{i_m,2},p^m_{i_m},q^m_1,q^m_2)\}_{m=1}^{\infty}$$ 
$2 \leq i_m \leq m$, converges into the same  quotient resolution, that is constructed according to the first step
of the sieve procedure [Se6].

\item"{(4)}" the quotient resolution that is constructed from the convergent sequences of  pairs of weakly strictly solid specializations,
is a proper quotient resolution of the graded resolutions that was constructed from the convergent sequences, 
$\{(x^m_{i_m,1},p^m_{i_m},q^m_1)\}_{m=1}^{\infty}$, in the first step of the iterative procedure, i.e., the quotient resolution that
is 
constructed in the second step of the iterative procedure is not a graded closure of the graded resolution that was constructed in the first step,
but rather it is a quotient resolution of strictly smaller complexity according to the sieve procedure [Se6].
\endroster

We continue iteratively. At each step $n$ we look for a further subtriangle, for which the corresponding sequences of
weakly strictly solid specializations (of $Sld(x,p,q)$), 
  $\{(x^m_{i_m,n},p^m_{i_m},q^m_n)\}_{m=1}^{\infty}$, 
$n \leq i_m \leq m$, converge into the same limit group, and the combined weakly strictly solid specializations converge into
the same quotient resolution, that is constructed according to the general step of the sieve procedure. Furthermore, the quotient resolution
that is constructed in the $n$-th step of the procedure is not a closure of the quotient resolution that was constructed in step $n-1$, but
rather it is a quotient resolution of strictly smaller complexity according to the sieve procedure [Se6].   

By proposition 8.3 this iterative procedure terminates after finitely many steps. When it terminates we are left with a triangle
of specializations, (still denoted)  
$\{p^m_i\}_{i=1}^m$ and $\{q^m_j\}_{j=1}^m$, which is a subtriangle of the original triangle of specializations in $G=A*B$, for
which  $(p_i,q_j) \in ES(p,q)$ if and only if $j \leq i$. Furthermore, suppose that the iterative procedure terminated at step $n_1$. Then:
\roster
\item"{(1)}" for each pair $(p^m_i,q^m_j)$, there exists a weakly strictly solid specialization of $Sld(x,p,q)$, $(x^m_{i,j},p^m_i,q^m_j)$, 
where $1 \leq j \leq n_1$ and $n_1+1 \leq i \leq m$.

\item"{(2)}"
 every sequence of tuples of weakly strictly solid specializations: 
$$\{(x^m_{i_m,1},\ldots,x^m_{i_m,n_1},p^m_{i_m},q^m_1,\ldots,q^m_{n_1})\}_{m=n_1+1}^{\infty}$$ 
$n_1+1 \leq i_m \leq m$, converges into the same  quotient resolution, that was constructed according to the general step
of the sieve procedure [Se6].
\endroster

We denote the graded resolution that is constructed in the $n_1$-th step of the iterative procedure, $GRes_{n_1}$. By theorem 1.21 we
can associate a cover graded resolution with $GRes_{n_1}$, that we denote, $CGRes_{n_1}$, that satisfies 
all the properties that are
listed in theorem 1.21, and after possibly passing to a further subtriangle, we may assume that all the tuples of specializations,
$$\{(x^m_{i_m,1},\ldots,x^m_{i_m,n_1},p^m_{i_m},q^m_1,\ldots,q^m_{n_1})\}_{m=n_1+1}^{\infty}$$
$n_1+1 \leq i_m \leq m$, factor through the cover graded resolution, $CGRes_{n_1}$.

At this point we swap the raws in the subtriangle of specializations. For each $i$, $n_1+1 \leq i \leq m$, we replace $p^m_i$
with $p^m_{m+n_1+1-i}$, and for each $j$,
$n_1+1 \leq j \leq m$, we replace $q^m_j$
with $q^m_{m+n_1+1-j}$. Note that after this swap, for each pair of indices, $n_1+1 \leq i,j \leq m$, $(p^m_i,q^m_j) \in ES(p,q)$
if and only if $i \leq j$.

Each specialization $p^m_i$, $n_1+1 \leq i \leq m$, belongs to some fiber of the graded resolution, $CGRes_{n_1}$. If there is 
a subtriangle of our given triangle, (still denoted) $\{p^m_i\}$ and $\{q^m_j\}$, $n_1+1 \leq i,j \leq m$, 
for which for every triple $(i,j,m)$,
$n_1+1 \leq j < i \leq m$, there is no sequence of weakly strictly solid specializations of $Sld(x,p,q)$, 
$\{(x^m_{i,j}(s),p^m_i(s),q^m_j)\}$, where the sequence, $\{p^m_i(s)\}$, is a restriction to the elements $p$ of a test sequence 
in the fiber that contains $p^m_i$, we reached a terminal point of the first procedure (for proving the stability of
$ES(p,q)$), and we continue with the subtriangle, and the graded resolutions, $CGRes_{n_1}$, and the graded resolution, 
$GRes_{n_1}$, to the second part of the proof, i.e., to a procedure that allows us to analyze ungraded resolutions and not
graded ones, as we did along the second iterative procedure in proving the stability of varieties and Diophantine sets
(theorems 8.2 and 8.5).

Otherwise, if there is no subtriangle with this property, there exists a subtriangle of the given triangle, 
$\{p^m_i\}$ and $\{q^m_j\}$, $n_1+1 \leq i,j \leq m$, 
for which for every triple $(i,j,m)$,
$n_1+1 \leq j < i \leq m$, there exists a sequence of weakly strictly solid specializations of $Sld(x,p,q)$, 
$\{(x^m_{i,j}(s),p^m_i(s),q^m_j)\}$, where the sequence, $\{p^m_i(s)\}$, is a restriction to the elements $p$ of a test sequence 
in the fiber that contains $p^m_i$. 

In this case there exists a  further subtriangle
of the triangle that is associated with $CGRes_{n_1}$, (that we still denote) 
$\{p^m_i\}$ and $\{q^m_j\}$, $n_1+1 \leq i,j \leq m$, 
for which:
\roster
\item"{(1)}" for every pair $(i,m)$,
$n_1+1  < i \leq m$, there exists a sequence of weakly strictly solid specializations of $Sld(x,p,q)$,
$\{(x^m_{i,j}(s),p^m_i(s),q^m_{n_1+1})\}$, and an associated test sequence of the cover graded resolution, $CGRes_{n_1}$:
$$(z^m_{i}(s),p^m_i(s),q^m_1,\ldots,q^m_{n_1})$$ that is all in the fiber that contains $p^m_i$. 

\item"{(2)}" every combined sequence,
$\{(z^m_{i_m}(s_m),x^m_{i_m,n_1+1}(s_m),p^m_{i_m}(s_m),q^m_1,\ldots,q^m_{n_1},q^m_{n_1+1})\}$, where $n_1+1  < i_m \leq m$, and
$m< s_m$, converges into the same graded resolution, $FGRes_{n_1+1}$ (it is graded with respect to the parameter subgroup
$<q_1,\ldots,q_{n_1+1},q_{n_1+1}>$. In particular, the sequence:
$$\{(z^m_{i_m}(s_m),p^m_{i_m}(s_m),q^m_1,\ldots,q^m_{n_1})\}$$ converges into the graded resolution, $GRes_{n_1}$.
\endroster

The graded resolution, $FGRes_{n_1+1}$, can be viewed as a formal graded resolution, and it has the same structure as the graded
resolution, $GRes_{n_1}$. With the graded resolution, $FGRes_{n_1+1}$, we associate a cover graded resolution according to theorem
1.21, that we denote, $CFGRes_{n_1+1}$. We can pass to a further subtriangle, so that the sequences of specializations:  
$\{(z^m_{i_m}(s_m),x^m_{i_m,n_1+1}(s_m),p^m_{i_m}(s_m),q^m_1,\ldots,q^m_{n_1},q^m_{n_1+1})\}$, where $n_1+1  < i_m \leq m$, and
$m< s_m$, actually factor through the graded cover resolution, $CFGRes_{n_1+1}$. 

As for $i,j$, $n_1 < i,j \leq m$, $(p^m_i,q^m_j) \in E(p,q)$ if and only if $i \leq j$, the tuple,
$\{(x^m_{i,n_1+1},p^m_{i},q^m_{n_1+1})\}$, that extends to a specialization of $CFGRes_{n_1+1}$, can not be weakly strictly solid
for $n_1+1  < i \leq m$. Hence, on the specializations that factor through $CFGRes_{n_1+1}$ we may further impose one
of finitely many Diophantine conditions that force the specializations, 
$\{(x^m_{i,n_1+1},p^m_{i},q^m_{n_1+1})\}$, that extends to a specialization of $CFGRes_{n_1+1}$, to be non weakly strictly solid
for $n_1+1  < i \leq m$. 

We further pass to a subtriangle, so that every sequence of specializations, 
$(x^m_{i_m,n_1+1},p^m_{i_m},q^m_{n_1+1})$, $n_1+1  < i_m \leq m$, 
together with the specialization of the Diophantine condition that 
forces it to be non weakly strictly solid, converges into a quotient resolution of $FGRes_{n_1+1}$, that we denote, $GRes_{n_1+1}$.
Since generic points in $FGRes_{n_1+1}$ restrict to weakly strictly solid specializations, $GRes_{n_1+1}$, is a proper
quotient resolution of $FGRes_{n_1+1}$, i.e., it is a quotient resolution of strictly smaller complexity than $GRes_{n_1}$ 
(in light of the sieve procedure).

We continue iteratively. At each step we first look for a subtriangle, 
(still denoted) $\{p^m_i\}$ and $\{q^m_j\}$, $ i,j \leq m$, 
for which for every triple $(i,j,m)$,
$ j < i \leq m$, there is no sequence of weakly strictly solid specializations of $Sld(x,p,q)$, 
$\{(x^m_{i,j}(s),p^m_i(s),q^m_j)\}$, where the sequence, $\{p^m_i(s)\}$, is a restriction to the elements $p$ of a test sequence 
in the fiber that contains $p^m_i$ in the cover of the graded resolution that was constructed in the previous step. 
If there is such a subtriangle, we reached a terminal point of this part of the procedure.

Otherwise, if there is no subtriangle with this property, there exists a subtriangle of the given triangle, 
$\{p^m_i\}$ and $q^m_j$, $ i,j \leq m$, 
for which for every triple $(i,j,m)$,
$ j < i \leq m$, there exists a sequence of weakly strictly solid specializations of $Sld(x,p,q)$, 
$\{(x^m_{i,j}(s),p^m_i(s),q^m_j)\}$, where the sequence, $\{p^m_i(s)\}$, is a restriction to the elements $p$ of a test sequence 
in the fiber that contains $p^m_i$ in the cover of the graded resolution that was constructed in the previous step. 

In this case we (possibly)  pass to a further subtriangle, and collect test sequences that extend to weakly strictly solid 
specializations. The combined specializations converge into a graded resolution that is similar to $FGRes_{n_1+1}$. On the
specializations in a further subtriangle, that factor through a cover graded resolution, we impose one of finitely many possible 
Diophantine conditions, and the combined specializations converge into a quotient resolution which is a proper quotient
resolution (of strictly smaller complexity) of the resolution that was constructed in the previous step. 

By the termination of the sieve procedure (theorem 22 in [Se6], see also the proof of theorem 4.4), this iterative procedure 
terminates after finitely many steps. We denote its terminating step, $n_2$, and clearly $n_2>n_1$. 
When it terminates there exists a subtriangle
for which for every triple $(i,j,m)$,
$ j < i \leq m$, there is no sequence of weakly strictly solid specializations of $Sld(x,p,q)$, 
$\{(x^m_{i,j}(s),p^m_i(s),q^m_j)\}$, where the sequence, $\{p^m_i(s)\}$, is a restriction to the elements $p$ of a test sequence 
in the fiber that contains $p^m_i$ in the cover of the graded resolution that was constructed in step $n_2$. 

At this point we once again swap the raws in the subtriangle of specializations. For each $i$, $n_2+1 \leq i \leq m$, 
we replace $p^m_i$
with $p^m_{m+n_2+1-i}$, and for each $j$,
$n_2+1 \leq j \leq m$, we replace $q^m_j$
with $q^m_{m+n_2+1-j}$. After this swap, for each pair of indices, $n_2+1 \leq i,j \leq m$, $(p^m_i,q^m_j) \in ES(p,q)$
if and only if $j \leq i$.

We continue according to the first iterative procedure, i.e., at each step we collect weakly strictly solid specializations,
that converge into a proper quotient resolution (a resolution of strictly smaller complexity) of the graded resolution that
was constructed in the previous step. This iterative procedure terminates after finitely many steps. After it terminates we
swap the raws in the triangle once again, and continue iteratively along the steps of the second iterative procedure. This
(second) iterative procedure terminates after finitely many steps, and we continue iteratively. At each step
we first swap the raws, and then continue along the first or second iterative procedure until they terminate (after
finitely many steps), swap, and continue along the other procedure.

This iterative procedure of swapping and continuing along the two iterative procedures until they terminate, terminates after finitely
many steps, by the proof of theorem 4.4. Once it terminates, we are left with a graded resolution, $GRes_{n_k}$, which is
graded with respect to the parameter subgroup, $<q_1,\ldots,q_{n_k}>$, a cover of this graded resolution that satisfies the 
properties of cover graded resolutions that are listed in theorem 1.21, $CGRes_{n_k}$, 
and a
subtriangle of specializations, 
(still denoted) $\{p^m_i\}$ and $\{q^m_j\}$, $n_k+1 \leq i,j \leq m$ that extend to specializations that factor through the cover graded
resolution, $CGRes_{n_k}$, for which $(p^m_i,q^m_j) \in ES(p,q)$ if and only if $n_k+1 \leq j \leq i \leq m$. Furthermore, each specialization
$p^m_i$ extends to a specialization that factors through $CGRes_{n_k}$. A test sequence of specializations in the fiber of
$CGRes_{n_k}$ that contains $p^m_i$, restricts to specializations, $\{p^m_i(s)\}$, that extend to weakly strictly solid specializations
of $Sld(x,p,q)$, $\{(x^m_{i,j}(s),p^m_i(s),q^m_j)\}$, if and only if $n_k+1 \leq j \leq i \leq m$.    

\medskip
Now, we can continue to the second part as we did in proving the stability of varieties and Diophantine sets (theorems 8.2 and 8.5).
The (cover) graded resolution, $CGRes_{n_k}$,
terminates in a free product of a rigid or a weakly solid factor with (possibly) finitely
many elliptic factors and (possibly) a free group. 
By theorems 1.14 and 1.15, with the rigid or solid factor, one can associate finitely many 
(combinatorial) configurations, so that  each configuration contains finitely many fractions, and finitely
many elliptic elements. A rigid or almost shortest (weakly) strictly solid specialization is given by
fixed words in the fractions and the elliptic elements. The value of these 
fractions depend only on the specialization of the defining parameters, and not on the (specific) rigid or 
almost shortest (weakly) strictly solid
specialization, whereas the elliptic elements do depend on the specific rigid or weakly strictly solid specialization
(and not only on the specialization of the defining parameters). See theorems 1.14 and 1.15 for the precise details.

We can pass to a further subtriangle of specializations, $\{p^m_i\}_{i=n_k+1}^m$ and $\{q^m_j\}_{j=n_k+1}^m$, so that
$(p^m_i,q^m_j) \in ES(p,q)$ if and only if $j \leq i$,   
and assume that the specializations of the rigid or weakly solid factor of the terminal limit group of $CGRes_{n_k}$ that are
associated with the specializations, $\{(p^m_i,q^m_j)\}$, from the subtriangle, are all associated with  one (fixed) combinatorial
configuration (out of  the finitely many 
combinatorial configurations) that is presented in theorems 1.14 and 1.15,  and is associated with the rigid or
weakly solid factor of the terminal limit group of $CGRes_{n_k}$.

We denote the fractions that appear in the (fixed)  combinatorial configuration (that get the same values for all 
the rigid or almost
shortest strictly solid specializations that are associated with the same specializations of the parameters ($q^m_1,\ldots,q^m_{n_k}$),
and with the fixed combinatorial configuration), $v^m_1,\ldots,v^m_f$, and for brevity $v^m$. 
We denote the elliptic elements that are associated with 
the values of the parameters ($q^m_1,\ldots,q^m_{n_k}$), $eq^m_1,\ldots,eq^m_d$, and for brevity $eq^m$,
and with the rigid or almost shortest (weakly) strictly 
solid specializations: 
$ep^m_1,\ldots,ep^m_g$, and for brevity $ep^m$ (see the statements of theorems 1.14 and 1.15 for these notions). 

At this point we start the second part of our argument that involves two interchanging iterative procedures like the first part,
and in which all the resolutions are ungraded, i.e., they terminate in  elliptic subgroups, in a similar way to the second part
of the argument that was used to prove stability of varieties and Diophantine sets.  

First, 
With each element, $p^m_i$, $n_k < i \leq m$,
from the triangle
of specializations that is associated with the cover graded resolution, $CGRes_{n_k}$, we associate a tuple of elliptic  
specializations, $ep^m_i$, that is
associated with it and with the fixed combinatorial configuration that is associated with the rigid or weakly solid factor
of the terminal limit group of $CGRes_{n_k}$, and with the fixed values of the elements, $v^m$. 

We continue with the
triangle of specializations, $\{ep^m_i\}$ and $\{q^m_j\}$, $n_k < i,j \leq m$. We swap each line of the triangle, i.e.,
we replace a pair, $(ep^m_i,q^m_i)$, with the pair, $(ep^m_{m-i+n_k+1},q^m_{m-i+n_k+1})$, for $n_k < i \leq m$. After this swap,
the corresponding pair, $(p^m_i,q^m_j) \in ES(p,q)$, $n_k < i,j \leq m$,  if and only if $i \leq j$. 

We apply the first iterative procedure for this triangle of specializations. 
Given the triangle of pairs, and the specialization, $v^m$ and $eq^m$, we can pass to a subtriangle (still denoted with the same 
indices), such that every sequence of pairs, $\{(ep^m_{n_k+1},q^m_{j_m},v^m,eq^m)\}_{m=1}^{\infty}$, $n_k+1 \leq j_m \leq m$,
converges into the same limit group
(over free products), $E$. By passing to a further subtriangle, we may assume that every such sequence,
converges into the same ungraded resolution (that can be viewed as multi-graded with respect to the elliptic tuples, 
$ep^m_{n_k+1}$
 and $eq^m$): 
$E=E_0 \to E_1 \to \ldots \to E_s$, where 
$E_s$ is a free product of finitely many f.g.\ elliptic factors and possibly a free group. We denote this ungraded resolution, $ERes_1$.

We continue iteratively. At each step  we look for a further subtriangle, for which the corresponding sequences of
weakly strictly solid specializations (of $Sld(x,p,q)$), 
  $\{(x^m_{n,i_m},p^m_{n},q^m_{j_m})\}_{m=1}^{\infty}$, 
$n \leq j_m \leq m$, converge into the same limit group, and the combined weakly strictly solid specializations converge into
the same quotient resolution, that is constructed according to the general step of the sieve procedure. Note that the
constructed quotient resolution is ungraded and terminates in a free product of elliptic factors and possibly a 
free group. Furthermore, the quotient resolution
that is constructed in the general step of the procedure is not a closure of the quotient resolution that was constructed in 
the previous step, but
rather it is a quotient resolution of strictly smaller complexity according to the sieve procedure [Se6].   

By proposition 8.3 this iterative procedure terminates after finitely many steps, that we denote $\ell_1$ ($\ell_1>n_k$). 
When it terminates we are left with an ungraded resolution, that we denote $ERes_{\ell_1}$. With the ungraded resolution,
$ERes_{\ell_1}$ we associate a cover resolution, $CERes_{\ell_1}$, that satisfies the properties of theorem 1.21. With
$CERes_{\ell_1}$, there is an associated triangle
of specializations, (still denoted)  
$\{ep^m_i\}_{i=\ell_1+1}^m$ and $\{q^m_j\}_{j=\ell_1+1}^m$, which is a subtriangle of the original triangle of specializations in $G=A*B$, for
which  the corresponding specialization, $(p^m_i,q^m_j) \in ES(p,q)$, if and only if $i \leq j$. Furthermore, for every pair, 
$(ep^m_i,q^m_j)$,
for which $\ell_1+1 \leq i \leq j \leq m$, the triple $(ep^m_i,v^m,q^m_j)$ extends to a specialization that factors through the cover resolution
$CERes_{\ell_1}$.

At this point we continue as we did in the first part of the argument. First, we 
swap the raws in the subtriangle of specializations. 
Note that after this swap, for each pair of indices, $\ell_1+1 \leq i,j \leq m$, the corresponding specializations, $(p^m_i,q^m_j) \in ES(p,q)$
if and only if $j \leq i$.

Each specialization $(v^m,q^m_j)$, $\ell_1+1 \leq j \leq m$, belongs to some fiber of the graded resolution, $CERes_{\ell_1}$. If there is 
a subtriangle of our given triangle, (still denoted) $\{ep^m_i\}$ and $\{q^m_j\}$, $\ell_1+1 \leq i,j \leq m$, 
for which for every triple $(i,j,m)$,
$\ell_1+1 \leq i < j \leq m$, there is no sequence of weakly strictly solid specializations of $Sld(x,p,q)$, 
$\{(x^m_{i,j}(s),p^m_i(s),q^m_j(s))\}$, where the sequence, $\{(p^m_i(s),q^m_j(s))\}$, is a restriction to the pair $(p,q)$ of a test sequence 
in the fiber that contains $q^m_j$, we reached a terminal point of the  procedure (i.e.,  a terminal point of the second part 
of the argument for proving stability of $ES(p,q)$). In this case we can extract a contradiction, and conclude stability, 
from this subtriangle, and the graded resolution
$CERes_{\ell_1+1}$.  

Otherwise, if there is no subtriangle with this property, there exists a subtriangle of the given triangle, 
$\{ep^m_i\}$ and $\{q^m_j\}$, $\ell_1+1 \leq i,j \leq m$, 
for which for every triple $(i,j,m)$,
$\ell_1+1 \leq i < j \leq m$, there exists a sequence of weakly strictly solid specializations of $Sld(x,p,q)$, 
$\{(x^m_{i,j}(s),p^m_i(s),q^m_j(s))\}$, where the sequence, $\{(p^m_i(s),q^m_j(s))\}$, is a restriction to the pair $(p,q)$ of a test sequence 
in the fiber that contains $q^m_j$. 

In this case there exists a  further subtriangle
of the triangle that is associated with $CERes_{\ell_1}$, (that we still denote) 
$\{p^m_i\}$ and $\{q^m_j\}$, $\ell_1+1 \leq i,j \leq m$, 
for which:
\roster
\item"{(1)}" for every pair $(j,m)$,
$\ell_1+1  < j \leq m$, there exists a sequence of weakly strictly solid specializations of $Sld(x,p,q)$, 
$\{(x^m_{\ell_1+1,j}(s),p^m_{\ell_1+1}(s),q^m_j(s))\}$, and an associated test sequence of the cover  resolution, $CERes_{\ell_1}$:
$$\{(z^m_{\ell_1+1,j}(s),v^m(s),p^m_{\ell_1+1}(s),q^m_j(s),eq^m_1,\ldots,eq^m_{n_k},ep^m_{n_k+1},\ldots,ep^m_{\ell_1})\}$$ that is all in the fiber that contains $(v^m,q^m_j)$. 

\item"{(2)}" every combined sequence:
$$\{(z^m_{\ell_1+1,j_m}(s_m),x^m_{\ell_1+1,j_m}(s_m),v^m(s_m),p^m_{\ell_1+1}(s_m),q^m_{j_m}(s_m),eq^m_1,\ldots,eq^m_{n_k},
ep^m_{n_k+1},\ldots,ep^m_{\ell_1},ep^m_{\ell_1+1}
)\}$$ 
where $\ell_1+1 < j_m \leq m$, and
$m< s_m$, converges into the same ungraded resolution, $FERes_{\ell_1+1}$.
In particular, the sequence:
$$\{(z^m_{\ell_1+1,j_m}(s_m),v^m(s_m),p^m_{\ell_1+1}(s_m),q^m_{j_m}(s_m),eq^m_1,\ldots,eq^m_{n_1},
ep^m_{n_k+1},\ldots,ep^m_{\ell_1}
)\}$$ converges into the ungraded resolution, $ERes_{\ell_1}$.
\endroster

The  resolution, $FERes_{\ell_1+1}$, can be viewed as a formal (ungraded) resolution, and it has the same structure as the 
resolution, $ERes_{\ell_1}$. With the  resolution, $FERes_{\ell_1+1}$, we associate a cover  resolution, with a f.p.\ completion and
terminal limit group, according to theorem
1.21, that we denote, $CFERes_{\ell_1+1}$. We can pass to a further subtriangle, so that the sequences of specializations:  
$$\{(z^m_{\ell_1+1,j_m}(s_m),x^m_{\ell_1+1,j_m}(s_m),v^m(s_m),p^m_{\ell_1+1}(s_m),q^m_{j_m}(s_m),eq^m_1,\ldots,eq^m_{n_k},
ep^m_{n_k+1},\ldots,ep^m_{\ell_1},ep^m_{\ell_1+1}
)\}$$ 
actually factor through the cover resolution, $CFERes_{\ell_1+1}$.

As for $i,j$, $\ell_1+1 \leq i,j \leq m$, $(p^m_i,q^m_j) \in E(p,q)$ if and only if $j \leq i$, the tuple,
$\{(x^m_{\ell_1+1,j},p^m_{\ell_1+1},q^m_j)\}$, $\ell_1+1 < j \leq m$, 
that extends to a specialization of $CFERes_{\ell_1+1}$, can not be weakly strictly solid.
Hence, on the specializations that factor through $CFERes_{\ell_1+1}$ we may further impose one
of finitely many Diophantine conditions that force the specializations, 
$\{(x^m_{\ell_1+1,j},p^m_{\ell_1+1},q^m_j)\}$, $\ell_1+1 < j \leq m$, that extends to a specialization of $CFERes_{\ell_1+1}$, 
to be non weakly strictly solid.

We further pass to a subtriangle, so that every sequence of specializations, 
$(x^m_{\ell_1+1,j_m},p^m_{\ell_1+1},q^m_{j_m})$, $\ell_1+1  < j_m \leq m$, and their extensions to specializations of $CFERes$, 
together with the specialization of the Diophantine condition that 
forces it to be non weakly strictly solid, converges into a quotient resolution of $FERes_{\ell_1+1}$, that we denote, $ERes_{\ell_1+1}$.
Since generic points in $FERes_{\ell+1_1+1}$ restrict to weakly strictly solid specializations, $ERes_{\ell_1+1}$, is a proper
quotient resolution of $FERes_{\ell_1+1}$, i.e., it is a quotient resolution of strictly smaller complexity than $ERes_{\ell_1}$ 
(in light of the sieve procedure).

We continue iteratively, precisely as we did in the iterative procedure that was used in the first part of th argument. 
At each step we first look for a subtriangle, 
(still denoted) $\{ep^m_i\}$ and $\{q^m_j\}$, $ i,j \leq m$, 
for which for every triple $(i,j,m)$,
$ i < j \leq m$, there is no sequence of weakly strictly solid specializations of $Sld(x,p,q)$, 
$\{(x^m_{i,j}(s),p^m_i(s),q^m_j(s))\}$, where the sequence, $\{(p^m_i(s),q^m_j(s))\}$, is a restriction to the elements $p,q$ of a test sequence 
in the fiber that contains $(v^m,q^m_j)$ in the cover of the ungraded resolution that was constructed in the previous step. 
If there is such a subtriangle, we reached a terminal point of this part of the procedure.

Otherwise, if there is no subtriangle with this property, there exists a subtriangle of the given triangle, 
$\{ep^m_i\}$ and $q^m_j$, $ i,j \leq m$, 
for which for every triple $(i,j,m)$,
$ i < j \leq m$, there exists a sequence of weakly strictly solid specializations of $Sld(x,p,q)$, 
$\{(x^m_{i,j}(s),p^m_i(s),q^m_j)\}$, where the sequence, $\{(p^m_i(s),q^m_j(s))\}$, is a restriction to the elements $p,q$ of a test sequence 
in the fiber that contains $(v^m,q^m_j)$ in the cover of the graded resolution that was constructed in the previous step. 

In this case we (possibly)  pass to a further subtriangle, and collect test sequences that extend to weakly strictly solid 
specializations. The combined specializations converge into a graded resolution that is similar to $FERes_{\ell_1+1}$. On the
specializations in a further subtriangle, that extend to specializations that 
factor through a cover  resolution we impose one of finitely many possible 
Diophantine conditions, and the combined specializations converge into a quotient resolution which is a proper quotient
resolution (of strictly smaller complexity) of the resolution that was constructed in the previous step. 

By the termination of the sieve procedure (theorem 22 in [Se6], see also the proof of theorem 4.4), this iterative procedure 
terminates after finitely many steps. We denote its terminating step, $\ell_2$, and clearly $\ell_2>\ell_1$. 
When it terminates there exists a subtriangle
for which for every triple $(i,j,m)$,
$ i < j \leq m$, there is no sequence of weakly strictly solid specializations of $Sld(x,p,q)$, 
$\{(x^m_{i,j}(s),p^m_i(s),q^m_j(s)\}$, where the sequence, $\{(p^m_i(s),q^m_j(s))\}$, is a restriction to the elements $p,q$ of a test sequence 
in the fiber that contains $(v^m,q^m_j)$ in the cover of the graded resolution that was constructed in step $\ell_2$. 

We continue as we did in the first part of the argument.
At this point we once again swap the raws in the subtriangle of specializations. 
 After this swap, for each pair of indices, $\ell_2+1 \leq i,j \leq m$, $(p^m_i,q^m_j) \in ES(p,q)$
if and only if $i \leq j$.

We continue according to the first iterative procedure that was used in the second part of the argument. i.e., 
at each step we collect weakly strictly solid specializations,
that converge into a proper quotient resolution (a resolution of strictly smaller complexity) of the graded resolution that
was constructed in the previous step. This iterative procedure terminates after finitely many steps. After it terminates we
swap the raws in the triangle once again, and continue iteratively along the steps of the second iterative procedure that was used
in the second part of the argument. This
(second) iterative procedure terminates after finitely many steps, and we continue iteratively. At each step
we first swap the raws, and then continue along the first or second iterative procedure until they terminate (after
finitely many steps), swap, and continue along the other procedure.

This iterative procedure of swapping and continuing along the two iterative procedures until they terminate, terminates after finitely
many steps, by the proof of theorem 4.4. Once it terminates, we are left with an ungraded resolution, $ERes_{\ell_t}$, 
 a cover of this  resolution that satisfies the 
properties of cover  resolutions that are listed in theorem 1.21 (and in particular has a f.p.\ completion and terminal limit group), $CERes_{\ell_t}$, 
and a
subtriangle of specializations, 
(still denoted) $\{ep^m_i\}$ and $\{q^m_j\}$, $\ell_t+1 \leq i,j \leq m$ that extend to specializations that factor through the cover 
resolution, $CERes_{\ell_t}$, for which (the corresponding specialization) $(p^m_i,q^m_j) \in ES(p,q)$ if and only if $\ell_t+1 \leq i \leq j \leq m$. 
Furthermore, each specialization
$(v^m,q^m_j)$ extends to a specialization that factors through $CERes_{\ell_t}$. A test sequence of specializations in the fiber of
$CERes_{\ell_t}$ that contains $(v^m,q^m_j)$, restricts to specializations, $\{(p^m_i(s),q^m_j(s))\}$, that extend to weakly strictly solid specializations
of $Sld(x,p,q)$, $\{(x^m_{i,j}(s),p^m_i(s),q^m_j(s))\}$, if and only if $\ell_t+1 \leq i \leq j \leq m$.    

We denote the ungraded resolution that is constructed in step $\ell_t$, $ERes$. With $ERes$ we associate a cover (ungraded) resolution,
$CERes$, according to theorem 1.21. As (by theorem 1.21) the completion and the terminal limit group of the cover
resolution, $CERes$, are finitely presented, we may pass to a further subtriangle of the tuples of specializations that were used to
construct the resolution $ERes$, such that all the specializations in this subtriangle factor through $CERes$.

Recall that the subgroup $<q>$ is embedded in the graded cover $CERes$, and the elements $p$ can be expressed as words in
the elements that are associated with $v^m$ in $CERes$, with the elliptic elements $ep$. Given the cover resolution, $CERes$,we look at 
all its test sequences, for which the restriction of the test sequences to the elements $p,q$, extend to 
weakly strictly solid specializations of the weakly solid limit group, $Sld(x,p,a)$. By the techniques that were used to construct the
formal Makanin-Razborov diagram, with the collection of these test sequences we can associate a formal Makanin-Razborov diagram (see
theorem 2.7), where each of the finitely many resolutions in this formal Makanin-Razborov diagram is a closure of
$CERes$, and each such resolution has a f.p.\ completion and terminal limit group.

A specialization of $Sld(x,p,q)$ is not weakly strictly solid if it satisfies one out of finitely many Diophantine conditions.
Given each of the resolutions in this (formal) Makanin-Razborov diagram, we look at all the test sequences that factor through it,
for which the elements that are associated with the weakly strictly solid specializations of $Sld(x,p,q)$, satisfy one of the finitely
many Diophantine conditions that force them to be non weakly strictly solid specializations. With the collection of these sequences we can
once again associate a formal Makanin-Razborov diagram, in which every resolution has a f.p.\ completion and terminal limit group.  

Now, the subtriangle of specializations that factor through the cover graded resolution, $CERes$, extend to specializations that factor through 
resolutions in the (formal) Makanin-Razborov diagram that is associated with $CERes$. Specialization in the subtriangle
that is associated with $CERes$ satisfy:
$(p^m_i,q^m_j) \in ES(p,q)$ if and only if $i \leq j$. By the iterative procedure that was applied to construct the terminal
ungraded resolution $ERes$, test sequences in $CERes$ that are contained in the fiber of of $(v^m,q^m_j)$, extend to weakly strictly solid specializations
of $Sld(x,p,q)$ if and only if $i \leq j$ as well.

The resolutions in the formal Makanin-Razborov diagrams that are associated with $CERes$ enable one to reduce the question of whether
a fiber in $CERes$ contains a test sequence that extends to weakly strictly solid specializations of $Sld(x,p,q)$ to a
(finite) disjunction of conjunctions of (fixed) AE predicates over  the factors of the free product $G=A*B$. Hence, the existence
of the (infinite) subtriangle of specializations that is associated with $CERes$, contradicts the stability of the constructed AE predicates
over the factors $A$ and $B$. Therefore, $ES(p,q)$ is stable, and theorem 8.6 follows.

\line{\hss$\qed$}

\medskip
Having proved that varieties, Diophantine sets, and sets of specializations of the defining parameters for which there 
exists a rigid or a weakly strictly solid specialization of a rigid or a weakly solid limit group, are stable over a
free product in which the factors are stable, we are finally ready to complete the proof of theorem 8.1, i.e., to prove that
the theory of a free product of stable groups is stable. 

The argument that we use to prove stability of a general definable set over a free product of stable groups,
is based  on the argument that was used
to prove that the set, $ES(p,q)$, is stable over such free products in theorem 8.6. To adapt this approach,
 we use the uniform geometric description of a definable set over free products, that is stated and proved in theorem 6.1.

Let $Def(p,q)$ be a coefficient-free definable set. Wlog we will assume that $Def(p,q)$ is an $E{(AE)}^k$ set.
Recall that theorem 6.1 associates with $Def(p,q)$ finitely
many graded resolutions, $DRes_1,\ldots,DRes_g$, and with each such graded resolution, there are finitely
many associated auxiliary resolutions. The graded resolutions, $DRes_i$,  have the properties of resolutions in the graded Makanin-Razborov
diagram (over free products) of a f.p.\ group (theorem 1.22), and in particular they terminate in a free product of a rigid or
a weakly solid limit group (that can be embedded in a f.p.\ ungraded completion) with finitely many f.p.\ elliptic factors. The auxiliary
resolutions have the same properties as resolutions in a graded formal Makanin-Razborov diagram (theorem 2.7), 
and their terminal limit groups have the same properties as the terminal limit groups of  the graded resolutions, $DRes_i$.

By theorem 6.1, the question of whether a specialization of the defining parameters is in $Def(p,q)$ or in its complement,
over an arbitrary non-trivial free product that is not isomorphic to $D_{\infty}$, depends entirely on the set of
all possible specializations of the terminal limit groups of the graded resolutions, $DRes_i$ (that are associated with $Def(p,q)$),
and their possible extensions to specializations of the terminal limit groups of the (finitely many) associated resolutions (see
theorem 6.1 for details). This enables one to reduce a sentence over free products to sentences over its factors (theorem 6.3), and
even to get a form of quantifier elimination for predicates over free products (theorem 6.4). Note that the graded resolutions that
are constructed in theorem 6.1 and are 
associated with $Def(p,q)$ are good for all free products, regardless of the number of factors,
whereas the reduction of a sentence to its factors, and the quantifier elimination, do depend on the number of factors in the given
free product.

By theorem 6.1 (see also the proof of theorem 6.3),  a specialization of the defining parameters $p,q$ is in $Def(p,q)$, if and
only if
specializations of the terminal limit groups of the graded resolutions, $DRes_i$ (that are associated with $Def(p,q)$), and their 
associated auxiliary resolutions, satisfy certain properties that are listed in theorem 6.1 (and correspondingly in theorem 6.3).

\noindent
Hence, our approach to prove stability of the set $Def(p,q)$ over free products consists of two parts, as in proving stability of the
set $ES(p,q)$ (theorem 8.6). In each part we apply a sequence of interchanging iterative procedures, where the iterative procedures are
similar to the ones that were used in the proof of theorem 8.6, and the way we interchange them depends and follows the requirements
on the terminal limit group of the graded resolutions, $DRes$, that are associated with $Def(p,q)$, and the terminal
limit groups of their associated auxiliary
resolutions, as these requirements appear in theorem 6.1 (and in its proof).  

The terminal limit group of any of those (finitely many) graded resolutions, is a free product of a rigid or a weakly solid factor 
(w.r.t.\ to the parameter
subgroup $<p,q>$) with (possibly) finitely many elliptic factors.
Furthermore, the rigid or weakly solid factor in such a terminal limit group can be embedded in an ungraded f.p.\ completion. 
Hence, each terminal limit group of one of the graded resolutions, $DRes_i$, or one of its associated auxiliary resolutions,
can be embedded into a f.p.\ group, which is the free product of the f.p.\ ungraded completion into which the rigid or
weakly solid factor embeds, free product with the (finitely many) f.p.\ elliptic factors of the terminal limit group.

Let $G=A*B$ be a non-trivial 
free product, that is not isomorphic to $D_{\infty}$, in which
both
$A$ and $B$ are stable, and suppose that the coefficient-free definable set, $Def(p,q)$, is unstable over $G$.
Since $Def(p,q)$ is unstable over $G$, for every positive integer $m$, there exists two sequences of tuples
with elements in $G$,
$\{p^m_i\}_{i=1}^m$ and $\{q^m_j\}_{j=1}^m$, such that $(p_i,q_j) \in Def(p,q)$ if and only if $j \leq i$.

Let $Def(p,q)$ be given by the coefficient-free predicate:
$$  Def(p,q) \ = \ \exists t \ \forall y_1 \ \exists x_1 \ \ldots \ \forall y_k \ \exists x_k \
 \Sigma(t,y_1,x_1,\ldots,y_k,x_k,p,q)=1 \, \wedge$$
$$\wedge \, \Psi(t,y_1,x_1,y_k,x_k,p,q) \neq 1$$ 
By theorem 6.1, given the triangle of pairs, $\{(p^m_i,q^m_j) \ 1 \leq i,j \leq m\}$,   
we can pass to a subtriangle (still denoted with the same 
indices), and with each pair, $(p^m_i,q^m_1)$, 
associate a specialization of the terminal limit group of one of the resolutions,
$DRes$, that are associated with $Def(p,q)$, that restricts to a rigid or a weakly strictly solid specialization of the rigid or 
weakly solid factor of that terminal limit group. Furthermore, these specializations of the terminal limit group of $DRes$, are
supposed to testify that the various specializations, $(p^m_i,q^m_1)$ are in the set $Def(p,q)$. we can pass to a further subtriangle
(still denoted with the same indices), such that these specializations of the terminal limit group of $DRes$,
$(t^m_{i,1},p^m_i,q^m_1)$, $1 \leq i \leq m$,  have the property that every sequence of the form: 
$\{(t^m_{i_m,1},p^m_{i_m},q^m_1)\}_{m=1}^{\infty}$, converges into the same limit group
(over free products), $L(t,p,q)$. By passing to a further subtriangle, we may assume that every such sequence,
$\{(t^m_{i_m,1},p^m_{i_m},q^m_1)\}_{m=1}^{\infty}$, 
converges into the same graded resolution, 
with respect to the parameter subgroup $<q>$: 
$L=L_0 \to L_1 \to \ldots \to L_s$, where 
$L_s$ is rigid or solid. We denote this graded resolution, $GRes_1$.

We continue as we did in the first part of proving that $ES(p,q)$ is stable (theorem 8.6). 
At this stage we look for a  further subtriangle (still denoted with the same indices), for which: 
\roster
\item"{(1)}" for each pair $(p^m_i,q^m_2)$, there exists a specialization of the terminal limit group of one of the graded
resolutions, $DRes$, that restricts to a rigid or a weakly  strictly solid specialization of the rigid or solid factor of that
terminal limit group, $(t^m_{i,2},p^m_i,q^m_2)$, 
$2 \leq i \leq m$. Furthermore, we may assume that these specializations of the terminal limit group of $DRes$ testify that the
specializations, $(p^m_i,q^m_2)$, $2 \leq i \leq m$, are in the set $Def(p,q)$.

\item"{(2)}" 
 every sequence of these specializations (of the terminal limit group of $DRes$), $\{(t^m_{i_m,2},p^m_{i_m},q^m_2)\}_{m=1}^{\infty}$, 
$2 \leq i_m \leq m$, converges into the same limit group, $U(t,p,q)$.

\item"{(3)}"
 every sequence of pairs of these specializations of the terminal limit group of $DRes$, 
$\{(t^m_{i_m,1},t^m_{i_m,2},p^m_{i_m},q^m_1,q^m_2)\}_{m=1}^{\infty}$, 
$2 \leq i_m \leq m$, converges into the same  quotient resolution, that is constructed according to the first step
of the sieve procedure [Se6].

\item"{(4)}" the quotient resolution that is constructed from the convergent sequences of  pairs of specializations of the terminal
limit group of $DRes$,
is a proper quotient resolution of the graded resolutions that was constructed from the convergent sequences, 
$\{(t^m_{i_m,1},p^m_{i_m},q^m_1)\}_{m=1}^{\infty}$, in the first step of the iterative procedure, i.e., the quotient resolution that
is 
constructed in the second step of the iterative procedure is not a graded closure of the graded resolution that was constructed in the first step,
but rather it is a quotient resolution of strictly smaller complexity according to the sieve procedure [Se6].
\endroster

We continue iteratively as we did in the first part of the proof of theorem 8.6. At each step $n$
 we look for a further subtriangle, for which the corresponding sequences of specializations of the terminal limit group of a
resolution $DRes$, that is associated with $Def(p,q)$, that restrict to rigid or weakly strictly solid specializations of the 
rigid or weakly solid factors of this terminal limit group, 
  $\{(t^m_{i_m,n},p^m_{i_m},q^m_n)\}_{m=1}^{\infty}$, 
$n \leq i_m \leq m$, converge into the same limit group, and the combined specializations 
specializations converge into
the same quotient resolution, that is constructed according to the general step of the sieve procedure. Furthermore, the quotient resolution
that is constructed in the $n$-th step of the procedure is not a closure of the quotient resolution that was constructed in step $n-1$, but
rather it is a quotient resolution of strictly smaller complexity according to the sieve procedure [Se6].   

By proposition 8.3 this iterative procedure terminates after finitely many steps. When it terminates we are left with a triangle
of specializations, (still denoted)  
$\{p^m_i\}_{i=1}^m$ and $\{q^m_j\}_{j=1}^m$, which is a subtriangle of the original triangle of specializations in $G=A*B$, for
which  $(p_i,q_j) \in Def(p,q)$ if and only if $j \leq i$. Furthermore, suppose that the iterative procedure terminated at step $n_1$. Then:
\roster
\item"{(1)}" for each pair $(p^m_i,q^m_j)$, there exists a specialization of a terminal limit group of a resolution, $DRes$, $(t^m_{i,j},p^m_i,q^m_j)$, 
where $1 \leq j \leq n_1$ and $n_1+1 \leq i \leq m$, that restricts to a rigid or a weakly strictly solid specialization of the
rigid or weakly solid factor of that terminal limit group.

\item"{(2)}"
 every sequence of these specializations: 
$$\{(t^m_{i_m,1},\ldots,t^m_{i_m,n_1},p^m_{i_m},q^m_1,\ldots,q^m_{n_1})\}_{m=n_1+1}^{\infty}$$ 
$n_1+1 \leq i_m \leq m$, converges into the same  quotient resolution, that was constructed according to the general step
of the sieve procedure [Se6].
%
\endroster

We denote the graded resolution that is constructed in the $n_1$-th step of the iterative procedure, $GRes_{n_1}$. 
By theorem 1.21 we
can associate a cover graded resolution with $GRes_{n_1}$, that we denote, $CGRes_{n_1}$, that satisfy all the properties that are
listed in theorem 1.21, and after possibly passing to a further subtriangle, we may assume that all the tuples of specializations,
$$\{(t^m_{i_m,1},\ldots,t^m_{i_m,n_1},p^m_{i_m},q^m_1,\ldots,q^m_{n_1})\}_{m=n_1+1}^{\infty}$$
$n_1+1 \leq i_m   \leq m$, factor through the cover graded resolution, $CGRes_{n_1}$.

Each specialization $p^m_i$, $n_1+1 \leq i \leq m$, belongs to some fiber of the graded resolution, $CGRes_{n_1}$. If there is 
a subtriangle of our given triangle, (still denoted) $\{p^m_i\}$ and $\{q^m_j\}$, $n_1+1 \leq i,j \leq m$, 
for which for every triple $(i,j,m)$,
$n_1+1 \leq j \leq i \leq m$, there is a sequence of  specializations of the terminal limit group of $DRes$, that restrict
to rigid or weakly strictly solid specializations of the rigid or weakly solid factor of that terminal limit group,
$\{(t^m_{i,j}(s),p^m_i(s),q^m_j)\}$, where the sequence, $\{p^m_i(s)\}$, is a restriction to the elements $p$ of a test sequence 
in the fiber that contains $p^m_i$, and these elements testify that the tuples, $(p^m_i(s),q^m_j)$, are in the definable set,
$Def(p,q)$, we reached a terminal point of the first procedure (for proving the stability of
$Def(p,q)$), and we continue with the subtriangle, and the graded resolutions, $CGRes_{n_1}$, and the graded resolution, 
$GRes_{n_1}$, to the next iterative procedure, that studies those pairs $(p^m_i,q^m_j)$, $n_1+1 \leq i < j \leq m$, that do not belong
to $Def(p,q)$.

Otherwise, if there is no subtriangle with this property, there exists a subtriangle of the given triangle, 
$\{p^m_i\}$ and $\{q^m_j\}$, $n_1+1 \leq i,j \leq m$, 
for which for every triple $(i,j,m)$,
$n_1+1 \leq j < i \leq m$, there exists a sequence of specializations of the terminal limit group of $DRes$, that restrict to
rigid or weakly strictly solid specializations of the rigid or weakly solid factor of that terminal limit group,
$\{(t^m_{i,j}(s),p^m_i(s),q^m_j)\}$, where the sequence, $\{p^m_i(s)\}$, is a restriction to the elements $p$ of a test sequence 
in the fiber that contains $p^m_i$, and:
\roster
\item"{(1)}" for every triple, $(i,j,m)$, $n_1 \leq j \leq i \leq m$, the sequence of specializations, 
$\{(t^m_{i,j}(s),p^m_i(s),q^m_j)\}$, and its associated test sequence of the cover graded resolution, $CGRes_{n_1}$, converge
into a closure of $CGRes_{n_1}$, and the given specialization of the terminal limit group of $DRes$, $(t^m_{i,j},p^m_i,q^m_j)$,
that testifies that $(p^m_i,q^m_j) \in Def(p,q)$, is a specialization of a cover of that closure of $CGRes_{n_1}$.

\item"{(2)}" the sequence of specializations of the terminal limit group of $DRes$:
$$\{(t^m_{i,j}(s),p^m_i(s),q^m_j)\}$$ do not testify that the tuples, $(p^m_i(s),q^m_j(s))$ are in the definable set, $Def(p,q)$.
\endroster

In this case there exists a  further subtriangle
of the triangle that is associated with $CGRes_{n_1}$, (that we still denote) 
$\{p^m_i\}$ and $\{q^m_j\}$, $n_1+1 \leq i,j \leq m$, 
for which:
\roster
\item"{(1)}" for every triple $(i,j,m)$,
$n_1+1 \leq j \leq i \leq m$, there exists a sequence of specializations of the terminal limit group of
an auxiliary resolution, $Y_1DRes$, that is associated with $DRes$, that restricts to rigid or weakly strictly solid specializations
of the rigid or weakly solid factor of that terminal limit group of $Y_1DRes$,  
$\{(y1^m_{i,j}(s),t^m_{i,j}(s),p^m_i(s),q^m_j)\}$, that extends the sequence of specializations of the terminal limit group of $DRes$,
$\{(t^m_{i,j}(s),p^m_i(s),q^m_j)\}$, that is associated with the pair, $(p^m_i,q^m_j)$. Furthermore, these specializations of
the terminal limit group of the auxiliary resolution $Y_1DRes$, testify that the corresponding specializations of $DRes$, 
$\{(t^m_{i,j}(s),p^m_i(s),q^m_j)\}$, do not testify that the tuples, $(p^m_i,q^m_j)$, are in $Def(p,q)$.

\item"{(2)}" every combined sequence:
$$\{(z^m_{i_m,n_1+1}(s_m),y1^m_{i_m,n_1+1}(s_m),t^m_{i_m,n_1+1}(s_m),p^m_{i_m}(s_m),q^m_1,\ldots,q^m_{n_1},q^m_{n_1+1})\}$$ 
where $n_1+1 \leq  i_m \leq 
m$ and
$m< s_m$, converges into the same graded resolution, $FGRes_{n_1+1}$ (it is graded with respect to the parameter subgroup
$<q_1,\ldots,q_{n_1+1}>$. In particular, the sequence:
$$\{(z^m_{i_m,n_1+1}(s_m),p^m_{i_m}(s_m),q^m_1,\ldots,q^m_{n_1})\}$$ converges into the graded resolution, $GRes_{n_1}$.
\endroster

The graded resolution, $FGRes_{n_1+1}$, can be viewed as a formal graded resolution, and it has the same structure as the graded
resolution, $GRes_{n_1}$. With the graded resolution, $FGRes_{n_1+1}$, we associate a cover graded resolution according to theorem
1.21, that we denote, $CFGRes_{n_1+1}$. We can pass to a further subtriangle, so that the sequences of specializations:  
$$\{(z^m_{i_m,n_1+1}(s_m),y1^m_{i_m,n_1+1}(s_m),t^m_{i_m,n_1+1}(s_m),p^m_{i_m}(s_m),q^m_1,\ldots,q^m_{n_1},q^m_{n_1+1})\}$$ 
where $n_1+1 \leq j  i_m \leq m$, and
$m< s_m$, actually factor through the graded cover resolution, $CFGRes_{n_1+1}$. 

The sequences,
$\{(t^m_{i,n_1+1}(s),p^m_i(s),q^m_{n_1+1})\}$, do not testify that the tuples, $(p^m_i,q^m_{n_1+1})$, are in $Def(p,q)$, and their extensions to
specializations of the terminal limit group of the auxiliary resolution, $Y_1DRes$, 
$\{(y1^m_{i,n_1+1}(s),t^m_{i,n_1+1}(s),p^m_i(s),q^m_{n_1+1})\}$, testify for that. Furthermore, the specializations, 
$(t^m_{i,n_1+1},p^m_i,q^m_{n_1+1})$,
$n_1+1  \leq i \leq m$, do testify that $(p^m_i,q^m_{n_1+1}) \in Def(p,q)$. Hence, on the specializations that factor 
through the cover, $CFGRes_{n_1+1}$, we can further impose either one of finitely many Diophantine conditions that forces
the specializations of the elements $y1$, that are specializations of the terminal limit group of an auxiliary resolution,
$Y_1DRes$, to restrict to non-rigid or non weakly strictly solid specialization of the rigid
or weakly solid specialization factor of that terminal limit group, or to impose the existence of extensions of these specializations
of the terminal limit group of $Y_1DRes$, to an auxiliary resolution, $X_1Y_1DRes$ (see theorem 6.1 for these auxiliary resolutions).

We further pass to a subtriangle, so that every sequence of specializations, 
$(y1^m_{i_m,n_1+1},t^m_{i_m,n_1+1},p^m_{i_m},q^m_{n_1+1})$, $n_1+1 \leq  i_m \leq m$, together with their extension to 
specializations of $CFGRes$, and together with the specialization of the 
Diophantine condition that 
forces them to restrict to non-rigid or non weakly strictly solid specializations of the rigid or weakly solid factor of the terminal
limit group of $Y_1DRes$, or together with the corresponding specializations of the terminal limit group of an auxiliary
resolution, $X_1Y_1DRes$, that restrict to rigid or weakly strictly solid specializations of the rigid or weakly solid factor of that
terminal limit group,  converges into a quotient resolution of $FGRes_{n_1+1}$, that we denote, $GRes_{n_1+1}$.

If the obtained quotient resolution is a proper quotient resolution of $GRes_{n_1}$, i.e., not a closure of it, we continue
iteratively, by looking for a subtriangle with test sequences that extend to specializations of the terminal
limit group of one of the auxiliary resolution, $Y_1DRes$, and then force either one of finitely many Diophantine conditions or
a further extension to the terminal limit group of one of the auxiliary resolutions, $X_1Y_1DRes$. 

If the obtained resolution is not a proper quotient resolution of $GRes_{n_1}$, we pass to a further subtriangle, in which test
sequences of a cover of the obtained quotient resolution (which is a closure of $GRes_{n_1}$), extend to specializations of
the terminal limit group of one the auxiliary resolutions, $Y_2X_1Y_1DRes$, that restrict to rigid or weakly strictly solid
specializations of the rigid or weakly solid factor of that terminal limit group.

We continue iteratively, according to the iterative procedure that was applied in the first $n_1$ steps, and according to the 
proof of theorem 6.1 that constructed the graded resolutions, $DRes$, and their associated auxiliary resolutions, that are
all associated with the definable set, $Def(p,q)$. This iterative procedure terminates after finitely many steps according to
proposition 8.3, and the proof of theorem 4.4. When the iterative procedure terminates, at step $n_k$, we are left
with a quotient resolution that we denote, $GRes_{n_k}$, that is graded with respect to the parameters,
$q^m_1,\ldots,q^m_{n_k}$. 

With $GRes_{n_k}$ we associate a cover resolution according to theorem 1.21, that we denote $CGRes_{n_k}$. With $GRes_{n_k}$
and its cover, $CGRes_{n_k}$,  there is an associated subtriangle of specializations,
$\{p^m_i\}$ and $\{q^m_j\}$, $n_k+1 \leq i,j \leq m$, so that every specialization, $p^m_i$ extends to a specialization that factors 
through the cover $CGRes_{n_k}$. By the properties of the iterative procedure that was used to construct the quotient resolution,
$GRes_{n_k}$, for any specialization, $q^m_j$, $n_k+1 \leq j \leq m$, and any specialization, $p^m_i$, $j \leq i \leq m$, there
exists a  sequence of specializations,
$\{(p^m_i(s),q^m_j)\}$, where the specializations, $\{p^m_i(s)\}$, are the restriction to the variables $p$, of a test sequence of
specializations in the fiber that contains $p^m_i$, so that the sequence,  
$\{(p^m_i(s),q^m_j)\}$, extend to specializations,
$\{(t^m_{i,j}(s),p^m_i(s),q^m_j)\}$, that are all specializations of the terminal limit group of a resolution $DRes$, that restrict
to rigid or weakly strictly solid specializations of the rigid or weakly solid factor of that terminal limit group, and each
of these specializations testify that the corresponding specialization,
$\{(p^m_i(s),q^m_j)\}$, is in the definable set $Def(p,q)$.

\medskip
At this point we swap the raws in the subtriangle of specializations, as we did in  proving the stability of the set $ES(p,q)$ 
(theorem 8.6). For each $i$, $n_k+1 \leq i \leq m$, we replace $p^m_i$
with $p^m_{m+n_k+1-i}$, and for each $j$,
$n_k+1 \leq j \leq m$, we replace $q^m_j$
with $q^m_{m+n_k+1-j}$. Note that after this swap, for each pair of indices, $n_k+1 \leq i,j \leq m$, $(p^m_i,q^m_j) \in Def(p,q)$
if and only if $i \leq j$.

Each specialization $p^m_i$, $n_k+1 \leq i \leq m$, belongs to some fiber of the graded resolution, $CGRes_{n_k}$. If there is 
a subtriangle of our given triangle, (still denoted) $\{p^m_i\}$ and $\{q^m_j\}$, $n_k+1 \leq i,j \leq m$, 
for which for every triple $(i,j,m)$,
$n_k+1 \leq j < i \leq m$, there is no sequence of specializations of the terminal limit group of one of the resolutions, $DRes$,
that restrict to rigid or weakly strictly solid specializations of the rigid or weakly solid factor of that terminal limit group, 
$\{(t^m_{i,j}(s),p^m_i(s),q^m_j)\}$, where the sequence, $\{p^m_i(s)\}$, is a restriction to the elements $p$ of a test sequence 
in the fiber that contains $p^m_i$ (in the variety that is associated with $CGRes_{n_k}$), and these specializations of the
terminal limit group of $DRes$, testify that $(p^m_i(s),q^m_j) \in Def(p,q)$,
we reached a terminal point of the first procedure (for proving the stability of
$Def(p,q)$), and we continue with the subtriangle, and the graded resolution, $CGRes_{n_k}$, and the graded resolution, 
$GRes_{n_k}$, to the second part of the proof (in which we analyze ungraded resolutions and not
graded ones, as we did along the second iterative procedure in proving the stability of $ES(p,q)$).

Otherwise,  there is a subtriangle with this property, and we pass to that subtriangle.
In this case there exists a  further subtriangle
of the triangle that is associated with $CGRes_{n_k}$, (that we still denote) 
$\{p^m_i\}$ and $\{q^m_j\}$, $n_k+1 \leq i,j \leq m$, 
for which:
\roster
\item"{(1)}" for every pair $(i,m)$,
$n_k+1  < i \leq m$, there exists a sequence of 
specializations of the terminal limit group of one of the resolutions, $DRes$,
that restrict to rigid or weakly strictly solid specializations of the rigid or weakly solid factor of that terminal limit group, 
$\{(t^m_{i,n_k+1}(s),p^m_i(s),q^m_{n_k+1})\}$, 
 and an associated test sequence of the cover graded resolution, $CGRes_{n_k}$,
$(z^m_{i,n_k+1}(s),p^m_i(s),q^m_1,\ldots,q^m_{n_k})$, that is all in the fiber that contains $p^m_i$. Furthermore, the specializations,
$\{(t^m_{i,n_k+1}(s),p^m_i(s),q^m_{n_k+1})\}$, testify that the tuples, $(p^m_i(s),q^m_{n_k+1}) \in Def(p,q)$.

\item"{(2)}" every combined sequence,
$\{(z^m_{i_m,n_k+1}(s_m),t^m_{i_m,n_k+1}(s_m),p^m_{i_m}(s_m),q^m_1,\ldots,q^m_{n_k},q^m_{n_k+1})\}$, where $n_k+1  < i_m \leq m$, and
$m< s_m$, converges into the same graded resolution, $FGRes_{n_k+1}$ (it is graded with respect to the parameter subgroup
$<q_1,\ldots,q_{n_k+1}$. In particular, the sequence:
$\{(z^m_{i_m,n_k+1}(s_m),p^m_{i_m}(s_m),q^m_1,\ldots,q^m_{n_k})\}$, converges into the graded resolution, $GRes_{n_k}$.
\endroster

With the graded resolution, $FGRes_{n_k+1}$, which is a graded formal closure of $GRes_{n_k}$, 
we associate a cover graded resolution according to theorem
1.21, that we denote, $CFGRes_{n_k+1}$. We can pass to a further subtriangle, so that the sequences of specializations:  
$\{(z^m_{i_m,n_k+1}(s_m),t^m_{i_m,n_k+1}(s_m),p^m_{i_m}(s_m),q^m_1,\ldots,q^m_{n_k},q^m_{n_k+1})\}$, where $n_k+1  < i_m \leq m$, and
$m< s_m$, actually factor through the graded cover resolution, $CFGRes_{n_k+1}$. 

As for $i,j$, $n_k < i,j \leq m$, $(p^m_i,q^m_j) \in Def(p,q)$ if and only if $i \leq j$, on the tuples, $(t^m_{i,j},p^m_i,q^m_j)$,
for $j < i$, that extend to specializations that factor through $CFGRes_{n_k+1}$, 
we may impose one of two conditions. Either they satisfy one of finitely many
Diophantine conditions that force these specializations to restrict to non-rigid or non weakly strictly solid specializations of
the rigid or weakly solid factor of the terminal limit group of $DRes$, or each of these specializations can be extended to
a specialization of the terminal limit group of an auxiliary resolution, $Y_1DRes$, a specialization that restricts to
a rigid or a weakly strictly  solid specialization of the rigid or weakly solid factor of that terminal limit group. Furthermore,
these specializations of the terminal limit group of $Y_1DRes$ testify that a specialization, $(t^m_{i,j},p^m_i,q^m_j)$, does not
testify that $(p^m_i,q^m_j) \in Def(p,q)$.

We further pass to a subtriangle, so that every sequence of specializations, 
$(t^m_{i_m,n_k+1},p^m_{i_m},q^m_{n_k+1})$, $n_k+1 < i_m \leq m$, together with the specialization of the 
Diophantine condition that 
forces it to restrict to non rigid or non weakly strictly solid specialization of the terminal limit group of the associated
graded resolution, $DRes$, or the specialization of the terminal limit group of an auxiliary resolution, $Y_1DRes$,
converges into a quotient resolution of $FGRes_{n_k+1}$, that we denote, $GRes_{n_k+1}$. 
$GRes_{n_k+1}$ is either a proper quotient 
resolution of $FGRes_{n_k+1}$ 
(which is a closure of $GRes_{n_k}$), 
i.e., a quotient resolution of strictly smaller complexity than $GRes_{n_k}$ according to the sieve procedure [Se6],
or it is a closure of $FGRes_{n_k+1}$. If $GRes_{n_k+1}$ is a proper quotient resolution of $GRes_{n_k}$ we reached a terminal
point of this step of the procedure, and we continue to the next step. 

If $GRes_{n_k+1}$ is a closure of $GRes_{n_k}$, we associate with it a cover, $CGRes_{n_k+1}$,
according to theorem 1.21, with a f.p.\ completion
and terminal limit group. We further pass to a subtriangle of specializations, $\{p^m_i\}$ and $\{q^m_j\}$, $n_k+1 \leq i,j \leq m$, 
for which, $(p^m_i,q^m_j) \in Def(p,q)$ if and only if $ i \leq j$, and  the 
associated specializations, $(y1^m_{i,n_k+1},t^m_{i,n_k+1},p^m_i,q^m_{n_k+1})$, $n_k+1  < i \leq m$, 
extend to specializations that 
factor through the cover $CGRes_{n_k+1}$.

In the variety that is associated with the cover, $CGRes_{n_k+1}$, generic points (test sequences) in the
fibers that contain the  specializations that restrict to the tuples, 
$(y1^m_{i,n_k+1},t^m_{i,n_k+1},p^m_i,q^m_{n_k+1})$, $n_k+1 < i \leq m$, restrict to specializations of the terminal limit group of the
auxiliary resolution, $Y_1DRes$, that extend to specializations of the terminal limit group of at least one of the auxiliary
resolutions, $X_1Y_1DRes$, that restrict to rigid or weakly strictly solid specializations of the rigid or weakly solid factor
of the terminal limit group of that auxiliary resolution. Furthermore, these specializations of the terminal limit group
of $X_1Y_1DRes$, testify that the corresponding specializations of $Y_1DRes$, do not form an obstacle for proving that
$(p^m_i(s),q^m_{n_k+1}) \in Def(p,q)$ for $n_k+1  < i \leq m$, and restrictions of test sequences of $CGRes_{n_k+1}$ to the
variables $p$, $\{p^m_i(s)\}$.

We pass to a further subtriangle of specializations, $\{p^m_i\}$ and $\{q^m_j\}$, and use the test sequences that extend to 
specializations of the terminal limit group of an auxiliary resolution, $X_1Y_1DRes$, to associate a graded closure with 
$CGRes_{n_k+1}$. With this graded closure we associate a cover with a f.p.\ completion and terminal limit group according
to theorem 1.21. 

\noindent
We  pass to a further subtriangle, for which the specializations,     
$$(x1^m_{i,n_k+1},y1^m_{i,n_k+1},t^m_{i,n_k+1},p^m_i,q^m_{n_k+1})$$ $n_k+1  < i \leq m$,
extend to specializations that factor through the constructed cover resolution. 
On the tuples, $(x1^m_{i,n_k+1},y1^m_{i,n_k+1},t^m_{i,n_k+1},p^m_i,q^m_{n_k+1})$,
 that extend to specializations that factor through the constructed cover resolution, 
we may impose one of two conditions, similar to the ones that we imposed on the specializations of the cover resolution,
$CGRes_{n_k+1}$. Either they satisfy one of finitely many
Diophantine conditions that force these specializations to restrict to non-rigid or non weakly strictly solid specializations of
the rigid or weakly solid factor of the terminal limit group of $X_1Y_1DRes$, or each of these specializations can be extended to
a specialization of the terminal limit group of an auxiliary resolution, $Y_2X_1Y_1DRes$, a specialization that restricts to
a rigid or a weakly strictly  solid specialization of the rigid or weakly solid factor of that terminal limit group. Furthermore,
these specializations of the terminal limit group of $Y_2X_1Y_1DRes$ testify that the specializations, 
$(x1^m_{i,n_k+1},y1^m_{i,n_k+1},t^m_{i,n_k+1},p^m_i,q^m_{n_k+1})$, do not
testify that $(p^m_i,q^m_j) \in Def(p,q)$ for $n_k+1  < i \leq m$.

We further pass to a subtriangle, so that every sequence of specializations, 
$(x1^m_{i,n_k+1},y1^m_{i,n_k+1},t^m_{i_m,n_k+1},p^m_{i_m},q^m_{n_k+1})$, $n_k+1  < i_m \leq m$, 
together with the specialization of the Diophantine condition that 
forces them to restrict to non rigid or  non weakly strictly solid specializations of the rigid or weakly solid factor of the
terminal limit group of the auxiliary resolution, $X_1Y_1DRes$, or the specialization of the terminal limit group 
of an auxiliary resolution, $Y_2X_1Y_1DRes$,
converges into a quotient resolution of $GRes_{n_k+1}$. This quotient resolution is either a proper quotient 
resolution of $GRes_{n_k+1}$ 
(which is a closure of $GRes_{n_k}$), 
i.e., a quotient resolution of strictly smaller complexity than $GRes_{n_k}$ according to the sieve procedure [Se6],
or it is a closure of $GRes_{n_k+1}$. If the constructed quotient resolution  is a proper quotient resolution of $GRes_{n_k}$ 
we reached a terminal
point of this step of the procedure, and we continue to the next step. 

If the constructed quotient resolution is a closure of $GRes_{n_k}$, we continue iteratively along the steps of the auxiliary
resolutions that are associated with the graded resolution, $DRes$, and the given definable set, $Def(p,q)$, in theorem 6.1.
At each step we first pass to a subtriangle and associate a closure with the previously constructed quotient resolution, and
then impose either a Diophantine condition that forces the collected specializations to restrict to non rigid or non weakly
strictly solid specializations of the rigid or weakly solid factor of the terminal limit group of the corresponding
auxiliary resolution, or force the corresponding specializations to extend to specializations of the terminal limit
group of an associated auxiliary resolution. Then we pass to a further subtriangle and use it to construct a new
quotient resolution that is either a closure of $GRes_{n_k}$ or it is a proper quotient resolution, i.e. a quotient
resolution of strictly smaller complexity according to the sieve procedure.

\noindent
By theorem 6.1, and our assumptions that $\{(p^m_i,q^m_j)\} \in Def(p,q)$ if and only if $n_k+1  < i \leq m$, after finitely
many such steps, we must reach a point in which the constructed quotient resolution is indeed a proper quotient resolution.

We continue iteratively, and this iterative procedure terminates after finitely many steps
by proposition 8.3 and the argument that was used to
prove theorem 4.4. When it terminates we are left with a subtriangle of specializations, $\{p^m_i\}$ and $\{q^m_j\}$, for
which $(p^m_i,q^m_j) \in Def(p,q)$ if and only if $i \leq j$. Furthermore, with the subtriangle there is an associated
quotient resolution, $GRes_{b_1}$, that is graded with respect to a tuple, $q_1,\ldots,q_{b_1}$. With $GRes_{b_1}$ there is
an associated cover (graded) resolution, $CGRes_{b_1}$, and all the specializations $p^m_i$ extend to specializations that
factor through the cover, $CGRes_{b_1}$. Furthermore, for every $b_1+1 \leq j < i \leq m$, there is no test sequence in the
fiber that contains $p^m_i$ (in the variety that is associated with the closure,
$CGRes_{b_1}$), that restrict to specializations, $\{p^m_i(s)\}$, and the tuples $(p^m_i(s),q^m_j) \in Def(p,q)$. 

At this point we once again swap the raws in the subtriangle of specializations, so that after the swap,
in the subtriangle of
specializations, $\{p^m_i\}$ and $\{q^m_j\}$, $b_1+1 \leq i,j \leq m$, $(p^m_i,q^m_j) \in Def(p,q)$ if and only if $j \leq i$. 
We continue according to the first iterative procedure, until we obtain a quotient resolution, with an associated
cover and subtriangle of specializations, so that for each pair in the remaining subtriangle of
specializations, for every $j \leq i$, and $p^m_i$, a  test sequence in the fiber that is associated with $p^m_i$ restricts
to a sequence of specializations, $\{p^m_i(s)\}$, so that $(p^m_i(s),q^m_j) \in Def(p,q)$.

After this first iterative procedure  terminates we
swap the raws in the remaining subtriangle triangle once again, and continue iteratively along the steps of the second iterative procedure. This
(second) iterative procedure terminates after finitely many steps, and we continue iteratively. At each step
we first swap the raws, and then continue along the first or second iterative procedure until they terminate (after
finitely many steps), swap, and continue along the other procedure.

As in the analysis of the set $ES(p,q)$ in theorem 8.6, this iterative procedure of swapping and continuing along the two iterative 
procedures until they terminate, terminates after finitely
many steps, by the proof of theorem 4.4. Once it terminates, we are left with a graded resolution, $GRes_f$, which is
graded with respect to the parameter subgroup, $<q_1,\ldots,q_f>$, a cover of this graded resolution that satisfies the 
properties of cover graded resolutions that are listed in theorem 1.21, $CGRes_f$, 
and a
subtriangle of specializations, 
(still denoted) $\{p^m_i\}$ and $\{q^m_j\}$, $f+1 \leq i,j \leq m$, that extend to specializations that factor through the cover graded
resolution, $CGRes_f$, for which $(p^m_i,q^m_j) \in Def(p,q)$ if and only if $f+1 \leq j \leq i \leq m$. Furthermore, each specialization
$p^m_i$ extends to a specialization that factors through $CGRes_{f}$. A test sequence of specializations in the fiber of
$CGRes_{f}$ that contains $p^m_i$, restricts to specializations, $\{p^m_i(s)\}$, that satisfy $(p^m_i(s),q^m_j) \in Def(p,q)$
 if and only if $f+1 \leq j \leq i \leq m$.

\medskip
Now, we can continue to the second part as we did in proving the stability of varieties, Diophantine sets, and the sets $ES(p,q)$
 (theorems 8.2, 8.5 and 8.6).
The (cover) graded resolution, $CGRes_{f}$,
terminates in a free product of a rigid or a weakly solid factor with (possibly) finitely
many elliptic factors. By theorems 1.14 and 1.15, with the rigid or solid factor, one can associate finitely many 
(combinatorial) configurations, so that  each configuration contains finitely many fractions, and finitely
many elliptic elements, and a rigid or almost shortest (weakly) strictly solid specialization is given by
fixed words in the fractions and the elliptic elements. The value of these 
fractions depend only on the specialization of the defining parameters, and not on the (specific) rigid or 
almost shortest (weakly) strictly solid
specialization, whereas the elliptic elements do depend on the specific rigid or weakly strictly solid specialization
(and not only on the specialization of the defining parameters). 

We can pass to a further subtriangle of specializations, $\{p^m_i\}_{i=f+1}^m$ and $\{q^m_j\}_{j=f+1}^m$, so that
$(p^m_i,q^m_j) \in Def(p,q)$ if and only if $j \leq i$,   
and assume that the specializations of the rigid or weakly solid factor of the terminal limit group of $CGRes_{f}$ that are
associated with the specializations, $\{(p^m_i,q^m_j)\}$, from the subtriangle, are all associated with  one (fixed) combinatorial
configuration (out of  the finitely many 
combinatorial configurations) that is presented in theorems 1.14 and 1.15,  and is associated with the rigid or
weakly solid factor of the terminal limit group of $CGRes_{f}$. 

We denote the fractions that appear in the (fixed)  combinatorial configuration (that get the same values for all 
the rigid or almost
shortest strictly solid specializations that are associated with the same specializations of the parameters ($q^m_1,\ldots,q^m_{f}$),
and with the fixed combinatorial configuration), $v^m_1,\ldots,v^m_e$, and for brevity $v^m$. 
We denote the elliptic elements that are associated with 
the values of the parameters ($q^m_1,\ldots,q^m_{f}$), $eq^m_1,\ldots,eq^m_d$, and for brevity $eq^m$,
and with the rigid or almost shortest (weakly) strictly 
solid specializations: 
$ep^m_1,\ldots,ep^m_g$, and for brevity $ep^m$ (see the statements of theorems 1.14 and 1.15 for these notions). 

At this point we start the second part of our argument that is similar to the second part of the argument that was used to prove
the stability of the set $ES(p,q)$ in theorem 8.6. This second part involves two interchanging iterative procedures, that are
similar to the ones that were used in the  first part,
and in which all the resolutions are ungraded, i.e., they terminate in  elliptic subgroups. 

First, 
with each element, $p^m_i$, $f < i \leq m$,
from the triangle
of specializations that is associated with the cover graded resolution, $CGRes_{f}$, we associate a tuple of elliptic  
specializations, $ep^m_i$, that is
associated with it and with the fixed combinatorial configuration that is associated with the rigid or weakly solid factor
of the terminal limit group of $CGRes_{f}$, and with the fixed values of the elements, $v^m$. 

We continue with the
triangle of specializations, $\{ep^m_i\}$ and $\{q^m_j\}$, $f < i,j \leq m$. We swap each line of the triangle. 
After this swap,
the corresponding pair, $(p^m_i,q^m_j) \in Def(p,q)$, $f < i,j \leq m$,  if and only if $i \leq j$. 

We apply the first iterative procedure for the triangle of specializations. At each step we pass to a further subtriangle, and
require that specializations of the terminal limit group of one of the graded resolution, $DRes$, or of the terminal limit
group of one of its associated auxiliary resolutions, and are associated with (or extend) 
the pairs, $(p^m_i,q^m_j)$, $i \leq j$, will converge into
the same limit group, and the specializations of the quotient resolution that was constructed in the previous step, that are
extended to specializations of these terminal limit group converge into a quotient resolution of the resolution that was constructed in
the previous step. 

Precisely as we argued in the first part of the procedure, after finitely many steps we obtain a quotient resolution, $GRes_{\ell_1}$,
and a subtriangle of specializations, $\{p^m_i\}$ and $\{q^m_j\}$, $\ell_1 +1 \leq i,j \leq m$ that are associated with it, with
associated elliptic specializations, $\{ep^m_i\}$.
With the quotient (ungraded) resolution, we associate a cover (ungraded) resolution, that terminates in a f.p.\ group, such that all the
subtriangle of specializations factor through that cover. Furthermore, for each triple of indices $i,j,m$, $\ell_1 \leq i \leq j 
\leq m$, there exists
a test sequence in the fiber that contains $(v^m,q^m_j)$ (in the variety that is associated with the given cover
of $GRes_{\ell_1+1}$), that restrict to specializations, $\{(v^m(s),ep^m_i,p^m_i(s),q^m_j(s))\}$, so that for every index $s$, 
$(p^m_i(s),q^m_j(s)) \in Def(p,q)$. 

We continue as we did in the first part of the argument. First, we 
swap the raws in the subtriangle of specializations. 
Note that after this swap, for each pair of indices, $\ell_1+1 \leq i,j \leq m$, the corresponding specializations, $(p^m_i,q^m_j) \in 
Def(p,q)$,
if and only if $j \leq i$.

We apply the second iterative procedure that was used in the first part of the argument,
for the subtriangle of specializations. 
After finitely many steps we obtain a quotient resolution, $GRes_{\ell_2}$,
and a subtriangle of specializations, $\{p^m_i\}$ and $\{q^m_j\}$, $\ell_2 +1 \leq i,j \leq m$, that are associated with it, with
associated elliptic specializations, $\{ep^m_i\}$.
With the quotient (ungraded) resolution, we associate a cover (ungraded) resolution, that terminates in a f.p.\ group, such that all the
subtriangle of specializations factor through that cover. Furthermore, for each triple of indices $i,j,m$, $\ell_2+1 < i < j \leq m$, 
there exists
a test sequence in the fiber that contains $q^m_j$ (in the variety that is associated with the given cover
of $GRes_{\ell_2+1}$), that restrict to specializations, $\{(v^m(s),ep^m_i,p^m_i(s),q^m_j(s))\}$, so that for every index $s$, 
$(p^m_i(s),q^m_j(s)) \notin Def(p,q)$. 

We continue as we did in the first part of the argument.
At this point we once again swap the raws in the subtriangle of specializations. 
Then we apply the first procedure that was used in the first part of the argument, until we construct a quotient (ungraded) resolution,
a cover of this quotient resolution, and a subtriangle of specializations, $\{p^m_i\}$ and $\{q^m_j\}$, with associated elliptic
specializations, $\{ep^m_i\}$, that factor through the cover resolution.  
Furthermore, for each pair of indices $i,m$, $i<m$, there exists
a test sequence in the fiber that contains $(v^m,q^m_j)$ 
(in the variety that is associated with the given cover of the quotient resolution),
that restrict to specializations, $\{(v^m(s),ep^i_m,p^m_i(s),q^m_j(s))\}$, so that for every index $s$, 
and every $j$, $i \leq j \leq m$,
$(p^m_i(s),q^m_j(s)) \in Def(p,q)$. 

Then we swap the raws of the subtriangle once again, and apply the second iterative procedure that was used in the first part
of the argument. We obtain a quotient (ungraded) resolution, a cover of this quotient resolution (with a f.p.\ terminal limit group),
and a subtriangle of specializations,
 $\{p^m_i\}$ and $\{q^m_j\}$, with associated elliptic
specializations, $\{ep^m_i\}$, that factor through the cover resolution.  
Furthermore, for each triple of indices $i,j,m$, $ i < j \leq m$, 
there exists
a test sequence in the fiber that contains $q(v^m,^m_j)$ 
(in the variety that is associated with the given cover of the quotient resolution),
$\{(v^m(s),ep^m_i,p^m_i(s),q^m_j(s))\}$, so that for every index $s$, 
$(p^m_i(s),q^m_j(s)) \notin Def(p,q)$.

As we argued in the proof of theorem 8.6, this iterative procedure of swapping and continuing along the two iterative procedures 
until they terminate, terminates after finitely
many steps, by the proof of theorem 4.4. Once it terminates, we are left with an ungraded (quotient) resolution, $ERes_{\ell_t}$, 
 a cover of this  resolution that satisfies the 
properties of cover  resolutions that are listed in theorem 1.21 (and in particular has a f.p.\ completion and terminal limit group), 
$CERes_{\ell_t}$, 
and a
subtriangle of specializations, 
(still denoted) $\{ep^m_i\}$ and $\{q^m_j\}$, $\ell_t+1 \leq i,j \leq m$ that extend to specializations that factor through the cover 
resolution, $CERes_{\ell_t}$, for which (the corresponding specialization) $(p^m_i,q^m_j) \in Def(p,q)$ if and only if $\ell_t+1 \leq i \leq j \leq m$. 
Furthermore, each specialization
$q^m_j$ extends to a specialization that factors through $CERes_{\ell_t}$. A test sequence of specializations in the fiber of
$CERes_{\ell_t}$ that contains $q^m_j$, restricts to specializations, $\{(v^m(s),ep^m_i,p^m_i(s),q^m_j(s))\}$, for which 
(for every index $s$),
$(p^m_i(s),q^m_j(s)) \in Def(p,q)$, if and only if $\ell_t+1 \leq i \leq j \leq m$.

For brevity, we denote the ungraded resolution that is constructed in step $\ell_t$, $ERes$, and its closure, $CERes$.
Recall that the subgroup $<q>$ is embedded in the graded cover $CERes$, and the elements $p$ can be expressed as words in
the elements that are associated with $v^m$ in $CERes$, with the elliptic elements $ep^m_i$. Given the cover resolution, $CERes$, we look at 
all its test sequences, for which the restriction of the test sequences to the elements $p,q$, extend to specializations of the various 
terminal limit groups of the resolutions, $DRes$, and of their associated auxiliary resolutions, that are
all associated with $Def(p,q)$ according to theorem 6.1. All these test sequences, together with their extended specializations
(to the terminal limit groups of the graded resolutions $DRes$ and their auxiliary resolutions), can be collected in finitely
many graded resolutions that have similar properties to the resolutions in the formal Makanin-Razborov diagram in theorems 2.6
and 2.7. In particular, all the resolutions that are associated with these test sequences and with the closure, $CERes$, do
all have f.p.\ completion and terminal limit groups. We denote these ungraded resolutions that are associated with the closure,
$CERes$, $ACERes$.  

The specializations in the subtriangle that is associated with the closure, $CERes$, have the property that a generic point
(test sequence)  in
the fiber that contains a pair $(p^m_i,q^m_j)$, restricts to specializations,  
$(p^m_i(s),q^m_j(s)) \in Def(p,q)$, if and only if $ i \leq j \leq m$. The properties of the
graded resolutions, $DRes$, and their associated auxiliary resolutions, as listed in the statement of theorem 6.1, 
enable one to reduce a sentence over free products to
a finite disjunction of conjunctions of sentences over the factors of the free product in theorem 6.3. 
In exactly the same way (as in the argument that was used to prove theorem 6.3), the constructed resolutions, $ACERes$, enable
one to reduce the question of whether or not a tuple, $(p^m_i,q^m_j)$, is in $Def(p,q)$, to a finite disjunction of
conjunctions  of predicates in the
free variables (parameters), $ep^m_i$ and $eq^m_j$, over the factors $A$ and $B$ of the ambient free product $G=A*B$.

Hence, the existence
of the (infinite) subtriangle of specializations that is associated with $CERes$, contradicts the stability of the constructed  predicates
over the factors $A$ and $B$. Therefore, $Def(p,q)$ is stable, and theorem 8.1 follows for coefficient-free definable sets (predicates).

The stability of all coefficient-free predicates over free products of stable groups implies the stability of all predicates
over such free products. Indeed, we can replace each predicate by a coefficient-free predicate, by replacing coefficients with free
variables (parameters). The stability of the obtained coefficient-free predicate clearly implies the stability of the
original predicate. This proves the stability of free products of two stable groups and concludes the proof of theorem 8.1.

\line{\hss$\qed$}

Theorem 8.1 proves the stability of free products of finitely many stable groups.
In case we look at a sequence of groups for which every sentence is uniformly stable for the sequence, we obtain stability of the free
product of the entire sequence.

\vglue 1pc
\proclaim{Theorem 8.7} Let $G_1,G_2,\ldots$ be a sequence of groups. Suppose that every sentence (over groups) $\Phi$ is uniformly stable over the
sequence $\{G_i\}$. Then the infinite free product, $G_1*G_2*\ldots$, is stable.
\endproclaim

\nfp Let $G=G_1*G_2*\ldots$, and let $Def(p,q)$ be a coefficient-free definable set over $G$. The graded resolutions, $DRes$, and their
associated auxiliary resolutions, that were associated with $Def(p,q)$ in theorem 6.1, are universal (although not canonical), hence,
the conclusion of theorem 6.1 is valid for all countable free products and not just for finite free products.

The iterative procedures that were used to analyze the set,
$Def(p,q)$, over a free product of the form $A*B$, in the argument that was used to prove theorem 8.1, work without a change over
a countable free product. Hence, if $Def(p,q)$ is not stable over $G$, we can associate with it a quotient ungraded resolution, $ERes$, a
cover of this ungraded resolution, $CERes$, and a triangle of specializations, $\{p^m_i\}$ and $\{q^m_j\}$, $i,j \leq m$, that factor
through that cover, $CERes$, precisely as we did with non-stable predicate over a finite free product. 

With the cover $CERes$, we associate finitely many ungraded resolutions, $ACERes$, precisely as we did in the proof of theorem
8.1. Now, the instability of $Def(p,q)$, that is translated into the subtriangle of specializations that factor through
the cover, $CERes$, implies the existence of a predicate over the various factors in the free product $G=G_1*G_2*\ldots$,
that is not uniformly stable. This contradicts our assumption on the sequence of groups, $\{G_i\}$.

As we noted in the proof of theorem 8.1, every predicate can be transformed into a coefficient-free predicate by replacing
the constants by free variables (parameters). Hence, if all coefficient-free predicates are stable over the countable
free product  $G$, all the predicates are stable.

\line{\hss$\qed$}

\vglue 1.5pc
\centerline{\bf{\S9. Noetherianity}} 
\medskip

In this section we prove another basic property that a free product inherits from its factors, the Noetherianity of varieties.
This problem was brought to our attention by Ilya Kazachkov [Ca-Ka]. Its proof doesn't really require what we already proved in this
paper, but it does rely on the techniques that were used in constructing the  Makanin-Razborov diagram for varieties over a free product [Ja-Se].
Abderezak Ould-Houcine has informed us that he has an alternative proof of the problem. 
 
\vglue 1pc
\proclaim{Theorem 9.1} Let $A,B$ be equationally Noetherian groups. Then $A*B$ is equationally Noetherian.
\endproclaim

\nfp As an infinite dihedral group is linear, the theorem follows for $D_{\infty}$ by Guba's theorem [Gu]. Hence, we may assume that
$A$ and $B$ are equationally Noetherian, and $A*B$ is a non-trivial free product that is not isomorphic to $D_{\infty}$.

Suppose that there exists a coefficient-free  infinite system, $\Sigma(x)=1$, where $x$ is a finite tuple of variables, 
that is not equivalent to a finite subsystem. We denote by $\Sigma_n(x)=1$ the finite subsystem of $\Sigma(x)=1$ that contains 
the first $n$
equations (in $\Sigma$).

Since $\Sigma$ is not equivalent to a finite subsystem, there must exist a sequence of integers: $n_1 <n_2 < \ldots$, for which there exist
specializations of the tuple $x$ in $A*B$: $x_1,x_2,\ldots$, such that for every index $k$, $\Sigma_{n_k}(x_k)=1$ and $\Sigma_{n_{k+1}}(x_k) \neq 1$,
in $A*B$.

By theorem 18 in [Ja-Se], from the sequence of specializations, $\{x_k\}$, it is possible to extract a subsequence that converges into a 
resolution $Res$ (over free products): $L_0 \to L_1 \to \ldots \to \L_s$, where the $L_i$'s
 are limit groups over free products, the initial limit group,
$L_0$, is naturally generated by the elements $x$, that correspond to the variables in the system of equations,
$\Sigma(x)=1$,
and $L_s$ is a free product of finitely many elliptic factors and a (possibly trivial) free group. Note that the elliptic factors in $L_s$ need
not be finitely presented, hence, we can not deduce from theorem 18 in [Ja-Se] that a subsequence of the specializations $\{x_k\}$ factor
through the resolution $Res$. 

Also, since for each index $k$, $\Sigma_{n_k}(x_k)=1$, and the resolution, $Res$, is obtained from a convergent subsequence of the
specializations, $\{x_k\}$, all the equations in the infinite system, $\Sigma(x)$, hold as relations of the initial limit group $L_0$ of
$Res$.   

$L_s$ is the terminal limit group of the resolution, $Res$. $L_s$ need not be finitely presented. Let $L_s$ be the free product,
$A_1*\ldots*A_{\ell}*F_m$, where $A_1,\ldots,A_{\ell}$ are elliptic, and $F_m$ is a (possibly trivial) free group.  Each of the elliptic
factors, $A_i$, $i=1,\ldots,\ell$, is finitely generated and not necessarily finitely presented. Since both factors $A$ and $B$ of
the free product $G=A*B$, are equationally Noetherian, for each factor $A_i$, there exists a f.p.\ group, $FPA_i$, 
and an epimorphism, $\tau_i:FPA_i \to A_i$,
so that every homomorphism  $h:FPA_i \to A$ or $f:FPA_i \to B$, factors as $h= \hat h \circ \tau_i$ or $f = \hat f \circ \tau_i$, 
where $\hat h$ and $\hat f$ are homomorphisms from $A_i$ to $A$ and $B$ in correspondence. 

With the (well-structured) resolution, $Res$, we associate its completion, $Comp(Res)$. We set the completion, $Comp(FPRes)$, to be the
completion, that is obtained from $Comp(Res)$ by replacing the terminal limit group  of $Res$, $L_s$, by the f.p.\ group,
$FPL_s=FPA_1*\ldots *FPA_{\ell}*F_m$. $L_0$ is naturally embedded in $Comp(Res)$, and is generated by the elements $x$ that correspond
to the variables in the system, $\Sigma(x)=1$. Similarly, there is a subgroup generated by the elements $x$ in the completion, 
$Comp(FPRes)$.

Since a subsequence of the specializations, $\{x_k\}$ converges into the resolution, $Res$, and $Comp(FPRes)$ (and $FPL_s$) are
finitely presented, there exists a subsequence of this subsequence (of $\{x_k\}$) that factors through $Comp(FPRes)$, i.e.,
each of these specializations can be extended to a specialization that factors through $Comp(FPRes)$. Since $Comp(FPRes)$ differs
from $Comp(Res)$ only in the terminal limit group ($FPL_s$ versus $L_s$), and since every homomorphism of the terminal limit group
of  $Comp(FPRes)$, $FPL_s$, factors through the epimorphism onto the limit group $L_s$, the subsequence of the sequence of specializations,
$\{x_k\}$, that factor through $Comp(FPRes)$, factor through the original completion, $Comp(Res)$ (and through the original resolution, $Res$).
 However, all the equations of the system $\Sigma(x)=1$, hold as relations for the limit group $L_0$, that is generated by the elements
$x$ that correspond to the variables of $\Sigma$. This
 contradicts the assumption that for each $k$, $\Sigma_{n_k}(x_k)=1$, but $\Sigma_{n_{k+1}}(x_k) \neq 1$.

This proves that every coefficient free system of equations over $A*B$ is equivalent to a finite subsystem. For systems of equations
with coefficients, we need to repeat the same argument, working over the specific free product $G=A*B$. In this case, we first assume that
both $A$ and $B$ are countable (and then generalize the argument for any $A$ and $B$). We start
with  an infinite system 
of equations with coefficients, $\Sigma(x,a_j,b_j)=1$, where $x$ is a finite tuple of variables, and $a_j$ and $b_j$ are (finite) tuples of 
elements in $A$ and $B$, in correspondence, for each equation $j$, and the system $\Sigma$ is not equivalent to a finite subsystem. 

We choose the specializations $\{x_k\}$, and indices $ n_1 < n_2 < \ldots$, in the same way as in the coefficient-free case. i.e., 
$\Sigma_{n_k}(x_k,a_j,b_j)=1$, and 
$\Sigma_{n_{k+1}}(x_k,a_j,b_j) \neq 1$, for each index $k$. By modifying the proof of theorem 18 in [Ja-Se] for systems of equations with
coefficients (over the same coefficient group), there exists a subsequence of the specializations, $\{x_k\}$, that converges into a 
resolution with coefficients, $L_0 \to \ldots \to L_s$, where $L_0$ is a countable limit group, that is generated by the coefficient group $A*B$ 
and elements $x$
that correspond to the variables in the system $\Sigma$, and $L_s$, the terminal limit group, is a free product:
$C*D*A_1*\ldots*A_{\ell}*F_m$, where $A<C$, $B<D$, the $A_i$'s are elliptic factors, and $F_m$ is a possibly trivial free group.
Furthermore, $C$ is f.g.\ relative to $A$, $D$ is f.g.\ relative to $B$, and the $A_i$'s are all finitely generated. Since the initial
limit group $L_0$ is generated by $A$, $B$, and elements $x$ that correspond to the variables in $\Sigma$, and since $L_0$ is obtained as a
limit from a subsequence of the specializations, $\{x_k\}$, all the equations in $\Sigma$ hold as relations for the limit group, $L_0$  (note that
$L_0$ is a countable, not necessarily f.g.\ limit group, that is f.g.\ relative to the elliptic subgroups 
$A$ and $B$. Still all the basic results on limit
groups over free products that were obtained in [Ja-Se], including theorem 18, hold for such limit groups. The modification of
the basic results of [Ja-Se] from f.g.\ groups to countable 
groups that are f.g.\ relative to finitely many elliptic (coefficient) groups,
requires (in addition to [Ja-Se]), the techniques of [Gu] and [Se11], that allow one to analyze super-stable actions of 
such groups on real trees, by presenting them as inductive limits of f.g. groups, where these f.g.\ groups are
defined iteratively by enlarging the stabilizer of a single vertex in each step).

Since both factors $A$ and $B$ are equationally Noetherian, we can replace the completion of the constructed resolution
(with coefficients), with a completion that is f.p.\ relative to $A$ and $B$. This means that we can replace the terminal limit group,
$L_s$, by a free product, $FPC*FPD*FPA_1*\ldots*FPA_{\ell}*F_m$, where:
\roster
\item"{(1)}" $FPC$ is f.p.\ relative to $A$, and $FPD$ is f.p.\ relative to $B$. There exist natural epimorphisms, $FPC \to C$, 
that restricts to the identity homomorphism on $A$, and $FPD \to D$, that restricts to the identity homomorphism on $B$, 
such that every homomorphism, $h:FPC \to A$ and $f:FPD \to B$, that restricts to the identity homomorphism on $A$ and $B$ in 
correspondence,
factors through the corresponding  natural epimorphism: $FPC \to C$ or $FPD \to D$.

\item"{(2)}" $A_1,\ldots,A_{\ell}$ are all finitely presented.
For every index $i$, $i=1,\ldots,\ell$, there exists a natural epimorphism, $FPA_i \to A_i$, 
such that every homomorphism, $h:FPA_i \to A$ and $f:FPA_i \to B$,
factors through  the natural epimorphism: $FPA_i \to A_i$.
\endroster

Since the modified completion is f.p.\ relative to $A$ and $B$, there exists a subsequence of the specializations, $\{x_k\}$, that
factor through this modified completion. 
Since the modified completion differs from the original completion  only 
 in the terminal limit group, and since every homomorphism of the terminal limit group
of the modified completion into $A*B$, that restricts to the identity on $A$ and $B$, factors through the epimorphism onto the limit group $L_s$, 
the terminal limit group of the original completion, the subsequence of the sequence of specializations,
$\{x_k\}$, that factor through the modified completion,  factor through the original completion. As in the
coefficient-free case, this contradicts the fact that all
  the equations in the system $\Sigma(x,a_j,b_j)=1$, hold as relations for the limit group $L_0$, the initial limit
group of the original resolution, that is generated by the elements
$x$ that correspond to the variables of $\Sigma$, together with $A$ and $B$, as we have assumed  
  that for each index $k$, $\Sigma_{n_k}(x_k,a_j,b_j)=1$, but $\Sigma_{n_{k+1}}(x_k,a_j,b_j) \neq 1$.

This proves that every  system of equations with coefficients over $A*B$, where both $A$ and $B$ are countable,  
is equivalent to a finite subsystem. For system of equations over free products $A*B$, for general equationally Noetherian $A$ and $B$, we
need to slightly modify our argument, and essentially reduce it to the countable case.

Once again,  we start
with  an infinite system 
of equations with coefficients, $\Sigma(x,a_j,b_j)=1$, where $x$ is a finite tuple of variables, and $a_j$ and $b_j$ are (finite) tuples of 
elements in $A$ and $B$, in correspondence, for each equation $i$, and the system $\Sigma$ is not equivalent to a finite subsystem ($A$ and $B$ are general groups).

We choose the specializations $\{x_k\}$, and indices $ n_1 < n_2 < \ldots$, in the same way as in the countable case. i.e., 
$\Sigma_{n_k}(x_k,a_j,b_j)=1$, and 
$\Sigma_{n_{k+1}}(x_k,a_j,b_j) \neq 1$, for each index $k$. 

\noindent
Each specialization $x_k$ is an element in the free product $A*B$. We set $A_c$ to be a countable subgroup of $A$ that is 
generated by the elements, $\{a_j\}$, that appear in the various equations of the system, $\Sigma(x,a_j,b_j)=1$, together
with all the elements of $A$ that appear in the words in $A*B$ that are associated with the countable collection of elements 
$\{x_k\}$. Similarly, we set $B_c$ to be a countable subgroup of $B$ that is generated by $\{b_j\}$, and the elements of $B$
that appear in the words in $A*B$ that are associated with the elements $\{x_k\}$.

Now,  both the coefficients in the equations of $\Sigma(x,a_j,b_j)=1$, and the specializations, $\{x_k\}$, are all contained 
in the countable subgroup, $A_c*B_c$. Hence, we repeat the construction of the resolution in the countable case, and extract
a subsequence of the specializations $\{x_k\}$, that converges into a resolution:
$L^c_0 \to \ldots \to L^c_s$ over the countable free product, $A_c*B_c$. The initial limit group, $L^c_0$ , is generated
by  $A_c$, $B_c$, and finitely many elements that are associated with the variables $x$. Since $L^c_0$ is obtained as a limit of
a subsequence of the specializations, $\{x_k\}$, all the equations in the system, $\Sigma(x,a_j,b_j)$, hold as relations
in the limit group, $L^c_s$.

$L^c_s$, the terminal limit group, is f.g.\ over the groups $A_c$ and $B_c$, and admit a free product:
$C_c*D_c*A_1*\ldots*A_{\ell}*F_m$, where $A_c<C_c$, $B_c<D_c$, the $A_i$'s are elliptic factors, and $F_m$ is a possibly trivial 
free group. We replace the subgroups $C_c$ and $D_c$, by the groups, $C=C_c*_{A_c} \, A$, and $D=D_c*_{B_c} \, B$. As $C_c$ and 
$D_c$ are f.g.\ over $A_c$ and $B_c$ in correspondence, $C$ and $D$ are f.g.\ over $A$ and $B$, in correspondence.
We set $L_s$ to be $L_s=C*D*A_1*\ldots *A_{\ell}*F_m$. $L_s$ is f.g.\ over $A$ and $B$.

With the resolution $Res_c$: $L^c_0 \to \ldots \to L^c_s$, we can associate its completion, $Comp(Res_c)$. With $Comp(Res_c)$ we
associate another completion, that we denote $Comp(Res)$, that is obtained from the completion, $Comp(Res_c)$, by replacing
the terminal limit group of the resolution $Res_c$, $L^c_s$, with the group, $L_s$. Note that $Comp(Res_c) < Comp(Res)$,
and since all the equations in
the system, $\Sigma$, hold as relations in the limit group $L^c_0$, and hence in $Comp(Res_c)$, they hold as relations in
the completion, $Comp(Res)$, as well.

Since both factors $A$ and $B$ are equationally Noetherian, we can replace the completion, $Comp(Res)$,
with a completion that is f.p.\ relative to $A$ and $B$. This means that we can replace the terminal limit group,
$L_s$, by a free product, $FPC*FPD*FPA_1*\ldots*FPA_{\ell}*F_m$, that satisfy properties (1) and (2) that were used in
modifying the completion in case $A$ and $B$ are countable. In particular, the modified free product is f.p.\ relative to 
$A$ and $B$, and any homomorphism of the modified free product into $A*B$, that restricts to the identity on $A$ and $B$,
factors through the completion, $Comp(Res)$.

Since the modified completion is f.p.\ relative to $A$ and $B$, there exists a subsequence of the specializations, $\{x_k\}$, that
factor through this modified completion. Hence, this subsequence factor (i.e., can be extended to specializations that factor)
through the completion, $Comp(Res)$.
As in the
coefficient-free and countable cases, this contradicts the fact that all
  the equations in the system $\Sigma(x,a_j,b_j)=1$, hold as relations in $Comp(Res)$, whereas we assumed
  that for each index $k$, $Sigma_{n_k}(x_k,a_j,b_j)=1$, but $\Sigma_{n_{k+1}}(x_k,a_j,b_j) \neq 1$.
This finally proves that $A*B$ is equationally Noetherian.

\line{\hss$\qed$}


\vskip 2em

\end